\documentclass{amsart}
\usepackage[utf8]{inputenc}
\usepackage[T1]{fontenc}
\usepackage{amsmath,amssymb,amsthm,dsfont}

\usepackage[english]{babel}
\usepackage[all]{xy}

\theoremstyle{plain}
\newtheorem{Th}{Theorem}[section]
\newtheorem{Lem}[Th]{Lemma}
\newtheorem{Cor}[Th]{Corollary}
\newtheorem{Prop}[Th]{Proposition}
\newtheorem{Def}[Th]{Definition}

\newtheorem{Rem}[Th]{Remark}

\newtheorem{Thintro}{Theorem}

\newcommand{\QQ}{\mathbb{Q}}
\newcommand{\KK}{\mathbb{K}}
\renewcommand{\AA}{\mathbb{A}}
\DeclareMathOperator{\im}{im}

\DeclareMathOperator{\rank}{rk}

\DeclareMathOperator{\GL}{GL}
\DeclareMathOperator{\proj}{P}
\DeclareMathOperator{\ELL}{L}

\DeclareMathOperator{\End}{End}
\DeclareMathOperator{\supp}{supp}
\newcommand{\ie}{\textit{i.e.}\;}

\newcommand{\cf}{\textit{c.f.}\;}
\newcommand{\dist}{\mathrm{d}}
\newcommand{\ZZ}{\mathbb{Z}}
\newcommand{\NN}{\mathbb{N}}
\newcommand{\FF}{\mathbb{F}}
\newcommand{\PP}{\mathbb{P}}
\newcommand{\EE}{\mathbb{E}}
\newcommand{\CC}{\mathbb{C}}
\newcommand{\RR}{\mathbb{R}}

\newcommand{\eps}{\varepsilon}
\newcommand{\essker}{\underline{\ker}}
\newcommand{\evrank}{\underline{\rank}}
\newcommand{\eigen}{\mathcal{E}}


\title[Products of random matrices]{Limit theorems for a strongly irreducible product of independent random matrices under optimal moment assumptions}

\author{Axel P\'eneau}
\address{Université de Tours, Université d’Orléans, CNRS, IDP, UMR 7013, Tours, France}

\email{axel.peneau@univ-tours.fr}
\date{\today}

\begin{document}

	\begin{abstract}
		Let $\nu$ be a probability distribution over the semi-group of square matrices of size $d \ge 2$ over a locally compact field $\mathbb{K}$, \textit{e.g.} $\mathbb{R}$.
		We consider the random walk $\overline{\gamma}_n := \gamma_0\cdots\gamma_{n-1}$ for $(\gamma_k)_{k \in \mathbb{N}}$ independent of law $\nu$.
		Let $s_1 \ge s_2 \ge \dots \ge s_d$ be the singular values given by the Cartan projection. 
		Under a contraction assumption on $\nu$, we show that $(\log\frac{s_1}{s_2}(\overline{\gamma}_n))_{n \in\NN}$, escapes to infinity linearly and satisfies exponential large deviations inequalities below its escape rate. 
		This extends the notion of simplicity of the top Lyapunov exponent.
		We also show that the image of a generic line by $\overline{\gamma}_n$ as well as its eigenspace of maximal eigenvalue both converge to the same random line $\ell^\infty$ at an exponential speed.

		If we moreover assume that $\nu$ is supported on the group of invertible matrices and that the push-forward distribution $N_*\nu$ is $\mathrm{L}^p$ for $N: g \mapsto\log\|g\|\|g^{-1}\|$ and for some $p > 0$, then we show that $- \log\mathrm{d}(\ell^\infty, H)$ is uniformly $\mathrm{L}^p$ for all proper subspace $H \subset \mathbb{R}^d$.
		For $p = 1$, we moreover show that the rescaled logarithm of each coefficient of $\overline{\gamma}_n$ almost surely converges to the top Lyapunov exponent.
		
		To prove these results, we do not rely on the existence of the stationary measure nor on the existence of the Lyapunov exponents.
		Instead we describe an effective way to group the i.i.d. factors into i.i.d. random words that are somehow aligned in the Cartan decomposition. 
		We moreover have an explicit control over the moments. 
	\end{abstract}
	
	\maketitle
	

	\section{Introduction}
	
	\subsection*{Motivation and main results}
	
	To highlight the strength and versatility the approach developed in the present article let us state two important results proven in the present article.
	Let $d \ge 2$ and let $E = \RR^d$ be the Euclidean space endowed with its Euclidean norm and with canonical basis $(e_1, \dots, e_d)$. 
	Let $(\gamma_n)_{n \ge 0} \in \GL(E)$ be a random independent and identically distributed sequence of invertible square matrices.
	In other words $\gamma : \Omega \to \GL(E)^\NN$ is a measurable map such that $\gamma_*\PP = \nu^{\otimes\NN}$ for $\nu$ a probability measure on $\GL(E)$ and $(\Omega, \PP)$ a probability space.
	To $\gamma$, we associate the random walk $(\overline{\gamma}_n)_{n \ge 0}$ defined as\footnote{Through the present paper, we sometimes write $A := B$ instead of $A = B$, to specify that we define a free expression $A$ as being equal to a given expression $B$, like here, or to remind that a given $A$ is by definition equal to $B$, like $e := \sum_{k = 0}^{+\infty}1/k!$. This is used to disambiguate with the use of the notation $A = B$ to claim that a given $A$ is equal to a given $B$, like $\ZZ^\times = \{-1, 1\}$ or to say that we now assume that two given $A$ and $B$ are equal, like "let $\ker(\phi) = \im(\phi)$".} $\overline{\gamma}_n:= \gamma_0 \cdots \gamma_{n-1}$ for all $n \ge 0$.
	The random matrix $\overline{\gamma}_n$ is then a product of $n$ independent random matrices.	
	
	This simple model of non commuting random walk has been intensely studied during the last century with applications to branching processes, quantum physics, number theory and geometric group theory, discussed in more details in Paragraph \ref{intro:background} of the present article.
	
	We write $^te_i \overline{\gamma}_n e_j$ for the entry of $\overline{\gamma}_n$ at the coordinates $(i,j)$.
	We write $\rho_1(\overline{\gamma}_n) \ge \rho_2(\overline{\gamma}_n) \ge \dots \ge \rho_d(\overline{\gamma}_n) \ge 0$ for the spectral values of $\overline{\gamma}_n$ \ie the moduli of its eigenvalues. 
	We say that a matrix $h$ is proximal when $\rho_1(h) > \rho_2(h)$, which means that the eigenspace of $h$ associated to its top eigenvalue is a line.
	
	We say that a semi-group $\Gamma < \GL(E)$ is strongly irreducible if there is no subspace $\{0\} \subsetneq V \subsetneq \RR^d$ that is virtually invariant, \ie there exists $F \subset \Gamma$ finite such that $FV = \Gamma V$.
	We say that a semi-group $\Gamma < \GL(E)$ is proximal if it contains a proximal element, this is equivalent to saying that the closure of $\RR \Gamma$ contains a rank one projection.

	Through the present paper, we assume that $\nu$ is strongly irreducible and proximal in the sense that the semi-group generated by its support is.
	This notion is discussed in more details in Section \ref{intro:def}.
	A notable qualitative result proven in the present article is the following.
	
	\begin{Thintro}[Law of large numbers for the coefficients]\label{th:intro-lln}
		Assume that the measure $\nu$ is strongly irreducible and proximal and let $(\gamma_k)_{k \ge 0} \sim \nu^{\otimes\NN}$.
		Assume that $\EE (\log^+\|\gamma_0\|) < +\infty$ and $\EE(\log^+\|\gamma_0^{-1}\|) < +\infty$.
		Then there exists a constant $\lambda_1 = \lambda_1(\nu)$ such that for all $1 \le i,j \le d$ and $\PP$-almost surely:
		\begin{equation*}
			\lim_{n \to +\infty} \frac{\log|^te_i\overline{\gamma}_n e_j|}{n} = \lim_{n \to +\infty} \frac{\log\rho_1(\overline{\gamma}_n)}{n} = \lambda_1.
		\end{equation*}
	\end{Thintro}
	
	Theorem \ref{th:intro-lln} is notable for answering the long-standing conjecture of convergence of the coefficients under $\ELL^1$ moment assumption.
	From the work of Furstenberg and Kesten \cite{porm}, it was known that $\frac{\log\|\overline{\gamma}_n\|}{n} \to \lambda_1$ almost surely, with $\lambda_1 = \lim \frac{\EE(\log\|\overline{\gamma}_n\|)}{n}$. 
	This is a particular case of the sub-additive ergodic Theorem, stated a few years later by Kingman in \cite{K68}.
	
	Weaker versions of Theorem \ref{th:intro-lln} were stated by Guivarc'h and Lepage \cite{GuivarchTCL}, \cite{lepage82} with an exponential moment condition \ie assuming that $\EE(\|\gamma_0\|^\beta) < +\infty$ and $\EE(\|\gamma_0^{-1}\|^\beta) < +\infty$ for some $\beta > 0$ and by Xiao, Grama and Liu \cite{xiao2021limit} using the work of Benoist and Quint \cite{Benoist2016CentralLT} with an $\ELL^2$ moment condition \ie assuming that $\EE(\log^2\|\gamma_0\|) < +\infty$ and $\EE(\log^2\|\gamma_0^{-1}\|) < +\infty$.
	
	Let us now state a notable quantitative result of the present article.
	Given $g \in\GL(E)$, we write  $s_1(g) \ge s_2(g) \ge \dots \ge s_d(g)$ for the singular values of $g$.
	It is known from Oseledets multiplicative ergodic Theorem that given $(\gamma_n )_{n \ge 0}$ i.i.d. that satisfy $\EE(\log\|\gamma_0\|) < +\infty$, there exists a Lyapunov spectrum $\lambda_1 \ge \lambda_2 \ge \dots \ge \lambda_d \ge 0$ such that $\frac{s_k(\overline{\gamma}_n)}{n} \to \lambda_k$ almost surely and for all $k \in\{1, \dots,d\}$.
	When $\lambda_2 > 0$, this implies that $\log\frac{s_1}{s_2}(\overline{\gamma}_n)/{n} \to \lambda_1 - \lambda_2$.
	In \cite{GR1985}, Guivarc'h and Raugi show that with the algebraic assumptions of Theorem \ref{th:intro-lln}, we moreover have $\lambda_1 > \lambda_2$, they also state a large deviations principle under an exponential moment assumption.

	Under the same algebraic assumptions as in \cite{GR1985}, we obtain the following remarkable result, which gives a one-sided large deviation principle and an exponential contraction result without any moment assumptions on $\nu$. 
	In this case, the Lyapunov spectrum may not be well defined.
	
	\begin{Thintro}[Exponential contraction without moment condition]\label{th:intro-contraction}
		Let $\nu$ be a strongly irreducible and proximal probability measure over $\GL(E)$ and let $\gamma := (\gamma_k)_{k \ge 0} \sim \nu^{\otimes\NN}$.
		Then there exists a constant $0 < \lambda_{1,2}(\nu) \le +\infty$ such that $\frac{1}{n}\log\frac{s_1}{s_2}(\overline{\gamma}_n) \to \lambda_{1,2}(\nu)$ almost surely.
		Moreover, for all $\alpha < \lambda_{1,2}(\nu)$, there exist constants $C, \beta > 0$ such that
		\begin{equation}
			\PP\left(\log \frac{s_1(\overline{\gamma}_n)}{s_2(\overline{\gamma}_n)} \le \alpha n \right) \le C e^{-\beta n}
		\end{equation}
		for all $n$ and there exists a random line $\ell^\infty(\gamma) \in \proj(E)$ such that for all non-random vector $v \in E\setminus\{0\}$:
		\begin{equation}\label{contraction-intro}
			\PP\left(\dist([\overline{\gamma}_n v], \ell^\infty(\gamma)) \ge  e^{-\alpha n}\right) \le C e^{-\beta n}.
		\end{equation}
	\end{Thintro}

	Here $[v] := \RR v \in \proj(E)$ denotes the projective class of $v \in E \setminus \{0\}$ and $\dist$ is the natural distance on the compact set $\proj(E)$.
	The almost sure convergence $[\overline{\gamma}_n v] \to \ell^\infty$ was proven in \cite{GR1985} via a martingale argument.
	However, the exponential contraction in probability given by the right hard member of \eqref{contraction-intro} was only proven under a finite exponential moment assumption on $\log\|\gamma_0^{-1}\|\|\gamma_0\|$.
	
	Theorem \ref{th:intro-contraction} tells that, contrary to classical large deviations inequalities for sums of random variables, this exponential contraction is not guaranteed by the tameness of $\nu$ but by the geometry of $\GL(E)$.
	Another strength of the framework developed in the present article is that small modification of the measure $\nu$ in the weak-$*$ topology (which may break the moment assumption) do not break the exponential contraction given by \eqref{contraction-intro}.
	
	Before giving more details on the techniques used in the present article, let us quickly remind the traditional framework used to study products of random matrices in order to put emphasis on the innovative nature of the present paper.
	
	An interesting observation of Furstenberg (which is the ground for boundary theory \cite{furstenberg1973boundary}) is that since $\GL(E)$ acts by homeomorphisms on the projective space $\proj(E)$, which implies that there exists a stationary probability measure $\xi$ on $\proj(E)$ (\ie given a random direction $[x] \sim \xi$, independent of $\gamma_0$, we have $[\gamma_0 x] \sim \xi$).
	Under the hypotheses of Theorem \ref{th:intro-contraction}, this stationary measure is unique and is the law of $\ell^\infty(\gamma)$.
	Given $[x_0]\sim \nu$ independent of $\gamma$, the sequence $\log\|\gamma_{n-1} \cdots \gamma_0 x\|$ is a stochastic process that can be studied with classical ergodic theoretic tools.
	
	With this kind of techniques, Guivarc'h, Raugi and Lepage only obtained the analogous to Theorem \ref{th:intro-lln}, under an exponential moment assumption on $\log^+\|\gamma_0\|+\log^+\|\gamma_0^{-1}\|$, that was upgraded three decades later to a second order moment assumption by Xiao, Grama and Liu \cite{xiao2021limit}, using the regularity results of Benoist and Quint \cite{CLT16}, (see Corollary \ref{cor:regxi} for an improvement of this regularity result).
	
	The issue with this approach is that the "moments" of $\nu$ appear to control the moment of each step of the stochastic process, to control its mixing speed and to control the regularity of the stationary measure.
	The approach proposed in the present work takes inspiration from the pivoting technique described in \cite[Section~2]{pivot}. 
	The technique used in section \ref{sec:pivot}
	Nevertheless, we do not rely on Gouëzel's results since Gromov's hyperbolicity axioms are not satisfied in our case.
	For a reader who is not familiar with traditional operator  techniques, the originality of this approach has the advantage of making the present paper self included, it also means that we can not use historical results to take shortcuts through technicities.
	
	Though it is done at a high technical cost, we do not simply extend known results using stronger ergodic theoretic tools but offer a new point of view on products of random matrices.
	We construct the invariant measure $\xi$ on $\proj(E)$ as the distribution of the almost sure limit of the sequence $([\overline{\gamma}_n x])_{n \ge 0}$.
	Without assuming the matrices to be invertible, the limit of $([\overline{\gamma}_n x])_{n \ge 0}$ exists almost surely for $x$ outside of a countable union of subspaces of $E$ that have dimension $k < d$, on which $\ell^\infty(\gamma)$ is the non-zero limit point of the sequence $([\overline{\gamma}_n x])_{n \ge 0}$ (\textit{c.f.}, Theorem \ref{th:cvspeed} and Lemma \ref{lem:descri-ker} and \ref{lem:essker}).
	Since the limit does not depend on $x$, we denote it by $\ell^\infty(\gamma)$.
	It follows that $\ell^\infty$ is an equivariant map so $\xi = \ell^\infty_*\nu^{\otimes\NN}$ is $\nu$-stationary.
	From \eqref{contraction-intro}, we deduce an exponential mixing result for the associated Markov chain (\textit{c.f.}, Corollary \ref{cor:ergodic} of Theorem \ref{th:cvspeed}).
	
	Assuming the matrices to be invertible, we moreover give a constructive and explicit link between the distribution of $\log\|\gamma_0\|+\log\|\gamma_0^{-1}\|$ and the regularity of $\xi$ (\textit{c.f.}, Corollary \ref{cor:regxi} of Theorem \ref{th:slln}). 
	Theorem \ref{th:intro-lln} follows from this result.

	In Paragraph \ref{intro:def} we remind some classical notations for semi-groups of matrices and give a proper definition of the assumptions made on the i.i.d. sequence of matrices we consider.
	The main technicality of this part is to properly define the notion of strong irreducibility and proximality for a semi-group of singular matrices, which is not present in the literature. 
	
	In Paragraph \ref{sec:intro-pivot} we introduce the notion of alignment, that is central in the present work and detailed in Section \ref{kak}, as well as a compact notation for grouping consecutive factors of a sequence, that is used in Sections \ref{markov} and \ref{sec:pivot}.
	We also state the main result of the present paper: Theorem \ref{th:pivot}, proven in Section \ref{sec:pivot}, which describes a coupling of our i.i.d. sequence of matrices with a random sequence of integers, such that the associated grouping is almost Markovian and satisfies nice alignment properties.
	Though central, Theorem \ref{th:pivot} is quite technical and the results of sections \ref{intro:non-invertible}, \ref{intro-flag} and \ref{intro:invertible}, which motivate the introduction of the pivotal times decomposition are stated as generalization of known results to make them more palatable.
	In Section \ref{results}, we detail how all the results stated in the introduction are deduced from Theorem \ref{th:pivot}.
	The reader who is not familiar with the classical theory of filtrations for Markov processes may skip Section \ref{sec:pivot} in a first read and simply assume Theorem \ref{th:pivot} to be true in Section \ref{results}.
	
	That being said, the real added value of the present work lies in Section \ref{sec:pivot}.
	Indeed, the constructive nature of the proof of Theorem \ref{th:ex-piv} opens for further development in applied mathematics and computation.
	Following the proof, the constants mentioned in the results of the present article can be effectively computed with enough information on $\nu$.
	Namely, we only need to know the values of the constants described in Definition \ref{def:quant-contraction}.
	
	Moreover, the axiomatic description of the pivoting technique given in Section \ref{sec:pivot} allows to translate the results of the present paper to random walks on groups and semi-groups that share large-scale geometrical properties with $\GL(E)$ or $\End(E)$, for example relatively hyperbolic groups, groups acting on Bruhat-Tits buildings or on more general non-flat $\mathrm{CAT(0)}$ spaces.

	\subsection{The rank one toy model}\label{intro:rank-one}

	We say that an endomorphism $h \in \End(E)$ of the Euclidean space $E =\RR^d$ has rank one if $h(E)$ is a line. 
	In Theorem \ref{th:intro-contraction}, we say that the ratio $\frac{s_2}{s_1}(\overline{\gamma}_n)$ has limit $0$. 
	Saying that an endomorphism $h$ has rank one is equivalent to saying that $\frac{s_2}{s_1}(h) = 0$. So morally, what Theorem \ref{th:intro-contraction} says is that looking at a product of a lot random matrices is like looking at a rank $1$ matrix.
	To get an intuition of the kind of contraction results we expect for products of random i.i.d. matrices, let us first look at products of rank one matrices.
	
	We remind that for all rank one endomorphism $h$, there exist a linear form $f$ and a vector $v$ such that $h(x) = f(x) \cdot v$ for all $x \in E$, in this case, we write\footnote{This notation can be confusing because the image of a vector $x$ is then written as $hx = v f x$ which is the product of $v$ on the left with the scalar $fx$ on the right, The order of this product doesn't matter  because $\KK$ is commutative.} $h = vf$.
	the direction $[v]$ determines the image of $h$ and the hyperplane $\ker(f)$ determines the Kernel of $h$. The operator norm $\|h\| = \max \|hx\| / \|x\|$ is then equal to $\|v\| \|f\|$.

	Let $\nu$ be a measure supported on the set of rank one matrices in $\End(E)$ and let $(\gamma_k)_{k \ge 0} \sim \nu^{\otimes\NN}$ be a random i.i.d. sequence.
	At all time $n$, assuming that the product $\overline\gamma_n = \gamma_0\cdots \gamma_{n-1}$ is not zero, it shares its image with $\gamma_0$ and its kernel with $\gamma_{n-1}$.
	It follows from Theorem \ref{th:limit-uv} that $[\overline{\gamma}_n v]$ converges exponentially fast to a random limit that does not depend on $v$. 
	In the rank one case, this image is in fact stationary, provided that $\overline{\gamma}_n v \neq 0$.
	
	On the other hand, the kernels (or by duality, the co-images) of $(\overline\gamma_n)$ form an i.i.d. sequence. 
	Moreover at all time $n \ge 2$, the kernel and the image of $\overline\gamma_n$ are independent.
	In the general case, we see that for a non random choice of $f \in E \setminus\{0\}$ and $v \in E \setminus\{0\}$ the correlation between $[\overline\gamma_n x]$ and $[f\overline\gamma_n]$ decays exponentially fast in the sense of Corollary \ref{cor:decay}.
	
	If we look at the norm of $\overline{\gamma}_n$, by a telescopic argument, we get:
	\begin{equation*}
		\|\overline\gamma_n\| = \prod_{k = 0}^{n-1} \|\gamma_k\| \frac{\|\overline\gamma_{k+1}\|}{\|\overline\gamma_k\|\|\gamma_k\|}.
	\end{equation*}
	Moreover, $\frac{\|\overline\gamma_{k+1}\|}{\|\overline\gamma_k\|\|\gamma_k\|} = \frac{\|\gamma_{k-1} \gamma_k\|}{\|\gamma_{k-1}\|\|\gamma_k\|}$ for all $k \ge 1$ so we can write:
	\begin{align*}
		\|\overline\gamma_n\| &= \prod_{k = 0}^{n-1} \|\gamma_k\| \prod_{k = 1}^{n-1}\frac{\|\gamma_{k-1} \gamma_k\|}{\|\gamma_{k-1}\|\|\gamma_k\|} \\
		& = \prod_{k = 0}^{n-1} \|\gamma_k\| \prod_{\substack{k = 1 \\ k \in 2\ZZ}}^{n-1}\frac{\|\gamma_{k-1} \gamma_k\|}{\|\gamma_{k-1}\|\|\gamma_k\|}\prod_{\substack{k = 1 \\ k \in 1+ 2\ZZ}}^{n-1}\frac{\|\gamma_{k-1} \gamma_k\|}{\|\gamma_{k-1}\|\|\gamma_k\|}.
	\end{align*}
	So $\|\overline{\gamma}_n\|$ is the product of three i.i.d. products. 
	Therefore, if $\PP(\gamma_0\gamma_1 = 0) = 0$, then $\PP(\overline\gamma_n = 0) = 0$ for all $n \ge 0$. Otherwise, $(\overline\gamma_n)$ is stationary to $0$ and the random first time $n$ such that $\overline{\gamma}_n = 0$ has a finite exponential moment.
	
	When $\PP(\gamma_0\gamma_1 = 0) = 0$,
	\begin{equation}\label{rankone:norm}
		\|\overline\gamma_n\| = \sum_{k = 0}^{n-1} \log\|\gamma_k\| + \sum_{\substack{k = 1 \\ k \in 2\ZZ}}^{n-1}\log\frac{\|\gamma_{k-1} \gamma_k\|}{\|\gamma_{k-1}\|\|\gamma_k\|} + \sum_{\substack{k = 1 \\ k \in 1+ 2\ZZ}}^{n-1}\log\frac{\|\gamma_{k-1} \gamma_k\|}{\|\gamma_{k-1}\|\|\gamma_k\|}.
	\end{equation}
	Therefore, we get a law of large numbers for $\log\|\overline\gamma_n\|$ when the random variables are integrable, \ie $\EE|\log\|\gamma_0\||$ and $\EE|\log\frac{\|\gamma_0\gamma_1\|}{\|\gamma_0\|\|\gamma_1\|}|$ are both finite. 
	Given a unit linear form $f$ and a unit vector $v$, such that $\EE|\log\frac{\|f\gamma_0\|}{\|f\|\|\gamma_0\|}|$ and $\EE|\log\frac{\|\gamma_0 v\|}{\|\gamma_0\|\|v\|}|$ are finite, we get a law of large numbers for $\log|f \overline\gamma_n v|$.
	More precisely,
	\begin{equation}\label{ranone:coef}
		\lim_n \frac{\log|f \overline\gamma_n v|}{n} = \EE\left(\log\frac{\|\gamma_0\gamma_1\|}{\|\gamma_1\|}\right) = \lambda_1(\nu).
	\end{equation}
	
	If we write $\sigma(g, v) = \log\|gv\|/\|v\|$ for the norm cocycle, seen as a measurable function $\End(E) \times \proj(E) \to [-\infty, + \infty)$, then we get that $\lambda_1(\nu) = \int \sigma(\gamma_0, v) d\nu(\gamma_0)d\xi[v]$, where $\xi$ is the distribution of the image of $\gamma_1$.
	The measure $\xi$ plays the role of a $\nu$-stationary measure in this case because $[\gamma v]$ is the image of $\gamma_0$, which has the same law as the image of $\gamma_1$.
	
	We remind that if we assume $\nu$ to be supported on the space of invertible matrices, the existence of an invariant measure is guaranteed by the fact that $\GL(E)$ acts by homeomorphisms on $\proj(E)$.
	For rank one matrices on the other hand, it is clear that the convolution $\nu* \xi$ either gives weight to the degenerate projective class, denoted by $[0]$, or it is the law of the image of $\gamma_0$.  
	Hence, the only candidate to be an invariant measure is the law of the image of $\gamma_0$ and it is an invariant measure when $\PP(\gamma_0 \gamma_1) = 0$.
	Otherwise, $\overline\gamma_n$ is stationary to $0$.

	\subsection{Notion of contraction}\label{intro:def}
	
	Let us properly define the contraction assumptions under which the results of the present paper hold.
	
	Let $\KK$ be a local field \ie a field endowed with an absolute value $|\cdot| = \KK \to \RR_{\ge 0}$ that makes it locally compact\footnote{We remind that given a topological field $\KK$ that is locally compact, one can construct an absolute value using the Haar measure on $(\KK, +)$. Given a Haar measure $\mu$, and $k \in \KK \setminus \{0\}$, we set $|k|$ to be the density of $\mu$ with respect to the push-forward of $\mu$ by $k$.}.
	Local fields are classified into three categories: the Archimedean fields, $\RR$ and $\CC$, in which $\ZZ$ embeds isometrically; the ultrametric fields of positive characteristic, which are the $\FF_q((t))$'s and their finite extensions and the ultrametric field of characteristic $0$, which are the $\QQ_p$'s and their finite extensions. For all prime number $p$, $\QQ_p = \{\sum_{k > v} a_k p^k\mid v \in \ZZ, 0 \le a_k < p, a_v \neq 0\}$ denotes the field of $p$-adic numbers, which is naturally endowed with the absolute value $|\sum_{k > v} a_k p^k| = p^{-v}$ and for $q = p^j$, $\FF_q$ denotes the finite field of cardinal $q$ and $\FF_q((t)) = \{\sum_{k > v} a_k t^k\mid v \in \ZZ, a_k \in \FF_q, a_v \neq 0\}$ denotes the field of Laurent series over $\FF_q$, which is naturally endowed with the absolute value $|\sum_{k > v} a_k t^k| = q^{-v}$.
	
	Let $d \ge 2$ be fixed and let $E = \KK^d$, when $\KK$ is Archimedean, we endow $E$ with the Euclidean (real) or Hermitian (complex) norm $\|\cdot\|: (x_1, \dots, x_d) \mapsto \sqrt{|x_1|^2 + \cdots +|x_d|^2}$, otherwise, we endow $E$ with an ultrametric norm\footnote{Not all ultrametric norms are isometric to the canonical ultrametric norm on $\KK^d$ but they are all equivalent \ie bounded by multiples of each other, which is all that matters for the present work. That being said, the geometric results of section \ref{section:contraction} hold for any choice of norm.}, for example $\|\cdot\|: (x_1, \dots, x_d) \mapsto \max\{|x_1|, \dots,|x_d|\}$.

	Given $1 \le k \le d$, we write $\bigwedge^k E$ for the $k$-th exterior product of $E$, which is simply the dual space of the space of $k$-linear alternate forms. 
	For all $j, k \ge 1$ such that $j + k \le d$, we have a natural bi-linear map $\wedge: \bigwedge^j E \times \bigwedge^k E \to \bigwedge^{j+k} E$.
	There is also a canonical\footnote{In the literature, the choice of this norm may vary up to a multiplicative constant.} norm $\|\cdot\|$ on $\bigwedge^k E$ for all $k$, which is characterized by the fact that the inequality $\|x \wedge y\| \le \|x\|\|y\|$ holds for all $j,k \ge 1$, all $x \in\bigwedge^j E$ and all $y \in \bigwedge^k E$ and is sharp in the sense that there exists $y' \in \bigwedge^k E$ such that $x \wedge y = x \wedge y'$ and such that $\|x \wedge y'\| = \|x\|\|y'\|$.
	
	To give a more explicit construction, the norm on $\bigwedge^k E$ is the Euclidean (resp. Hermitian, resp. ultrametric) norm associated to the basis $\left(e_{i_1} \wedge \dots \wedge e_{i_k}\right)_{1 \le i_1 < \dots < i_k \le d}$, where $\left(e_i\right)_{i \le i\le d}$ is the canonical basis of $E$.
	
	We write $\End(E)$ for the semi-group of linear endomorphisms of $E$, which is naturally identified with the set of square matrices of size $d$ with coefficients in $\KK$.
	We write $\GL(E)$ for the group of linear automorphisms of $E$.
	
	Given $g \in\End(E)$, we write $s_1(g) \ge s_2(g) \ge \dots \ge s_d(g) \ge 0$ for the singular values of $g$ counted with multiplicity \ie the absolute values of the coefficients of the diagonal matrix in the polar decomposition. We remind that the Cartan decomposition Theorem tells us that for all matrix $g \in \End(\KK^d)$, there exists a decomposition $g = k a k'$, where $k$ and $k'$ are isometries and $a$ is diagonal. 
	Moreover, the multi-set\footnote{A multi-set is a set $S$ with a multiplicity measure $m : S \to \NN_{\ge 1}$. Here, we see a finite multi-set as the equivalence class of a finite sequences under permutation of the indices.} $(s_1, \dots, s_d)$ of the absolute values of the coefficients of $a$ does not depend on the choice of $k$ and $k'$.
	
	The first singular value $s_1(g)$ is simply the operator norm of $g$, that we also denote by $\|g\| = \max_{v \in E \setminus \{0\}} \|g v\|/\|v\|$ and the following characterization gives the higher order singular values:
	\begin{equation}\label{def-sv}
		\forall g \in\End(E),\; \forall 2 \le k \le d, \; 
		s_k(g) = \max_{\substack{V \subset E, \\ \dim(V) = k}} \min_{v \in V \setminus\{0\}}\frac{\|g v\|}{\|v\|} = \frac{\|\bigwedge^k g\|}{\|\bigwedge^{k-1} g\|},
	\end{equation}
	where $\bigwedge^k g \in\End(\bigwedge^k E)$ is the linear map characterized by $\bigwedge^k g: x_1 \wedge \cdots \wedge x_k \mapsto gx_1 \wedge \cdots \wedge gx_k$. 
	Note that $\bigwedge^k : \End(E) \to \End(\bigwedge^k E)$ is a semi-group morphism.
	Formula \eqref{def-sv} gives two alternate definitions for the singular values that do not rely on the construction of a Cartan decomposition. 
	
	In practice, seeing singular values as ratio of sub-multiplicative quantities is what allows us to define the Lyapunov spectrum.
	Namely, for all measure $\nu$ on $\End(E)$, such that $\int \log\|g\|d\nu < + \infty$, the sequence $\int \log\|\bigwedge^k g\| d\nu^{*n}g$ is well defined and sub-additive.
	Therefore $\frac1n \int \log\|\bigwedge^k g\| d\nu^{*n}g$ has a limit, denoted by $\lambda_1(\bigwedge^k_*\nu)$ and by Kingman's sub-additive Ergodic Theorem, $\log\|\bigwedge^k \overline\gamma_n\| \to \lambda_1(\bigwedge^k_*\nu)$ for $\nu^{\otimes\NN}$ almost all sequence $(\gamma_n)_{n \ge 0}$. 
	Then for all $k > 2$, if $\lambda_1(\bigwedge^{k - 1}_*\nu) >  -\infty$, we write $\lambda_k(\nu) := \lambda_1(\bigwedge^k_*\nu) - \lambda_1(\bigwedge^{k - 1}_*\nu)$.
	Then $\log s_k(\overline\gamma_n) \to \lambda_k(\nu)$ for all $k$ and for $\nu^{\otimes\NN}$ almost all sequence $(\gamma_n)_{n \ge 0}$.

	Given $g \in\End(E)$, we write $\rho_1(g) \ge \rho_2(g) \ge \dots \ge \rho_d(g)$ for the spectral values of $g$ \ie the moduli of its eigenvalues counted with multiplicity.
	Dynamically speaking, $\rho_k(g) := \lim_n \sqrt[n]{s_k(g^n)}$. 
	We say that an endomorphism $g$ is proximal if $\rho_1(g) > \rho_2(g)$. 
	Dynamically speaking, $g$ is proximal if and only if the sequence of the projective classes $([g^n])_n$ converges to the class of a rank one projection, \ie the class of a rank one map whose image is not included in its kernel.
	
	By the dynamic characterization of the spectrum, $\log\rho_k(g) = \lambda_k(\delta_g)$ for all $k \ge 1$ and for all $g$. 
	Let us simply write $\lambda_k(g) := \log\rho_k(g)$ for the $k$-th spectral exponent of $g$.
	
	In the present article, we extend the notion of strong irreducibility to semi-groups of non-invertible matrices in the following way.

	\begin{Def}[Strong irreducibility]\label{def:sr}
		We say that a set $S \subset \End(E)$ is irreducible if there is no subspace $0 \subsetneq V \subsetneq E$ such that $S V \subset V$, \ie $sv \in V$ for all $s \in S$ and all $v \in V$.
		We say that a set $S \subset \End(E)$ is strongly irreducible or acts strongly irreducibly on $E$ if there is no finite union of subspaces $\{0\} \subsetneq \bigcup_{i \in I} V_i \subsetneq E$ such that $\bigcup_{i \in I, s \in S} s V_i \subset \bigcup_{i \in I, s \in S} s V_i$ and no finite union of subspaces $\{0\} \subsetneq \bigcup_{j \in J} W_j \subsetneq E^*$ such that $\bigcup_{j\in J, s \in S} W_j s \subset \bigcup_{s,j} s W_j$.
	\end{Def}
	
	\begin{Def}
		We say that a measure is strongly irreducible when its support is.
	\end{Def}

	Note that if there exists $0 \subsetneq V \subsetneq E$ such that $S V \subset V$, then for $W = V^\bot = \{w \in E^* \mid \forall v \in V,\,wv = 0\}$, we have $W S \subset W$ and $\{0\} \subsetneq W \subsetneq E^*$. 
	Therefore, the notion of irreducibility is invariant by transposition in the sense that a set $S \subset \End(E)$ is irreducible if and only if its transpose $^tS \subset \End(E^*)$ is.
	
	To get a symmetric notion of strong irreducibility, we need to have one condition on $S$ and one on $^tS$.
	Indeed, for all finite union of subspaces $\{0\} \subsetneq \bigcup_{i \in I} V_i \subsetneq E$, such that $\bigcup_{i \in I} V_i$ spans $E$, the set $S := \{ vf \mid v \in \bigcup V_i, f \in E^*\}$ of all rank-one matrices whose image is included in $\bigcup V_i$ stabilizes $\bigcup V_i$.
	On the other hand, for all $w \in E^* \setminus\{0\}$, there exists a vector $v \in  \bigcup_{i \in I} V_i$ such that $wv \neq 0$. Given such a $v$, we have $w S \supset \{w vf \mid f \in E^*\} = E^*$ so there is no $\{0\} \subsetneq A \subsetneq E^*$ such that $AS \subset A$.
	
	For $S \subset \GL(E)$, saying that $S$ is strongly irreducible is equivalent to saying that it is virtually irreducible in the sense that there exists a subgroup $S \subset G \subset \GL(E)$ that admits a finite index subgroup $H \triangleleft G$ that is not irreducible. 
	That way, for $0 \subsetneq V \subsetneq E$ a $H$-invariant subspace, the union $\bigcup_{g \in G} g V$, which is naturally $S$-invariant, is equal to $\bigcup_{g \in F} g V$ for $F$ a finite subset of $G$ such that $G = FH$.

	\begin{Def}[Proximality]\label{def:prox}
		Let $\Gamma \subset \End(E)$ be a semi-group. 
		We say that $\Gamma$ is proximal if one of the following equivalent conditions are satisfied:
		\begin{enumerate}
			\item There exists a rank $1$ projection in the closure of $\KK\Gamma$.
			\item There exists an element $g \in \Gamma \setminus\{0\}$ that is proximal in the sense that $\lambda_1(g) > \lambda_2(g)$.
		\end{enumerate}
	\end{Def}
	
	Indeed, given a rank $1$ projection and a sequence $g_n$ such that $[g_n] \to [\pi]$, by continuity, $\lambda_2(g_n) / \lambda_1(g_n) \to 0$.
	Conversely, given $g$ such that $\lambda_1(g) > \lambda_2(g)$, we have $[g^n] \to [\pi_g]$, where $\pi_g$ is the projection onto the eigenspace of $g$ associated to the eigenvalue of maximal modulus along the sum of all other characteristic spaces.
	
	\begin{Def}\label{def:proximal}
		Let $\nu$ be a probability measure over $\End(E)$.
		We say that $\Gamma$ is proximal if one of the following equivalent conditions are satisfied:
		\begin{enumerate}
			\item All measurable semi-group $\Gamma \subset \End(E)$ such that $\nu(\Gamma) = 1$ is proximal.
			\item The semi-group $\Gamma_\nu = \bigsqcup_{n \ge 0} \supp(\nu^{*n})$ is proximal.
			\item There exists an integer $n$ such that $\nu^{*n}\{g \mid \lambda_1(g) > \lambda_2(g)\} > 0$.
		\end{enumerate}
	\end{Def}

	One can check whether a measure is proximal or not by looking at the powers of its support. 
	More precisely, we can guarantee that a measure $\nu$ is proximal by finding a finite sequence in the support of $\nu$ whose product is proximal.
	Now to check that a measure $\nu$ is not proximal, we need to check that there is no proximal element in any power of the support of $\nu$, which is not a first order property of the support of $\nu$.
	The issue with Definition \ref{def:proximal} is that it does not give an algorithm to determine whether a finitely supported measure $\nu$ is proximal or not.
	
	Indeed, there is no universal lower bound on the smallest integer $n$ for which there exists a proximal element in $S^n$ for all sets $S \subset\GL(E)$ that generate a strongly irreducible and proximal semi-group.
	For example, we take $S = \{R_{\theta}, R'_{\rho}\}$ where $R_\theta$ is the rotation of angle $\theta$ in the basis $(e_1, e_2)$ and $R'_\rho$ is the rotation of angle $\rho$ in the basis $(2e_1, e_2/2)$, we pick $0 < \theta < \rho$ to be irrational multiples of $\pi$ (to guarantee strong irreducibility) and arbitrarily small. 
	Note that an element $g$ of $\GL(\RR^2)$ is proximal if and only if it has exactly two fixed points on $\proj(E)$.
	Let us view $\proj(\RR^2)$ as $\RR / \pi \ZZ$ endowed with the cyclic order.
	For all $x \in \proj(\RR^2)$, we have $R_\theta x = x + \theta$ and $x + \rho / 2 \le R'_{\rho} x \le x + 2\rho$.
	Moreover, the bounds of the inequality are sharp up to a correcting term of order $\rho^2$ so $h := R_\theta^{\lfloor(\pi - \rho / 2)/\theta \rfloor}R'_\rho$ is proximal for $\rho$ small enough (because $h[e_1] < [e_1] < [e_2] < h[e_2]$).
	On the other hand, given $n$ and $g \in S^{\cdot n}$, it holds for all $x$ that $x + n \min\{\eta, \rho / 2\} \le g(x) \le x + n \max\{\eta, 2\rho\}$, so we need $n \ge \pi / \max\{\eta, \rho \sqrt{2}\}$ for $g$ to be proximal.

	A celebrated of Goldshied and Margulis \cite{GM1987} gives the following algebraic characterization of the contraction property for groups of matrices over Archimedean fields.
	
	\begin{Prop}[Algebraic characterization of contraction \cite{GM1987}]\label{prop:zariski}
		Let $E$ be a Euclidean or Hermitian vector space.
		Let $\nu$ be a probability measure on $\GL(E)$.
		Let $\Gamma_\nu$ be the semi-group generated by the support of $\nu$.
		Let $G_\nu$ be the Zariski closure of $\Gamma_\nu$.
		Then $G_\nu$ is a group and the following propositions are equivalent:
		\begin{enumerate}
			\item The measure $\nu$ is strongly irreducible and proximal.
			\item The group $G_\nu$ is strongly irreducible and proximal.
			\item The connected component of the identity in $G_\nu$ is irreducible and proximal. 
		\end{enumerate}
	\end{Prop}
	
	If $E$ is ultrametric, the set $\{g \in \GL(E)\mid s_1(g) = s_d(g)\}$ is an open (and therefore Zariski dense) subgroup of $\GL(E)$ that is not proximal so Proposition \ref{prop:zariski} does not hold.
	This subject is discussed in more details in \cite{cone-limite}.
		
	A notable consequence of Proposition \ref{prop:zariski} is that the law of a random variable $\gamma \in \GL(E)$ is strongly irreducible and proximal if and only if the law of $\gamma^{-1}$ is.
	In the present paper, we drop the assumption that $\nu$ is supported on $\GL(E)$, which breaks this symmetry.
	
	In the beginning of section \ref{sec:kernel}, we prove that the following characterization of strong irreducibility and proximality (similar to the one used by by Guivarc'h and Raugi) still holds when the matrices are not invertible. 
	
	\begin{Prop}[Characterization of strong irreducibility and proximality via a contraction property for semi-groups of matrices]\label{prop:boundary}
		Let $\nu$ be a probability measure on $\End(E)$.
		Let $\Gamma_\nu$ be the semi-group generated by the support of $\nu$.
		The following propositions are equivalent:
		\begin{enumerate}
			\item The measure $\nu$ is strongly irreducible and proximal.
			\item The space $\partial'_1(\nu)$ of rank one matrices in the closure of $\KK \Gamma_\nu$ in $\End(E)$ is a product $\partial'_1(\nu) = \Lambda'_1(\nu) \Lambda'_1(^t\nu)$, where $\Lambda'_1(\nu) \subset E$ and $\Lambda'_1(^t\nu) \subset E^*$ are non-empty $\KK$-homogeneous (in the sense: stable by multiplication by elements of $\KK$) and are not included in a finite union of subspaces.
			\item For all finite families $(v_i)_{i \in I} \in E\setminus\{0\}$ and $(f_j)_{j \in J} \in E^*\setminus\{0\}$, there exists a rank one endomorphism $\pi \in \partial'_1(\nu)$ such that $f_j \pi v_i \neq 0$ for all $i \in I$ and all $j \in J$.
		\end{enumerate}
	\end{Prop}
	
	The homogeneity of $\Lambda'_1(\nu) := \bigcup_{\pi \in \partial'_1(\nu)} \pi(E) \subset E$ comes from the fact that the closure of $\KK\Gamma_\nu$ is homogeneous.
	The fact that it is non-empty comes from the fact that $\nu$ is proximal.
	The factorization $\partial'_1(\nu) = \Lambda'_1(\nu) \Lambda'_1(^t\nu)$ comes from the irreducibility of $\nu$.
	
	When the measure is not supported on $\GL(E)$, the knowledge of the support of $\nu$ is not sufficient to decide whether or not sequence $(\overline\gamma_n)$ is stationary to $0$.
	For example on $\End(\RR^2)$, take $\nu$ to be the law of $R \begin{pmatrix} 0 & 1 \\ 0 & 0 \end{pmatrix} R^{-1}$ with $R$ uniformly distributed along the Haar measure on $\mathrm{SO}(\RR^2) \simeq \RR/ 2\pi \ZZ$.
	Then $\nu$ is supported on a set of nilpotent matrices so all the semi-groups of full $\nu$-measure contain $0$. 
	However, for $(\gamma_n)\sim \nu^{*n}$, we get $\|\overline\gamma_n\| = \prod_{k = 1}^{n-1} |\sin(\theta_k)|$, where $\theta_k$ is the angle associated to the rotation $R_{k-1}^{-1} R_k$. 
	By invariance of the Haar measure, the sequence $(\theta_k)_{k \ge 1}$ is i.i.d. along the Haar measure so $\PP(|\sin \theta_k| = 0) = 0$.
	Therefore $\overline\gamma_n \neq 0$, almost surely and for all $n$.
	On the other hand if the distribution of the $R_k$'s has any atom, then $\overline\gamma_n$ is almost surely stationary to $0$ and, in this simple model, this is an equivalence.
	
	In higher dimension the $0-1$ law for the stationarity at $0$ yields a result on stationarity of the rank.
	In Section \ref{sec:kernel}, we prove the following result.
	
	\begin{Prop}[Eventual rank]\label{prop:rank}
		Let $\nu$ be a probability measure on $\End(E)$ and let $\gamma \sim \nu^{\otimes\NN}$.
		There exists an integer $0 \le r \le \dim(E)$ such that $\PP(\rank(\overline\gamma_n) < r) = 0$ for all $n$ and $n_r(\gamma) := \min\{n \in \NN \mid \rank(\overline\gamma_n) = r\}$ is almost surely finite and has a finite exponential moment.
		We write $r = \evrank(\nu)$ and call it the eventual rank of $\nu$.
	\end{Prop}
	
	Since $\rank(\overline\gamma_n)$ is non-increasing, it follows that $\rank(\overline\gamma_n) = \evrank(\nu)$ for all $n$ larger than $n_{\evrank(\nu)}$, which justifies the name eventual rank.
	
	The objects of interest in Theorem \ref{th:intro-contraction} are the columns $\overline\gamma_n v$. Knowing that $\overline\gamma_n \neq 0$, almost surely and for all $n$ does not mean that $\overline\gamma_n v \neq 0$ for all vector $v$.
	However, the following result tells us that $\PP(\overline\gamma_n v = 0)$ is bounded away from $1$ and equal to $0$ as soon as $v$ does not belong to a countable union of proper subspaces of $E$.
	
	\begin{Prop}[Essential kernel]\label{prop:essker}
		Let $\nu$ be an irreducible probability measure on $\End(E)$ and let $(\gamma_n)\sim \nu^{\otimes\NN}$.
		The set $\{v \in E \mid \exists n, \PP(\overline\gamma_n v = 0) > 0\}$ is a countable union of subspaces of $E$ that all have dimension at most $\dim(E) - \evrank(\nu)$.
		Moreover, there exists a constant $\alpha < 1$ such that $\PP(\overline\gamma_n v = 0) \le \alpha$ for all $v \in E \setminus\{0\}$ and for all $n \in \NN$.
	\end{Prop}

	\subsection{Alignment and pivotal extraction}\label{sec:intro-pivot}
	
	First let us introduce the notion of alignment for endomorphisms.
	
	\begin{Def}[Coarse alignment of matrices]\label{def:ali}
		Let $g \in E^* \cup \End(E)$ and let $h \in E \cup \End(E)$ be non-trivial. Let $0 \le \eps \le 1$, when the inequality
		\begin{equation}\label{def-ali}
			\frac{\|g h\|}{\|g\|\|h\|} > \eps
		\end{equation} 
		holds, we say that $g$ is aligned with $h$ and write $g \AA^\eps h$.
		By convention, we also write $g \AA^\eps h$ when $g = 0$ or $h = 0$.
	\end{Def}
	
	\begin{Def}[Contraction coefficient]
		Let $g \in \End(E) \setminus\{0\}$. 
		We write:
		\begin{equation}
			\sigma(g) := \frac{s_2(g)}{s_1(g)} = \frac{\|g\wedge g\|}{\|g\|^2}.
		\end{equation} 
	\end{Def}
	
	We remind that by the Cartan decomposition, the endomorphism $g$ is at distance $s_2(g)$ from the closest rank one endomorphism.
	Therefore, the quantity $\sigma(g)$ is the projective distance between $[g]$ and the closest rank one matrix. 
	
	\begin{Def}[Schottky measure]\label{def:schottky}
		let $\mu$ be a probability measure over $\End(E) \setminus\{0\}$ and let $\rho > 0$. 
		We say that $\mu$ is $\rho$-Schottky for $\AA^\eps$ if the inequalities $\mu\{\gamma\mid g \AA^\eps \gamma\} > 1 - \rho$ and $\mu\{\gamma\mid \gamma \AA^\eps g\} > 1 - \rho$ hold for all $g \in \End(E) \setminus\{0\}$.
	\end{Def}
	
	The easy part of our result is to show that a strongly irreducible and proximal measure is contracting in the following sense.
	
	\begin{Def}[Quantification of contraction]\label{def:quant-contraction}
		Let $\nu$ be a probability measure on $\End(E)$.
		Let $\alpha, \eps > 0$ and let $m \in \NN$.
		We say that $\nu$ is $(\alpha, \eps, m)$-contracting if there exists a probability measure $\mu$ over $\End(E)$, that is $1/6$-Schottky for $\AA^\eps$, such that $\sigma < \eps^6/48$ on the support of $\mu$ and such that $\alpha \mu < \nu^{*m}$.
	\end{Def}

	Note that saying that $\alpha \mu < \nu^{*m}$ implies that there exists a probability measure $\tilde\mu$ on $\End(E)^m$ such that $\alpha \tilde\mu < \nu^{\otimes m}$.	
	This notion of quantitative contraction is a weak-$*$ open condition in the following sense.

	\begin{Prop}\label{prop:weak-cont}
		Let $\alpha, \eps > 0$ and let $m \in \NN$. 
		The space of $(\alpha, \eps, m)$-contracting probability measures on $\End(E)$ is open in the weak-$*$ topology.
	\end{Prop}
	
	It will work in our advantage for applications as a statistical sample of arbitrarily large size, with arbitrarily small experimental uncertainty gives an arbitrarily good weak-$*$ approximation of any probability measure $\nu$.
	
	\begin{Lem}\label{lem:contraction}
		Let $\nu$ be a strongly irreducible and proximal probability measure on $\End(E)$. 
		There exist $\alpha, \eps > 0$ and $m \in \NN$ such that $\nu$ is $(\alpha, \eps, m)$-contracting.
	\end{Lem}
	
	Conversely, if a measure $\nu$ supported on $\GL(E)$ is $(\alpha, \eps, m)$-contracting for some $\alpha, \eps, m > 0$, then it is strongly irreducible and proximal.
	This fact is known and used in \cite{AMS-esmigroup} for example.
	However, if we assume $\nu$ to be $1/6$-Schottky for $\AA^\eps$ and supported over a finite set of rank one endomorphisms, then $\nu$ preserves a finite bundle of lines so it is not strongly irreducible but it is $(1, \eps, 1)$-contracting.
	For $\nu$ supported on $\End(E)$, the main difference between strong irreducibility and contraction is that for a strongly irreducible $\nu$ and for all finite family $v_1, \dots, v_N \in \essker(\nu) \setminus\{0\}$, the probability $\PP(\forall i, \overline\gamma_n v_i = 0)$ is bounded below by a positive constant while if we only assume contraction, it may not hold when $N > 6$.
	That being said, we introduce the notion of quantitative contraction because it allows us to link the speed of contraction given by Theorem \ref{th:intro-contraction} and its corollaries with weak-$*$ continuous quantities. The side effect is that it relaxes the strong irreducibility condition for products of random singular matrices.

	Before giving the statement of our main technical Theorem, we need to introduce the following notation, used through the rest of the present paper.
	
	\begin{Def}[Power notation for extractions]\label{def:power-extract}
		Let $\Gamma$ be a semi-group. 
		Let $(\gamma_n)_{n\ge 0} \in \Gamma^\NN$ and let $(p_k)_{k \ge 0} \in \NN^\NN$.
		For all $k \ge 0$, we write:
		\begin{align*}
			\widetilde{\gamma}^p_k &:= \widetilde\gamma_{\overline{p}_k, \overline{p}_{k+1}} = (\gamma_{\overline{p}_k}, \gamma_{\overline{p}_k+1}, \dots, \gamma_{\overline{p}_{k+1}-1}) \\
			\gamma^p_k &:= \gamma_{\overline{p}_k, \overline{p}_{k+1}} = \gamma_{\overline{p}_k} \cdots \gamma_{\overline{p}_{k+1}-1}.
		\end{align*}
	\end{Def}	
	
	It is important to take some time to understand these notations as they will be useful in the proofs of the results. 
	If we look at the walk, we get $\overline\gamma^p_n = \overline\gamma_{\overline{p}_n}$ for all $n$ so when $\Gamma$ is a group, the data of $(\overline\gamma^p_n)$ is equivalent to the data of a sub-sequence of $(\overline\gamma_n)_{n \ge 0}$. 
	On the other hand, the data of the sequence $(\widetilde\gamma^w_k)_{k \ge 0}$ is equivalent to the joint data of the sequences $(\gamma_n)_{n \ge 0}$ and $(p_k)_{k \ge 0}$.
	
	The $\widetilde{\cdot}$ on $\widetilde{\gamma}^p_k$ is used to specify that we interpret it as a word on the alphabet $\End(E)$. 
	Through the present paper, removing the $\widetilde{\cdot}$ on a word amounts to taking the ordered product of its letters, we denote by $\Pi$ the product map $\Pi: \bigcup_{n = 0}^{+\infty} \Gamma^n \to \Gamma\, ; \,(g_0, \cdots, g_{l-1}) \mapsto g_0 \cdots g_{l-1}$, which is a projection. 
	
	The purpose of the following Theorem is to find an extracting sequence $(p_n)$ such that $-\log\sigma(\overline\gamma^p_n)$ escapes at least linearly. 
	To get the exponential right member in \eqref{contraction-intro}, we need the $p_n$'s to have a finite exponential moment.
	
	\begin{Th}[Pivoting Technique]\label{th:pivot}
		Let $\nu$ be a probability measure on $\End(E)$ and let $(\gamma_n)\sim \nu^{\otimes\NN}$.
		Assume that $\nu$ is $(\alpha, \eps, m)$-contracting for some constants $\alpha, \eps > 0$ and $m \in \NN$.
		Let $\tilde{\mu} \le \nu^{\otimes m} / \alpha$ be a $1/6$-Schottky probability measure for $\AA^\eps$ such that $\sigma < \eps^6/48$ on the support of $\mu$.
		Then there exists\footnote{Unless specified otherwise, when we say that there exists a random variable, we assume it to be defined on a large enough standard probability space so that we can draw a uniform random variable in $[0,1]$ independently of all previously defined random variables.} a sequence $(p_n)_{n \ge 0}$ of random independent integers whose laws have finite exponential moments and depend only on $\alpha$ and $m$ and such that the following points hold :
		\begin{enumerate}
			\item For all $k$, $p_{2k + 1} = m$ and the conditional distribution of $\widetilde{\gamma}^p_{2k+1}$ with respect to $(\widetilde\gamma^p_{k'})_{k' \neq 2k + 1}$ is\footnote{Since the conditional distribution is only defined up to equality on a set of probability $1$, we could say say that is is almost surely bounded. To avoid repeating "almost surely" everywhere, an assertion involving random variables to be understood as holding on an event of probability $1$.} bounded above by $\frac{3}{2}\tilde\mu$ and therefore projects onto a $1/4$-Schottky measure for $\AA^\eps$.\label{pivot:schottky}
			\item The sequence $(\widetilde\gamma^p_{2k + 2})_{k \ge 0}$ is i.i.d. and independent of $\widetilde{\gamma}^p_{0}$ and for all $k \ge 1$, the conditional distribution of $\widetilde{\gamma}^p_{2k+1}$ with respect to $(\widetilde\gamma^p_{k'})_{k' \neq 2k + 1}$ is given by a measurable function $\eta(\gamma^p_{2k},\widetilde\gamma^p_{2k + 2})$, where $\eta$ is a measurable map that takes values in the space of probability measures and does not depend on $k$.\label{pivot:renewal}
			\item For all $0 < i < j < k$ and for all $f \AA^{\eps / 2} \gamma^p_{i}$ and $\gamma^p_{k - 1} \AA^{\eps/2} v$,\label{pivot:herali}
			\begin{equation*}
				f \gamma^p_{i} \cdots \gamma^p_{j-1} \AA^{\eps / 2} \gamma^p_{j} \cdots \gamma^p_{k-1} v.
 			\end{equation*}
 			\item \label{pivot:indep-tail} For all measurable set $S$ such that $\tilde{\mu}(S^m) = 1$ and $\nu(S) < 1$, the data of $(\gamma_n)_{n,\gamma_n \notin S}$ and $(p_n)_{n \ge 0}$ are independent relatively to the data of $\{n \mid \gamma_n \in S\}$. 
 			It means that there exists a random i.i.d. sequence $(g_n) \sim \nu_{S^c}^{\otimes\NN}$, where $\nu_{S^c}$ denotes the normalized restriction of $\nu$ to $\End(E) \setminus S$, that is independent of the data of $(p_n)$ and such that for all $n$, either $\gamma_n = g_n$ or $\gamma_n \in S$.
 			Moreover, for all event $A \subset \End(E)\setminus S$,
 			\begin{equation}\label{dom-tail}
 				\forall n \in \NN, \; \PP(\gamma_n \in A \mid (p_k)_{k \ge 0}) \le \nu(A) / (1- \alpha).
 			\end{equation}
		\end{enumerate}
	\end{Th}
		
	Point \eqref{pivot:renewal} implies that the times $(\overline{p}_{2k + 2})_{k \ge 0}$ act like renewal times for the random walk, in the sense that the distribution of $(\gamma_j)_{j \ge \overline{p}_{2k + 2}}$ does not depend on $k$ and the data of the \textit{future} $(\gamma_j)_{j \ge \overline{p}_{2k + 2}}$ is independent of the \textit{past} $(\gamma_j)_{j < \overline{p}_{2k + 1}}$ (and $\overline{p}_{2k + 1} = \overline{p}_{2k + 2}- m$).
	
	Let us quickly illustrate how we use Theorem \ref{th:pivot} by giving a sketch of the proof of Theorem \ref{th:intro-lln}.	
	Point \eqref{pivot:herali}, becomes really useful when combined with point \eqref{pivot:schottky}.
	Indeed, given a non-random $f \in E^*\setminus\{0\}$, by point \eqref{pivot:schottky}, $\PP(f_0\overline\gamma^p_i \AA^\eps \gamma^p_i) > 3/4$ for all $i \in 2\NN + 1$.
	Therefore, the first integer $i$ such that $f_0\overline\gamma^p_l \AA^{\eps / 2} \gamma^p_l$ has a finite exponential moment and since the $p_k$'s, are independent and have bounded exponential moment, the sum $\overline{p}_l$ also has a finite exponential moment.
	On the other hand, for all non-random $v \in E \setminus\{0\}$ and all $n$, we find an odd integer $r$ such that $\gamma^p_r \AA^{\eps / 2} \gamma{\overline{p}_{r+1}} \cdots \gamma_{n-1} v$ and get a bound of the exponential moment of $n - \overline{p}_{r+1}$.
	
	Then, using point \eqref{pivot:herali} and the basic results of Section \ref{kak}, we can bound the log-ratio $-\log \frac{|f \overline\gamma_n v|}{\|\overline\gamma_n\|}$ by the sum of the $N(\gamma_k)$'s for all $0 \le k < \overline{p}_l$ and all $\overline{p}_{r + 1} \le k < n$.
	Assuming $N $ to be bounded by a constant $B$ on the set $S$ of point \eqref{pivot:indep-tail}, this sum in now bounded by the sum of the $N(g_k)$ for $k$ in the same intervals, to which we add the random variable $B \overline{p}_l + n - \overline{p}_{r + 1}$ (which has a bounded exponential moment).
	With the assumptions of Theorem \ref{th:intro-lln}, we get that $\sum N(g_k)$ is a sum of exponentially many integrable random variables with the same law, this sum is bounded in law by an integrable random variable. We conclude via Borel-Cantelli's Lemma.
	
	Note that from \eqref{pivot:indep-tail}, without \eqref{dom-tail}, we get that $\PP(\gamma_n \in A \mid (p_k)_{k \ge 0}) \le \PP(g_k \in A) = \nu(A) / (1- \nu(S))$.
	However, we may assume without restriction that $\alpha$ is arbitrarily small but we may not assume that $\nu(S)$ is far from $1$ without putting restrictions on $\nu$.

	Since the space of $1/6$-Schottky measures is weak-$*$ open, for $S$ of large enough $\nu$-measure, the normalized restriction of $\tilde\mu$ to $S^m$ stays Schottky and bounded by $\nu^{\otimes m} / \alpha$.
	That way, we may assume $S$ to be compact; for $\nu$ supported on $\GL(E)$, it means that the map $N : g \mapsto \log\|g\|\|g^-1\|$ is bounded on $S$. 
	This trick is crucial in the proof of Theorem \ref{th:intro-lln}.

	We moreover have an explicit construction for the sequence $(p_n)_n$.
	Let us describe the law of $(p_n)_n$ as we construct it in Theorem \ref{th:ex-piv} (from which Theorem \ref{th:pivot} follows).
	Note that by point \eqref{pivot:schottky}, the $p_{2k+1}$ are non-random and equal to $m$ and by point \eqref{pivot:renewal}, the $(p_n)$'s are independent and the $p_{2k + 2}$ all have the same law so the joint law of $(p_n)$ is determined by law of $p_0$ and the law of $p_2$.
	Let us now describe the law of $p_0$ and $p_2$.
	
	For all $0 < \rho < 1$, we write $\mathcal{G}_\rho = \sum \delta_k \rho^k(1-\rho)$ for the geometric distribution of scale factor $\rho$ (which has expectation $\rho/(1-\rho)$) and we write $\mathcal{B}_\rho = \rho \delta_1 + (1-\rho)\delta_0$ for the Bernoulli distribution of parameter $\rho$.
	Let $(w_k)_{k \ge 0}\sim \mathcal{G}_{1-\alpha}^{\otimes\NN}$, let $(v_k)_{k \ge 0} \sim \mathcal{G}_{1/3}^{\otimes \NN}$, let $(b_k) \sim \mathcal{G}_{1/4}^{\otimes\NN}$ and let $(s_k)_{k \ge 0} \sim \mathcal{B}_{2/3}^{\otimes \NN}$ be independent.
	Let $m_0 = 0$ and for all $k \ge 0$, let $m_{k + 1} = (m_k + s_k + (s_k - 1) (b_k + 1))_+$.
	For all $k$, let $l_k = \max\{j\mid m_j = k\}$ be the time of last return in $k$.
	Let $j_0 = l_0$ and $j_2 = l_1 - l_0 - 1$.
	For $i = 0,2$, let $q_i = \sum_{k = 0}^{j_i} (v_k + 2)$ and let $p_i = - m + m \sum_{k = 0}^{q_i - 1} (w_k + 1)$.
	
	This gives explicit constructions for the laws of $p_0$ and $p_2$. 
	We can compute the expectation of $p_2$ as follows.
	By the law of large numbers for absorbed random walks, the sequence $(m_k)$ has drift $\lim m_k / k = \EE(s_k + (s_k - 1) b_k)_+ = 2/3 - (1/3)(4/3) = 2/9$ so $\EE(l_{k + 1}- l_k) = 9/2$.
	Therefore $\EE(q_2) = (9/2)(5/2)$ and $\EE(p_2) = \frac{45m}{4\alpha}  - m$.
	More details are given in Section \ref{sec:pivot}.

	\subsection{Contraction results without moment assumptions}\label{intro:non-invertible}

	Through the present paper, $E$ denotes a Euclidean, Hermitian or ultrametric vector space of dimension $d \ge 2$.
	We fix a probability measure $\nu$ that we assume to be strongly irreducible and proximal and to have positive eventual rank but we may have $\nu(\GL(E)) < 1$.
	Since all the results of the present section are derived from Theorem \ref{th:pivot} applied to a Schottky measure given by Lemma \ref{lem:contraction}, we may replace the assumption that $\nu$ is strongly irreducible and proximal by the assumption that it is $(\alpha, \eps, m)$-contracting for some $\alpha, \eps, m$, in the sense of Definition \ref{def:quant-contraction}.
	
	\begin{Th}[Contraction]\label{th:escspeed}
		Let $\nu$ be a proximal and strongly irreducible probability distribution on $\End(E)$ with positive eventual rank. 
		Let $(\gamma_n) \sim \nu^{\otimes\NN}$. 
		Write $\overline{\gamma}_n:= \gamma_0 \cdots \gamma_{n-1}$ for all $n$.
		Then there exists a positive constant $\lambda_{1,2}(\nu) \in (0,+\infty]$ such that $\frac{- \log\sigma(\overline{\gamma}_n)}{n} \to \lambda_{1,2}(\nu)$ almost surely. 
		Moreover, the following large deviations inequalities hold:
		\begin{equation}\label{ldsqz}
			\forall \alpha < \lambda_{1,2}(\nu),\;\exists C, \beta > 0 ,\; \forall n\in\NN,\; \PP(- \log\sigma(\overline{\gamma}_n) \le \alpha n) \le C e^{-\beta n}.
		\end{equation}
	\end{Th}
	
	\begin{Rem}[On the order of the product]
		Given a measure $\nu$, we write $^t\nu$ for the push-forward of $\nu$ by the transpose map that we call transpose of $\nu$.
		An endomorphism shares its singular values with its transpose. 
		Moreover, for a product $g = \gamma_0 \cdots \gamma_{n-1}$, we have $^tg = ^t\gamma_{n-1} \cdots ^t\gamma_0$.
		We conclude that the convergence in probability $\frac{- \log\sigma(\overline{\gamma}_n)}{n} \to \lambda_{1,2}(\nu)$ also holds when we replace the left-to-right product $\overline{\gamma}_n = \gamma_0\cdots \gamma_{n-1}$ with the right-to-left product $\gamma_{n-1} \cdots \gamma_0$ so $\lambda_{1,2}(\nu) = \lambda_{1,2}(^t\nu)$, just like in the $\ELL^1$ case. 
	\end{Rem}
	
	\begin{Rem}[Computing the value of $\lambda_{1,2}(\nu)$ with the pivotal times decomposition]
		In the proof of Theorem \ref{th:escspeed}, we compute the value of $\lambda_{1,2}(\nu)$ with the following formula, for $(p_n)$ and $m$ as in Theorem \ref{th:pivot}:
		\begin{equation}\label{lyap-gap-pivot}
			\lambda_{1,2}(\nu) = \lim_{n \to + \infty} \frac{\lambda_{1,2}(\gamma^p_{2} \cdots \gamma^p_{2n})}{n \EE(p_2 + m)}
		\end{equation}
		In particular, when $\nu$ is $(\alpha, \eps, m)$-contracting for some constants $\alpha, \eps, m > 0$, it holds that $\lambda_{1,2} \ge \frac{4\alpha}{45 m}|\log(\eps^6 / 48)|$.
		This gives a lower bound on $\lambda_{1,2}(\nu)$ that is stable under small perturbations in the weak-$*$ topology.
	\end{Rem}
	
	The quantitative estimates of \eqref{ldsqz} in Theorem \ref{th:escspeed} are new even in the setting of~\cite{GR86}, where the almost sure convergence $\frac{- \log\sigma(\overline{\gamma}_n)}{n} \to \lambda_{1,2}(\nu)$ was proven under a first moment assumption.
	In Theorem \ref{th:escspeed}, we not only prove the almost sure convergence without assuming $\log\|\gamma_0\|^{-1}$ to be integrable but we also state exponential large deviations inequalities below the escape speed.
	
	We remind that $\proj(E)$ denotes projective space associated to $E$ and $[\cdot]: E\setminus\{0\} \to \proj(E)$ denotes the projection map. 
	Let us endow $\proj(E)$ with the metric:
	\begin{equation}\label{distance-def}
		\dist:([x],[y]) \mapsto \frac{\|x\wedge y\|}{\|x\|\|y\|} = \min_{c \in\KK} \frac{\|x - cy\|}{\|x\|}.
	\end{equation}

	\begin{Th}[Quantitative convergence of the image]\label{th:cvspeed}
		Let $\nu$ be a strongly irreducible and proximal probability distribution on $\End(E)$ that has positive eventual rank. 
		Let $(\gamma_n) \sim \nu^{\otimes\NN}$. 
		There exists a random line $\ell^\infty(\gamma) \in \proj(E)$ such that for all $\alpha < \lambda_{1,2}(\nu)$, there exist constants $C,\beta > 0$, such that:
		\begin{equation}\label{erc}
			\forall v \in E \setminus \{0\},\,
			\forall n\in\NN,\,
			\PP\left(\dist([\overline{\gamma}_n v], \ell^\infty) \ge e^{-\alpha n}\mid \overline{\gamma}_n v \neq 0\right) \le C e^{-\beta n}.
		\end{equation}
	\end{Th}

	Note that if we drop the assumption that $\nu$ has positive eventual rank, Inequality \eqref{erc} is non-trivial as long as $\PP(\overline{\gamma}_n = 0) < 1$ for all $n$ and in this case it is an interesting question to ask whether or not it remains true or not with the same algebraic assumptions on $\nu$. 
	
	When $\sup_n \PP(\overline{\gamma}_n v = 0) = 0$ (which holds for all $v$ outside of a countable union of proper subspaces by Proposition \ref{prop:essker}) \eqref{erc} implies that $[\overline{\gamma}_n v] \to \ell^\infty(\gamma)$ almost surely.
	Otherwise the sequence $([\overline{\gamma}_n v])$ has two limit points in $\proj(E)\sqcup \{[0]\}$, $\ell^\infty(\gamma)$ and $[0]$.
	
	For $\nu$ supported on $\GL(E)$, the existence of the limit line $\ell^\infty(\gamma)$, distributed along the $\nu$-stationary measure, and the almost sure convergence $[\overline\gamma_n v]\to \ell^\infty(\gamma)$ were proven by Guivarc'h and Raugi in \cite{GR1985}, via a martingale argument on $\mathrm{Prob}(\proj(E))$.
	In fact, in the invertible case, we know that there exists a random limit line $\ell^\infty$ such that $[\overline\gamma_n v]\to \ell^\infty$ for all $v$ if and only if the measure $\nu$ is strongly irreducible and proximal.
	Since Theorem \ref{th:cvspeed} holds for any contracting measure, we deduce that contracting measure are strongly irreducible and proximal and therefore Lemma \ref{lem:contraction} is an equivalence for measures supported on $\GL(E)$.
	
	Even though it is interesting to be able to generalize this result for measures that give non-zero weight to singular matrices, the real strength of Theorem \ref{th:cvspeed} lies in the exponential convergence in probability given by \eqref{erc}.
	Indeed, the same formula was obtained by Guivarc'h and Raugi but only under an exponential moment assumption on $N_*\nu$, which has not been relaxed since then.

	We define $\xi_\nu$ to be the distribution of $\ell^\infty(\gamma)$ for $(\gamma_n)_n \sim \nu^{\otimes\NN}$. 
	Then $\xi_\nu$ is the only $\nu$-stationary measure on $\proj(E)$ in the sense that $\nu * \xi_\nu = \xi_\nu$.
	Indeed, if we write $\ell_1^\infty(\gamma)$ for the almost sure limit of $[\gamma_1\cdots \gamma_n v]$ for a given $v$ outside of the essential kernel of $\nu$, then $\ell_1^\infty(\gamma)$ has law $\xi_\nu$ and is independent of $\gamma_0$.
	Moreover, $\gamma_0 \ell_1^\infty(\gamma) = \ell^\infty(\gamma)$ with probability $1$, which proves that $\nu * \xi_\nu = \xi_\nu$.
	
	In Theorem \ref{th:limit-uv}, we prove that moreover, the limit line is almost surely given by a partially defined function $\ell^\infty: \End(E)^\NN \to \proj(E)$ that is shift-equivariant in the sense of Definition \ref{def:l-infty} and that does not depend on the choice of $\nu$.
	
	A direct consequence of Theorem \ref{th:cvspeed} is the following result, which we can see as an exponential mixing result for the Wasserstein metric on the projective space.
	We remind that given two probability measures $\eta$ and $\mu$ on a metric space $(X, \dist)$, the Wasserstein distance $W(\eta, \mu)$ is the infimum of $\EE(\dist(x,y))$ for all coupling $x \sim \eta$, $y\sim \mu$, it is also the maximum of $\int \phi d\eta - \int \phi d\mu$, for all $1$-Lipschitz function $\phi : X \to \RR$.
	
	\begin{Cor}[Exponential mixing]\label{cor:ergodic}
		Let $\nu$ be a strongly irreducible and proximal probability distribution on $\End(E)$. 
		There exists a unique $\nu$-stationary probability distribution $\xi_\nu$ on $\proj(E)$. Moreover, there exist constants $C, \beta$ such that for all probability distribution $\xi$ on $\proj(E)\setminus\essker(\nu)$ and for all Lipschitz function $\phi: \proj(E) \to \RR$ with Lipschitz constant $\lambda(\phi)$,
		\begin{equation}
			\forall n\in\NN,\; \left|\int_{\proj(E)}\phi\mathrm{d}\xi_\nu - \int_{\proj(E)}\phi\mathrm{d}\nu^{*n}*\xi\right|\le \lambda(\phi) C e^{-\beta n}.
		\end{equation}
	\end{Cor}
	
	\begin{proof}
		Let $0 < \alpha < \lambda_{1,2}(\nu)$ and let $C, \beta >0$ be as in in Theorem \ref{th:cvspeed}.
		Let $[x], ({\gamma}_n)_n \sim \xi \otimes\nu^{\otimes\NN}$.
		Since $\PP(x \in \essker(\nu)) = 0$, $\overline{\gamma}_n x \neq 0$ almost surely and for all $n$. 
		Let $\ell^\infty$ be as in Theorem \ref{th:cvspeed}.
		Then $\EE(f(\ell^\infty)) = \int f\mathrm{d}\xi_\nu$ for all $f$.
		Moreover, $\PP\left(\dist([\overline{\gamma}_n x], \ell^\infty) \le e^{-\alpha n}\right) \le Ce^{-\beta n}$ for all $n$. 
		Moreover the distance map takes values in $[0,1]$ so $\EE\left(\dist([\overline{\gamma}_n x], \ell^\infty)\right) \le C e^{-\beta n} + e^{-\alpha n} \le C'e^{-\beta' n}$ for $C' = C+1$ and $\beta' = \min\{\beta, \alpha\}$.
		Then, by convexity of the absolute value, $|\EE\left(\phi[\overline{\gamma}_n x]\right) - \EE\left(\phi(\ell^\infty)\right)| \le \EE\left(|\phi[\overline{\gamma}_n x] - \phi(\ell^\infty)|\right) \le \EE\left(\lambda(\phi)\dist([\overline{\gamma}_n x], \ell^\infty)\right) \le \lambda(\phi) C'e^{-\beta' n}$.
	\end{proof}
	
	Note that saying that $\xi$ is supported on $\proj(E)\setminus\essker(\nu)$ is not very restrictive because any measure that gives measure $0$ to all hyperplanes would satisfy that condition. 
	However, $\xi_\nu$ itself may give positive measure to some hyperplanes. 
	For example if $\nu$ is the barycentre of the Haar measure on the group of isometries and a Dirac mass $\delta_\pi$ at a projection endomorphism $\pi = \pi^2$, then $\nu$ is strongly irreducible and proximal with positive rank and $\xi_\nu$ is the average of the isometry-stationary measure and of the Dirac mass on the image of $\pi$. In particular $\xi^\infty_\nu$ gives positive measure to any hyperplane that contains the image of $\pi$. 
	
	We can also interpret Theorem \ref{th:intro-contraction} as a result on exponential decay of correlation.
	
	\begin{Cor}[Exponential decay of correlation] \label{cor:decay}
		Let $\nu$ be a strongly irreducible and proximal probability distribution on $\End(E)$ and let $(\gamma_n) \sim \nu^{\otimes\NN}$.
		For all $\alpha < \lambda_{1,2}(\nu)$, there exist constants $C,\beta > 0$, such that for all $n$, there exist two independent random variables $w \in \proj(E^*)$ and $u \in \proj(E)$ such that for all $f \in E^*\setminus\{0\}$ and all $v \in E \setminus\{0\}$,
		\begin{equation*}
			\PP(\dist([f \overline\gamma_n], w) + \dist([\overline\gamma_n v], u) \ge 2e^{-\alpha n}) \le 2C e^{-\beta (n-1) / 2}.
		\end{equation*}
	\end{Cor}
	
	\begin{proof}
		Let $C, \beta > 0$ be as in Theorem \ref{th:cvspeed} applied to $\nu$ and $^t\nu$ \ie such that \eqref{erc} holds for $\gamma \sim \nu^{\otimes\NN}$ and for $\gamma \sim ^t\nu^{\otimes\NN}$.
		Let $(\gamma'_k)_k\sim \nu^{\otimes\ZZ}$ be independent of the joint data of $(\gamma_k)$.
		Let $n$ be fixed.
		For all $j \le n/2 - 1$, set $g_j := \gamma_j$ and for all $j > n/ 2 - 1$, set $g_j := \gamma'_j$.
		We denote by $u$ the almost sure limit of $([\overline{g}_k v])$.
		That way $u$ is independent of the data of $(\gamma_k)_{n / 2 \le k < n}$ and $(\gamma'_k)_{k < n/2}$.
		Moreover, by Theorem \ref{th:cvspeed},
		\begin{equation}
			\PP(\dist([\overline\gamma_n v], u) \ge e^{-\alpha n}) \le C e^{-\beta \lfloor n / 2 \rfloor}.
		\end{equation}
		Then we define a sequence $(h_j)_{j \ge 0} \sim ^t\nu^{\otimes\NN}$ where $h_j := ^t\gamma_{n - j - 1}$ for all $j < n/2 $ and $h_j := ^t\gamma'_{n - j - 1}$ for all $j \ge n/2$.
		Let $w$ be the almost sure limit of $([\overline{h} f])$ 
		That way, the joint data of $(h_j)_j$ is independent of the joint data of $(g_j)_j$ and therefore $w$ and $u$ are independent.
		Moreover, by Theorem \ref{th:cvspeed} again
		\begin{equation}
			\PP(\dist([f \overline\gamma_n], w) \ge e^{-\alpha n}) \le C e^{- \beta n / 2}.
		\end{equation}
		We conclude by triangular inequality.
	\end{proof}
	
	Let $g$ be a square matrix such that $\rho_1(g) > \rho_2(g)$. 
	Then the top eigenvalue of $g$ is simple and belongs to $\KK$. 
	We write $\eigen^+(g)$ for the associated eigenspace, which is a line and therefore an element of $\proj(E)$.
	For all non-nilpotent matrix $g$, we write $\lambda_{1,2}(g) = \lim_n \frac{-\log\sigma(g^n)}{n} = \log\frac{\rho_1(g)}{\rho_2(g)} \in [0, + \infty]$

	\begin{Th}[Contraction properties for the dominant eigenspace]\label{th:eigenspace}
		Let $E$ be a Euclidean space and let $\nu$ be a proximal and strongly irreducible probability distribution on $\End(E)$. 
		Let $(\gamma_n)\sim\nu^{\otimes\NN}$, let $\lambda_{1,2}(\nu)$ be as in Theorem \ref{th:escspeed} and let $\ell^\infty$ be as in Theorem \ref{th:cvspeed}.
		Then for all $\alpha < \lambda_{1,2}(\nu)$, there exist constants $C, \beta >0$ such that for all $n > 0$,
		\begin{equation}\label{eigenspace}
			\PP\left(\lambda_{1,2} (\overline{\gamma}_n) \le\alpha n) \cup \dist(\eigen^+(\overline{\gamma}_n),\ell^\infty) \ge e^{-\alpha n}\right) \le C e^{-\beta n}.
		\end{equation}
	\end{Th}
	
	Theorem \ref{th:eigenspace} implies that $\eigen^+(\overline{\gamma}_n) \to \ell^\infty$ almost surely and that $\liminf_n \lambda_{1,2}(\overline\gamma_n) \ge \lambda_{1,2}(\nu)$ almost surely.
	In dimension $2$, $\lambda_{1,2} \le -\log\sigma$ so $\frac{1}{n}\log\frac{\rho_1}{\rho_2}(\overline{\gamma}_n) \to \lambda_{1,2}(\nu)$ almost surely by Theorem \ref{th:eigenspace} combined with Theorem \ref{th:escspeed}.
	In dimension $d \ge 3$, the question whether $\lambda_{1,2}(\nu)$ is actually the limit of $\lambda_{1,2}(\overline\gamma_n)$ remains open when $\lambda_{1,2}(\nu)$ is finite and $\nu$ is not in $\ELL^1$ (see corollary \ref{cor:gen-eigen} for the $\ELL^1$ case).

	\subsection{Convergence in the space of flags}\label{intro-flag}
	
	In \cite{GR86} \cite{GR89}, Guivarc'h and Raugi derive contraction results on the flag space for totally strongly irreducible measures (\ie the push-forward $\bigwedge^k_*\nu$ is strongly irreducible for all $k < d$).
	In the present paragraph, we point out that the same strategy works without moment conditions.
	Note that the Zariski density of $\Gamma_\nu$ guarantees total irreducibility (in both the Archimedean and ultrametric cases).
	
	We write $\mathrm{Gr}(E)$ for the set of subspaces of $E$.
	For all $k$, we write $\mathrm{Gr}_k(E)$ for the set of subspaces of $E$ that have dimension $k$.
	Given $\Theta \subset \{1, \dots, d-1\}$, we denote by $\mathrm{Fl}_{\Theta}(E)$ the set of partial flags of type $\Theta$.
	This is the space of totally ordered (for the inclusion order relation) collections of subspaces of $E$ whose set of dimensions (\ie its image by $\dim: \mathrm{Gr}(E) \to \{1, \dots, d\}$) is $\Theta$.
	Given $F$ a partial flag of type $\Theta$ and $k \in\Theta$, we write $F_k \in\mathrm{Gr}_k(E)$ for the (unique by definition) element of $F$ that has dimension $k$.
	Given two flags $F$ and $F'$ of type $\Theta$, we define:
	\begin{equation}
		\dist(F, F') = \max_{k \in \Theta}\dist(\textstyle\bigwedge^k F_k, \textstyle\bigwedge^k F'_k),
	\end{equation}
	where $\textstyle\bigwedge^k F_k$ is seen as an element of $\proj(\bigwedge^k E)$ and $\dist$ is defined by \eqref{distance-def}.
	For all $g \in \End(E)$, we write $\mathcal{F}^+(g)$ for the collection for $\lambda \ge 0$ of the sums of the characteristic spaces of $g$ associated to an eigenvalue of modulus at least $\lambda$.
	Write $\Theta(g) := \{k\mid  \rho_k > \rho_{k+1}\}$.
	
	The following straightforward corollary of Theorems \ref{th:cvspeed} and \ref{th:eigenspace} can be interpreted as an extension of Oseledets multiplicative ergodic theorem in the i.i.d. case without the $\ELL^1$ moment condition.
	
	\begin{Cor}[Contraction property in the flag space]
		Let $\nu$ be a probability measure on $\End(E)$ and let $(\gamma_n)_{n \ge 0} \sim \nu^{\otimes\NN}$.
		Let $q \le \min\{\evrank(\nu), d-1\}$ be such that $\bigwedge^k_*\nu$ is strongly irreducible for all $k \le q$.
		Denote by $\Theta(\nu)$ the set of indices $1 \le k \le q$ such that $\bigwedge^k_*\nu$ is proximal.
		Then there exists a measurable map $\mathcal{F}^\infty$ that is defined $\nu^{\otimes\NN}$ almost everywhere and such that the following points hold:
		\begin{enumerate}
			\item For all $1 \le i < j \le q+1$, there exists a non-random constant $\lambda_{i,j}(\nu)$ such that $\lim_{n} \frac1n \log\frac{s_i}{s_j}(\overline{\gamma}_n) = \lambda_{i,j}(\nu)$ almost surely and $\lambda_{i,j}(\nu) > 0$ when $i \in\Theta(\nu)$.\label{flag:spectrum}
			\item For all $\alpha < \min_{i \in \Theta(\nu)} \lambda_{i, i+1}(\nu)$, there exists constants $C, \beta > 0$ such that:
			\begin{equation}\label{eigen-flag}
				\PP(\dist(\mathcal{F}^+(\overline{\gamma}_n), \mathcal{F}^\infty(\gamma)) \ge e^{-\alpha}) \le C e^{-\beta n}
			\end{equation}
			and for all non random flag $F \in \mathrm{Fl}_{\Theta(\nu)}(E)$:
			\begin{equation}\label{limit-flag}
				\PP(\dist(\overline{\gamma}_n F, \mathcal{F}^\infty(\gamma)) \ge e^{-\alpha}\mid \overline{\gamma}_n F \in \mathrm{Fl}_{\Theta(\nu)}) \le C e^{-\beta n}.
			\end{equation}
		\end{enumerate} 
	\end{Cor}

	\begin{proof}
		We apply Theorems \ref{th:escspeed}, \ref{th:cvspeed} and \ref{th:eigenspace} to the random sequence $(\bigwedge^k \overline{\gamma}_n)_{n \ge 0}$ for all $k \in \Theta(\nu)$.
		
		For all $i \in \Theta(\nu)$, we set $\lambda_{i, i+1}(\nu) := \lambda_{1,2}(\bigwedge^i(\nu))$.
		By property of the singular values, $\frac{s_{1}}{s_2}(\bigwedge^i g) = \frac{s_i}{s_{i+1}}(g)$ for all endomorphism $g$ of rank at least $i$ so $\lim_{n} \frac1n \log\frac{s_i}{s_{i+1}}(\overline{\gamma}_n) = \lambda_{i,i+1}(\nu)$.
		For all $i \in\{1, \dots, q\}\setminus\Theta(\nu)$, we set $\lambda_{i, i+1}(\nu) = 0$
		Since $\bigwedge^k_*\nu$ is strongly irreducible but not proximal, $s_1 / s_2$ is bounded away from $0$ on $\Gamma_{\bigwedge^k_* \nu}$ (\textit{c.f.} \cite[Proposition (1.6)]{GR89} or Lemma \ref{lem:exiprox} in the present article for the non-invertible case) therefore $\frac1n \log\frac{s_i}{s_{i+1}}(\overline{\gamma}_n) \to 0$ almost surely.
		This yield point \eqref{flag:spectrum} for $\lambda_{i,j}(\nu) = \sum_{k = i}^{j-1} \lambda_{k,k+1}(\nu)$.		
		
		By Theorem \ref{th:eigenspace}, $\eigen^+(\bigwedge^k \overline{\gamma}_n)$ converges almost surely to $\ell^\infty(\bigwedge^k \gamma)$ in $\proj(\bigwedge^k E)$ for all $k \in \Theta(\nu)$.
		We need to check that for almost every realisation of $\gamma$, the $k$-line $\ell^\infty(\bigwedge^k \gamma)$ is pure in the sense that it is of the form $\bigwedge^k V$ for a given $V \in \mathrm{Gr}_k(E)$.
		Since $\overline{\gamma}_n \in \End(E)$, $\eigen^+(\bigwedge^k \overline{\gamma}_n)$ is pure and the space of pure $k$-line is closed.
		Let us write $\mathcal{F}^\infty_k(\gamma) \in \mathrm{Gr}_k(E)$ for the unique subspace of $E$ such that $\ell^\infty(\bigwedge^k \gamma) = \bigwedge^k \mathcal{F}^\infty_k(\gamma)$.
		
		By Theorem \ref{th:eigenspace}, for all $k \in \Theta(\nu)$, there exists constants $C_k, \beta_k > 0$ such that $\PP(\dist(\bigwedge^kF_k^{\mathrm{eigen}}(\overline{\gamma}_n), \bigwedge^k F_k^\infty(\gamma)) \ge e^{-\alpha}) \le C_k e^{-\beta_k n}$. This yields \eqref{eigen-flag} for $C := \sum_{k \in \Theta(\nu)} C_k$ and $\beta := \min_{k \in \Theta(\nu)} \beta_k$.
		
		Theorem \eqref{th:cvspeed}, tells us that for all $k \in \Theta(\nu)$, there exists constants $C_k, \beta_k > 0$ such that $\PP(\dist(\bigwedge^k\overline{\gamma}_n F, \mathcal{F}^\infty(\gamma)) \ge e^{-\alpha}\mid \bigwedge^k \overline{\gamma}_n F_k \neq 0) \le C_k e^{-\beta_k n}$. 
		Since the $F_k$'s are nested, the events $(\bigwedge^k \overline{\gamma}_n F_k = 0)$ also are.
		Moreover, $(\bigwedge^{\max\Theta(\nu)} \overline{\gamma}_n F_{\max \Theta(\nu)} = 0)$ is bounded away from $1$ by Lemma \ref{lem:descri-ker}, we denote by $p$ its upper bound.
		It follows that \eqref{eigen-flag} holds for $C := \sum C_k /(1-p)$ and $\beta := \min \beta_k$.
	\end{proof}
	
	The above corollary implies the existence of a unique $\nu$ stationary measure on $\mathrm{Fl}_{\Theta(\nu)}(E)$.
	Moreover, it tells us that this is the measure of a limit flag, which is given by a partially defined and shift-equivariant measurable map, like the Oseledets space decomposition.
	This gives a generalization for the specific i.i.d. case of Oseledets Theorem.

	\subsection{Regularity results with optimal moment assumptions}\label{intro:invertible}
	
	In the present paragraph, $E$ is a Euclidean, Hermitian or ultrametric vector space of dimension $d \ge 2$ and $\nu$ is a probability distribution over $\GL(E)$.
	For all $g \in \GL(E)$, we write $N(g) = \log(\|g\|\|g^{-1}\|)$.
	We assume that $\nu$ is strongly irreducible and proximal as characterized by Proposition \ref{prop:zariski}.
	
	Under these purely algebraic assumption, we give concentration results for the distribution of the coefficients and of the spectral radius that depend on the tail of $N_*\nu$.
	From these results, combined with Theorem \ref{th:cvspeed}, we deduce regularity results for the stationary measure.
	
	\begin{Th}[Probabilistic estimates for the coefficients and the spectral radius]\label{th:slln}
		Let $\nu$ be strongly irreducible and proximal over $\GL(E)$.
		There exist constants $C, \beta > 0$ such that for all $f \in E^* \setminus \{0\}$, all $v \in E \setminus \{0\}$, for all $n\in\NN$, for $\overline{\gamma}_n \sim \nu^{*n}$ and for all $t \in \RR_{\ge 0}$,
		\begin{equation}\label{dom-coef}
			\PP\left(\log\frac{\|f\|\|\overline{\gamma}_n\|\|v\|}{|f \overline{\gamma}_n v|} > t\right) \le C e^{-\beta n} + \sum_{k=1}^{+\infty} C e^{-\beta k} N_*\nu\left( t / k , +\infty \right).
		\end{equation}
		Moreover:
		\begin{equation}\label{dom-radius}
			\forall t \ge 0, \; \PP\left(\log\frac{\|\overline{\gamma}_n\|}{\rho_1(\overline{\gamma}_n)} > t\right)\le \sum_{k=1}^{+\infty} C e^{-\beta k} N_*\nu\left( t / k , +\infty \right).
		\end{equation}
	\end{Th}
	
	We remind that, Theorem \ref{th:slln} is derived from Theorem \ref{th:pivot} and for that reason, we can express fitting constants $C$ and $\beta$ in terms of the constants $\alpha, \eps, m$ for which $\nu$ in $(\alpha, \eps, m)$-contracting in the sense of Definition \ref{def:quant-contraction}.
	
	The term $C e^{- \beta n}$, that is present in \eqref{dom-coef} but not in \eqref{dom-radius} naturally appears in the proof of Theorem \ref{th:coef-control}, from which \eqref{dom-coef} is derived. 
	This term is annoying because it prevents us from giving a uniform bound on the moments of $- \log\frac{|f \overline{\gamma}_n v|}{\|f\| \|\overline{\gamma}_n\| \|v\|}$.
	However it is not a mere proof artefact, to give an example, let $\nu$ be a lazy measure, \ie that gives measure $1/2$ to the identity matrix, let $f \in E^* \setminus \{0\}$ and $v \in \ker(f) \setminus \{0\}$, then $\PP(|f \overline{\gamma}_n v| = 0) \ge 2^{-n}$ for all $n$.
	
	The result of Theorem \ref{th:slln} may seem technical at first glance, let us illustrate the scope of the result by using it to prove more digestible results, Theorem \ref{th:intro-lln} on almost sure convergence of the coefficients and Corollary \ref{cor:regxi} on regularity of the stationary measure.
	
	Given $C, \beta > 0$ and $\nu$ a probability measure on $\GL(E)$, we denote by $\zeta^{C,\beta}_\nu$ the probability distribution on $\RR_{\ge 0}$ characterized by the fact that for all $t \ge 0$,
	\begin{equation}\label{def-zetanu}
		\zeta^{C,\beta}_\nu(t,+\infty) := \min \left\{ 1,\sum_{k=1}^{+\infty} C e^{-\beta k}N_*{\nu}\left( t / k, +\infty\right) \right\}.
	\end{equation} 
	Note that $\sum_{k=1}^{+\infty} C e^{-\beta k} < + \infty$ so $\sum_{k=1}^{+\infty} C e^{-\beta k} N_*{\nu} (t / k, +\infty)$ goes to $0$ as $t \to +\infty$ by monotonous convergence.
	So $\zeta^{C,\beta}_\nu$ is indeed a probability distribution on $[0, + \infty)$ and since everything is non-negative,
	\begin{align}
		\EE\left(\zeta^{C,\beta}_\nu\right) & =  \int_{0}^{\infty} \zeta^{C,\beta}_\nu(t,+\infty) dt \notag\\
		& \le \int_{0}^{\infty} \sum_{k  = 0}^{+ \infty} Ce^{-\beta k} N_*\nu(t / k,+\infty) dt \notag \\
		& \le \sum_{k  = 0}^{+ \infty} Ce^{-\beta k} \int_{0}^{\infty} k N_*\nu(t, +\infty) dt \notag \\
		\EE\left(\zeta^{C,\beta}_\nu\right) & \le \frac{C}{e^\beta - 2 + e^{- \beta}} \EE(N_*\nu). \label{expectation}
	\end{align}
	
	In particular, $\EE\left(\zeta^{C,\beta}_\nu\right)$ is finite when $\EE(N_*\nu)$ is.
	However, Theorem \ref{th:slln} does not give a satisfying result in case $N_*\nu$ has a finite exponential moment because $\zeta^{C,\beta}_\nu$ does not have a finite exponential moment in this case.
	This is not an issue because we already know by \cite{GR86} that when $N_*\nu$ has a finite exponential moment, the logarithm of the coefficients and of the spectral radius both satisfy large deviations inequalities above and below the top Lyapunov exponent and there exists a measure $\zeta$ that also has a finite exponential moment and plays the same role as the measure $\zeta_\nu^{C, \beta}$.
	See \cite[Chapter~14]{livreBQ} for a precise description of the behaviour of products of random matrices under exponential moment assumptions.
	
	Let us now derive Theorem \ref{th:intro-lln} from Theorem \ref{th:slln}.
	
	\begin{proof}[Proof of Theorem \ref{th:intro-lln}]
		Let $C, \beta > 0$ be as in Theorem \ref{th:slln}.
		By \eqref{expectation}, $\EE(\zeta_\nu^{C,\beta}) < +\infty$. 
		Then for all $\eps > 0$ and by \eqref{dom-coef},
		\begin{align*}
			\sum_{n = 1}^{+ \infty} \PP\left( \log\frac{\|f\| \|\overline{\gamma}_n\| \|v\|}{|f \overline{\gamma}_n v|} \ge n \eps \right) & \le \sum_{n = 1}^{+\infty} C e^{-\beta n} + \sum_{n\in\NN} \zeta_\nu^{C, \beta}(n \eps, + \infty) \\
			& \le \frac{C}{1 - e^{-\beta}} + \eps^{-1}\EE(\zeta_\nu^{C,\beta}) < +\infty.
		\end{align*}
		Then by Borel-Cantelli's Lemma, $n^{-1}\log\frac{\|f\|\|\overline{\gamma}_n\|\|v\|}{|f \overline{\gamma}_n v|} \to 0$ almost surely. 
		Moreover, by Furstenberg and Kesten~\cite[Theorem~1]{porm}, $n^{-1}\log{\|\overline{\gamma}_n\|} \to \lambda_1(\nu)$ almost surely.
		By \eqref{dom-radius} and by the above computation, \[\sum_{n = 1}^{+ \infty} \PP\left( \log\frac{\|\overline{\gamma}_n\|}{\rho_1(\overline{\gamma}_n)} \ge n \eps \right) \le \eps^{-1}\EE(\zeta_\nu^{C,\beta})\] for all $\eps > 0$.
		Then by the above argument, $n^{-1}\log{\|\overline{\gamma}_n\|} \to \lambda_1(\nu)$ almost surely.
	\end{proof}
	
	Note also that \eqref{dom-coef} implies that for all non-random sequence $\alpha_n \to 0$ and for all $1 \le i,j \le \dim(E)$, the sequence $\left( \frac{\left| (\overline{\gamma}_n)_{i,j} \right|^{\alpha_n}}{\left\| \overline{\gamma}_n \right\|^{\alpha_n}} \right)_{n \in \NN}$ converges in probability to $1$, \ie its distribution converges to the Dirac measure at $1$ in the weak-$*$ topology.
	That holds for any strongly irreducible and proximal probability measure over $\GL(E)$ without moment assumption. 
	The same holds for $\rho_1$ by \eqref{dom-radius}.
	This observation implies that any weak convergence result (\textit{e.g.} a central limit theorem) that holds for $\log \left\| \overline{\gamma}_n \right\|$ also holds for $\log \left|(\overline{\gamma}_n)_{i,j} \right|$ and $\lambda_1(\overline{\gamma}_n)$.
	
	Let us now deduce a result on the regularity of the stationary measure from Theorem \ref{th:slln}.
	The formulation \eqref{regxislp} is analogous to the regularity result for the stationary measure on hyperbolic groups \cite[Proposition~5.1]{Benoist2016CentralLT}. 
	This is also an improvement of the regularity results for the stationary measure on the projective space stated in \cite[Proposition~4.5]{CLT16}.
	
	Given $V \in \mathrm{Gr}(E)$, and $l \in\proj(E)$, we write $\dist (\proj(V), l):= \min_{ v \,\in \,V\setminus\{0\}}\dist([v],l)$ and we set $\mathcal{N}_r(V):= \left\{l \in\proj(E)\,\middle|\,\dist(\proj(V),l) < r\right\}$ for $0 < r \le 1$.
	Given $p >0$, we say that a random variable $x \in\RR_{\ge 0}$ has a finite $\ELL^p$ moment if $\EE(x^p) < +\infty$ and we say that $x$ has a weak $\ELL^p$ moment if $\sup_{t \ge 0} t^p \PP(x > t) < + \infty$.
	By integration by parts, $\EE(x^p) = \int_0^\infty pt^{p-1} \PP(x > t) dt$ for all $p > 0$.
	So with the same computations as for \eqref{expectation}, for all $p > 0$, the (weak) $\ELL^p$ moment of $\zeta_\nu^{C,\beta}$ is bounded above by the (weak) $\ELL^p$ moment of $N_*\nu$ times a constant that depends on $C, \beta$ and $p$.

	\begin{Cor}[Regularity of the measure]\label{cor:regxi}
		Let $E$ be a Euclidean vector space.
		Let $\nu$ be a strongly irreducible and proximal probability measure on $\GL(E)$. Let $C, \beta$ be as in Theorem \ref{th:slln} and let $\xi_\nu$ be the only $\nu$-stationary measure on $\proj(E)$. Then,
		\begin{equation}\label{regxi}
			\forall V \in\mathrm{Gr}(E) \setminus\{E\}, \forall 0 < r \le 1,\; \xi_\nu (\mathcal{N}_r(V)) \le \zeta_\nu^{C,\beta} (|\log(r)|, +\infty).
		\end{equation}
		Let $p >0$.
		If we assume that $N_*\nu$ has a finite $\ELL^p$ moment then there exists a constant $C'$ such that:
		\begin{equation}\label{regxislp}
			\forall V \in\mathrm{Gr}(E) \setminus\{E\},\; \int_{l \in \proj(E)}|\log\dist(\proj(V), l)|^p d\xi_\nu(l) \le C'.
		\end{equation}
		If we assume that $N_*\nu$ has a weak $\ELL^p$ moment then there exists a constant $C'$ such that:
		\begin{equation}\label{regxiwlp}
			\forall V \in\mathrm{Gr}(E) \setminus\{E\}, \forall 0 < r < 1,\; \xi_\nu (\mathcal{N}_r(V)) \le C'|\log(r)|^{-p}.
		\end{equation}
	\end{Cor}
	
	\begin{proof}
		First note that \eqref{regxislp} and \eqref{regxiwlp} follow directly from \eqref{regxi}.
		Let us prove \eqref{regxi}.
		Let $f \in E^* \setminus \{0\}$ be such that $f V = \{0\}$ and let $u \in E\setminus\{0\}$ be a constant.
		Let $\gamma = (\gamma_n)_{n \ge 0} \sim \nu^{\otimes \NN}$.
		We claim that:
		\begin{equation}\label{reg}
			\xi_\nu (\mathcal{N}_r(V))  \le \limsup_{n \to + \infty} \PP\left(|f \overline{\gamma}_n u| < r \|f\| \|\overline{\gamma}_n\| \|u\|\right).
		\end{equation}
		Let $\ell^\infty(\gamma)$ be as in Theorem \ref{th:cvspeed}.
		We have $\ell^\infty(\gamma) \sim \xi_\nu$ so $\xi_\nu (\mathcal{N}_r(V))  = \PP(\ell^\infty(\gamma)) \in \mathcal{N}_r(V))$.
		Moreover, by Theorem \ref{th:cvspeed},$[\overline{\gamma}_n u] \to \ell^\infty$ almost surely and $\mathcal{N}_r(V)$ is open, so $\ell^\infty(\gamma) \in \mathcal{N}_r(V)$ almost surely implies that there exists $n_0$ such that $[\overline{\gamma}_n u] \in \mathcal{N}_r(V)$ for all $n \ge n_0$.
		Therefore, $\xi_\nu (\mathcal{N}_r(V)) = \lim_n \PP([\overline{\gamma}_n u] \in \mathcal{N}_r(V))$.
		
		Moreover, $\mathcal{N}_r(V) \subset \left\{[v]\mid |f v| < r \|f\|\|v\| \right\}$, which yields $([\overline{\gamma}_n u] \in \mathcal{N}_r(V)) \allowbreak \subset (|f \overline{\gamma}_n u| < r \|f\| \|\overline{\gamma}_n u\|)$.
		Moreover $\|\overline{\gamma}_n u\| \le \|\overline{\gamma}_n\| \|u\|$ so \eqref{reg} holds.
		Moreover, by \eqref{dom-coef}, $\limsup_{n} \PP\left(|f \overline{\gamma}_n u| < r \|f\| \|\overline{\gamma}_n\| \|u\|\right) \le \zeta_\nu^{C,\beta} (|\log(r)|, +\infty)$, which yields \eqref{regxi}.
	\end{proof}

	By reducing to the strongly irreducible and proximal case, we can get rid of the algebraic assumptions for the law of large numbers on the spectral radius.
	
	\begin{Cor}[Law of large numbers for the spectral radius without algebraic assumptions]\label{cor:gen-eigen}
		Let $\nu$ be any measure on $\GL(E)$.
		Let $\left(\gamma_n\right)_{n \ge 0} \sim \nu^{\otimes\NN}$.
		Assume that $\EE(\log^+\|\gamma_0\|) < +\infty$ and that $\EE(\log^+\|\gamma_0^{-1}\|) < +\infty$.
		Then for all $1 \le k \le d$, it holds almost surely that:
		\begin{equation}
			\lim_{n \to + \infty} \frac{\lambda_k(\overline{\gamma}_n)}{n} = \lambda_k(\nu).\label{gen-eigen}
		\end{equation}
	\end{Cor}
	
	\begin{proof}
		We remind that by definition $\lambda_k(\nu) = \lim_{n \to + \infty} \frac{\log s_k(\overline{\gamma}_n)}{n}$.
		Let $\Gamma_\nu$ be the closed semi-group generated by $\nu$.
		
		Assume that $\Gamma_\nu$ is irreducible and that $\lambda_1(\nu) > \lambda_2(\nu)$.
		Then by \cite{GR1985}, $\Gamma_\nu$ is actually strongly irreducible and proximal so $\frac{\lambda_1(\overline{\gamma}_n)}{n} \to \lambda_1(\nu)$ almost surely by Theorem \ref{th:intro-lln}.
		
		Now let us reduce the general case to the irreducible case with $\lambda_1 > \lambda_2$.
		First note that for all $g$, by sub-multiplicativity of the norm, $\rho_1(g) \le s_1(g) = \|g\|$; the same holds for $g^{-1}$ so $s_d(g) \le \rho_d(g)$ and therefore $s_d(g) \le \rho_k(g)\le s_1(g)$ for all $1 \le k \le d$.
		Hence, if we assume that $\lambda_1(\nu) = \lambda_d(\nu)$, then \eqref{gen-eigen} holds.
		
		Assume that $\lambda_1(\nu) > \lambda_2(\nu)$ and that $\Gamma_\nu$ is not irreducible.  
		Then $\Gamma_\nu$ is conjugated to a semi-group of block-wise upper triangular matrices whose diagonal blocks are all irreducible.
		Moreover, we know that the Lyapunov spectrum of $\nu$ is the union of the Lyapunov spectrums of the diagonal blocks and the spectrum of each element $g \in \Gamma_\nu$ is the union of the spectrums of the diagonal blocks of $g$. 
		
		Then for all $g \in\Gamma_\nu$, we write $\pi(g)$ for the restriction of $g$ to the diagonal block on which $\lambda_1(\nu)$ appears. 
		Then $\pi_*\nu$ is strongly irreducible by assumption and $\lambda_1(\pi_* \nu) = \lambda_1(\nu)$ by assumption, moreover $\lambda_2(\pi_*\nu) = \lambda_j(\nu)$ for some $j \ge 2$ so $\lambda_2(\pi_*\nu) < \lambda_1(\pi_*\nu)$. 
		Therefore, $\frac{\lambda_1\pi(\overline{\gamma}_n)}{n} \to \lambda_1(\nu)$ almost surely. 
		Moreover, $\rho_1\pi(g) \le \rho_1(g)$ for all $g \in\Gamma_\nu$ so $\frac{\lambda_1(\overline{\gamma}_n)}{n} \to \lambda_1(\nu)$ almost surely.
		
		Now assume that $\lambda_1(\nu) = \lambda_2(\nu) > \lambda_d(\nu)$.
		Let $j$ be the smallest integers such that $\lambda_{j}(\nu) > \lambda_{j+1}(\nu)$.
		Then $\lambda_1 (\bigwedge^j \nu ) = j\lambda_j(\nu) > \lambda_2(\bigwedge^j \nu) = (j-1) \lambda_j(\nu) + \lambda_{j+1}(\nu)$.
		Therefore and by the above argument $\frac{\lambda_1\left(\bigwedge^j \overline{\gamma}_n\right)}{n} \to j\lambda_j(\nu)$ almost surely.
		Hence $\frac{\lambda_1\left(\overline{\gamma}_n\right) + \cdots +\lambda_j\left(\overline{\gamma}_n\right)}{n} \to j\lambda_1(\nu)$ almost surely.
		Moreover $\limsup_n \frac{\lambda_k\left(\overline{\gamma}_n\right)}{n} \le \lambda_1(\nu)$ for all $k$. 
		It yields that $\frac{\lambda_1(\overline{\gamma}_n)}{n} \to \lambda_1(\nu)$ almost surely.
		
		In conclusion, the convergence $\lambda_1(\overline{\gamma}_n) / n \to \lambda_1(\nu)$ holds almost surely without any algebraic assumptions on $\nu$.
		For $k \ge 2$, \eqref{gen-eigen} follows from the fact that $\lambda_k(g) = \lambda_1(\bigwedge^k g) - \lambda_1(\bigwedge^{k-1} g)$.
	\end{proof}

	One may wonder whether the reduction used in The proof of Corollary \ref{cor:gen-eigen} can be used to prove the almost sure convergence of $\log|f\overline\gamma_n v| / n$ without algebraic assumptions.
	The short answer is no and a counterexample goes as follows.
	Let $\nu$ be the distribution or a random rotation $R_{2^{-x} \pi}\in \GL(\RR^2)$, where $x$ is a random integer that takes value $k$ with probability $2^{-k!}$ for all $k \ge 1$. 
	The support of $\nu$ generates a dense subgroup of $\mathrm{SO}(\RR^2)$ so it admits the Lebesgue measure as its unique stationary measure on $\proj(\RR^2)$ but it is not contracting.
	Let $(\gamma_n) \sim \nu^{\otimes \NN}$. We claim that the random walk $\overline\gamma_n$ is recurrent in the sense that it returns infinitely many times at $\mathrm{Id}$ and therefore the anti-diagonal coefficients of $\overline\gamma_n$ do not satisfy a law of large numbers.
	
	Let us give a sketch of the proof of that fact.
	For all $k$, we write $t_k$ for the first time at which $\gamma_{t_k} = R_{2^{-k} \pi}$.
	Then $\PP(t_k < 2^{-k! - k}) \le 2^{-k}$ for all $k$ so by Borel Cantelli's lemma, there exists a random integer $k_0$ such that $t_k \ge 2^{k! - k}$ for all $k \ge k_0$.
	Let us fix an integer $k$ and assume that $k \ge k_0$, then the conditional law of $(\gamma_n)_{0 \le n < 2^{k! - k}}$ with respect to the event $(k \ge k_0)$ is the independent coupling of $2^{k! - k}$ copies of the restriction of $\nu$ to the group generated by $R_{\pi 2^{-k-1}}$, which has order $2^{k}$. 
	
	This random walk is recurrent and returns at $\mathrm{Id}$ a random mixing time of order $2^{(k-1)!}$. 
	Indeed, it returns in the sub-group generated by $R_{\pi 2^{-k-2}}$ after a random geometrically distributed time with expectation $2^{(k-1)!}$.
	By induction it returns to $\mathrm{Id}$ after a time of order $2^{(k-2)!} << 2^{(k-1)!}$, so with large probability (of order $1 - 2^{-(k-1)}$), it does so before going out of $R_{\pi 2^{-k-1}}$ again.
	Therefore, for infinitely many values of $k$, we find an integer $2^{(k-1)!} \le n < 2^{k! - k}$ such that $\overline\gamma_n = \mathrm{Id}$.

	\subsection{Background and example of application}\label{intro:background}
	
	The study or products of random matrices bloomed with the eponym article~\cite{porm} where Furstenberg and Kesten construct an escape speed for the logarithm of the norm using the sub-additivity. This proof was generalized by Kingman's sub-additive ergodic Theorem~\cite{K68}. This article followed the works of Bellman~\cite{Bellman} who showed the almost sure convergence of the rescaled logarithms of coefficients as well as a central limit theorem for one specific example.
	In~\cite{porm} Furstenberg and Kesten exhibit a law of large numbers for the norm under a strong $\ELL^1$ moment condition for $\log\|\cdot\|$.
	For matrices that have positive entries and under a mixed $\ELL^1$-$\ELL^\infty$ moment condition, they prove a law of large numbers for the coefficients (entries) and under an additional $\ELL^{2+\delta}$ moment assumption, they exhibit a central limit Theorem.
	These works on matrices inspired the theory of measurable boundary theory for random walks on groups~\cite{furstenberg1973boundary}. 
	
	Under a finite exponential moment assumption, in~\cite{lepage82} Le Page exhibits a spectral gap for the action of $\nu$ on a well chosen space of Hölder regular function over the projective space.
	From that follows a Central limit Theorem \cite{lepageTCL} and large deviations inequalities.
	A more quantitative approach to the large deviation inequalities is developed in the works of Sert and Aoun~\cite{sert2018large}~\cite{aoun2022random}.

	In~\cite{livrebl}, Lacroix gives an example of application of the theory of products of random matrices to quantum physics. 
	Namely, in the second part of \cite{livrebl}, the theory of products of random matrices is applied to the study of the Schrödinger operator associated to one particle evolving in a random (stochastic) one dimensional environment with bounded width. In this case, the eigenfunctions of that random Schrödinger operator are associated to the dynamic of a random sequence of matrices is $\mathrm{Sp}(\CC^{2w})$, where $w$ denotes the width of the "one dimensional" environment. 
	In this case, the positivity of the top Lyapunov exponent prevents the existence of normalizable eigenfunction, which in turn implies diffusion. 
	
	The framework developed in the present paper would allow to drop the finite first moment assumption in this case. 
	In applications, moment conditions come naturally. 
	However, one may want to study a random renewed environment (\ie a concatenation of random i.i.d. sections of random size) for which the renewal length (\ie the size of the independent sections) has infinite expectation.
	This kind of random processes that renew after a heavy tailed number of steps appear naturally in Hawkes processes with excitation. 
	Random renewed environments may not be stochastic when the renewal time has infinite expectation so the results of \cite{livrebl} do not apply directly.
	The methods developed in the present article allow to extend some results of Lacroix to this non-stochastic setting.

	In~\cite{livreBQ} Yves Benoist and Jean-François Quint give an extensive state of the art overview of the field of study with an emphasis on the algebraic properties of semi-groups. Later, in~\cite{xiao2021limit} Xiao, Grama and Liu use~\cite{CLT16} to show that coefficients satisfy a law of large numbers under some technical $ \ELL^2 $ moment assumption. 
	We can also mention~\cite{grama2020zeroone} and~\cite{xiao2022edgeworth} that give other probabilistic estimates for the distribution of the coefficients. 
	The strong law of large numbers and central limit-theorem for the spectral radius were proven by Aoun and Sert in~\cite{aoun2020central} and in~\cite{aoun2021law} under an $\ELL^2$ moment assumption. 
	
	The importance of alignment of matrices was first noted in~\cite{AMS-esmigroup} along with the importance of Schottky sets. 
	These notions were then used by Aoun~\cite{RA2011} to show that outside of an exponentially rare event, independent draws of an irreducible random walk that has finite exponential moment generate a free group. 
	The results of the present article allow us to drop the moment assumption. 
	In~\cite{cuny2016limit} and~\cite{cuny2017komlos}, Cuny, Dedecker, Jan and Merlevède give KMT estimates for the behaviour of $\left(\log\|\overline{\gamma}_n\|\right)_{n\in\NN}$ under $\ELL^p$ moment assumptions for $p > 2$.
	
	The main difference between these previous works and this paper is that the measure $\nu$ has to be supported on the General Linear group $\GL(E)$ for the historical methods to work.
	Indeed, they rely of the existence of the stationary measure $\xi_\nu$ on $\proj(E)$, which is a consequence of the fact that $\GL(E)$ acts continuously on $\proj(E)$, which is compact.
	Some work has been done to study non-invertible matrices in the specific case of matrices that have real positive coefficients. In~\cite{porm}, Furstenberg and Kesten show limit laws for the coefficients under an $\ELL^\infty$ moment assumption, in~\cite{AM87} and~\cite{KS87_semigroup} Mukherjea, Kesten and Spitzer show some limit theorems for matrices with non-negative entries that are later improved by Hennion in~\cite{HH97} and more recently improved by Cuny, Dedecker and Merlevède in~\cite{cuny2023limit}. 
	
	It is worth mentioning that when studying products of random matrices with non-negative coefficients, we may look at the action of the random walk on the image of the non-negative quadrant in $\proj(E)$.
	This set is a compact simplex of dimension $d-1$, on which the semi-group of non-negative matrices acts continuously. 
	Classical tools of ergodic Theory apply rather well to this setting.
	Note that the relation $g \AA^\eps h$ holds as long as $\min\frac{^te_i g e_j}{\|g\|} > \eps$ when $g$ and $h$ have non-negative coefficients. 
	Therefore, applying Theorem \ref{th:pivot} to a non-degenerate product of random matrices with non-negative coefficients, we find a measure $\mu \le \nu^{\otimes m} /\alpha$ that is $0$-Schottky.
	In this setting Lemma \ref{lem:first-part} implies Theorem \ref{th:ex-piv} and we do not need the pivot algorithm described in Definition \ref{def:pivot-alg}.

	\subsection{About the pivoting technique}
	
	While it is important to mention the work 
	In Section \ref{sec:pivot} we mainly use the tools introduced in~\cite{pivot}, some of them having been introduced or used in former works like~\cite{bmss20} where Adrien Boulanger, Pierre Mathieu, Cagri Sert and Alessandro Sisto introduce the notion of large deviations inequalities below a given speed for random walks in discrete hyperbolic groups or~\cite{devin} where Mathieu and Sisto show some bi-lateral large deviations inequalities in the context of distributions that have a finite exponential moment.
	
	In~\cite{pivot} Sébastien Gouëzel uses the pivoting technique in the setting of hyperbolic groups to get large deviations estimates below the escape speed and to show the continuity of the escape speed. 
	For us, the most interesting part of Gouëzel's work is the "toy model" described in section 2. 
	In~\cite{choi1}~\cite{choi2}~\cite{choi3} Inhyeok Choi applies the pivoting technique to show results that are analogous to the ones of Gouëzel for the mapping class group of an hyperbolic surface as well as other groups that satisfy some contraction property. 
	
	In~\cite{CFFT22}, Chawla, Forghani, Frisch and Tiozzo use another view of the pivoting technique and the results of~\cite{pivot} to show that the Poisson boundary of random walk with finite entropy on a group that has an acylindrical action on an hyperbolic space is in fact the Gromov Boundary of said space. 
	In \cite{chawla2025poissonboundarydiscretesubgroups}, the same authors show that the Poisson boundary of random walks with finite entropy on discrete Zariski dense subgroups of real Lie groups is the Furstenberg boundary, answering a question asked in a previous version of the present article.
	
	\section{Local-to-global properties for the alignment of matrices}\label{kak}

	In this section, we describe the geometry of the semi-group $\Gamma:= \End(E)$ for $E$ a Euclidean, Hermitian or ultrametric vector space of dimension $d \ge 2$. 
	To get a good intuition, we may think of $\KK = \RR$ since the inequalities used in the proofs hold in the complex and ultrametric cases.
	
	We want to translate ideas of hyperbolic geometry into the language of products of endomorphisms. 
	The idea is to exhibit a local-to global property in the fashion of~\cite[Theorem~4]{Cannon1984}. 
	That way we can adapt the construction of pivotal times in~\cite[Section~2]{pivot} to the setting of products of random matrices.
	
	This strategy is similar to the ones used by Yves Benoist in \cite{Benoist1997} and Abels, Margulis and Soifert in \cite{AMS-esmigroup} to study asymptotic properties of linear groups.

	\subsection{Alignment and contraction coefficient}\label{sec:ali-sqz}
	
	Let us fix $0 < \eps < 1$, given two matrices or homomorphisms\footnote{Contrary to the Archimedean case Archimedean case, two isomorphic ultrametric spaces may not be isometric so we can not identify homomorphisms with matrices by the simple choice of a basis} $g$ and $h$, we write $g \AA^\eps h$ when the product $gh$ is well defined (\ie $g$ has as many columns as $h$ has rows) and $\|gh\| > \eps \|g\|\|h\|$.
	This extends the notation introduced Definition \ref{def:ali}.
	
	We remind that given $x,y$ two vectors, the norm of their exterior product is characterized by $\|x \wedge y\| =  \min_{a\in\KK} \|x - ay\|\|y\|$. 
	That way, 
	\begin{equation*}
		\|h\wedge h\| = \max_{x,y \in E\setminus\{0\}}\min_{a\in\KK} \frac{\|h(x - ay)\|\|h(y)\|}{\|x\|\|y\|}
	\end{equation*}
	for all linear map $h \in \mathrm{Hom}(E,F)$.

	\begin{Def}[Singular gap]\label{def:sqz}
		Let $E,F$ be Euclidean, Hermitian or ultrametric vector spaces and $h \in \mathrm{Hom}(E,F)\setminus\{0\}$. We define the first singular gap of $h$ as:
		\begin{equation*}
			\sigma(h) := \frac{\|h \wedge h\|}{\|h\|^2}.
		\end{equation*}
	\end{Def}
	
	Given a matrix $h$ of size $d \times d$, the quantity $\sigma(h)$ is in the interval $[0,1]$ and it is equal to the ratio $\frac{s_2(h)}{s_1(h)}$ where $s_1(h) \ge s_2(h) \ge \dots \ge s_d(h) \ge 0$ are the singular values of $h$.
	
	\begin{Def}[Distance between projective classes]\label{def:distance}
		Let $E$ be a Euclidean, Hermitian or ultrametric space. We denote by $\proj(E)$ the projective space of $E$ \ie the set of lines in $E$, endowed with the following distance map:
		\begin{equation}\label{distanceprojective}
			\forall x, y \in E \setminus \{0\},\; \dist([x],[y]):= \frac{\|x\wedge y\|}{\|x\|\|y\|} = \min_{a \in \KK}\frac{\|x - ay\|}{\|x\|}.
		\end{equation}
	\end{Def}
	
	\begin{Lem}[Lipschitz property for the norm cocycle]\label{lem:product-is-lipschitz}
		Let $E$ and $F$ be two Euclidean, Hermitian or ultrametric spaces and let $f \in\End(E)\setminus\{0\}$. 
		For all $x, y \in E \setminus\{0\}$,
		\begin{equation}\label{prod-is-lip}
			\left|\frac{\|fx\|}{\|f\|\|x\|} - \frac{\|fy\|}{\|f\|\|y\|}\right| \le \dist([x],[y]).
		\end{equation}
	\end{Lem}
	
	\begin{proof}
		Let $f \in\mathrm{Hom}(E,F)$ and let $x,y \in E \setminus\{0\}$. 
		We show that $\|fx\| \le \|fy\| + \|f\|\|x\|\dist([x],[y])$, which implies \eqref{prod-is-lip} since $x$ and $y$ play symmetric roles.
		Let $c \in\KK$ be such that $\|x - y c\| = \min_{a \in \KK} \|x - y a\|$.
		Then by Definition \ref{def:distance}, $\dist([x],[y]) = \|x - y c\| / \|x\|$. 
		Moreover $\|yc\| \le \|x\|$ by property of the orthogonal projection in the Euclidean and Hermitian cases and by ultrametric inequality otherwise.
		By triangular inequality and by homogeneity of the norm,
		\begin{equation*}
			\|f x\|  \le \|f y c\| + \|f (x - y c)\| 
			\le \|f y\|+ \|f\|\|x - y c\|. \qedhere
		\end{equation*}
	\end{proof}

	\begin{Def}\label{def:dom}
		Let $E,F$ be Euclidean, Hermitian or ultrametric spaces and let $h\in\mathrm{Hom}(E,F)\setminus\{0\}$.
		For all $0 < \eps \le 1$, we define $V^\eps(h):= \{x\in E;\|h x\| \ge \eps \|h\|\|x\|\}$ and $U^\eps(h):= h(V^\eps(h))$ and $W^\eps(h):= U^\eps(^th)$.
	\end{Def}
	
	Note that for all $h$ and all $0 \le \eps \le 1$ the sets $V^\eps(h)$, $U^\eps(h)$ and $W^\eps(h)$, are homogeneous in the sense that they are invariant under scalar multiplication. 	
	Note also that for $h$ an endomorphism of rank one, and for all $0 < \eps < 1$, the cone $U^\eps(h)$ is the image of $h$ so it has diameter $0$ in the projective space. 
	The idea to have in mind is that, given $h$ a matrix that has a large singular gap, the set $U^\eps(h)$ has a small diameter in the following sense.
	
	\begin{Lem}\label{lem:cont-prop}
		Let $E,F$ be Euclidean spaces, let $h\in\mathrm{Hom}(E,F)\setminus\{0\}$. Let $u \in U^1(h)\setminus\{0\}$ and let $u' \in U^\eps(h)\setminus\{0\}$. Then:
		\begin{equation}\label{cont-prop}
			\dist([u], [u']) \le \eps^{-1}{\sigma(h)}.
		\end{equation}
	\end{Lem}
	
	\begin{proof}
		Let $v \in V^1(h)$ and let $v'\in V^\eps(h)$ be such that $u = hv$ and $u' = hv'$.
		Then $u\wedge u' = \bigwedge^2 h(v \wedge v')$ so:
		\begin{equation*}
			\|u\wedge u'\|< \left\|\textstyle\bigwedge^2 h\right\| \|v \wedge v'\|.
		\end{equation*}
		Now saying that $v \in V^1(h)$ and $v'\in V^\eps(h)$ means that $\|u\| = \|h\| \|v\|$ and $\|u'\| \ge \eps \|h\| \|v'\|$. 
		Hence:
		\begin{equation*}
			\|u\| \|u'\| \ge \eps \|h\|^2\|v\|\|v'\|.
		\end{equation*}
		Then by taking the quotient of the above inequalities,
		\begin{equation*}
			\frac{\|u\wedge u'\|}{\|u\| \|u'\|} < \frac{\left\|h \wedge h\right\|}{\eps \|h\|^2}\frac{\|v\wedge v'\|}{\|v\| \|v'\|}\le \frac{\left\|h \wedge h\right\|}{\eps \|h\|^2}.
		\end{equation*}
		By definition, the term on the left is $\dist([u], [u'])$ and the member on the right is $\eps^{-1}\sigma(h)$.
	\end{proof}

	Lemma \ref{lem:cont-prop} tells us that the projective image of $U^\eps(h)$ has diameter at most $\eps$ as long as $\sigma(h) \le \eps^2 $.
	With the toy model analogy, the condition $\sigma(h) \le \eps^2$ will play the role of the condition for  word to be non-trivial. 
	We will extensively use the following simple remarks.
	
	\begin{Lem}
		Let $g$ and $h$ be non-zero matrices such that the product $gh$ is well defined. We have $g\AA^\eps h$ if and only if ${^th}\AA^\eps {^tg}$. Moreover $\sigma({^th})=\sigma(h)$.
	\end{Lem}
	
	\begin{proof}
		This well known fact is the consequence if three simple points. 
		One is that $\|h\| = \|{^th}\| = \max_{f,v} |fhv|/(\|f\|\|v\|)$ for all matrix $h$. 
		The second is that $^t(gh) = ^th ^tg$ so $\|gh\| = \|^th ^t g\|$.
		The third point is that $^th\wedge ^th = ^t(h\wedge h)$ so $\|^th \wedge ^th\| = \|h \wedge h\|$. 
	\end{proof}
	
	\begin{Lem}\label{lem:sharp-ali}
		Let $g$ and $h$ be non-zero matrices such that the product $gh$ is well defined and let $0 < \eps < 1$. 
		If there exist $u \in U^1(h)\setminus\{0\}$ and $w \in W^1(g)\setminus\{0\}$ such that $\frac{|wu|}{\|w\|\|u\|} > \eps$, then $g\AA^\eps h$.
		If $g\AA^\eps h$, then there exist $u \in U^\eps(h)\setminus\{0\}$ and $w \in W^\eps(g)\setminus\{0\}$ such that $w \AA^\eps u$, $g \AA^\eps u$, $w \AA^\eps h$.
	\end{Lem}
	
	\begin{proof}
		Let $u \in U^1(h)\setminus\{0\}$ and $w \in W^1(g)\setminus\{0\}$.
		Assume that $\frac{|wu|}{\|w\|\|u\|} > \eps$.
		Let $f \in V^1({^tg})$ and $v \in V^1(h)$ be such that $w = fg$ and $u = hv$.
		We have $\frac{|fghv|}{\|fg\|\|hv\|} > \eps$, therefore $\frac{|fghv|}{\|f\|\|g\|\|h\|\|v\|} > \eps$, so $\frac{\|gh\|}{\|g\|\|h\|} > \eps$, which means that $g \AA^\eps h$. 
		This proves the first implication of Lemma \ref{lem:sharp-ali}.
		
		Now assume that $g \AA^\eps h$.
		Let $f \in E^*\setminus\{0\}$ and let $v \in E\setminus\{0\}$ be such that $|fghv| = \|f\|\|gh\|\|v\|$. 
		Then $\frac{|fghv|}{\|f\|\|g\|\|h\|\|v\|} = \frac{\|gh\|}{\|g\|\|h\|} > \eps$. 
		Let $u:= hv$ and let $w:= fg$. 
		Then,
		\begin{equation*}
			\eps < \frac{|wu|}{\|w\|\|u\|} \cdot \frac{\|fg\|\|hv\|}{\|f\|\|g\|\|h\|\|v\|} = \frac{|fgu|}{\|f\|\|g\|\|u\|} \cdot \frac{\|hv\|}{\|h\|\|v\|}= \frac{|whv|}{\|w\|\|h\|\|v\|} \cdot \frac{\|fg\|}{\|f\|\|g\|}.
		\end{equation*}
		All the factors are in $[0, 1]$ so $\frac{|w u|}{\|w\|\|u\|} > \eps$ hence $w \AA^\eps u$, and $\frac{\|g u\|}{\|g\| \|u\|} > \frac{|f g u|}{\|f\| \|g\| \|u\|} \ge \eps$ hence $g \AA^\eps u$, and $\frac{\|w h\|}{\|w\| \|h\|} >  \frac{|whv|}{\|w\|\|h\|\|v\|}\ge \eps$ hence $w \AA^\eps h$.
		Moreover $\frac{\|h v\|}{\|h\|\|v\|} > \eps$, hence $v \in V^\eps(h)$ so $u \in U^\eps(h)$, and $\frac{\|fg\|}{\|f\|\|g\|} > \eps$ so $w \in W^\eps(g)$. 
		This proves the second implication  of Lemma \ref{lem:sharp-ali}.
	\end{proof}

	\begin{Lem}\label{lem:c-prod}
		Let $g$ and $h$ be non-zero matrices. Assume that $g \AA^\eps h$. Then one has:
		\begin{align}
			\sigma(gh) & < \eps^{-2} \sigma(g) \sigma(h). \label{lenali}
		\end{align}
		Moreover, for every non-zero vectors $u\in U^1(g)\setminus\{0\}$, and $u'\in U^1(gh)\setminus\{0\}$,
		\begin{equation}\label{limbd}
			\dist([u], [u']) < \frac{\sigma(g)}{\eps}.
		\end{equation}
	\end{Lem}
	
	\begin{proof}
		Note that the norm of the $\wedge$ product is sub-multiplicative because it is a norm so:
		\begin{equation}\label{m2}
			\|g h \wedge g h\| \le \|g \wedge g\|\|h \wedge h\|.
		\end{equation}
		Moreover, $g \AA^\eps h$ so $\|gh\|^{-2} < \eps^{-2} \|g\|^{-2}\|h\|^{-2}$, we combine this inequation with \eqref{m2}, using Definition \ref{def:sqz} for $\sigma$ and we find \eqref{lenali}.
		
		Now to prove \eqref{limbd}, we only need to show that $U^1(gh) \subset U^\eps(g)$ and use \eqref{cont-prop} from Lemma \ref{lem:cont-prop}. Indeed, for $v \in V^1(gh)$, one has $\|g h v\| \ge \eps \|g\| \|h\| \|v\| > \eps \|g\| \|h v\|$ which means that $h v \in V^\eps(g)$, so $g h v \in U^\eps(g)$.
		Hence $U^1(gh) \subset U^\eps(g)$. 
		Moreover, by \eqref{cont-prop} in Lemma \ref{lem:cont-prop} $\dist([u],[u']) < \frac{\sigma(g)}{\eps}$ for all $u \in U^1(g)$ and all $u' \in U^\eps(g)$, which proves \eqref{limbd}.
	\end{proof}

	\subsection{Local-to-global properties for the alignment}
	
	In this paragraph, we prove local to global properties for the alignment relation under a contraction condition. 
	We remind that if $g_0, \dots, g_n$ is a family of matrices of rank one (\ie $\sigma = 0$) such that $g_1 \cdots g_n \neq 0$, then for all $k$ and for all $0 < \eps \le 1$, the local and global alignment relations $(g_{k-1} \AA^\eps g_k)$, $(g_0 \cdots g_{k-1} \AA^\eps g_k \cdots g_n)$, $(g_0 \cdots g_{k-1} \AA^\eps g_k)$ and $(g_{k-1} \AA^\eps g_k \cdots g_n)$ are equivalent.
	
	Another analogy is the free group $\FF$, used as toy model in \cite[Section~2]{pivot}.
	We say that two words $g,h$ in $\FF$ are aligned when the last letter of the reduced form of $g$ is not the inverse of the first letter in the reduced form of $h$. 
	That way, given a family of non-trivial words $g_0, \dots, g_n$, saying that $g_{k-1} \AA g_k$ for all $1 \le k\le n$ is equivalent to saying that $g_{0} \cdots g_{k - 1}\AA g_{k} \cdots g_n$ for all $k$. It is also equivalent to saying that for all $k$, there exists $0 \le l_k \le k-1$ and $k \le r_k \le n$ such that $g_{l_k}\cdots g_{k-1} \AA g_{k}\cdots g_{r_k}$.
	
	For product of matrices that have rank more than one, we only assume that $\sigma(g_n) \le \eps^6/48$. In this case, instead of equivalences, we get that the local $\eps$-alignment condition implies a global $\eps/2$ alignment.
	This is the reason why we get $\AA^{\eps / 2}$ in point \eqref{pivot:herali} of Theorem \ref{th:pivot} while we assumed $(\alpha, \eps, m)$-contraction.

	\begin{Lem}[Transmission of the alignment]\label{lem:herali}
		Let $f,g,h$ be non-zero matrices such that the product $fgh$ is well defined. Assume that $\sigma(g) \le \eps^2/4$ and that $f \AA^{\eps} g \AA^{\frac{\eps}{2}} h$. Then $f \AA^{\frac{\eps}{2}} gh$.
	\end{Lem}

	\begin{proof}
		Let $u \in U^1(gh)$.
		Then $\|ghu\| = \|gh\|\|u\| > \frac{\eps}{2} \|g\|\|h\|\|u\|$ so $\|ghu\| > \frac{\eps}{2} \|g\|\|hu\|$ and therefore $u \in U^{\frac{\eps}{2}} (g)$.
		Let $v \in V^1(fg)\setminus\{0\}$. 
		Then $\|fgv\| > \|f\|\|g\|\|v\|\eps$ so $g v \in U^{\eps}(g)\setminus\{0\}$.
		Let $u'\in U^1(g)\setminus\{0\}$. 
		By Lemma \ref{lem:cont-prop} $\dist([gv],[u']) \le \frac{\eps}{4}$ and by Lemma \ref{lem:c-prod} $\dist([u], [u'])\le \frac{\eps}{4}$ so by triangular inequality $\dist([u],[gv])\le \frac{\eps}{2}$. 
		It follows from Lemma \ref{lem:product-is-lipschitz} that $\frac{\|fu\|}{\|f\|\|u\|} \ge \frac{\|fgv\|}{\|f\|\|gv\|} - \frac{\eps}{2}$. 
		Moreover $\frac{\|fgv\|}{\|f\|\|gv\|} \ge \frac{\|fgv\|}{\|f\|\|g\|\|v\|} > \eps$.
		Therefore $\frac{\|fu\|}{\|f\|\|u\|} > \frac{\eps}{2}$ and $u \in U^1(gh)$ so $f \AA^{\frac{\eps}{2}} gh$ by Lemma \ref{lem:sharp-ali}.
	\end{proof}
	
	\begin{Rem}[The ultrametric case is easier]
		Let $\KK$ be a ultrametric locally compact field and let $f,g,h$ be matrices with entries in $\KK$ such that the product $f,g,h$ is non trivial. 
		Let $0 < \eps \le 1$ be so that $\sigma(g) < \eps^2$ and $g \AA^\eps h$.
		Then $f \AA^\eps gh$ if, and only if, $f \AA^\eps g$. 
		This fact allows us to get rid of the multiplicative constants in the results of the present section in the ultrametric case.
	\end{Rem}
	
	\begin{Lem}\label{lem:triple-ali}
		Let $f$, $g$ and $h$ be non-zero matrices such that the product $fgh$ is well defined and let $0 < \eps \le 1$. Assume that $f \AA^\eps g \AA^\eps h$ and that $\sigma(g) < \eps^4/4$. Then $\sigma(fgh) < 4 \eps^{-4} \sigma(f) \sigma(g) \sigma(h) < 4 \eps^{-4} \sigma(g)$.
	\end{Lem}
	
	\begin{proof}
		By \eqref{lenali} in Lemma \ref{lem:c-prod}, $\sigma(gh) < \eps^{-2}\sigma(g) \sigma(h)$. 
		Moreover, by Lemma \ref{lem:herali}, $f \AA^{\frac{\eps}{2}} gh$ so $\sigma(fgh) < 4 \eps^{-2} \sigma(f) \sigma(gh) \le 4 \eps^{-4} \sigma(f) \sigma(g) \sigma(h)$.
		We conclude using the fact that $\sigma$ takes values in $[0,1]$.
	\end{proof}
	
	The following Lemma and Corollary allow us to derive Theorems \ref{th:escspeed} and \ref{th:cvspeed} from Theorem \ref{th:pivot}.
	
	\begin{Lem}[Contraction property for aligned chains]\label{lem:c-chain}
		Let $n \ge 1$ and let $g_0, \dots, g_n$ be non-zero matrices. Assume that $\sigma(g_k) < \eps^2/8$ and $g_{k} \AA^\eps g_{k+1}$ for all $k \in\{0, \dots, n-1\}$. Then $g_0 \AA^{\frac{\eps}{2}} g_{1} \cdots g_n$ and
		\begin{gather}
			\| g_0\cdots g_n \| \ge \left(\frac{\eps}{2}\right)^n\prod_{j = 0}^n \|g_j\| \quad\text{and} \label{alinorm-chaine}\\
			\sigma(g_0\cdots g_n) \le 2^{2n}\eps^{-2n}\prod_{j = 0}^n \sigma(g_j).\label{lenali-chaine}
		\end{gather}
	\end{Lem}
	
	\begin{proof}
		We prove by induction on $0 \le k < n$ that $g_{k} \AA^{\frac{\eps}{2}} g_{k + 1}\cdots g_n$. 
		This holds for $k = n-1$ because $g_{n-1} \AA^\eps g_n$. 
		Let $1 \le k < n$ and assume by induction that $g_{k} \AA^{\frac{\eps}{2}} g_{k+1}\cdots g_n$. 
		Lemma \ref{lem:herali}, applied to $f = g_{k-1}$, $g:= g_k$ and $h:= g_{k+1}\cdots g_n$, yields that $g_{k-1} \AA^{\frac{\eps}{2}} g_{k} \cdots g_n$, which concludes the induction. 
		
		As a consequence, $\|g_k \cdots g_n\| \ge \frac{\eps}{2} \|g_k\| \|g_{k+1} \cdots g_n\|$ for all $0 \le k < n$.
		Therefore, $\|g_k\cdots g_n\| \ge \left(\frac{\eps}{2}\right)^{n-k}\|g_k\| \allowbreak \|g_{k+1}\|\cdots \|g_n\|$, for all $0 \le k <n$, which yields \eqref{alinorm-chaine}.
		
		Let us apply Lemma \ref{lem:c-prod} to $g_{k} \AA^{\frac{\eps}{2}} g_{k + 1}\cdots g_n$. 
		Inequality \eqref{lenali} yields $\sigma(g_k \cdots g_n) < 4 \eps^{-2}\sigma(g_k)\sigma(g_{k+1} \cdots g_n)$ for all $0 \le k < n$.
		This implies \eqref{lenali-chaine} by a telescopic argument.
	\end{proof}

	\begin{Cor}[Limit line]\label{cor:limit-line}
		Let $(g_n)_{n\in\NN}$ be a sequence of matrices.
		Assume that for all $n \in \NN$, one has $g_n \AA^\eps g_{n+1}$ and $\sigma(g_{n+1}) \le \eps^2/8$.
		Then there exists $\ell^\infty \in \proj(E)$ such that:
		\begin{equation}\label{limit-line}
			\forall n\in\NN,\; \forall u_n \in U^1(g_0 \cdots g_{n-1}) \setminus \{0\},\; \dist([u_n], \ell^\infty) \le \frac{2 \sigma(g_0 \cdots g_{n-1})}{\eps}.
		\end{equation}
	\end{Cor}
	
	\begin{proof}
		For all $k > 0$, let $u_k \in U^1(g_0\cdots g_{k-1}) \setminus \{0\}$.
		Let $0 < k < m$.
		Lemma \ref{lem:alipart} applied to $(g_0, \dots, g_{m-1})$ yields $(g_0 \cdots g_{k-1}) \AA^\frac{\eps}{2} (g_k \cdots g_{m-1})$.
		Hence, Lemma \ref{lem:c-prod} implies:
		\begin{equation}\label{limite-mn}
			\dist([u_k],[u_m]) <  \frac{2 \sigma(g_0 \cdots g_{k-1})}{\eps}
		\end{equation} 
		and \eqref{lenali-chaine} yields $\sigma(g_0 \cdots g_{k-1}) < 2^{-k-1}$.
		Hence, the sequence $([u_k])_{k \ge 1}$ is a Cauchy sequence in $\proj(E)$, therefore it has a limit. 
		Let $\ell^\infty := \lim_m [u_m]$.
		For all $k$, taking $m \to +\infty$ in \eqref{limite-mn} yields
		\begin{equation*}\label{limit-uk}
			\dist([u_k],\ell^\infty) \le \frac{2 \sigma(g_0 \cdots g_{k-1})}{\eps}.
		\end{equation*}
		To conclude, note that for all $n \in \NN$, if we only change the value of $u_n$ in the sequence $(u_k)_{k \ge 0}$ then the limit $l_\infty$ does not change and \eqref{limit-uk} still holds for $k = n$ as long as $u_n \in U^1(g_0 \cdots g_{n-1}) \setminus \{0\}$.
	\end{proof}

	Let us now state more technical results on local to global properties fo the alignment.

	\begin{Lem}[Rigidity of the alignment]\label{lem:rig-ali}
		Let $f,g_1, g_2,h$ be non-zero matrices. Assume that $\sigma(g_i) \le \eps^2/12$ for $i \in\{1,2\}$ and that $f \AA^{\frac{\eps}{2}} g_1 \AA^\eps g_2  \AA^{\frac{\eps}{2}} h$. Then $f g_1 \AA^{\frac{\eps}{2}} g_2 h$.
	\end{Lem}
	
	\begin{proof}
		The proof is similar to the proof of Lemma \ref{lem:herali}.
		Let $w \in W^1(fg_1)$ and $u \in U^1(g_2 h)$. 
		Then $w \in W^{\frac{\eps}{2}} (f g_1)$ and $u \in U^{\frac{\eps}{2}}(g_2)$.
		Let $w' \in W^\eps(g_1)$ and $u' \in U^\eps(g_2)$ be such that $w' \AA^\eps u'$.
		Such $u', w'$ exist by Lemma \ref{lem:sharp-ali}.
		
		By Lemma \ref{lem:cont-prop} applied to $^tg_1$ and by triangular inequality, we get $\dist([w], [w']) < \frac{3\sigma(g_1)}{\eps} \le 3\eps/12 = \eps/4$. 
		With the same argument $\dist([u], [u']) < \eps/4$.
		By Lemma \ref{lem:product-is-lipschitz} we get $w \AA^{\eps - \frac{\eps}{4}} u'$ and $w \AA^{\frac{\eps}{2}} u$.
		Then by Lemma \ref{lem:sharp-ali}, we get $f g_1 \AA^{\frac{\eps}{2}} g_2 h$.
	\end{proof}

	\begin{Lem}[local-to-global property for the alignment]\label{lem:alipart}
		Let $g_0,\dots,g_n$ be non-zero matrices such that the product $g_0\cdots g_n$ is well defined. 
		Assume that $\sigma(g_i) \le \eps^2/12$ for every $k \in \{1,\dots,n-1\}$.
		Assume also that $g_0\AA^\eps g_1 \AA^\eps \cdots \AA^\eps g_n$, \ie $g_k \AA^\eps g_{k+1}$ for all $k \in \{0, \dots, n-1\}$. 
		Then for all $k \in \{1, \dots, n\}$,
		\begin{equation}\label{jsgh}
			g_0 \cdots g_{k-1} \AA^\frac{\eps}{2} g_{k} \cdots g_{n}.
		\end{equation}
	\end{Lem}
	
	\begin{proof}
		The case $k =1$ is a reformulation of Lemma \ref{lem:c-chain} and the case $k = n$ is the transpose of the case $k = 1$.
		
		Let $2 \le k \le n-1$.
		Lemma \ref{lem:herali} applied to $(g_k, \dots, g_n)$ yields $g_k \AA^\frac{\eps}{2} g_{k+1} \cdots g_n$. 
		The same Lemma \ref{lem:herali} applied to $(^t g_{k-1}, \dots, ^t g_0)$ yields $g_0 \cdots g_{k-2} \AA^{\frac{\eps}{2}} g_{k-1}$.
		Moreover, $g_{k-1} \AA^\eps g_k$ and $\sigma(g_{k-1}), \sigma(g_k) \le \eps^2/12$.
		Combine all that above and apply Lemma \ref{lem:rig-ali} to obtain \eqref{jsgh}.
	\end{proof}

	Lemmas \ref{lem:alipart} and Corollary \ref{cor:limit-line} give a good intuition on how local alignment properties of a sequence of matrices are linked to contraction properties for the action on the projective space.
	
	To prove Theorem \ref{th:pivot}, we need to deal with a notion of half local half global alignment condition.
	Let us give a more complete version of Lemma \ref{lem:alipart}, which is exactly what we need for the pivoting technique to work.
	
	\begin{Lem}\label{lem:Atilde}
		Let $n \ge 0$ and $g_{-1}, g_0, g_1, \dots, g_{2n}$ be non-zero matrices. 
		Assume that $g_{-1} \AA^\eps g_0$, that $\sigma(g_0) \le \eps^2/12$ and that for all $0 \le i < n$, we have $\sigma(g_{2i + 1}) \le \eps^6/48$ and:
		\begin{equation}\label{hyp-lem-atilde}
			g_0\cdots g_{2i}\AA^\eps g_{2i + 1} \AA^\eps g_{2i + 2}.
		\end{equation}
		Then $g_{-1}\AA^\frac{\eps}{2}(g_0\cdots g_{2n})$ and $\sigma(g_0 \cdots g_{2n}) \le \sigma(g_0)$.
		
		Let $f, h$ be non-zero matrices such that the products $f g_{-1}$ and $g_n h$ are well defined.
		Assume moreover that $\sigma(g_{-1}) \le \eps^6/48$, that $f \AA^\frac{\eps}{2} g_{-1}$ and that $g_0 \cdots g_n \AA^\frac{\eps}{2} h$.
		Then $fg_{-1} \AA^\frac{\eps}{2} g_0\cdots g_{2n} h$.
	\end{Lem}
	
	\begin{proof}
		Let $0 \le i < n$. 
		By Lemma \ref{lem:triple-ali} and by \eqref{hyp-lem-atilde},
		\begin{equation*}
			\sigma(g_0\cdots g_{2i+2}) < 4\eps^{-4} \sigma(g_0\cdots g_{2i})\sigma(g_{2i+1}) \le \sigma(g_0\cdots g_{2i}) \eps^2/12.
		\end{equation*} 
		By induction on $i$, we obtain $\sigma(g_0\cdots g_{2i}) \le \sigma(g_0)(\eps^2/12)^{i} \le (\eps^2/12)^{i + 1}$.
		Moreover, by \eqref{hyp-lem-atilde} and Lemma \ref{lem:herali}, $g_0\cdots g_{2i} \AA^{\frac{\eps}{2}} g_{2i + 1} g_{2i + 2}$ for all $i$.
		Given $u_i \in U^1(g_0\cdots  g_{2i}) \setminus \{0\}$ and $u_{i+1} \in U^1(g_0\cdots  g_{2i +2}) \setminus \{0\}$, \eqref{limbd}:
		\begin{equation*}
			\dist([u_i], [u_{i+1}]) < \frac{2\sigma(g_0\cdots g_{2i})}{\eps} \le \frac{\eps}{6} (\eps^2/12)^{i}.
		\end{equation*}
		For all $u_0 \in U^1(g_0)$ the triangular inequality yields:
		\begin{equation*}
			\dist([u_0], [u_i]) \le \frac{\eps}{4} \sum_{k = 0}^{i-1} (\eps^2/12)^{k} \le \frac{\eps/6}{1-\eps^2/12} \le \frac{2 \eps}{11}.
		\end{equation*}
		Let $u \in U^\eps(u_0)$ be such that $g_{-1}\AA^\eps u$, such a $u$ exists by Lemma \ref{lem:sharp-ali} applied to $g_{-1} \AA^\eps g_0$.
		Then by Lemma \ref{lem:cont-prop}, $\dist([u], [u_0]) \le \eps/8$ and by triangular inequality, $\dist([u], [u_n]) \le \eps/12 + 2\eps/11 \le 3\eps/11 \le \eps/2$.
		By Lemma \ref{lem:product-is-lipschitz}, $g_{-1} \AA^\frac{\eps}{2} u_n$.
		Moreover $u_n \in U^1(g_0\cdots g_{2n}) \setminus\{0\}$ so $g_{-1}\AA^\frac{\eps}{2}(g_0\cdots g_{2n})$ by Lemma \ref{lem:sharp-ali}.
		
		Now assume that $\sigma(g_{-1}) \le \eps^6/48$ and let $f,h$ be such that $f \AA^\frac{\eps}{2} g_{-1}$ and $g_0 \cdots g_n \AA^{\frac{\eps}{2}} h$. 
		If $n = 0$, then $fg_{-1} \AA^\frac{\eps}{2} g_0\cdots g_{2n} h$ by Lemma \ref{lem:rig-ali}.
		Assume that $n \ge 1$.
		Then, $\sigma(g_0\cdots g_n) \le (\eps^2 / 12)^{n+1} \le \eps^4/144$.
		For all $u' \in U^1(g_0\cdots g_n h)$, Lemma \ref{lem:c-prod} yields
		\begin{equation*}
			\dist([u_n], [u']) < \frac{2 \sigma(g_0\cdots g_n)}{\eps} \le \frac{\eps^3}{72}.
		\end{equation*}
		Therefore $\dist(u_n, u) \le \frac{\eps^3}{72} + \frac{\eps}{12} + \frac{2\eps}{11}$ by triangular inequality.
		Let $w \in W^{\eps}(g_{-1})$ be such that $w \AA^\eps u$, let $w'\in W^1(fg_{-1})$ and $w'' \in W^1(g_{-1})$.
		By Lemma \ref{lem:cont-prop}, $\dist([w],[w'']) \le \frac{\eps^5}{48}$.
		By Lemma \ref{lem:c-prod}, $\dist([w'],[w'']) \le \frac{2\eps^5}{48}$ so by triangular inequality, $\dist([w],[w']) \le \frac{3\eps^5}{48}$. 
		If follows from Lemma \ref{lem:product-is-lipschitz} that $\frac{|w'u|}{\|w'\|\|u\|} \ge \eps - \frac{\eps^5}{16}$.
		Applying \ref{lem:product-is-lipschitz} again yields $\frac{|w'u'|}{\|w'\|\|u'\|} \ge \eps - \frac{\eps^5}{16} - \frac{\eps^3}{72} - \frac{\eps}{12} - \frac{2\eps}{11} > \eps/2$, which concludes the proof.
	\end{proof}

	\subsection{Link between singular values and eigenvalues}
	
	The present paragraph is dedicated to the proof of the following elementary result. 
	
	\begin{Lem}\label{lem:eigen-align}
		Let $g$ be an endomorphism such that $\sigma(g) \le \eps^2/8$ and $g \AA^\eps g$. 
		Then $g$ is proximal and the following inequalities hold:
		\begin{gather}
			\rho_1(g) \ge \frac{\eps}{2} \|g\|,\label{ev-sv} \\
			\frac{\rho_2(g)}{\rho_1(g)} \le  \frac{4 \sigma(g)}{\eps^2},\label{sqz-prox} \\
			\forall u \in U^1(g), \dist([u], \eigen^+(g)) \le \frac{2\sigma(g)}{\eps}. \label{dom-eigen}
		\end{gather}
	\end{Lem}
	
	\begin{proof}
		Consider $(g_k)_{k\ge 0}$ to be the constant sequence equal to $g$.
		First, by Lemma \ref{lem:c-chain}, $\sigma(g^n) \le \left(\frac{4\sigma(g)}{\eps^2}\right)^n$ and $\|g^n\| \ge \left(\frac{\eps}{2}\|g\|\right)^n$ for all $n$.
		Therefore, $\liminf_n \|g^n\|^{1/n} \ge \frac{\eps}{2}\|g\|$.
		Moreover, by the spectral theorem, this inferior limit is a limit  equal to ${\rho_1(g)}$. This proves \eqref{ev-sv}.
		
		To get \eqref{dom-eigen}, we apply Corollary \ref{cor:limit-line}. 
		We get a line $\ell^\infty$ such that $[u_n] \to \ell^\infty$ for any choice of sequence $u_n \in U^1(g^n)\setminus\{0\}$. 
		Let $\ell^\infty$ be such a line and $(u_n)_{n \ge 0}$ be such a sequence.
		by Corollary \ref{cor:limit-line}, $\dist([u], \ell^\infty) \le \frac{2 \sigma(g)}{\eps}$ for all $u \in U^1(g)$. 
		Therefore, we only need to show that $\ell^\infty = \eigen^+(g)$. 
		Let $e \in \eigen^+(g)\setminus \{0\}$. 
		Then $\|ge\| = \rho_1(g) \|e\|$ so $e \in V^{\eps /2}(g)$  by \eqref{ev-sv}, as a consequence $e \in U^\frac{\eps}{2}(g)$. 
		Moreover, $e$ is an eigenvector (\ie $e \in \KK g e$) so $e \in U^{\eps /2}(g)$. 
		This reasoning holds for all powers of $g$ so $e \in U^{\eps /2}(g^n)$ for all $n$. 
		Therefore $\ell^\infty = [e]$.
	\end{proof}

	\section{Random products and extractions}\label{markov}
	
	This section is dedicated to the proof of basic results for products of random matrices.
	In Paragraph \ref{sec:kernel}, we describe the probabilistic behaviour of the kernel $\ker(\overline{\gamma}_n)$ in the case $\nu(\GL(E)) < 1$.
	In the second paragraph, we give a proof of Lemma \ref{lem:contraction}.
	
	We use the following notations. 
	Given $\Gamma$ a metric semi-group with a bi-lateral unit element $\mathbf{1}_\Gamma$, we write $\widetilde{\Gamma}:= \bigsqcup_{n \ge 0} \Gamma^n$ for the set of words on the alphabet $\Gamma$ and write $\odot$ for the concatenation product. 
	We write $()$ for the empty word, which is the unit of $\widetilde{\Gamma}$.
	
	As a disjoint union of finite products of metric space, $\widetilde{\Gamma}$ is endowed with a natural metrizable topology\footnote{Several metric classes give the same topology on $\widetilde{\Gamma}$, given a sequence $(\gamma_n)_{n \ge 0}$, the sequence $((\gamma_0, \dots, \gamma_n))_{n \in\NN}$ may or may not be a Cauchy sequence depending on the choice of metric.  
	For the sake of intuition, it is better endow $\widetilde\Gamma$ with a complete locally compact metric. 
	Given $k \le n$ and given two words $\tilde{\gamma}:= (\gamma_0,\dots,\gamma_{k-1})$ and $\tilde{\gamma}':= (\gamma'_0,\dots, \gamma_{n-1})$, we set $\dist_{\widetilde{\Gamma}}\left(\tilde{\gamma},\tilde{\gamma}'\right) = \sum_{i = 0}^{k-1} \dist_\Gamma(\gamma_i, \gamma'_i) + \sum_{i = k}^{n-1} \left(1 + \dist_\Gamma(\mathbf{1}_\Gamma, \gamma'_i)\right)$, that way closed bounded sets are compact.} and a Borel $\sigma$-algebra.
	We also endow $\widetilde{\Gamma}$ with two continuous semi-group morphisms:
	\begin{align*}
		L :  &\left(\widetilde{\Gamma}, \odot\right) \to \left(\NN, +\right); \;(\gamma_0, \dots, \gamma_{n-1})  \mapsto  n, \\
		\Pi: &\left(\widetilde{\Gamma}, \odot\right) \rightarrow \left(\Gamma, \cdot\right); \; (\gamma_0, \dots, \gamma_{n-1}) \mapsto \gamma_0 \cdots \gamma_{n-1}.
	\end{align*}

	Given $(\tilde{\gamma}_i)_{i \ge 0} \in \widetilde{\Gamma}^\NN$ a sequence which is not eventually stationary to the empty word, we write $\bigodot_{i = 0}^{+\infty} \tilde{\gamma}_i \in \Gamma^\NN$ for the infinite concatenation.
	In other words, for all $k \in \NN$ and for $0 \le n < \sum_{i = 0}^k L(\tilde{\gamma}_i)$, the entry of index $n$ of the sequence $\bigodot_{i = 0}^{+\infty} \tilde{\gamma}_i$ is the $(1+n)$-th letter of the word $\tilde{\gamma}_0 \odot \cdots \odot \tilde{\gamma}_{k}$.
	One may be tempted to see $\Gamma^\NN$ as a boundary for $\widetilde{\Gamma}$ in one of the following senses:
	\begin{enumerate}
		\item It is a Furstenberg boundary for the action of $\widetilde{\Gamma}$ on $\Gamma^\NN$ by concatenation on the left, which is continuous for the product topology. Given a word $\tilde{\gamma} = (\gamma_0, \dots, \gamma_{l-1})\in \widetilde{\Gamma}$ and a sequence $(g_n)_{n \ge 0} \in \Gamma^\NN$ and given an integer $j$, we define the element at index $j$ in the sequence $\tilde{\gamma} \odot (g_n)_{n \ge 0} \in \Gamma^\NN$ to be $\gamma_j$ when $j < l$ and $g_{j-l}$ when $j \ge l$.
		\item Assuming $\Gamma$ to be complete, it is the visual boundary (\ie the points of the completion that aren't represented by stationary sequences) for the product metric on $\widetilde{\Gamma}$, which is generated by the distance $\dist_{\widetilde{\Gamma}}(\tilde{\gamma}, \tilde{g}) = 2^{-l} + \sum_{k = 0}^{l - 1} \min\{2^{-k}, \dist_\Gamma(\gamma_k, g_k)\}$, where $l = \min\{L(\tilde{\gamma}), L(\tilde{g})\}$.
	\end{enumerate}  
	With this vision, given a sequence $(\tilde{\gamma}_i)_{i \ge 0}$ of words, the infinite concatenation $\bigodot_{i = 0}^{+\infty} \tilde{\gamma}_i$ is the limit of the sequence $(\tilde{\gamma}_0 \odot \cdots \odot \tilde{\gamma}_{k})_{k \ge 0}$ both for the completion (because $(\bigodot_{i = 0}^k\tilde{\gamma}_i)_{k \ge 0}$ is a Cauchy sequence for $\dist_{\widetilde{\Gamma}}$) and in the sense of Furstenberg (Because in the product topology, $\bigodot_{i = 0}^k\tilde{\gamma}_i \odot g \to \bigodot_{i = 0}^{+\infty} \tilde{\gamma}_i$ for all sequence $g = (g_n)_{n \ge 0} \in \Gamma^\NN$).
	
	A small technicity arises here since the distance $d_{\widetilde{\Gamma}}$ for which $\Gamma^\NN$ can be identified with the visual boundary of $\widetilde{\Gamma}$ is not coherent with the metric structure induced by the notation $\widetilde{\Gamma} = \bigsqcup_{n \ge 0} \Gamma^n$, which would make $\widetilde{\Gamma}$ complete (again, assuming $\Gamma$ to be complete).
	Though interesting and deeply intertwined with the foundations of Furstenberg's Boundary Theory, these technicities are of little importance in the present work since both metrics generate the same Borel $\sigma$-algebra.
	
	Given $(w_k)_{k \ge 0} \in \NN^\NN$, we write $\overline{w}_k:= w_0 + \cdots + w_{k-1}$ for all $k \in\NN$ and given $(\gamma_n)_{n \ge 0} \in \Gamma^\NN$, we write $\gamma^w_k:= \gamma_{\overline{w}_k} \cdots \gamma_{\overline{w}_{k+1} - 1}$ and $\widetilde{\gamma}^w_k:= \left(\gamma_{\overline{w}_k}, \dots, \gamma_{\overline{w}_{k+1} - 1}\right)$ for all $k \in\NN$.
	That way, $\bigodot_{k = 0}^{+ \infty} \widetilde{\gamma}^w_k = (\gamma_n)_{n \ge 0}$ when $(w_k)_{k}$ is not eventually stationary to $0$.
	
	Given $\tilde\kappa_1$ and $\tilde\kappa_2$ two probability distributions on $\widetilde{\Gamma}$, we write $\tilde\kappa_1 \odot \tilde\kappa_2$ for their convolution product, \ie the distribution of $\tilde{\gamma}_1 \odot \tilde{\gamma}_2$ when $\left(\tilde{\gamma}_1, \tilde{\gamma}_2\right) \sim \tilde\kappa_1 \otimes \tilde\kappa_2$.
	Given $(\tilde\kappa_i)_{i \ge 0}$ a family of probability distributions not eventually stationary to $\delta_{()}$, we write $\bigodot_{i = 0}^{+\infty} \tilde\kappa_i$ for the infinite convolution, \ie the distribution of $\bigodot_{i = 0}^{+\infty} \tilde{\gamma}_i \in \Gamma^\NN$, when $\left(\tilde{\gamma}_n\right) \sim \bigotimes \tilde\kappa_i$.
	If we set $\tilde\kappa_i := \tilde\kappa$ for all $i$ and for a given measure $\kappa \neq \delta_{()}$ then we write $\tilde\kappa^{\odot \NN}$ for $\bigodot_{i = 0}^{+\infty} \tilde\kappa_i$.
	Using the Furstenberg boundary language, $\tilde\kappa^{\odot \NN}$ is characterized as the unique $\tilde\kappa$-stationary measure on $\Gamma^\NN$.

	\subsection{Rank and essential kernel of a probability distribution}\label{sec:kernel}
	
	In this paragraph, we describe the probabilistic behaviour of the kernel of a product of i.i.d. random matrices and prove Propositions \ref{prop:rank} and \ref{prop:essker}.
	This allows us to study products of random matrices that are not almost surely invertible and not asymptotically almost surely zero.
	Given $h$ a linear map, we denote by $\rank(h)$ the rank of $h$ \ie the dimension of the image of $h$.
	
	We say that a probability measure $\nu$ is supported on a set $S$ if $\nu(S) = 1$, we write $\supp(\nu)$ for the smallest closed set of full $\nu$-measure.
	Let us empathize that we may have $\nu(\GL(E)) = 1$ but $\supp(\nu) \not\subset \GL(E)$
	
	In this paragraph $E$ still denotes a finite dimensional vector space of dimension $d \ge 2$ over an unspecified local field but the results hold for any measurable field, like $\CC[X]$ or the algebraic closure of $\QQ_p$ for example.

	\begin{Def}[Rank of a distribution]\label{def:rank}
		Let $\nu$ be a probability measure on $\End(E)$. We define the eventual rank of $\nu$ as the largest integer $\evrank(\nu)$ such that:
		\begin{equation}
			\forall n \ge 0,\, \nu^{*n}\left\{\gamma \in \End(E) \; \middle| \; \rank(\gamma) < \evrank(\nu)\right\} = 0.
		\end{equation}
	\end{Def}

	\begin{Lem}[Eventual rank of a distribution]\label{lem:rank}
		Let $\Gamma:= \End(E)$. 
		Let $\nu$ be a probability measure on $\Gamma$ such that $\evrank(\nu) < d$. 
		There exists a probability measure $\tilde\kappa$ on $\widetilde{\Gamma}$ such that $\tilde{\kappa}^{\odot \NN} = \nu^{\otimes\NN}$ and $\Pi_*\tilde\kappa$ is supported on the set of endomorphisms of rank equal to $\evrank(\nu)$.
		Moreover $L_*\tilde\kappa$ has finite exponential moment.
	\end{Lem}
	
	\begin{proof}
		Given a non-random sequence $\tilde{\gamma} = (\gamma_n)_{n \in \NN} \in \Gamma^\NN$, the sequence $\left(\rank(\overline{\gamma}_n)\right)_{n \in \NN}$ is a non-increasing sequence of non-negative integers so it is eventually stationary.
		Write $r_\gamma$ for the limit of $\left(\rank(\overline{\gamma}_n)\right)_{n\in\NN}$.
		Then for all sequence $\gamma \in\Gamma^\NN$, there exists an integer $n' \ge 1$ such that $\rank\left(\overline{\gamma}_n\right) = r_\gamma$ for all $n \ge n'$.
		Write $n_\gamma$ for the minimal such $n'$. 
		Note that $\gamma \mapsto r_\gamma$ and $\gamma \mapsto n_\gamma$ are measurable maps.
		
		Now let $(\gamma_n)_{n\in\NN} \sim \nu^{\otimes \NN}$ be a random sequence. 
		We define $\tilde\kappa$ to be the distribution of $(\gamma_0, \dots, \gamma_{n_\gamma - 1})$. 
		Note that for all $k$, the event $(n_{\gamma} = k)$ is measurable for the $\sigma$-algebra generated by $(\gamma_0, \dots, \gamma_{k-1})$, which is independent of $\left(\gamma_{n+k}\right)_{n \ge 0}$.
		Therefore, the conditional distribution of $(\gamma_{n + n_{\gamma}})_{n\in\NN}$ with respect to $(\gamma_0, \dots, \gamma_{n_\gamma - 1})$ is $\nu^{\otimes\NN}$.
		Hence, $\tilde\kappa \odot \nu^{\otimes\NN} = \nu^{\otimes\NN}$ (and $\tilde\kappa \neq \delta_{()}$) so $\tilde\kappa^{\odot \NN} = \nu^{\otimes\NN}$.
		
		Moreover $\Pi_*\tilde\kappa$ is the distribution of $\overline{\gamma}_{n_\gamma}$ which has rank $r_\gamma$ and $L_*\tilde\kappa$ is the distribution of $n_\gamma$.	
		Therefore, we only need to show that $r_\gamma$ is almost surely constant and that $n_\gamma$ has finite exponential moment. 
		
		Let $r_0$ be the essential lower bound of $r_\gamma$ \ie the largest integer such that $\PP(r_\gamma \ge r_0) =1$. 
		Let $n_0 \in\NN$ be such that $\PP(\rank(\overline{\gamma}_{n_0}) = r_0) > 0$ and write $\alpha:= \PP(\rank(\overline{\gamma}_{n_0}) = r_0)$. 
		We claim that such an integer $n_0$ exists. 
		Indeed, by minimality, $\PP(r_\gamma > r_0) < 1$ so $\PP(r_\gamma = r_0) > 0$, which means that $\PP(\rank(\overline{\gamma}_{n_\gamma}) = r_0) > 0$.
		Let $n_0$ to be such that $\PP(n_\gamma \le n_0 \cap r_\gamma = r_0) > 0$. 
		Such an $n_0$ exists, otherwise $n_0$ would be almost surely infinite, which is absurd.
		
		Furthermore, $\PP(\rank(\gamma_{kn_0}\cdots \gamma_{(k+1)n_0-1}) = r_0) > 0$ for all $k \in\NN$ because the sequence $(\gamma_n)$ is i.i.d. 
		Moreover these events are independents so for all $k \in\NN$,
		\begin{align*}
			\PP(\forall k' < k, \rank(\gamma_{k'n_0}\cdots\gamma_{(k'+1)n_0-1}) > r_0 ) & = (1-\alpha)^k.
		\end{align*}
		The rank of a product is bounded above by the rank of each of its factors so:
		\begin{align*}
			\forall k\in\NN,\;\PP\left( \rank(\overline{\gamma}_{kn_0}) > r_0 \right) & \le (1-\alpha)^k, \\
			\text{hence, }\forall n\in\NN, \; \PP\left( \rank\left(\overline{\gamma}_{\lfloor\frac{n}{n_0}\rfloor n_0}\right) > r_0 \right) & \le (1-\alpha)^{\lfloor\frac{n}{n_0}\rfloor}.
		\end{align*}
		For all $n\in\NN$, we have $\lfloor\frac{n}{n_0}\rfloor \ge \frac{n}{n_0}-1$ and $\lfloor\frac{n}{n_0}\rfloor n_0 \le n$ so:
		\begin{align*}
			\forall n\in\NN, \; \PP\left(\rank(\overline{\gamma}_{n}) > r_0\right) & \le (1-\alpha)^{\frac{n}{n_0}-1}.
		\end{align*}
		Let $C = \frac{1}{1-\alpha}$ and $\beta = \frac{-\log{(1-\alpha)}}{n_0 }> 0$. 
		Then $\PP(\rank(\overline{\gamma}_{n}) > r_0 ) \le C e^{-\beta n}$ for all $n$. 
		Note that $ \PP\left(\rank(\overline{\gamma}_{n}) > r_0\right) \ge \PP(r_\gamma > r_0)$ for all $n$ so $\PP(r_\gamma > r_0) = 0$, which implies that $r_0 = r_\gamma$. 
		Hence $\PP\left(\rank(\overline{\gamma}_{n}) > r_0\right) = \PP( n < n_\gamma)$ for all $n$.
		It follows that $\PP(n_\gamma > n)\le C e^{-\beta n}$, which means that $n_\gamma$ has a finite exponential moment.
	\end{proof}

	\begin{Def}[Essential kernel]\label{def:essker}
		Let $\nu$ be a probability distribution on $\End(E)$. 
		We define the essential kernel of $\nu$ as:
		\begin{equation}\label{essker}
			\essker(\nu):= \left\{v \in E \,\middle|\,\exists n\in \NN, \nu^{*n}\{h\in\End(E)\mid hv = 0\} > 0 \right\}.
		\end{equation}
	\end{Def}

	\begin{Lem}\label{lem:descri-ker}
		Let $E$ be a Euclidean space of dimension $d \ge 2$ and let $\nu$ be a probability distribution on $\End(E)$. 
		There exists a probability distribution $\kappa$ on $\End(E)$ which is supported on the set of rank $\evrank(\nu)$ endomorphisms and such that:
		\begin{equation}\label{essker-kappa}
			\essker(\nu) = \essker(\kappa)  \{v \in E \mid \kappa\{h\in\End(E)\mid hv = 0\} > 0\}.
		\end{equation}
		Moreover, for all $v \in E$, the sequence $(\nu^{*n}\{h \mid hv = 0\})_n$ is non-decreasing and
		\begin{equation}
			\lim_{n \to +\infty}\nu^{*n}\{h \mid hv = 0\} = \kappa\{h \mid  hv = 0\}.\label{limker}
		\end{equation}
		Hence, there exists a constant $\alpha < 1$ such that:
		\begin{equation}\label{max-ker}
			\forall v\in E, \,\sup_{n\in\NN}\nu^{*n}\{h\in\End(E)\mid hv = 0\} \in [0, \alpha] \cup \{1\}.
		\end{equation}
		At last, the set:
		\begin{equation}\label{ker-is-inv}
			\left\{v \in E\,\middle|\, \sup_{n\in\NN}\nu^{*n}\{h\in\End(E)\mid hv = 0\}  = 1\right\}
		\end{equation}
		is a subspace of $E$ which is $\nu$-almost surely invariant.
	\end{Lem}
	
	\begin{proof}
		Let $(\gamma_n)_{n \ge 0} \sim \nu^{\otimes\NN}$. We define the random integer:
		\begin{equation*}
			n_0:= \min \left\{n \in\NN\,\middle|\,\rank(\gamma_{n-1} \cdots \gamma_0) = \evrank(\nu) \right\}.
		\end{equation*} 
		Let $g:= \gamma_{n_0 - 1} \cdots \gamma_0$ and let $\kappa$ be the distribution of $g$.
		Then by Lemma \ref{lem:rank} applied to the transpose of $\nu$, the random integer $n_0$ has finite exponential moment. 
		Let $v \in E$, and let $n \in\NN$.
		We claim that:
		\begin{equation}\label{ker-n}
			\nu^{*n}\{h \in \End(E)\mid  hv = 0\} = \PP(\gamma_{n-1} \cdots \gamma_0 v = 0) \le \PP(g v = 0).
		\end{equation}
		When $n < n_0$, we factorize $g = (\gamma_{n_0 - 1} \cdots \gamma_{n}) (\gamma_{n-1} \cdots \gamma_0)$, from which we deduce that $\ker(\gamma_{n-1} \cdots \gamma_{0}) \subset \ker(g)$.
		Therefore $(\gamma_{n-1} \cdots \gamma_0 v = 0,\, n < n_0) \subset (gv = 0, \, n< n_0)$.
		When $n \ge n_0$, $\ker(\gamma_{n-1} \cdots \gamma_{0}) \supset \ker(g)$. 
		By definition of $n_0$, $\rank(\gamma_{n-1} \cdots \gamma_0) = \rank(g)$ so $\ker(\gamma_{n-1} \cdots \gamma_{0}) = \ker(g)$. 
		It follows that:
		\begin{equation}\label{ker-egal}
			(\gamma_{n-1} \cdots \gamma_0 v = 0,\, n \ge n_0) = (gv = 0, \, n \ge n_0).
		\end{equation}
		Hence $(\gamma_{n-1} \cdots \gamma_0 = 0) \subset (gv = 0)$, which proves \eqref{ker-n}.
		A direct consequence of \eqref{ker-n} is that:
		\begin{equation}\label{ker-m}
			\essker(\nu) \subset \left\{v \in E \,\middle|\, \PP(gv = 0) > 0 \right\}.
		\end{equation}
		It follows from \eqref{ker-egal} that
		\begin{equation*}
			\forall n,\;\PP(\gamma_{n-1} \cdots \gamma_0 v = 0) \ge \PP(g v = 0) - \PP(n_0 > n).
		\end{equation*}
		Moreover $\PP(n_0 > n) \to 0$, which yields \eqref{limker} by double inclusion.
		Therefore, for all $v \in \essker(\kappa)$ (\ie such that $\PP(g v = 0)> 0$), there exists an integer $n \in\NN$ such that $\PP(\gamma_{n-1} \cdots \gamma_0 v = 0) > 0$ so $\essker(\kappa) \subset \essker(\nu)$.
		Combined with \eqref{ker-m}, this proves \eqref{essker-kappa}.
		
		Let us prove \eqref{max-ker}.
		Let $V = \{v \in E \mid  \PP(gv = 0) = 1\}$.
		First note that $V$ is a subspace of $E$ and that $gV = \{0\}$ with probability $1$. 
		Let
		\begin{equation}
			\alpha:= \sup_{n\in\NN, v \in E \setminus V}\nu^{*n}\{h\in\End(E)\mid hv = 0\}.
		\end{equation}
		Assume by contradiction that $\alpha = 1$.
		Let $(v_n)_{n \ge 0}$ be a non-random sequence in $E \setminus V$ such that $\PP(g v_n = 0) \ge 1 - 2^{-n}$ for all $n$. 
		Then $\sum_{n = 0}^{+\infty} \PP\left(g v_n \neq 0\right) < +\infty$. 
		Therefore, by Borel-Cantelli's Lemma, the random integer $m_0:= \max\{n \in\NN \mid  g v_n \neq 0\}$ is almost surely finite.
		Set $V':= \bigcap_{m \to +\infty} \sum_{n\ge m} \KK v_n$.
		Since $E$ is finite dimensional, the non-increasing sequence of subspaces $(\sum_{n\ge m}\KK v_n)_{m \ge 0}$ is eventually stationary, \ie there exists an integer $m \in\NN$ such that $v_n \in V'$ for all $n \ge m$.
		Therefore $V' \not\subset V$.
		Moreover, $g(V') \subset \sum_{n \ge m_0} \KK g v_n = \{0\}$ on the set $(m_0 < +\infty)$, which has probability $1$.
		It follows that $V' \subset V$, which is a contradiction so $\alpha < 1$.
		
		Let us prove that the set described in \eqref{ker-is-inv} is $\nu$-almost-surely invariant.
		Assume by contradiction that $\PP(\gamma_0(V) = V) \neq 1$. 
		Let $v \in V$ be such that $\PP(\gamma_0 v \notin V) > 0$. Then for all $n\in\NN$,
		\begin{align*}
			\PP(\gamma_{n}\cdots \gamma_0 v \neq 0) 
			& \ge \PP(\gamma_0 v \notin V) \PP(\gamma_n\cdots \gamma_1 \gamma_0 v \neq 0\mid \gamma_0 v \notin V) 
			\\
			&\ge \PP(\gamma_0 v \notin V) (1-\alpha)> 0,
		\end{align*}
		which is absurd because $\PP(\gamma_{n} \cdots \gamma_0 v \neq 0) \to 0$.
	\end{proof}

	Let un now prove Proposition \ref{prop:essker} which tells that $\essker(\nu)$ is a finite union of subspaces of $E$ of dimension at most $d-\evrank(\nu)$.
	Note that in dimension $d = 2$ or more generally when $\evrank(\nu) = d-1$, this follows from \eqref{essker-kappa}.
	Indeed, for a random matrix $h \sim \kappa$ of rank $d-1$, the set $\{[v] \in \proj(E)\mid \PP(hv = 0) > 0\}$ is the collection of atoms of the random variable $\ker(h) \in \proj(E)$ and we know that a random variable has countably may atoms.
	Let us prove the following Lemma which directly implies Proposition \ref{prop:essker}.
	
	\begin{Lem}\label{lem:essker}
		Let $\nu$ be a probability distribution on $\End(E)$. The set $\essker(\nu)$ is a countable union of subspaces of $E$ that each have dimension at most $\dim(E) - \evrank(\nu)$. 
	\end{Lem}
	
	\begin{proof}
		Let $d':= \dim(E) - \evrank(\nu)$. For all $k \in\{0, \dots, \dim(E)\}$, we denote by $\mathrm{Gr}_{k}(E)$ the set of subspaces of $E$ of dimension $k$. Let $\kappa$ be as in Lemma \ref{lem:descri-ker}
		First we show that $\essker(\nu)$ is included in a countable union of subspaces of dimension exactly $d'$.
		Given $n \in \NN$ and $\alpha > 0$, we define:
		\begin{equation*}
			K_\alpha:= \{x\in E\mid \kappa\{h\in\End(E)\mid h x =0\} \ge \alpha\}
		\end{equation*}
		Note that $\essker(\nu) = \bigcup_{\alpha > 0} K_{\alpha} = \bigcup_{m \in\NN} K_{2^{-m}}$ so we only need to show that $K_{2^{-m}}$ is included in a countable union of subspaces for all $m \in \NN$.
		
		Let $m \in\NN$, we claim that $K_{2^{-m}}$ is included in a union of at most $\binom{d'2^m}{d'}$ subspaces of $E$ of dimension $d'$.
		Let $g \sim \kappa$, write $\alpha:= 2^{-m}$ and assume that $K_\alpha \neq \{0\}$.
		
		Let $N$ be an integer and let $(x_1, \dots, x_N) \in K_\alpha$. 
		We say that the family $(x_i)_{1\le i\le N}$ is in general position up to $d'$ if for all $1 \le i_1 < \dots < i_k \le N$ with $k \le d'+1$, the space $\langle x_{i_1}, \dots, x_{i_k} \rangle$ has dimension exactly $k$.
		Let us show that in this case:
		\begin{equation}
			N \le \frac{d'}{\alpha}. \label{n-le-d'/alpha}
		\end{equation}
		To all index $i \in \{1, \dots, N\}$, we associate a random integer variable $a_i:= \mathds{1}_{g x_i = 0} \in \{0,1\}$ \ie such that $a_i = 1$ when $g(x_i) = 0$ and $a_i = 0$ otherwise.
		Note that $\ker(g)$ has dimension at most $d'$ almost surely. 
		As a consequence, for all $\le i_1 < \dots < i_{d'+1} \le N$, we have $\dim\langle x_{i_j} \rangle_{1\le j \le d'+1} > \dim(\ker(g))$ almost surely. Hence $\PP\left(\langle x_{i_j} \rangle \subset \ker(g)\right) = 0$ so $g x_{i_j} \neq 0$ for at least one index $j$.
		This means that, with probability $1$, the random set of indices $\{1\le i \le N\mid g x_i =0\}$ does not admit any subset of size $d'+1$ so its cardinal is at most $d'$. 
		In other words, $\sum_{i=1}^N a_i\le d'$ almost surely.
		Now, note that by definition of $\KK_\alpha$, $\EE(a_i) = \PP(g x_i = 0) \ge \alpha$ for all $i \in\{1,\dots, N\}$. Hence $N \alpha \le \sum_{i=1}^N E(a_i) \le d'$, which proves \eqref{n-le-d'/alpha}.
		
		Now we want to construct an integer $N \le d'/\alpha$ and a family $(x_1, \dots, x_N) \in K_\alpha$ that is in general position up to $d'$ and such that:
		\begin{equation}\label{k-alpha}
			K_\alpha \subset \bigcup_{1\le i_1 < \dots < i_{d'} \le N} \langle x_{i_j} \rangle_{1\le j \le d'}.
		\end{equation}
		We do it by induction. 
		Since we assumed that $K_\alpha \neq \{0\}$, there exists a non-zero vector $x_1 \in K_\alpha$. 
		Let $j \ge 1$ and assume that we have constructed a sequence $(x_1, \dots, x_j) \in K_\alpha^j$ that is in general position up to $d'$. 
		If
		\begin{equation*}
			K_\alpha \subset \bigcup_{1 \le i_1 < \dots < i_{d'} \le j} \langle x_{i_j} \rangle_{1\le j \le d'},
		\end{equation*}
		then we set $N:= j$ and the algorithm ends as \eqref{k-alpha} is satisfied.
		Otherwise, we take:
		\begin{equation*}
			x_{j+1} \in K_\alpha \setminus\left(\bigcup_{1\le i_1 < \dots < i_{d'} \le j} \langle x_{i_1},\dots, x_{i_{d'}} \rangle\right).
		\end{equation*}
		In this case, we have constructed a family $(x_1, \dots, x_{j+1}) \in K_\alpha$ that is in general position up to $d'$ and we can proceed to the next step. 
		This process reaches the  stopping condition after at most $\lfloor d' /\alpha \rfloor$ steps by \eqref{n-le-d'/alpha}.
		This proves that $K_{2^{-m}}$ is included in a countable union of subspaces of dimension exactly $d'$ for all $m$ so $\essker(\nu)$ also is.
		
		To conclude, we show that $\essker(\nu)$ is in fact equal to a countable union of subspaces of dimension at most $d'$. 
		Let $(V_k)_{k\in\NN} \in \mathrm{Gr}_{d'}(E)^{\NN}$ be such that $\essker(\nu) \subset \bigcup V_k$. 
		Let us construct a family $(V_{k_0, \dots, k_{j}})_{j \le d',(k_0, \dots, k_j) \in \NN^{j+1}}$ such that $(V_{k_0})_{k_0 \in\NN} = (V_k)_{k\in\NN}$ and:
		\begin{itemize}
			\item for all $j \le d'$, $\essker(\nu) \subset \bigcup_{(k_0, \dots, k_j) \in \NN^{j+1}} V_{k_0, \dots, k_{j}}$,
			\item for all $j \le d'$ and all multi-index $(k_0, \dots, k_{j})$, $V_{k_0, \dots, k_{j}} \subset V_{k_0, \dots, k_{j-1}}$, with equality if and only if $V_{k_0, \dots, k_{j-1}} \subset \essker(\nu)$. 
		\end{itemize}
		Given such a family,
		\begin{equation}
			\essker(\nu) = \bigcup_{(k_0, \dots, k_{d'}) \in \NN^{d'+1}} V_{k_0, \dots, k_{d'+1}}.
		\end{equation}
		Indeed, for all multi-index $(k_0, \dots, k_{d'}) \in \NN^{d+1}$, such that $V_{k_0, \dots, k_{d'}} \not\subset \essker(\nu)$, we have $V_{k_0} \supsetneq V_{k_0, k_1} \supsetneq \dots \supsetneq V_{k_0, \dots, k_{d'}}$ by assumption. 
		Therefore $\dim(V_{k_0, \dots, k_{d'}}) \le 0$, which is a contradiction.
		
		We construct the family $(V_{k_0, \dots, k_{j}})$ by induction. 
		Let $0 \le c \le d'$. 
		Assume that we have constructed a family $\left(V_{k_0, \dots, k_{j}}\right)_{j \le c, (k_0, \dots, k_j) \in \NN^{j+1}}$, such that for all $j \le c$,
		\begin{equation}
			\essker(\nu) \subset \bigcup_{(k_0, \dots, k_j)\in\NN^j} V_{k_1, \dots, k_{j}},
		\end{equation}
		and such that $V_{k_0, \dots, k_{j}} \subset V_{k_0, \dots, k_{j-1}}$ for all $j \le c$ and for all multi-index $(k_0, \dots, k_j)$, with equality if and only if $V_{k_0, \dots, k_{j-1}} \subset K$.
		
		Let $(k_0, \dots, k_c)$ be a multi-index such that $V_{k_1, \dots, k_{c}} \not \subset \essker(\nu) =\essker(\kappa)$ and let $h \sim \kappa$. 
		Then $h(V_{k_0, \dots, k_{c}}) \neq \{0\}$ with probability $1$ so the restriction of $h$ to $V_{k_0, \dots, k_{c}}$ has rank at least $1$.
		By the previous argument, the set:
		\begin{equation*}
			\essker(\nu) \cap V_{k_0, \dots, k_{c}} = \{x\in V_{k_0, \dots, k_{c}}\mid  \PP(h x = 0) > 0\}
		\end{equation*}
		is included in the union of a countable family of subspaces of $V_{k_0, \dots, k_{c}}$ that have dimension $\dim\left(V_{k_0, \dots, k_{c}}\right) - 1$.
		For all multi-index $(k_0, \dots, k_c)\in\NN^{c+1}$, if $V_{k_1, \dots, k_{c}} \not\subset K$, we set $(V_{k_0, \dots, k_{c+1}} )_{k_{c+1} \in \NN}$ to be such a family; otherwise, we set $V_{k_0, \dots, k_{c+1}}:= V_{k_0, \dots, k_{c}}$.
	\end{proof}

	\subsection{Rank one boundary of a semi-group}
	
	In the present paragraph, we prove Proposition \ref{prop:boundary} and Lemma \ref{lem:contraction}.
	Given a subset $A$ of a topological space $X$, we denote by $\mathbf{cl}_X(A)$ the closure of $A$ in $X$.	
	The results of the present paragraph are well known in the invertible case and used in \cite{AMS-esmigroup} or \cite{aoun2022random} for example.

	\begin{Def}[Rank one boundary and limit sets]
		Let $\Gamma < \End(E)$ be a sub-semi-group. Let $\overline\Gamma:= \mathbf{cl}_{\End(E)}(\KK \Gamma)$. We denote by $\partial_1\Gamma$ the rank-one boundary of $\Gamma$, defined as:
		\begin{equation}
			\partial_1\Gamma:=  \{[\gamma]\mid \gamma\in\overline\Gamma,\; \rank(\gamma) = 1\} 
		\end{equation}
		We define the left and right limit sets of $\Gamma$ as:
		\begin{align*}
			\Lambda_1\Gamma & := \left\{[hv]\,\middle|\,[h] \in\partial_1 \Gamma,\; v \in E \setminus \ker(h)\right\} \subset \proj(E) \\
			\Lambda_1^t\Gamma & := \left\{[fh]\,\middle|\,[h] \in\partial_1 \Gamma,\; f \in E^* \setminus \ker({^th})\right\} \subset \proj(E^*).
		\end{align*}
	\end{Def}
	
	Note that $\overline{\Gamma}$ is also a semi-group because the product is bilinear and continuous.
	
	\begin{Def}[Range and boundary of a distribution]
		Let $E$ be a Euclidean space and let $\nu$ be a probability measure on $\End(E)$. 
		We denote by $\Gamma_\nu$ the range of $\nu$ defined as the smallest closed sub-semi-group of $\End(E)$ that has measure $1$ for $\nu$. 
		We define $\partial_1(\nu):= \partial_1\Gamma_\nu$, $\Lambda_1(\nu):= \Lambda_1\Gamma_\nu$ and $\Lambda_1(^t\nu):= \Lambda_1^t\Gamma_\nu$.
	\end{Def}
	
	With the notations of Proposition \ref{prop:boundary}, $\partial_(\nu) = \proj(\partial_1'(\nu))$, $\Lambda_1(\nu) = \proj(\Lambda_1'(\nu))$ and $\Lambda_1(^t\nu) = \proj(\Lambda_1'(^t\nu))$
	
	Let us now break down the proof of Proposition \ref{prop:boundary} as announced in the introduction.
	
	\begin{Lem}\label{lem:boundary}
		Let $E$ be a Euclidean space and let $\Gamma < \End(E)$ be an irreducible semi-group. Then:
		\begin{equation}
			\partial_1\Gamma = \left\{[uw] \,\middle|\, [u] \in\lambda_1\Gamma, [w] \in\Lambda_1^t\Gamma\right\}.   
		\end{equation}
	\end{Lem}

	\begin{proof}
		Note that the space $\proj\End(E)$ is metrizable so the closure is characterized by sequences. Let $\pi$ be a rank one endomorphism. Saying that $[\pi] \in \partial \Gamma$ is equivalent to saying that there exists a sequence $(\gamma_n) \in\left(\Gamma \setminus \{0\}\right)^\NN$ such that $[\gamma_n] \to [\pi]$. 
		Therefore, for all $[\pi_1], [\pi_2] \in \partial \Gamma$, and for all $g \in \Gamma$ such that $\pi_1 g \pi_2 \neq 0$, we have $[\pi_1 g \pi_2] \in\partial \Gamma$.
		
		Let $v_1 \in\lambda_1\Gamma$ and let $f_2\in\partial_v \Gamma$. Let $v_2 \in E$ and $f_1 \in E^*$ be such that $[v_1 f_1] \in \partial \Gamma$ and $[v_2 f_2] \in \partial \Gamma$. By definition of the irreducibility, there exists an element $g \in \Gamma$ such that $f_1 g v_2 \neq 0$. 
		Given such a $g$, $v_1 f_1 g v_2 f_2 \neq 0$ and $f_1 g v_2$ is a scalar so $[v_1 f_2] = [v_1 f_1 g v_2 f_2] \in \partial \Gamma$.
	\end{proof}

	In the end of the present section, we rely on the elementary fact that $\proj\End(E)$ is metrizable and compact. 
	For the sake of clarity, let us explicit a distance map on $\proj\End(E)$.
	Given $x, y \in \proj\End(E)$, we define:
	\begin{equation}\label{dist-proj-endo}
		\dist(x,y) = \min_{\substack{[g] = x, [h] = y, \\ \|g\| = \|h\|}}\frac{\|g - h\|}{\|g\|}.
	\end{equation}
	
	We remind that in Definition \ref{def:distance}, we characterised the distance on $\proj(E)$ as $\dist(x, y) = \inf_{[u] = x, [v] = y} \frac{\|u - v\|}{\|v\|}$, without the condition $\|u\| = \|v\|$. 
	With this characterisation, the distance is not obviously symmetric and it only is because the Euclidean, Hermitian and ultrametric norms lift canonically to the exterior algebra. 
	
	It turns out that when $E$ is ultrametric the inequality $\frac{\|u - v\|}{\|v\|} < 1$ implies $\|u\|=\|v\|$. 
	Moreover, the space $\End(E)$ is also ultrametric so the distance defined in \eqref{dist-proj-endo} is the one of Definition \ref{def:distance}.
	
	On the other hand, when $E$ is Euclidean or Hermitian, the condition $\|g\| = \|h\|$ in \eqref{distanceprojective} is not redundant and removing it breaks the symmetry.
	In fact, in dimension $d \ge 2$, the space $(\End(E), \|\cdot\|)$ is not Euclidean or Hermitian because its balls are not strictly convex.
	
	\begin{Lem}[Characterisation of proximality]\label{lem:exiprox}
		Let $\nu$ be a probability measure on $\End(E)$. 
		Let $(\gamma_n)\sim \nu^{\otimes\NN}$. 
		Assume that $\nu$ is irreducible. 
		Then the following assertions are equivalent:
		\begin{enumerate}
			\item\label{exiprox:sqz} There exists a constant $b > 0$ such that $\sigma(\overline{\gamma}_n) \ge b$ almost surely and for all $n \ge 0$.
			\item\label{exiprox:empty} $\partial\nu = \emptyset$.
			\item\label{exiprox:noprox} We have $\rho_1(\overline{\gamma}_n) = \rho_2(\overline{\gamma}_n)$ almost surely and for all $n \in\NN$.
		\end{enumerate}
	\end{Lem}
	
	\begin{proof}
		We first prove that point \eqref{exiprox:sqz} implies point \eqref{exiprox:noprox}.
		Assume that there exist $g \in \overline{\Gamma}_\nu$ such that $\rho_1(g) > \rho_2(g)$, which contradicts point \eqref{exiprox:noprox}.
		Then $\sigma(g^k) \to 0$, which contradicts point \eqref{exiprox:sqz}.
		
		We now prove that point \eqref{exiprox:noprox} implies point \eqref{exiprox:empty}.
		The proof relies on the fact that the set $\{g \mid \rho_1(g) = \rho_2(g)\}$ is closed in $\mathrm{End}(E)$ and homogeneous.
		Hence \eqref{exiprox:noprox} implies that $\rho_1(g) = \rho_2(g)$ for all $g \in \overline{\Gamma}_\nu$.
		For all rank one matrix $\pi \in \overline{\Gamma}_\nu$, we have $\rho_1(\pi) = \rho_2(\pi) = 0$ so $\im(\pi) \subset \ker(\pi)$.
		If such a matrix exists, then $\rho_1(g \pi) = 0$ for all $g \in \overline{\Gamma}_\nu$ and therefore $g(\im(\pi)) \subset \ker(\pi)$, which contradicts the irreducibility of $\overline\Gamma_\nu$. 
		It follows that $\partial\nu = \emptyset$.
		
		Finally, we prove that point \eqref{exiprox:empty} implies point \eqref{exiprox:sqz}.
		The set of rank one endomorphisms and $\overline\Gamma_\nu \setminus\{0\}$ are both closed and homogeneous in $\End(E) \setminus \{0\}$ and point \eqref{exiprox:empty} says that they do not intersect, which implies that the projective distance between these sets is positive.
		We remind that $\sigma$ is the distance to the nearest rank one endomorphism in $\proj\End(E)$ so point \eqref{exiprox:empty} implies that it is bounded below by a constant on $\overline\Gamma_\nu \setminus\{0\}$.
		Point \eqref{exiprox:sqz} follows from the fact that $\overline{\gamma}_n \in \overline\Gamma_\nu$ almost surely and for all $n$.
	\end{proof}

	\subsection{Construction of the Schottky measure}
	
	In the present paragraph, we prove Proposition \ref{prop:weak-cont} and Lemma \ref{lem:contraction}.
	We recall that given $0 < \eps \le 1$ and $g, h$ two matrices, we write $g \AA^\eps h$ when $\|gh\| > \eps \|g\|\|h\|$.

	Formally, a measurable binary relation over a measurable space $\Gamma$ is a measurable subset of $\Gamma \times \Gamma$. 
	Given $g,h \in\Gamma$, we write $g \AA h$ to say that $(g,h) \in \AA$. 
	Given $S, T \subset \Gamma$, we write $S \AA T$ to say that $S \times T \subset \AA$.
	Given $n \in\NN$ and $(g_0, \dots, g_n )\in\Gamma^{n+1}$, we write $g_0 \AA \dots \AA g_n$ to say that $g_i \AA g_{i+1}$ for all $i \in\{0, \dots, n-1\}$.
	
	\begin{Def}
		Let $\Gamma$ be a measurable space, let $\AA$ be a measurable binary relation on $\Gamma$ and let $0 \le \rho < 1$. 
		Let $\mu$ be a measure on $\Gamma$. 
		We say that $\mu$ is $\rho$-Schottky for $\AA$ if:
		\begin{equation}\notag
			\forall h \in \Gamma, \mu\{\gamma \in \Gamma \mid  \lnot h \AA \gamma\} \le \rho \quad\text{and}\quad \mu\{\gamma \in \Gamma \mid  \lnot \gamma \AA h\} \le \rho.
		\end{equation}
	\end{Def}
	
	\begin{Rem}
		Let $\AA \subset \AA'$ be two measurable binary relations and let $0 \le \rho \le \rho' < 1$. If a measure $\mu$ is $\rho$-Schottky for $\AA$, then it is also $\rho'$-Schottky for $\AA'$.
	\end{Rem}
	
	To construct a Schottky measure, we approximate the uniform distribution over a set of rank-one matrices in $\partial'_1(\nu)$ that are in general position in the sense that the sum of $k$ distinct matrices in this set has rank $k$ for all $k \le d$.
	
	\begin{Lem}\label{lem:pi-schottky}
		Let $\Gamma < \End(E)$ be a strongly irreducible and proximal semi-group and let $0 < \rho \le 1$. There exist an integer $M\in\NN$, a constant $\eps > 0$ and a family $([\pi_1], \dots, [\pi_M])\in\partial\Gamma^M$ such that:
		\begin{equation}\label{schottky-set}
			\forall h \in\End(E)\setminus\{0\},\; \#\{k\mid \pi_k\AA^\eps h\} \ge (1-\rho) M
			\text{ and }\#\{k\mid h\AA^\eps \pi_k\} \ge (1-\rho) M.
		\end{equation}
	\end{Lem}
	
	\begin{proof}
		Let $d := \dim(E)$. Let $m \in \NN$ and assume that we have constructed a family $([u_1], \dots, [u_m]) \in \Lambda_1\Gamma^{m}$ that is in general position in the sense that $\dim\left(\langle u_{i_1}, \dots, u_{i_k}\rangle \right) = k$ for all $k \le d$, and for all $1 \le i_1 < \dots < i_k \le m$. 
		Let:
		\begin{equation*}
			u_{m+1} \in \{u\in E\mid [u]\in \Lambda_1\Gamma\}\setminus \bigcup_{\substack{k \le d-1,\\ 1\le i_1 < \dots < i_k \le m}} \langle u_{i_j}\rangle_{1 \le j\le k}.
		\end{equation*}
		Such a $u_{m+1}$ exists because the set $\{u\in E\mid [u]\in\Lambda_1\Gamma\}\cup\{0\}$ is $\Gamma$-invariant by definition and not reduced to $\{0\}$ by Lemma \ref{lem:exiprox}. 
		Hence $\{u\in E\mid [u] \in \Lambda_1 \Gamma\}$ cannot be included in $\bigcup_{\substack{i_1<\dots <i_k}} \langle u_{i_j}\rangle_{1 \le j\le k}$, which is a finite union of proper subspaces. 
		Then $[u_{m+1}] \in\partial \Gamma$ by construction and we can easily check that $([u_1], \dots, [u_{m+1}])$ is in general position.
		
		Let $M := \lceil (d-1)/{\rho} \rceil$.
		Let $([u_1], \dots, [u_{M}]) \in \Lambda_1\Gamma^{M}$ and $([w_1], \dots, [w_{M}]) \in \Lambda_1^t \Gamma^{N}$ be in general position. 
		We can construct the family $(u_i)_{1  \le i \le N}$ by induction using the above construction.
		To construct $(w_i)_{1  \le i \le N}$, note that the transpose semi-group $^t\Gamma:= \{{^t\gamma}\mid \gamma \in \Gamma\} \subset \End(E^*)$ is also strongly irreducible and proximal and $\Lambda_1^t \Gamma = \Lambda_1{^t\Gamma}$.
		
		For all $k \in\{1, \dots,M\}$, we define $\pi_{k}:= u_k w_k$. 
		Then $([\pi_1], \dots, [\pi_M]) \in \partial \Gamma^N$ by Lemma \ref{lem:boundary}. 
		Since the $\pi_k$'s have rank one (\cf Section \ref{intro:rank-one}), we check that $g \AA^\eps \pi_{k}$ if and only if $g \AA^\eps u_k$ and that $\pi_{k} \AA^\eps h$ if and only if $w_k \AA^\eps h$. 
		
		Let $h \in\End(E) \setminus\{0\}$. 
		Let $I^h:= \{i\mid  u_i \in\ker(h)\}$. 
		The family $([u_i])_{1 \le i \le M}$ is in general position and $\langle u_i\rangle_{i\in I} \subset \ker(h)$ so $\# I^h \le d-1$. 
		Let $J^h := \{j \mid  w_j \in \ker({}^t h)\}$. 
		By the same argument, $\# J^h \le d - 1$.
		
		For all $h \in \End(E) \setminus \{0\}$, let:
		\begin{equation*}
			\psi(h):= \max_{\substack{I, J \subset \{1, \dots, M\},\\\#I\le d-1, \#J \le d-1}} \min\left\{\min_{i \notin I} \frac{\|h u_i\|}{\|h\|\|u_i\|},\, \min_{j \notin J} \frac{\|w_j h\|}{\|w_j\|\|h\|}\right\}.
		\end{equation*}
		By the previous argument, $\psi(h) > 0$ for all $h$. 
		Moreover the maps $h\mapsto \frac{\|h v\|}{\|h\|\|v\|}$ and $\frac{\|f h\|}{\|f\|\|h\|}$ are continuous for all $v \in E\setminus\{0\}$ and all $f\in E^* \setminus\{0\}$. 
		It follows that $\psi$ is continuous. 
		Moreover, $\psi$ is invariant by scalar multiplication so there exists a unique continuous map $\phi: \proj\End(E) \to (0,1]$ such that $\phi([h]) = \psi(h)$ for all $h \in\End(E) \setminus\{0\}$. 
		The projective space $\proj\End(E)$ is compact so $\phi (\proj\End(E))$ is compact in $(0, 1]$. Let $\eps > 0$ be its lower bound.
		
		Let $h \in\End(E)\setminus\{0\}$. 
		The set of indices $\{k \mid  \pi_k \AA^{\eps} h\}$ is $\{k \mid  w_k \AA^{\eps} h\}$, which has cardinality at least $M - d + 1 \ge (1- \rho) M$. 
		Indeed, given a subset $J$ that realises the maximum in the definition of $\psi(h)$, since $\psi(h) \ge \eps$, we have $w_j \AA^{\eps} h$ for all $j \in \{1, \dots, M\} \setminus J$.
		By the same argument $\#\{k \mid  h\AA^{\eps} \pi_k\} \ge M - d + 1 \ge (1- \rho) M$, which yields \eqref{schottky-set}.
	\end{proof}

	\begin{Cor}\label{cor:S'-schottky} 
		Let $\nu$ be a strongly irreducible and proximal probability measure on $\End(E)$, let $0 < \rho \le 1$ and let $\delta: (0,1) \to (0,1)$. There exist an integer $N$, two constants $\alpha', \eps \in (0,1)$, a family $(n_{k})_{1\le k \le N} \in\NN^{N}$ and a family $(S'_{k})_{1\le k \le N}$ of measurable subsets of $\Gamma$ such that $\nu^{*n_{k}} (S'_{k}) \ge \alpha'$ for all $1\le k \le N$ and such that:
		\begin{gather}
			\forall h \in\End(E),\; \#\{k\mid S'_k\AA^{2\eps} h\} \ge (1-\rho) N
			\text{ and }\#\{k\mid h\AA^{2\eps} S'_k\} \ge (1-\rho) N, \label{item:align}
			\\
			\forall 1 \le j \le N,\; \#\{k\mid S'_k\AA^{\eps} S'_j\} \ge (1-\rho) N
			\text{ and }\#\{k\mid S'_j\AA^{\eps} S'_k\} \ge (1-\rho) N, \label{item:self-align}
			\\
			\forall 1 \le k \le N,\; \forall s' \in S'_k, \; \sigma(s') \le \delta(\eps).\label{sigma-delta}
		\end{gather}
	\end{Cor}

	\begin{proof}
		Let $N \in \NN$, let $\eps > 0$ and let $([\pi_1], \dots, [\pi_N])\in\partial\nu^N$ be such that:
		\begin{equation*}
			\forall h \in\End(E)\setminus\{0\},\; \#\{k\mid \pi_k\AA^{3\eps} h\} \ge (1-\rho) N
			\text{ and }\#\{k\mid h\AA^{3\eps} \pi_k\} \ge (1-\rho) N.
		\end{equation*}
		Such a family exists by Lemma \ref{lem:pi-schottky}. 
		To all $k \in\{1, \dots, N\}$, we associate the open set:
		\begin{equation*}
			S'_k:= \left\{ s \in\End(E)\setminus\{0\} \,\middle|\, \dist([s], [\pi_k]) < \eps ,\;\sigma(s) < \delta(\eps)\right\},
		\end{equation*}
		where $\dist$ is defined by \eqref{dist-proj-endo}.
		By definition of $S'_k$, \eqref{sigma-delta} holds.
		
		To prove \eqref{item:align} let $k \in\{1, \dots, N\}$ and $h \in \End(E) \setminus\{0\}$. 
		Assume that $h \AA^{3\eps} \pi_k$.
		We claim that $h \AA^{2 \eps} s'$ for all $s' \in S'_k$. 
		By homogeneity, we may assume that $\|s' - \pi_k\| < \eps \|s'\| = \eps \|\pi_k\|$ so $\|hs' - h\pi_k\| < \eps \|h\|\|\pi_k\|$. 
		By assumption, $\left\|h \pi_k\right\| > 3 \eps \|h\|\|\pi_k\|$ and by triangular inequality $\left\|h s\right\| > 2 \eps \|h\|\|\pi_k\| = 2 \eps \|h\|\|s\|$. 
		By duality, the same reasoning works for the left alignment, which yields \eqref{item:align}.
		
		Now let $j,k \in\{1, \dots, N\}$ be such that $\pi_j \AA^{3\eps}\pi_k$ and let $s'_j \in S'_j$ and $s'_k \in S'_k$. 
		By the above argument, $s'_j \AA^{2\eps} \pi_k$ and $s'_j \AA^\eps s'_k$, which yields \eqref{item:self-align}.
		
		For all $k \in \{1, \dots, N\}$, the set $S'_k$ is an open neighbourhood of $\pi_k$ so it intersects $\overline\Gamma_\nu = \mathrm{cl}\left(\bigcup_{n=0}^\infty \KK \supp(\nu^{*n})\right)$. 
		By characterisation of the closure, $S'_k$ intersects $\bigcup_{n=0}^\infty \KK \supp(\nu^{*n})$
		Hence, there exists an integer $n_k$ such that $S'_k$ intersects $\KK \supp(\nu^{*n_k})$.
		Moreover $S'_k$ is invariant under non-trivial scalar multiplication so it intersects $\supp(\nu^{*n_k})$ and $\nu^{n_k}(S'_k) > 0$ by characterisation of the support. 
		Finally, pick such a family $(n_k) \in \NN^N$ and set $\alpha' := \min_k \nu^{*n_k}(S'_k)$. 
	\end{proof}

	Let us now use the local-to-global transmission of the alignment to prove that the $(n_k)$'s of Corollary \ref{cor:S'-schottky} can be picked to be all the same.

	\begin{Lem}\label{lem:schottky}
		Let $\nu$ be a strongly irreducible and proximal probability measure on $\End(E)$, let $0 < \rho < 1$ and $f = (0,1) \to (0,1)$. 
		There exist an integer $N$, two constants $\alpha, \eps \in (0,1)$, an integer $m$ and a family $(S_{k})_{1\le k \le N}$ of measurable subsets of $\End(E)$ such that $\nu^{*m} (S_{k}) \ge \alpha$ for all $1\le k \le N$ and such that:
		\begin{gather}
			\forall h \in\End(E),\; \#\{k\mid S_k\AA^{\eps} h\} \ge (1-\rho) N
			\text{ and }\#\{k\mid h\AA^{\eps} S_k\} \ge (1-\rho) N,\label{item:Schottky}\\
			\forall 1 \le k \le N,\;\forall s \in S_k, \; \sigma(s) \le \delta(\eps)\label{item:contract}.
		\end{gather}
	\end{Lem}
	
	\begin{proof}
		Without loss of generality, we may assume that $\rho < \frac{1}{3}$ and that $\delta(\eps) \le \eps^2/12$.
		Let $N \in\NN$, let $\alpha', \eps \in (0,1)$, let $(n_k)$ and let $(S'_k)$ be as in Corollary \ref{cor:S'-schottky}. To all index $k\in\{1, \dots, N\}$, we associate two indices $i_k, j_k \in \{1, \dots, N\}$ such that $S'_k \AA^\eps S'_{i_k}$, $S'_{i_k} \AA^\eps S'_{j_k}$ and $S'_{j_k} \AA^\eps S'_{k}$. 
		By \eqref{item:self-align} in Corollary \ref{cor:S'-schottky}, such indices $i_k,j_k$ exist because:
		\begin{equation*}
			\#\{(i,j) \,| \,S'_{i}\AA^\eps S'_k \AA^\eps S'_{i}\AA^\eps S'_{j} \AA^\eps S'_{k}\AA^\eps S'_{j}\}\ge (1 - 2\rho) N (1-3\rho) N > 0.
		\end{equation*}
		Hence the set of all possible values for $i_k, j_k$ is non-empty. 
		Let $m'$ be the smallest common multiple of $\{n_k + n_{i_k}\mid  1 \le k \le N\}$, let $m''$ be the smallest common multiple of $\{n_k + n_{j_k}\mid  1 \le k \le N\}$ and let $m = m' + m''$. Let $k \in \{1, \dots, N\}$. We define $p_k:= \frac{m'}{n_k + n_{i_k}}$ and $q_k:= \frac{m''}{n_k + n_{j_k}}$ and:
		\begin{equation*}
			S_k:= (S'_k \cdot S'_{i_k})^{\cdot p_k}\cdot(S'_{j_k} \cdot S'_k)^{\cdot q_k} = \Pi\left( (S'_k \times S'_{j_k})^{p_k}\times(S'_{j_k} \times S'_k)^{q_k}\right).
		\end{equation*}
		Therefore, 
		\begin{equation*}
			\nu^{*m}(S_k) \ge \left(\nu^{n_k}(S'_k)\nu^{n_{i_k}}(S'_{i_k})\right)^{p_k}\left(\nu^{n_{i_k}}(S'_{i_k})\nu^{n_{j_k}}(S'_{k})\right)^{q_k} \ge \alpha'^{2 p_k + 2 q_k} \ge \alpha'^{2m} =: \alpha.
		\end{equation*}
		
		Let $h \in \End(E) \setminus \{0\}$ and $k\in\{1, \dots, N\}$ and assume that $h \AA^{2\eps} S'_k$.
		Since $\delta(\eps) \le \eps^2/12$, it follows from Lemma \ref{lem:alipart} that $h \AA^{\eps} S_k$. 
		If we instead pick $h$ such that $S'_k \AA^{2\eps} h$, the same argument yields that $S_k \AA^\eps h$. 
		Hence, we can derive \eqref{item:Schottky} from \eqref{item:align}.
		
		Moreover, by Lemma \ref{lem:c-prod} and because we assumed that $\delta(\eps) \le \eps^2 / 8$, for all $k$ and for all $s \in S_k$, it holds that $\sigma(s) \le \delta(\eps)^{2 p_k + 2 q_k} 2^{4 p_k + 4 q_k - 4} \eps^{- 4 p_k - 4 q_k +2} \le \delta(\eps)$, which proves \eqref{item:contract}.
	\end{proof}
	
	Now let us prove the main result of this section.
	In fact, for the rest of the discussion, we will not directly use the formal definitions of strong irreducibility and proximality.
	Instead, we simply think of a strongly irreducible and proximal measure as a measure that satisfies the conclusions of Corollary \ref{cor:schottky} for $\rho = \frac{1}{6}$ and $\delta: \eps \mapsto \eps^6/48$.
	
	\begin{Cor}\label{cor:schottky}
		Let $E$ be a Euclidean space, let $\Gamma = \End(E)$ and let $\nu$ be a strongly irreducible and proximal probability measure on $\Gamma$, let $0 < \rho < 1$ and let $\delta: (0,1) \to (0,1)$. 
		There exists an integer $m$, two constants $\alpha, \eps \in (0,1)$ and a probability measure $\tilde\mu$ on $\Gamma^m$ such that:
		\begin{enumerate}
			\item The measure $\Pi_*\tilde\mu$ is $\rho$-Schottky for $\AA^\eps$.
			\item The measure $\tilde\mu$ is absolutely continuous with respect to $\nu^{\otimes m}$ in the sense that $\alpha \tilde\mu \le \nu^{\otimes m}$.
			\item We have $\sigma_*\Pi_*\tilde\mu\left[0, \delta(\eps)\right] = 1$ \ie $\sigma(s_1 \cdots s_m) \ge \delta(\eps)$ for all $(s_1, \dots, s_m)$ in the support of $\tilde\mu$.
		\end{enumerate}
	\end{Cor}
	
	\begin{proof}
		Let $m,N\in\NN$, $\alpha,\eps\in(0,1)$ and $(S_k)_{1 \le k\le N}$ be as in Lemma \ref{lem:schottky}. 
		Define $f: \Gamma^m \to \RR_{\ge 0}$ as:
		\begin{equation*}
			f:= \sum_{k = 1}^{N}\frac{\mathds{1}_{S_k} \circ \Pi}{\nu^{*m}(S_k)}.
		\end{equation*}
		Then $f \le \frac{N}{\alpha}$ because we assumed that $\nu^{*m}(S_k) \ge \alpha$ for all $k$. 
		Moreover, we can check that $\int_{\Gamma^m} f d\nu^{\otimes m} = N$ and that $\int_{\Pi^{-1}(S_k)} f d\nu^{\otimes m} \ge 1$ for all $k \in\{1, \dots, N\}$. 
		Let:
		\begin{equation*}
			\tilde\mu:= \frac{f \nu^{\otimes m}}{\int_{\Gamma^m} f d\nu^{\otimes m}}\quad\text{and}\quad \mu:= \Pi_*\tilde\mu.
		\end{equation*} 
		Then $\frac{\alpha}{N} \mu \le \nu^{*m}$ by definition. 
		Moreover $\mu\left(\bigcup_{i\in I}S_i\right)\ge \frac{\# I}{N}$, for all $I \subset \{1, \dots, N\}$.
		It follows from \eqref{item:Schottky} that $\mu$ is $\rho$-Schottky. 
		Moreover $\mu$ is supported on $\bigcup S_k$ and $\sigma(S_k)\subset \left[0, \delta(\eps)\right]$ for all $k$.
		As a consequence $\sigma_*\mu\left[0, \delta(\eps)\right] = 1$. 
	\end{proof}

	Corollary \ref{cor:schottky} directly implies Lemma \ref{lem:contraction}, its strength is that is gives us a compact set $S$ such that $\tilde\mu(S^m) = 1$, this allows us to capitalize on point \eqref{pivot:indep-tail} in Theorem \ref{th:pivot}. 
	
	Let us now prove Proposition \ref{prop:weak-cont}, which states that the space of $(\alpha, \eps, m)$ contracting probability measures is weak-$*$ open.
	We remind that the weak-$*$ topology on a metrizable space $X$ is the dual of the $\ell^\infty$ topology on continuous functions with compact support.
	On $\mathrm{Prob}(X)$, this topology is the the dual of the $\ell^\infty$ topology on bounded continuous functions, which is equivalent to the bounded Wasserstein topology, generated by the distance $\dist(\mu, \nu) = \inf_{x \sim \mu, y\sim \nu} \EE(\dist(x, y))$ for $\dist : X \times X \to [0,1]$ any bounded distance that generates the topology of $X$.
	
	Visually, a small transformation in the weak-$*$ topology consists in slightly moving all of the mass, which is what a numerical approximation does for example, and moving a small amount of mass arbitrarily far.  
	Statistically, given $(\gamma_n)_{n \ge 0} \sim \nu^{\otimes\NN}$, the sampling measure $\sum_{k = 0}^n \delta_{\gamma_k} / n$ converges almost surely to $\nu$ in the weak-$*$ topology so it is nice for contraction result to be robust under small transformations in this topology.
	Moments on the other hand are not weak-$*$ continuous so predictions for the law of large numbers made at a time $n$ for example, may not hold at a time much larger than $n$ if there are rare events (in the sense $\PP << 1/n$) that change the value of the expectation.

	An example is the random walk on $\ZZ$ with steps of law $\eta$, where $\eta$ gives weight $2/n!$ to $(-1)^n n!$ for $n > 2$, gives weight $1/2$ to $2$ and gives the remaining weight to $0$. 
	For $n$ large and at a time $n! << k << (n+1) !$, we do not see any steps of size larger that $n!$ and see many copies of each step of size at most $n!$. 
	So we can safely affirm that we are looking at a random walk with bounded steps and with bias equal to $(-1)^n$. 
	However, once we zoom out by a factor $n$, we see a random walk with opposite bias.
	
	The sequence $(\log\|\overline\gamma_n\|)_n$ may exhibit this type of multi-scale behaviour.
	For $(-\log\sigma(\overline\gamma_n))_n$ on the other hand, a change of scale can only increase the perceived escape speed.

	\begin{proof}[Proof of Proposition \ref{prop:weak-cont}]
		Let us fix $\alpha, \eps > 0$ and $m \in \NN$ and pick a measure $\nu$ on $\End(E)$ that is $(\alpha, \eps, m)$-contracting.
		We look for a weak-$*$ neighbourhood of $\nu$ that only contains $(\alpha, \eps, m)$-contracting measures.
		
		First, let us show that the space of measures that are $1/6$-Schottky for $\AA^\eps$ is weak-$*$ open.
		First, note that for a given $g \in \End(E)$ and for $A_g = \{\gamma\mid g \AA^\eps \gamma\}$, the condition $\mu(A_g) > 5/6$ is equivalent to saying that there exists a continuous map $\psi_g$ that takes values in $[0,1]$ and whose support is compact and included in the open set $\{\gamma\mid g \AA^\eps \gamma\}$ such that $\int \psi_g d \mu > 5/6$ (which is a weak-$*$ open condition).
		Given $\mu$ such that $\mu(A_g) > 5/6$ by $\sigma$-compactness of $A_g$, there exists a compact $K \subset A_g$ such that $\mu(K)> 5/6$ and since $A_g$ is open, there exists a compact $S \subset A_g$ such that the interior of $S$ contains $K$.
		Let us set $\psi_g$ to be equal to $1$ on $K$ and supported on $S$, for example $\psi_g(x) = 1- \dist(x, K \cup S^c)/\dist(x, S^c)$. 
		Then $\int\psi_g d\mu > \mu(K) > 5/6$ and $\mu'(S) > 5/6$ for all measure $\mu'$ such that $\int\psi_g d\mu' > 5/6$. Therefore $\mu'(A_g) > 0$.
		
		Given such a compact $S \subset A_g$, we have $\inf_{s \in S} \frac{\|gs\|}{\|g\|\|s\|} = \min \inf_{s \in S} \frac{\|gs\|}{\|g\|\|s\|} > \eps$ by continuity of the map $s \mapsto \frac{\|gs\|}{\|g\|\|s\|}$ on $\End(E) \setminus\{0\}$.
		It follows from the equi-continuity of the family of maps $g' \mapsto  \frac{\|g's\|}{\|g'\|\|s\|}$ for $s \in S$ that $S \subset A_{g'}$ for all $g'$ in a neighbourhood of $g$.
		By homogeneity of the alignment relation, $S \subset A_{g'}$ for all $g'$ in an homogeneous neighbourhood of $g$.
		By compactness of $\proj\End(E)$, there exists a finite family $\psi_1, \dots, \psi_N$ of continuous and compactly supported maps such that for all $g \in \End(E)\setminus\{0\}$, there exists an index $i$ such that $\supp(g_i) \subset A_g$ and $\int \psi_i d\mu > 5/6$ for all $i$.
		
		Then all measure $\mu'$ such that $\int \psi_i d\mu' > 5/6$ for all $i$ is $1/6$-Schottky for $\AA^\eps$.
		Moreover, the set of such measures is a finite intersection of weak-$*$ open sets and is therefore weak-$*$ open.
		
		Let us now fix a measure $\mu < \nu^{*m}$ that is $1/6$-Schottky.
		Let $(\gamma_0,\dots, \gamma_{m-1}) \sim \nu^{\otimes m}$ and let $A$ be an event such that $\PP(A) > \alpha$ and the conditional distribution of $\overline\gamma_m$ with respect to $A$ is $\mu$. 
		Let $(\gamma'_0,\dots, \gamma'_{m-1})$ be a random sequence.
		Let $\delta > 0$ be such that $\PP(\|\gamma_k - \gamma'_k\| \ge \delta) < \delta$, we want to exhibit a condition on $\delta$ that guarantees that the conditional law of $\overline\gamma'_m$ with respect to $A$ is $1/6$-Schottky for $\AA^\eps$.
		
		Let $\delta'$ be such that for $s \in \mu$ and for all $s'$ such that $\PP(\|s' - s\| \ge \delta') < \delta$, the law of $s'$ is $1/6$-Schottky for $\AA^\eps$. 
		Such a $\delta'$ exists because the set of measures that are $1/6$-Schottky for $\AA^\eps$ is weak-$*$ open.
		
		With probability at least $1 - m \delta$ the inequality $\|\gamma_k - \gamma'_k\| < \delta$ holds for all $k$.
		Therefore, $\|\overline\gamma_m - \overline\gamma'_m\| < m \delta \prod_{k = 0}^{m-1}(\|\gamma_k\| + \delta)$.
		Let $B$ be such that $\PP(\|\gamma_k\| \ge B) < \alpha \delta' / 2m$. 
		Then with probability at least $1- m \delta - \alpha\delta'/2$, we have $\|\overline\gamma_m - \overline\gamma'_m\| \le m \delta (B + \delta)^m$.
		Therefore,
		\begin{equation*}
			\PP(\|\overline\gamma_m - \overline\gamma'_m\| > m \delta (B + \delta)^m\mid A ) < m \delta / \alpha + \delta'/2.
		\end{equation*}
		Let us now assume that $0 < \delta \le \alpha \delta' / (2m)$ and that $m \delta (B + \delta)^m \le \delta'$, then we get that $\PP(\|\overline\gamma_m - \overline\gamma'_m\| \le \delta'\mid A ) < \delta'$ so the conditional distribution of $\overline\gamma'_m$ with respect to $A$ is $1/6$-Schottky for $\AA^\eps$.
		
		We deduce that all measure $\nu'$ that admits a coupling $\gamma \sim \nu$, $\gamma'\sim \nu'$ such that $\PP(\|\gamma- \gamma'\| \ge \delta) < \delta$, is also $(\alpha, \eps, m)$-contracting.
		So the space of $(\alpha, \eps, m)$-contracting measures contains a weak-$*$ neighbourhood of $\nu$.
		This reasoning holds for all $(\alpha, \eps, m)$-contracting measure $\nu$ so the space of $(\alpha, \eps, m)$-contracting measures is weak-$*$ open in $\mathrm{Prob}(\End(E))$.
	\end{proof}
	
	With slightly more advanced probabilistic constructions, we can push the proof a bit further and show that:
	 
	\begin{Prop}\label{prop:weak-bounded}
		Let $\alpha, \eps > 0$, let $m \in \NN$ and let $B > 0$.
		The space of measures $\nu$, supported on $\End(E)$ such that there exists a measure $\tilde\mu < \nu^{\otimes m} / \alpha$ satisfying:
		\begin{itemize}
			\item $\mu = \Pi_*\tilde\mu$ is $1/6$-Schottky for $\AA^\eps$ and
			\item $\max N(s_k) < B$ for almost every $(s_0, \cdots,s_{m-1}) \sim \mu$,
		\end{itemize}
		is open in the weak-$*$ topology.
	\end{Prop} 
	
	\begin{proof}
		With the same notations as in the proof of Proposition \ref{prop:weak-cont}, we replace $A$ with the event:
		\begin{equation*}
			A' = A \cap (\forall 0 \le k < m, N(\gamma_k) < B - \delta).
		\end{equation*}
		By assumption, $N(\gamma_k) < B$ almost surely on $A$, so for $\delta$ small enough, we can assume that $\PP(A') > \alpha$ and we get that:
		\begin{equation*}
			\PP(\|\overline\gamma_m - \overline\gamma'_m\| > \delta' \mid A' ) < \delta'.
		\end{equation*} 
		Now the conditional distribution $\mu'$ of $\overline\gamma_m$ with respect to $A'$ is at distance at most $1- \PP(A') / \PP(A)$ from $\mu$.
		To construct a coupling, set $g := \overline\gamma_m$ on $A'$ and set $g \sim \mu'$ on $A \setminus A'$.
		Then $\PP(\|g - \overline\gamma_m\| > 0 \mid A) \le 1- \PP(A') / \PP(A)$ and the conditional law of $g$ with respect to $A$ is equal to the conditional law of $\overline\gamma_m$ with respect to $A'$.
		Let $h$ be a random variable defined on $A'$ such that $\PP(\|h - \overline\gamma_m\| > 0 \mid A') \le 1 - \PP(A') / \PP(A)$ and such that the conditional law of $h$ with respect to $A'$ is equal to conditional law of $\overline\gamma_m$ with respect to $A$.
		Such a $h$ can be constructed by taking a measure preserving map $H : A \to A'$ such that $H_* g = \overline \gamma_n$ almost surely and setting $h = H_* \overline\gamma_n$.
		So if we assume that $1-\PP(A') / \PP(A) < \delta' - \PP(\|\overline\gamma_m - \overline\gamma'_m\| > \delta' \mid A' )$, we get:
		\begin{equation*}
			\PP(\|h - \overline\gamma'_m\| > \delta' \mid A' ) < \delta'.
		\end{equation*}
		Moreover $h \sim \mu$ so the conditional distribution of $\overline\gamma'_m$ with respect to $h$ is $1/6$-Schottky for $\AA^\eps$ and by construction $\max N(\gamma'_k) < B$ on $A'$.
	\end{proof}
	
	The space of probability measures that have positive rank which is closed and has empty interior in $\mathrm{Prob}(\End(E))$ so the assumption $\evrank(\nu) > 0$, though crucial in applications, is not statistically robust.

	\section{Construction of a globally aligned extraction}\label{sec:pivot}

	The present section is dedicated to the proof of Theorem \ref{th:pivot}.
	Let $\alpha$, $m$, $\rho$, $\eps$ and $\mu$ be as in Corollary \ref{cor:schottky} for $\delta: \eps \mapsto \eps^6/32$ and let $\AA = \AA^\eps$.
	
	The idea starts with the following observation.
	Since we have constructed a measure $\tilde\mu\le \alpha^{-1} \nu^{\otimes m}$, then $\nu^{\otimes m}$ is a convex combination of $\tilde\mu$ and another probability measure, that we denote by $\tilde\kappa$.
	This means that we can construct a random sequence $(\gamma_n)\sim \nu^{\otimes\NN}$ by alternating randomly between independent $\tilde\mu$-distributed blocks of length $m$ and $\tilde\kappa$-distributed blocks of length $m$.
	
	From such a decomposition, we construct another decomposition that alternates non-randomly between random blocks of random length and random $\tilde\mu$ distributed blocks.
	These blocks are globally independent. 
	Using the properties of $\mu:= \Pi_*\tilde\mu$, we also make sure that each block have the inner structure of being a concatenation of three blocs $\tilde{u}_0 \odot \tilde{s}_1 \odot \tilde{u}_2$ with $u_0 \AA s_1 \AA u_2$ (removing the $\tilde{}$ means taking the product) and $s_1$ is in the support of $\mu$, which means that $\sigma(s_1)$ is small enough to apply Lemma \ref{lem:triple-ali}.
	
	Starting from this sequence alternating between words whose product are close to rank one matrices and $\tilde\mu$ distributed words, we use an inductive construction to create a sequence of blocks whose products form an aligned sequence. 
	The issue is that with this construction the lengths of the blocks are not stopping times and the blocks are not independent.
	
	A crucial technical point is to make sure that each step of the construction is independent of the contents of each $\tilde\kappa$ distributed block. 
	Allowing us to have control over the tail of the distribution of each block and therefore over the moments.
	
	In the present section, we use the following notations for usual probability distributions.
	Given a parameter $\alpha \in (0,1)$, we write $\mathcal{B}_\alpha:= \alpha \delta_1 + (1-\alpha) \delta_0$ for the Bernoulli law of parameter $\alpha$ and $\mathcal{G}_\alpha:= \sum_{k\ge 0}{\alpha^k}(1-\alpha)\delta_k$.
	We write $\mathcal{U}$ for the uniform probability distribution on $[0,1]$ \ie the Lebesgue measure.
	Given an integer $m$, we will write $m\times_*$ and $m+_*$ for the push forward operators associated to the multiplication by $m$ and the addition of $m$ respectively.
	
	\begin{Def}[Power convolution]\label{def:power-convol}
		Let $\Gamma$ be a metric semi-group.
		We write $*$ for the convolution product of probability distributions on $\Gamma$.
		Let $\nu$ be a probability distribution on $\Gamma$. 
		Let $\eta$ be a probability distribution on $\NN$.
		We define the probability distributions $\nu^{*\eta}$ on $\Gamma$ and $\nu^{\otimes\eta}$ on $\widetilde{\Gamma}$ by:
		\begin{equation}
			\nu^{*\eta} = \sum_{k\in\NN} \eta\{k\} \nu^{*k} \quad\text{and}\quad\nu^{\otimes \eta} = \sum_{k\in\NN} \eta\{k\} \nu^{\otimes k}.
		\end{equation}
	\end{Def}
	
	In other words, given an i.i.d. sequence of integers, matrices or words $(\gamma_n)_{n\in\NN} \sim \nu^{\otimes\NN}$ and a random sequence of integers $(w_k)_{k \ge 0} \sim \eta^{\otimes \NN}$ that are independent, and using the notation of Definition \ref{def:power-extract}, we have $\left(\widetilde{\gamma}^w_k\right)_{k \ge 0} \sim \left(\nu^{\otimes \eta}\right)^{\otimes\NN}$ and $\left({\gamma}^w_k\right)_{k \ge 0} \sim \left(\nu^{*\eta}\right)^{\otimes\NN}$.
	Let us state a few straightforward results to understand the notion of power convolution.
	
	\begin{Lem}\label{lem:power-compare}
		Let $\nu$ be a probability distribution son $(\RR_{\ge 0}, +)$ and $\eta$ be a probability distributions on $(\NN, +)$.
		Let $(x_n)_{n \ge 0}$ be a random sequence and let $w$ be a random non-negative integer.
		Assume that $\PP(w \ge k) \le \eta(\NN_{\ge k})$ for all $k \ge 0$ and assume that $\PP(x_n > t\mid x_0, \dots, x_{n-1}, w) \le \nu(t, + \infty)$ almost surely, for all $t \ge 0$ and for all $n$.
		Then $\PP\left(\overline{x}_{w} > t\right) \le \nu^{*\eta}(t, + \infty)$ for all $t \ge 0$.
	\end{Lem}
	
	\begin{proof}
		Here is an example of proof using a monotonous coupling.
		Let $\left(u_k\right)_{k \ge -1} \sim \mathcal{U}^{\otimes\NN}$. 
		Let $p = \min\{k \in \NN_{\ge 0}\mid \eta(\NN_{\ge k}) \le u_{-1}\}$.
		For all $n \ge 0$, write $y_n  = \min\{t \ge 0\mid \nu(t, +\infty) \le u_k\}$.
		That way $\overline{y}_p \sim \nu^{*\eta}$.
		Let $w' = \min\{k \in \NN_{\ge 0}\mid  \PP(w \ge k) \le u_{-1}\}$. 
		Then $w'$ has the same distribution law as $w$ and $w' \le p$ for all values of $\left(u_k\right)_{k \ge -1} \in (0,1)^\NN$. 
		We define a sequence $(x'_n)_{n \ge 0}$ that has the same law as $(x_n)_{n \ge 0}$ by induction on $n$ with the following formula:
		\begin{equation*}
			\forall n \in\NN,\, x'_n = \min\left\{t \ge 0\,\middle|\, \PP(x_n > t\mid \forall k < n, x_{k} = x'_{k} \cap w = w') \le u_k\right\}.
		\end{equation*} 
		that way, $\PP(x'_n > t\mid x'_0, \dots, x'_{n-1}, w') = \PP(x_n > t\mid x_0, \dots, x_{n-1}, w)$ almost surely and for all $t \ge 0$ so the joint distribution of $\left(w',(x'_n)_{n \ge 0}\right)$ is equal to the joint distribution of $\left(w,(x_n)_{n \ge 0}\right)$.
		
		Moreover, $x'_n \le y_n$ for all $n$ so $\overline{x}'_n \le \overline{y}_n$ for all $n \ge 0$. 
		Moreover the sequences  $(\overline{x}'_n)$ and $(\overline{y}_n)$ are both non-decreasing so $\overline{x}'_{w'} \le \overline{y}_p$. 
	\end{proof}
	
	The following Lemmas illustrate how the notion of finite exponential moment behaves under power convolution.
	
	\begin{Lem}
		For all $0 < \alpha, \beta < 1$, $(1+_*\mathcal{G}_{1-\alpha})^{*(1+_*\mathcal{G}_{1-\beta})} = (1+_*\mathcal{G}_{1-\alpha \beta})$.
	\end{Lem}
	
	\begin{proof}
		Let $(a_n)_{n \ge 1},(b_n)_{n \ge 1}\sim\mathcal{B}_\alpha^{\otimes\NN} \otimes \mathcal{B}_\beta^{\otimes\NN}$. Then $(a_n b_n) \sim \mathcal{B}_{\alpha \beta}^{\otimes\NN}$
		For all $n \ge 1$, let $\overline{x}_n$ be the $n$-the element of $\{n\mid  a_n = 1\}$. 
		Then $(x_n)_{n \ge 0} \sim (1+_*\mathcal{G}_{1-\alpha})^{\otimes \NN}$. 
		Moreover, by independence, $(b_{\overline{x}_n})\sim \mathcal{B}_\beta^{\otimes\NN}$ and $(b_{\overline{x}_n})_{n \ge 1}$ is independent of $(x_n)_{n \ge 0}$. 
		Let $y \ge 1$ be the smallest integer such that $b_{\overline x_y } = 1$.
		Then $y \sim (1+_*\mathcal{G}_{1 - \beta})$ and $y$ is independent of $(x_n)$.
		So $\overline x_y \sim (1+_*\mathcal{G}_{1 - \alpha})^{*(1+_*\mathcal{G}_{1 - \beta})}$ and $\overline x_y $ is also the first index $n \ge 1$ such that $a_n b_n = 1$ so $\overline x_y \sim (1+_*\mathcal{G}_{1- \alpha \beta})$.
	\end{proof}
	
	\begin{Lem}\label{lem:sumexp}
		Let $\nu$ be a probability distribution on $\RR$ and $\eta$ be a probability distribution on $\NN$. Assume that $\nu$ and $\eta$ both have a finite positive exponential moment, then $\nu^{*\eta}$ also has.
		Let $w \sim \eta$, $x\sim \nu$ and $y \sim \nu^{*\eta}$. 
		Given $\alpha, \beta> 1$ such that $\EE(\beta^x) \le \alpha$, it holds that $\EE(\beta^{y})\le \EE(\alpha^w)$.
		Hence, there exists $\beta > 1$ such that $\EE(\beta^{y}) < +\infty$.
	\end{Lem}
	
	\begin{proof}
		We first show the second result. 
		Let $(w,(x_n)_{n \ge 0}) \sim \eta\otimes\nu^{\otimes\NN}$ and let $\alpha, \beta > 1$.
		Assume that $\EE(\beta^{x_n}) \le \alpha$ for all $n$.
		Then by independence, $\EE(\beta^{\overline{x}_w}) = \sum_k \PP(w = k) \EE(\beta^{\overline{x}_k}) \le \sum_k \PP(w = k) \alpha^k = \EE(\alpha^w)$.
		
		To conclude, note that the fact that $\eta$ and $\nu$ have a finite exponential moment means that there exists $\alpha, \beta_0 > 1$ such that $\EE(\alpha^w) < + \infty$ and $\EE\left(\beta_0^{x_0}\right) < +\infty$.
		By Jensen's inequality and by concavity of the power function, $\EE(\beta^{x_0}) \le \EE(\beta_0^{x_0})^\frac{\log(\beta)}{\log(\beta_0)}$ for all $1 < \beta < \beta_0$ so there exists $\beta > 1$ such that $\EE\left(\beta^{x_0}\right) \le \alpha$.
	\end{proof}

	\subsection{Statement of the result and motivation}
	
	In this paragraph, we give more detail on the probabilistic idea behind the construction.
	The following lemma is simple and not really new but central in the pivoting technique, we state it using the power notation for extractions in order to get familiar with this notation.
	
	\begin{Lem}\label{lem:random-altern}
		Let $\Gamma$ be a metric space.
		Let $\kappa$, $\eta$ be probability distributions on $\Gamma$ and let $\alpha \in (0,1)$. 
		Let $\nu:= \alpha \eta + (1-\alpha) \kappa$.
		Let $(u_n)_{n} \sim \kappa^{\otimes\NN}$ and let $(s_n)_n \sim \eta^{\otimes\NN}$ be independent of the data of $(u_n)_n$.
		Let $(w_{2n})_{n\in\NN} \sim \mathcal{G}_\alpha^{\otimes\NN}$ be independent of the joint data of $(u_n)_n$ and $(s_n)_n$ and let $w_{2k+1} = 1$ for all $k \in\NN$.
		For all $k \in\NN$, we define:
		\begin{align*}
			\tilde{g}_{2k}:= \widetilde{u}^w_{2k} & = \left(u_{\overline{w}_{2k}}, \dots, u_{\overline{w}_{2k+1} - 1}\right) = \left(u_{\overline{w}_{2k}}, \dots, u_{\overline{w}_{2k} + w_{2k} - 1}\right)\\
			\tilde{g}_{2k + 1}:= \widetilde{s}^w_{2k +1} & = \left(s_{\overline{w}_{2k + 1}}, \dots, s_{\overline{w}_{2k+2} - 1}\right) = \left(s_{\overline{w}_{2k}}\right).
		\end{align*}
		Let $(\gamma_n)_n = \bigodot_{k=0}^{+\infty} \tilde{g}_k$ \ie $\widetilde{\gamma}^w_k = \tilde{g}_k$ for all $k \in\NN$. 
		Then $\gamma_n \sim \nu^{\otimes \NN}$.
	\end{Lem}
	
	\begin{proof}
		Let $(u_n)_n$, $(s_n)_n$ and $(w_k)_k$ be as in Lemma \ref{lem:random-altern}.
		To all integer $n$, we associate a random number $i_n$, we write $i_n = 1$ when there exists an integer $k$ such that $\overline{w}_{2k} = n$ and $0$ otherwise.
		We first show that $(i_n)_n\sim \mathcal{B}_\alpha^{\otimes\NN}$.
		For that, we first observe that by definition of $\mathcal{G}_\alpha$,
		\begin{equation*}
			\forall t, k \in\NN, \PP(w_{2k} = t\mid w_{2k} \ge t) = \alpha.
		\end{equation*}
		Moreover, the $(w_{2k})_k$ are independent so:
		\begin{equation}\label{geom-alpha}
			\forall t, k \in\NN, \PP(w_{2k} = t\mid w_{2k} \ge t,(w_{2k'})_{k' < k}) = \alpha.
		\end{equation}
		Now for all $n \in\NN$, we write $q_n:= \max\{k \in\NN\mid \overline{w}_{2k} < n\}$.
		In other words, for all $k \in\NN$; $q_n = k$ if and only if $\overline{w}_{2k} < n$ and $\overline{w}_{2k+2} \ge n$.
		Note that $q_n$ is well defined for all $n$ because we did set $w_{2k + 1} := 1$ for all $k$.
		
		By construction, $\overline{w}_{2k+2} = \overline{w}_{2k} + w_{2k} + 1$ so we can decompose the event $(q_n = k)$ into the intersection of $\overline{w}_{2k} \le n$ (which is in the $\sigma$-algebra generated by $(w_{2k'})_{k' < k}$) and $w_{2k} \ge n - \overline{w}_{2k} - 1$.
		Moreover, $i_n = 1$ if and only if $\overline{w}_{2q_n + 2} = n$, which is equivalent to saying that $w_{2q_n} = n - \overline{w}_{2q_n} - 1$.
		Note also that in \eqref{geom-alpha} one can replace the constant parameter $t$ by a random parameter $t_n^k:= n - \overline{w}_{2k} - 1$ which can be expressed in terms of $(w_{2k'})_{k' < k}$ and then use \eqref{geom-alpha} on the (non-negative) level sets of $t_n^k$.
		Therefore \eqref{geom-alpha} yields
		\begin{equation}\label{ssfkqghv}
			\forall k, n \in\NN, \PP(w_{2k} = n - \overline{w}_{2k} - 1 \mid  w_{2k} \ge n - \overline{w}_{2k} - 1 \ge 0,(w_{2k'})_{k' < k}) = \alpha.
		\end{equation}
		By definition, $(w_{2k} \ge n - \overline{w}_{2k} - 1 \ge 0) = (q_n = k) = (t_n^k \ge 0) \cap (w_{2k} \ge t_n^k)$.
		The event $(t_n^k \ge 0)$ belongs to the $\sigma$-algebra generated by $(w_{2k'})_{k' < k}$
		Under the condition $(q_n = k)$, the event $(w_{2k} = n - \overline{w}_{2k} - 1)$ is equivalent to $(i_n = 1)$ so \eqref{ssfkqghv} yields $\PP(i_n = 1 \mid q_n = k, (w_{2k'})_{k' < k}) = \alpha$.
		This is true for all $k$ in the support of $q_n$, so we rephrase it as:
		\begin{equation*}
			\forall n \in\NN, \PP(i_n = 1 \mid  (w_{2k'})_{k' < q_n}) = \alpha.
		\end{equation*}
		Moreover, the data of $(i_{n'})_{n' < n}$ is determined by (and actually equivalent to) the joint data of $q_n$ and $(w_{2k'})_{k' < q_n}$.
		In conclusion, $\PP(i_n = 1\mid i_0,\dots i_{n-1}) = \alpha$ for all $n \in\NN$, which characterizes the fact that $(i_n)_n\sim \mathcal{B}_\alpha^{\otimes\NN}$.
		
		Note moreover that the data of $(i_n)_n$ is determined by the data of $(w_k)_k$ so it is independent of the joint data of $(u_n)_n$ and $(s_n)_n$. 
		Since all three sequences are i.i.d. and independent, the sequence $(u_n,s_n,i_n)_{n \in \NN}$ is i.i.d. in $(\Gamma \times \Gamma \times \{0,1\})^\NN$.
		For all $n \in\NN$, we have $\gamma_n = s_n$ when $i_n = 1$ and $\gamma_n = u_n$ when $i_n = 0$.
		In other words, $\gamma_n$ is given by the image of $(u_n, s_n, i_n)$ by a measurable function, so $(\gamma_n)$ is i.i.d. and there remains to check that $\gamma_0 \sim \nu$.
		Let $A \subset \Gamma$ be measurable.
		Since $(i_0 = 1)$ and $(i_1 = 0)$ are disjoint events,
		\begin{equation*}
			\PP(\gamma_0 \in A)  = \PP(\gamma_0 \in A \cap i_0 = 1) + \PP(\gamma_0 \in A \cap i_0 = 0)
		\end{equation*}
		so by independence
		\begin{align*}
			\PP(\gamma_0 \in A) & = \PP(s_0 \in A) \PP(i_0 = 1) + \PP(u_0 \in A) \PP(i_0 = 0)\\
			& = \alpha\eta(A) + (1-\alpha) \kappa(A) = \nu(A). \qedhere
		\end{align*}
	\end{proof}
	
	\begin{Cor}
		Let $E$ be a Euclidean vector space and let $\nu$ be a probability distribution on $\End(E)$. 
		Let $\alpha \in (0,1)$ and let $m \in\NN$.  
		Let $\tilde\mu$ be a probability distribution on $\End(E)^m$ such that $\alpha \tilde\mu \le \nu^{\otimes m}$ and let $\tilde\kappa = \frac{1}{1-\alpha}(\nu^{\otimes m} - \alpha \tilde\mu)$.
		Let $(u_k)_k\sim \tilde\kappa^{\odot\NN}$, let $(s_k)_k\sim \tilde\mu^{\odot\NN}$ and let $(w_{2k})_{k\in\NN} \sim \mathcal{G}_{\alpha}^{\otimes\NN}$ be globally independent random sequences defined on the same probability space.
		Let $w_{2k+1} = m$ for all $k \in\NN$ and let $(\tilde{g}_n)$ be the random sequence defined as $\tilde{g}_{2k} = \tilde{u}^w_{2k}$ and $\tilde{g}_{2k + 1} = \tilde{s}^w_{2k + 1}$ for all $k \in\NN$.
		Let $(\gamma_n)_{n\in\NN} \sim \bigodot_{k = 0}^{+\infty}\tilde{g}_k$.
		Then $\tilde{g}_k  = \gamma^{mw}_k$ for all $k \in\NN$ and $(\gamma_n)_n\sim \nu^{\otimes\NN}$.
	\end{Cor}
	
	\begin{proof}
		We apply Lemma \ref{lem:random-altern} to $\nu^{\otimes m} = \alpha \tilde\mu + (1-\alpha) \tilde{\kappa}$.
		Note that $(\widetilde{u}^m_k)_k\sim \tilde\kappa^{\otimes\NN}$ and $(\widetilde{s}^m_k)_k\sim \tilde\mu^{\otimes\NN}$.
		Let $\tilde{\tilde{g}}_{2k} = {\widetilde{{\widetilde{u}}^m}}^w_{2k}$ and $\tilde{\tilde{g}}_{2k+1} = {\widetilde{{\widetilde{s}}^m}}^w_{2k +1}$ for all $k \in\NN$.
		Then $\tilde{g}_k = \bigodot \tilde{\tilde{g}}_{k}$ for all $k$ and $\bigodot_{k = 0}^{+\infty}\tilde{\tilde{g}}_{k} = (\widetilde{\gamma}^m_k)_k$ because each letter of each $\tilde{\tilde{g}}_{k}$ is a word of length $m$.
		So Lemma \ref{lem:random-altern} tells us that $(\widetilde{\gamma}^m_k)_k \sim (\nu^{\otimes m})^{\otimes\NN}$ and $(\widetilde{\gamma}^m)^w_k = \tilde{\tilde{g}}_{k}$ for all $k \in\NN$ and as a consequence $(\gamma_n)_n\sim \nu^{\otimes \NN}$ and $\tilde{g}_k  = \gamma^{mw}_k$ for all $k$.
	\end{proof}
	
	We call the $\tilde{\kappa}^{\odot\NN}$-distributed sequence $u$ for "unknown" because we do not assume anything about $\tilde\kappa$. 
	We call the $\tilde\mu^{\odot\NN}$-distributed sequence $s$ for "Schottky" because in practice we need the product $\mu = \Pi_*\tilde\mu$ to be $\rho$-Schottky for a binary relation $\AA$ and for a constant $\rho < 1/5$. 	
	
	With the notation of Definition \ref{def:power-convol}, the words $\tilde{g}_{2k}$'s defined in Lemma \ref{lem:random-altern} have distribution law $\tilde\kappa^{\otimes \mathcal{G}_\alpha}$.
	Let us now give a fully detailed statement of the main result of this section.
	We give an abstract statement for topological semi-groups so that we may apply Theorem \ref{th:ex-piv} to random walks on relatively hyperbolic groups, semi-simple lie groups, groups acting on $Cat(0)$ buildings/
	In the present paper, we only use it for $\Gamma = \End(E)$ and for $\AA = \AA^\eps$ as in Definition \ref{def:ali}.
	
	\begin{Th}[Pivot extraction]\label{th:ex-piv}
		Let $\Gamma$ be a metric semi-group endowed with a measurable binary relation $\AA$.
		Let $\alpha\in (0,1)$, let $\rho \in(0,1/5)$ and let $m \in\NN$.  
		Let $\tilde\mu$ be a probability distribution on $\Gamma^m$ such that $\alpha \tilde\mu \le \nu^{\otimes m}$ and let $\tilde\kappa:= \frac{1}{1-\alpha}(\nu^{\otimes m} - \alpha \tilde\mu)$.
		Let $\mu = \Pi_*\tilde\mu$ and assume that $\mu$ is $\rho$-Schottky for $\AA$.
		Then there exist random sequences $({u}_k)_k, ({s}_k)_k, (w_{2k})_k \sim \tilde\kappa^{\odot\NN}\otimes \tilde\mu^{\odot\NN} \otimes \mathcal{G}_{\alpha}^{\otimes\NN}$ and two random sequences of integers $(v_k)_k$ and $(p_k)_k$ all defined on the same probability space such that with the notations $\check{w} := mw$, $(\gamma_n)_{n\ge 0}:= \bigodot_{k=0}^\infty \widetilde{u}^{\check w}_{2k} \widetilde{s}^{\check w}_{2k+1}$, $\check{v}:= \check{w}^v$, $\hat{v}_{2k}:= 2 + \check{v}_{4k}$, $\hat{v}_{2k+1} = 1$ for all $k$ and $\check{p}:= \hat{v}^p$, the following assertions hold:
		\begin{enumerate}
			\item \label{indep-v}
			The data of $(v_k)_k$ is independent of the joint data of $(\tilde{u}_k)_k$ and $(w_{k})_k$. 
			\item \label{law-of-v}
			For all $k \in\NN$, $v_{4k+1} = v_{4k+2} = v_{4k+3} = 1$ and $\left(
			\frac{v_{4k}-1}{2}\right)_{k\in\NN} \sim \mathcal{G}_{2\rho}^{\otimes\NN}$. 
			\item \label{item:zouli-alignemon}
			For all $k \in\NN$, $\gamma^{\check{v}}_{4k}\AA \gamma^{\check{v}}_{4k+1} \AA \gamma^{\check{v}}_{4k+2}$. 
			\item \label{item:hat-v-is-stopping-time}
			We have $(\widetilde{\gamma}^{\check{v}}_{4k+3})_{k\in\NN} \sim \tilde\mu^{\otimes \NN}$ and the sequence $((\widetilde{\gamma}^{\check{v}}_{4k},\widetilde{\gamma}^{\check{v}}_{4k+1}, \widetilde{\gamma}^{\check{v}}_{4k+2}))_{k\in\NN}$ is i.i.d. and independent of $(\widetilde{\gamma}^{\check{v}}_{4k+3})_{k\in\NN}$.
			\item \label{item:tail}
			The data of $(p_k)_k$ is independent of the joint data of $(\tilde{u}_k)_k$ and $(w_{k})_k$ and $(v_k)_k$ and $(\widetilde{\gamma}^{\hat{v}}_{k})_{k\in \{0,1,2\}+4\NN}$.
			\item \label{item:indep-p}
			Each $p_k$ is a positive odd integer, $p_{2k+1} = 1$ for all $k$ and $(p_{2k+2})_k$ is i.i.d. and independent of $p_0$.
			\item \label{item:exp-p}
			The distribution laws of $p_0$ and $p_2$ only depend on $\rho$ and have a finite exponential moment.
			\item \label{item:indep}
			The sequence $(\widetilde{\gamma}^{\check{p}}_{2k+2})_{k \ge 0}$ is i.i.d. and independent of $\widetilde{\gamma}^{\check{p}}_{0}$.
			\item \label{item:ali}
			For all $k\in\NN$, we have $\gamma^{\check{p}}_{2k} \AA \gamma^{\check{p}}_{2k+1}$.
			\item \label{item:ali-rec}
			For all $k \in\NN$, $\gamma^{\check{p}}_{2k+1} \AA \gamma^{\hat{v}}_{\overline{p}_{2k+2}}$ and there is a family of odd integers $1 = c_1^k < c_2^k < \dots < c_{j_k}^k = p_{2k+2}$ such that for all $1\le i < j_k$,
			\begin{equation}\label{ali-rec}
				\gamma^{\hat{v}}_{\overline{p}_{2k+2}}\cdots \gamma^{\hat{v}}_{\overline{p}_{2k+2}+c_i^k-1} \AA \gamma^{\hat{v}}_{\overline{p}_{2k+2}+c_i^k} \AA \gamma^{\hat{v}}_{\overline{p}_{2k+2}+c_i^k + 1} \cdots \gamma^{\hat{v}}_{\overline{p}_{2k+2}+c_{i+1}^k - 1}.
			\end{equation}
			\item For all $k \in\NN$, the conditional distribution of $\widetilde{\gamma}^{\check{p}}_{2k+1}$ with respect to the joint data of $(\widetilde{\gamma}^{\check{p}}_{k'})_{k'\neq 2k+1}$, $(\tilde{u}_k)_{k \ge 0}$, $(w_k)_{k \ge 0}$ and $(v_k)_{k \ge 0}$ is the normalized restriction of $\tilde\mu$ to the measurable set \label{item:tjr-piv}
			\begin{equation}\label{ceka}
				C_k:= \Pi^{-1}\left\{\gamma \in\Gamma\,\middle|\, \gamma^{\check{p}}_{2k} \AA \gamma^{\check{p}}_{2k+1} \AA \gamma^{\hat{v}}_{\overline{p}_{2k+2}} \right\}.		
			\end{equation}
		\end{enumerate} 
	\end{Th}
	
	Let us quickly break down what Theorem \ref{th:ex-piv} says.
	The notation $(\gamma_n)_{n\in\NN}:= \bigodot_{k=0}^\infty \tilde{u}^{w}_{2k} \tilde{s}^{w}_{2k+1}$ is simply a compact way to say that $(\gamma_n)_{n\in\NN}$ is characterized by the fact that $\widetilde{\gamma}^{mw}_{2k} = \tilde{u}^w_{2k}$ and $\widetilde{\gamma}^{mw}_{2k+1} = \tilde{s}^w_{2k+1}$ for all $k \in\NN$.
	The notation $\check{v}:= \check{w}^v$ is a compact way to say that we set $\check{v}_k = m(w_{\overline{v}_k} + \dots + w_{\overline{v}_{k+1}-1})$ for all $k \ge 0$.
	
	Let us point out a few facts that follow directly from Theorem \ref{th:ex-piv}.
	For all $k \in\NN$, we have $\widetilde{\gamma}^{\hat{v}}_{2k} = \widetilde{\gamma}^{\check{v}}_{4k} \odot \widetilde{\gamma}^{\check{v}}_{4k+1} \odot \widetilde{\gamma}^{\check{v}}_{4k+2}$ and $\widetilde{\gamma}^{\hat{v}}_{2k+1} = \widetilde{\gamma}^{\check{v}}_{4k+3}$, by definition of $\hat{v}$ and by the first part of \eqref{law-of-v}.
	So \eqref{item:hat-v-is-stopping-time} implies that the sequence $(\widetilde{\gamma}^{\hat{v}}_{2k})_k$ is i.i.d. and independent of $(\widetilde{\gamma}^{\hat{v}}_{2k+1})_{k\in\NN} \sim \tilde\mu^{\otimes \NN}$.
	
	It is also important to note that the set $C_k$ defined in \eqref{item:tjr-piv} has $\tilde\mu$ measure at least $1-2\rho$ by the Schottky property of $\mu$.
	So the product of the normalised restriction of $\tilde\mu$ to $C_k$ is $\frac{\rho}{1-2\rho}$-Schottky.
	
	Let us now detail how Theorem \ref{th:ex-piv} implies Theorem \ref{th:pivot}.
	
	\begin{proof}[Proof of Theorem \ref{th:pivot}]
		In Theorem \ref{th:pivot}, we are given a probability measure $\nu$ supported on the semi-group $\End(E)$.
		We are also given constants $\alpha, \eps > 0$ and $m \in \NN$ such that $\nu$ is $(\alpha, \eps, m)$-contracting, in the sense that there exists a measure $\mu \le \nu^{* m}/\alpha$ that is $1/6$-Schottky for the alignment relation $\AA^\eps$.
		Let us lift  $\mu$ to a probability measure $\tilde\mu \le \nu^{\otimes m}/\alpha$ such that $\Pi_* \tilde\mu = \mu$.
		
		We apply Theorem \ref{th:ex-piv} and get random sequences $(u_k)$ and $(s_k)$ in $\End(E)$ and $(w_k)$, $(v_k)$, $(p_k)$ in $\NN$ and the auxiliary sequences $(\check{w}_k)$, $(\check{v}_k)$, $(\hat{v}_k)$ and $(\check{p}_k)$ also in $\NN$.
		By Lemma \ref{lem:random-altern}, the sequence $(\gamma_n)_{n\ge 0}:= \bigodot_{k=0}^\infty \widetilde{u}^{\check w}_{2k} \widetilde{s}^{\check w}_{2k+1}$ has law $\nu^{\otimes\NN}$.
		Let us prove that the sequence $\check{p}_k$ satisfies the conclusions of Theorem \ref{th:pivot}.
		
		Before starting, let us check that the $\check{p}_k$'s all have a finite exponential moment using Lemmas \ref{lem:power-compare} and \ref{lem:sumexp}.
		By construction, the $(w_k)$'s are independent and have a bounded exponential moment and so to the $(\check{w}_k)$'s.
		By point \eqref{law-of-v} in Theorem \ref{th:ex-piv}, the $(v_k)$'s are independent and have a bounded exponential moment and py point \eqref{indep-v} the data of $(v_k)_k$ is independent of the data of $(\check{w}_k)_k$.
		Then by Lemma \ref{lem:sumexp}, the $\check{v}_k$'s have a bounded exponential moment.
		By points \eqref{item:tail}, \eqref{item:indep-p} and \eqref{item:exp-p} in Theorem \ref{th:ex-piv}, the $(p_k)$'s are independent, all have a finite exponential moment and their joint data is independent of the data of $(\check{v}_k)_{k \ge 0}$.
		Therefore, the $\check{p}_k$'s are independent and have a finite exponential moment.
		
		In point \eqref{pivot:schottky} of Theorem \ref{th:pivot}, we state that for all $k \ge 0$, the conditional distribution of $\widetilde{\gamma}^{\check{p}}_{2k + 1}$ with respect to the joint data of $(\widetilde{\gamma}^{\check{p}}_{k'})_{k'\neq 2k+1}$ is bounded above by $\frac{3}{2} \tilde\mu$, from which we deduce that this conditional distribution projects onto a $1/4$-Schottky measure nad form which the facts that $\check{p}_{2k + 1} = m$ follows.
		This follows from Point \eqref{item:tjr-piv} in Theorem \ref{th:ex-piv}.
		Indeed, the set $C_k$ described in \eqref{ceka} has $\tilde\mu$ measure at least $2/3$ by the Schottky property.
		
		In point \eqref{pivot:renewal} of Theorem \ref{th:pivot}, we state that the sequence $(\widetilde{\gamma}^{\check{p}}_{2k + 2})_{k \ge 0}$ is i.i.d. and independent of $\widetilde{\gamma}^{\check{p}}_{0}$, which is exactly what point \eqref{item:indep} says in Theorem \ref{th:ex-piv}, and we also state that for all $k$, the conditional distribution of $\widetilde{\gamma}^{\check{p}}_{2k + 1}$ with respect to the joint data of $(\widetilde{\gamma}^{\check{p}}_{k'})_{k'\neq 2k+1}$ is given by the evaluation of a measurable function at $(\gamma^{\check{p}}_{2k}, \widetilde{\gamma}^{\check{p}}_{2k + 2})$, this means that $\widetilde{\gamma}^{\check{p}}_{2k + 1}$ is independent of $(\widetilde{\gamma}^{\check{p}}_{k'})_{k'\neq 2k+1}$, relatively to the data of $\gamma^{\check{p}}_{2k}$ and $\widetilde{\gamma}^{\check{p}}_{2k + 2}$, from which we deduce a renewal property.
		The relative independence comes from point \eqref{item:tjr-piv} in Theorem \ref{th:ex-piv} again.
		Naively, $\eta$ returns the normalized restriction of $\tilde\mu$ to $C_k$. 
		Technically, $C_k$ depends on $\widetilde\gamma^{\hat{v}}_{\overline{p}_{2k+2}}$ which is a random prefix of $\widetilde{\gamma}^{\check{p}}_{2k + 2}$ so we need to show that this random prefix is chosen independently of the data of $(\widetilde{\gamma}^{\check{p}}_{k'})_{k'\notin \{2k+1, 2k + 2\}}$, which comes from point \eqref{item:indep-p} in Theorem \ref{th:ex-piv}.
		
		In point \eqref{pivot:herali} of Theorem \ref{th:pivot}, we state that for all triple of indices $0 < i' < j' < k'$ and for all $f \AA^{\eps/2} \gamma^{\check{p}}_{i'}$ and all $\gamma^{\check{p}}_{k'-1} \AA^{\eps / 2} h$, we have $f \gamma^{\check{p}}_{i'} \cdots \gamma^{\check{p}}_{j'-1} \AA^{\eps/2} \gamma^{\check{p}}_{j'} \cdots \gamma^{\check{p}}_{k'-1} h$.
		This statement is non-probabilistic and we show that it holds any realization of $\gamma^{\check{p}}$ that satisfies points \eqref{item:ali} and \eqref{item:ali-rec} of Theorem \ref{th:ex-piv}.
		We start the induction by proving that point \eqref{pivot:herali} of Theorem \ref{th:pivot} holds when $k' = j' + 1 = i' + 2$.
		
		Assume that $k' = i' + 2$.
		If $k'$ and $i'$ are even, we write $k := i'/2$.
		Let $c_1^k < \dots < c_{j_k}^k$ be as in point \eqref{item:ali-rec} in Theorem \ref{th:ex-piv}.
		We want to apply Lemma \ref{lem:Atilde} to $n := j_k + 1$, $g_{-1} := \gamma^{\check{p}}_{2k + 1}$, $g_0 := \gamma^{\hat{v}}_{2k + 2}$, $g_{2i + 1} := \gamma^{\hat{v}}_{\overline{p}_{2k+2}+c_i^k}$ and $g_{2 i + 2} := \gamma^{\hat{v}}_{\overline{p}_{2k+2}+c_i^k + 1} \cdots \gamma^{\hat{v}}_{\overline{p}_{2k+2}+c_{i+1}^k - 1}$ for all $0 \le i < n$.
		The alignment condition \eqref{hyp-lem-atilde} follows from \eqref{ali-rec} and $\sigma(g_{2i - 1}) < \eps^6/48$ because $g_{2i - 1}$ is in the support of $\mu$ by point \eqref{pivot:schottky} from Theorem \ref{th:pivot}, proven above.
		By assumption, we also have $f\AA^{\eps / 2}g_{-1}$ and $g_{0} \cdots g_n \AA^{\eps/2} h$.
		
 		All that is left to check is that $\sigma(g_0) < \eps^2 / 12$.
 		This comes from point \eqref{item:zouli-alignemon} in Theorem \ref{th:ex-piv} and Lemma \ref{lem:triple-ali}.
 		Indeed, $\gamma^{\hat{v}}_{2k + 2} = \gamma^{\check{v}}_{4k+4}\gamma^{\check{v}}_{4k+5}\gamma^{\check{v}}_{4k+6}$ by construction, $\overline{v}_{4k+5}$ is odd and $v_{4k + 5} = 1$ so $\gamma^{\check{v}_{4k+5}} = s^{\check{w}}_{\overline{v}_{4k+5}}$ is in the support of $\mu$ and therefore $\sigma(\gamma^{\check{v}}_{4k+5})< \eps^6/48$. 
 		Moreover, by point \eqref{item:zouli-alignemon} in Theorem \ref{th:ex-piv}, $\gamma^{\check{v}}_{4k+4} \AA^\eps \gamma^{\check{v}}_{4k+5} \AA^\eps \gamma^{\check{v}}_{4k+6}$ so by Lemma \ref{lem:triple-ali}, $\sigma(\gamma^{\hat{v}}_{2k + 2}) < \eps^2/12$.
 		
 		Therefore, we can apply Lemma \ref{lem:Atilde}. 
 		We obtain that $f \gamma^{\check{p}}_{i'} \AA^{\eps/2} \gamma^{\check{p}}_{k' -1} h$ and that $\sigma(\gamma^{\check{p}}_{2k + 2}) < \eps^2 / 12$.

 		Still assuming that $k' = i' + 2$, when $k'$ and $i'$ are odd, we set $i' := 2k + 1$.
 		We want to apply Lemma \ref{lem:rig-ali} to $g_1 := \gamma^{\check{p}}_{2k+1}$ and $g_2 := \gamma^{\check{p}}_{2k+2}$. 
 		By the above arguments, $\sigma(g_{1}) < \eps^6/48$ and $\sigma(g_{2}) < \eps^2 / 12$.
 		So it follows from Lemma \ref{lem:rig-ali} that $fg_1 \AA^{\eps/2} g_2h$.
 		
 		We have proven that point \eqref{pivot:herali} of Theorem \ref{th:pivot} holds when $k = i + 2$.
 		Let $i < j < k$ and assume that $k - i \ge 3$, let $f \AA^{\eps /2} \gamma_i^{\check{p}}$ and $\gamma_{k-1}^{\check{p}} \AA^{\eps / 2} h$.
 		If $k \ge j + 2$, since $\gamma^{\check{p}}_{k - 1} \AA^{\eps / 2} h$, it follows from point \eqref{pivot:herali} from Theorem \ref{th:pivot} applied to $j' = k-1$ $i' := k - 2$ and $f' := \mathrm{Id}$ (and to $k', h' = k, h$) that $\gamma^{\check{p}}_{k' - 2} \AA^{\eps / 2} \gamma^{\check{p}}_{k' - 1} h$.
 		In this case, the alignment $g \gamma_i^{\check{p}} \cdots \gamma_{j-1}^{\check{p}} \AA^{\eps / 2} \gamma_j^{\check{p}} \cdots \gamma_{k-1}^{\check{p}} h$ follows from point \eqref{pivot:herali} applied to $k' = k - 1$, to $h' = \gamma_{k-1}^{\check{p}} h$ and to the same $i,j,f$.
 		
 		If $k = j+1$, then $i \le j-2$. 
 		In this case, we apply point \eqref{pivot:herali} from Theorem \ref{th:pivot} to $j'= i+1$, $' = i + 2$ and $h' := \mathrm{Id}$, and we obtain that $f\gamma^{\check{p}}_{i'} \AA^{\eps/2} \gamma^{\check{p}}_{i' + 1}$.
 		Therefore, the alignment relation follows from point \eqref{pivot:herali} applied to $i' = i + 1$, to $f' = f\gamma_{i}^{\check{p}}$ and to the same $j,k,h$.
 		By induction, point \eqref{pivot:herali} of Theorem \ref{th:pivot} holds for all $0 < i' < j' < k'$.
 		
		Finally, point \eqref{pivot:indep-tail} of Theorem \ref{th:pivot} states that for all $S$ such that $\tilde\mu(S^m) = 1$ (\ie $s_k \in S$ almost surely and for all $s$) but $\nu(S) > 0$, there exists a sequence $(g_n) \sim \nu_{S^c}^{\otimes \NN}$ that is independent of $(\check{p}_n)$ and such that $\gamma_n = g_n$ whenever $\gamma_n \notin S$. That comes from the fact that the sequences $(u_n)_n$ and $(\check{p}_n)_n$ are independent and $\gamma_n = u_n$ whenever $\gamma_n \notin S$.
		
		Let us give an explicit construction of $(g_n)$.
		We start with a sequence $(g'_n) \sim \nu_{S^c}^{\otimes \NN}$ that is i.i.d. and independent of all the random variables constructed in Theorem \ref{th:ex-piv}.
		Then we set $g_n = g'_n$ when $u_n \in S$ and $g_n = u_n$ when $u_n \notin S$. 
		Since the blocks of $m$ consecutive elements of $(u_n)$, denoted by $(\widetilde{u}^m_k)_k$, are i.i.d., the blocks $(\widetilde{g}^m_k)$ also are.
		
		We simply need to show that $\widetilde{g}^m_0 := (g_0, \dots, g_{m-1}) \sim \nu_{S^c}^{\otimes m}$.
		Let $I =  \{k < m \mid u_k \notin S\}$ and let $I^c = \{0, \dots, m-1\} \setminus I$. 
		We claim that the relative distribution of $(g_0, \dots, g_{m-1})$ with respect to $I$ is constant and equal to $\nu_{S^c}^{\otimes m}$.
		The data of $I$ is independent of the data of $(g'_n)_n$ by definition so the relative distribution of $(g_k)_{k \in I^c} = (g'_k)_{k \in I^c}$ with respect to the data of $I$ is $\nu_{S^c}^{\otimes I^c}$.
		Moreover, the date of $I$ is independent of the data of $w_0$ and by definition of $(\gamma_n)_{n\ge 0}:= \bigodot_{k=0}^\infty \widetilde{u}^{\check w}_{2k} \widetilde{s}^{\check w}_{2k+1}$ when $w_0 \ge 1$, we have $(u_i)_{i \in I} = (\gamma_i)_{i \in I}$ and $I = I' := \{k < m \mid \gamma_k \notin S\}$ and $I' = \emptyset$ when $w_0 = 0$.
		Moreover, the distribution of $(\gamma_i)_{i \in I'}$ with respect to $I'$ is $\nu_{S^c}^{\otimes I'}$, so once we restrict to the event $(I' \neq \emptyset) \subset (I' = I)$, we get $(g_i)_{i \in I} = \nu_{S^c}^{\otimes I}$.
		
		Let $A \subset \End(E) \setminus S$ and let $n \ge 0$.
		We use the fact that $(\gamma_n \in A) \subset (u_n \in A)$ and $_n$ is independent of the joint data of $(p_k)_{k \ge 0}$ so $\PP(\gamma_n \in A \mid (p_k)) \le \PP(u_n \in A \mid (p_k)) = \PP(u_n \in A)$.
		Moreover, $\PP(u_n = \gamma_n \mid u_n \in A) = 1 - \alpha$ so $\PP(u_n \in A) \le  \nu(A) / (1 - \alpha)$, which yields \eqref{dom-tail}.
	\end{proof}

	\subsection{First step of the proof of Theorem \ref{th:ex-piv}}
	
	In this paragraph, we will use the notations of Theorem \ref{th:ex-piv}. 
	The first step of the proof is the construction of the sequence $(v_k)_{k\in\NN}$ that satisfies points \eqref{indep-v} to \eqref{item:hat-v-is-stopping-time}.
	
	The naive idea is to construct it using a stopping time.
	We know that the $\left(\gamma^{mw}_{2k+1}\right)$ are i.i.d. $\rho$-Schottky and their joint data is independent of $\left(\gamma^{mw}_{2k}\right)_{k\in\NN}$ so for all $n \in\NN$,
	\begin{equation}\label{xkvvjlqshfl}
		\PP\left(\gamma^{mw}_{0}\cdots \gamma^{mw}_{2n} \AA \gamma^{mw}_{2n+1} \AA \gamma^{mw}_{2n+2}\mid \gamma^{mw}_{0},\dots,\gamma^{mw}_{2n}\right) \ge 1-2\rho.
	\end{equation}
	Let $n_0$ to be the smallest integer such that $\gamma^{mw}_{0}\cdots \gamma^{mw}_{2{n_0}} \AA \gamma^{mw}_{2n_0+1} \AA \gamma^{mw}_{2n_0+2}$.
	We could define $v_0:= 2n_0+1$, that way $\gamma^{\check{v}}_0 = \gamma^{mw}_{0}\cdots \gamma^{mw}_{2{n_0}}$ and $\gamma^{\check{v}}_1 = \gamma^{mw}_{2n_0+1}$ and $\gamma^{\check{v}}_2 = \gamma^{mw}_{2n_0+1}$ so \eqref{item:zouli-alignemon} holds.
	
	In this case $v_0$ would have a finite exponential moment by \eqref{xkvvjlqshfl} and by a stopping time argument, the data of $\left(\gamma^{mw}_{2k + v_0 + 2}\right)_{k\in\NN}$ would be independent of $\left(\widetilde{\gamma}^{mw}_{0},\dots,\widetilde{\gamma}^{mw}_{v_0+1}\right)$ so $\widetilde{\gamma}^{\check{v}}_3 = \widetilde{\gamma}^{mw}_{2k + v_0 + 2}$ has law $\tilde\mu$, and is independent of the data of $\left(\widetilde{\gamma}^{\check{v}}_0, \widetilde{\gamma}^{\check{v}}_1, \widetilde{\gamma}^{\check{v}}_2\right)$ and the sequence $\left(\gamma^{mw}_{2k \overline{v}_4}\right)_{k\in\NN}$ is independent of $\left(\widetilde{\gamma}^{\check{v}}_0, \dots, \widetilde{\gamma}^{\check{v}}_3\right)$ and has the same distribution as $\left(\gamma^{mw}_{2k}\right)_{k\in\NN}$ so we can reiterate the process.
	
	The advantage of this naive approach is that the sequence $(v_k)_k$ can be computed explicitly from the data of $\left(\gamma^{mw}_{2k}\right)_{k\in\NN}$.
	The issue is that the data of $(v_k)_{k\in\NN}$ would not be independent of the data of $\left(\widetilde{\gamma}^{mw}_{2k}\right)_{k\in\NN}$ so we would not have the same control over the tail of the distribution of $\gamma^{\check{v}}_0$.
	
	Luckily, this correlation between $(v_k)_{k\in\NN}$ and $\left(\widetilde{\gamma}^{mw}_{2k}\right)_{k\in\NN}$ is easy to get rid of.
	A quick observation is that if we had:
	\begin{equation}\label{ihopeso}
		\PP\left(\gamma^{mw}_{0}\cdots \gamma^{mw}_{2n} \AA \gamma^{mw}_{2n+1} \AA \gamma^{mw}_{2n+2}\mid (\gamma^{mw}_{k})_{k \neq 2n+1}\right) = 1-2\rho,
	\end{equation}
	almost surely and for all $n$, then the conditional distribution of $n_0$ with respect to the data of $\left(\widetilde{\gamma}^{mw}_{2k}\right)_{k\in\NN}$ would be almost surely equal to $\mathcal{G}_{2\rho}$, giving us the independence of $v_0$ with the data of $\left(\widetilde{\gamma}^{mw}_{2k}\right)_{k\in\NN}$.
	However \eqref{ihopeso} has no reason to be satisfied as the Schottky property only gives us an inequality.
	
	Write $A_n:= \gamma^{mw}_{0}\cdots \gamma^{mw}_{2n} \AA \gamma^{mw}_{2n+1} \AA \gamma^{mw}_{2n+2}$.
	A sure-fire way to make \eqref{xkvvjlqshfl} into an equality is to artificially remove mass from $A_n$.
	Denote by $\mathcal{U}$ the (uniform) Lebesgue measure on $[0,1]$. 
	Consider a random uniformly distributed sequence $(\tau_n)_{n\in\NN} \sim \mathcal{U}^{\otimes\NN}$ that is independent of the joint data of all the formerly defined random variables.
	We then replace the event $A_n$ with:
	\begin{equation}
		A'_n := \left(\tau_n < \frac{(1-2\rho) \mathds{1}_{A_n}}{\PP\left(A_n\,\middle|\,(\gamma^{mw}_{k})_{k \neq 2n+1}\right)}\right).
	\end{equation}
	By construction $A'_n \subset A_n$ and $\PP\left(A_n\mid (\gamma^{mw}_{k})_{k \neq 2n+1}\right) = 1-2\rho$ almost surely and for all $n \in\NN$.
	In some sense, asking moreover that $\tau_n < \frac{(1-2\rho)}{\PP\left(A_n\,\middle|\,(\gamma^{mw}_{k})_{k \neq 2n+1}\right)}$ when $A_n$ holds is like inflicting a penalty to $A_n$ to lower its probability of success just enough to make the $(A'_j)_j$ i.i.d. and independent of the data of $\left(\widetilde{\gamma}^{mw}_{2k}\right)_k$.
	Let us now construct the sequence $v$ of Theorem \ref{th:ex-piv}.

	\begin{Lem}\label{lem:first-part}
		Let $\Gamma$ be a metric semi-group endowed with a measurable binary relation $\AA$.
		Let $\alpha\in (0,1)$, let $\rho \in(0,1/5)$ and let $m \in\NN$.  
		Let $\tilde\mu$ be a probability distribution on $\Gamma^m$ such that $\alpha \tilde\mu \le \nu^{\otimes m}$ and let $\tilde\kappa = \frac{1}{1-\alpha} (\nu^{\otimes m} - \alpha \tilde\mu)$.
		Let $\mu = \Pi_*\tilde\mu$ and assume that $\mu$ is $\rho$-Schottky for $\AA$.
		Then there exist random sequences $(u_k)_k, (s_k)_k, (w_{2k})_k \sim \tilde\kappa^{\odot \NN}\otimes \tilde\mu^{\odot\NN} \otimes \mathcal{G}_{\alpha}^{\otimes\NN}$ and a random sequence of integers $(v_k)_k$ such that if we write $(\gamma_n)_{n\in\NN}:= \bigodot_{k=0}^\infty \widetilde{u}^{w}_{2k} \tilde{s}^{w}_{2k+1}$ and $\check{w}:= mw$ and $\check{v}:= \check{w}^v$, then the following assertions hold:
		\begin{enumerate}
			\item The data of $(v_k)_k$ is independent of the joint data of $(u_k)_k$ and $(w_{k})_k$.
			\item For all $k \in\NN$, we have $v_{4k+1} = v_{4k+2} = v_{4k+3} = 1$ and $\left(
			\frac{v_{4k}-1}{2}\right)_{k\in\NN} \sim \mathcal{G}_{2\rho}^{\otimes\NN}$.
			\item For all $k \in\NN$, we have $\gamma^{\check{v}}_{4k}\AA \gamma^{\check{v}}_{4k+1} \AA \gamma^{\check{v}}_{4k+2}$.
			\item We have $(\widetilde{\gamma}^{\check{v}}_{4k+3})_{k\in\NN} \sim \tilde\mu^{\otimes \NN}$ and the sequence $\left((\widetilde{\gamma}^{\check{v}}_{4k},\widetilde{\gamma}^{\check{v}}_{4k+1}, \widetilde{\gamma}^{\check{v}}_{4k+2})\right)_{k\in\NN}$ is i.i.d. and independent of $(\widetilde{\gamma}^{\check{v}}_{4k+3})_{k\in\NN}$.
		\end{enumerate}
	\end{Lem}
	
	\begin{proof}
		Let $\Omega := \Gamma^\NN \times \Gamma^\NN \times \NN^\NN \times [0,1]^\NN$ and $\PP := \tilde\kappa^{\odot\NN}\otimes \tilde\mu^{\odot\NN} \otimes \left(\mathcal{G}_{\alpha}\otimes\delta_{1}\right)^{\odot\NN} \otimes \mathcal{U}^{\otimes \NN}$.
		For all $j \in\NN$, we denote by $\mathcal{C}_j$ the pull-back of the Borel $\sigma$-algebra by the map:
		\begin{equation}\label{Cj-filtration}
			(u_n)_n, (s_n)_n,(w_k)_k, (\tau_k)_k \longmapsto (u_n)_{n < m\overline{w}_{2j}}, (s_n)_{n < m\overline{w}_{2j}},(w_k)_{k < 2j},(\tau_k)_{k < j}
		\end{equation}
		Let $T: \Omega \to \Omega$ be the transformation
		\begin{equation*}
			T : (u_n)_n, (s_n)_n,(w_k)_k, (\tau_k)_k \longmapsto (u_{n+m \overline{w}_2})_n, (s_{n + m \overline{w}_2})_n,(w_{k+2})_k, (\tau_{k+1})_k.
		\end{equation*}
		We claim that $T$ is measure preserving and that for all $j \in\NN$, the push-forward of the Borel $\sigma$-algebra by $T^j$ is independent of $\mathcal{C}_j$.
		Note that for all $j \ge 0$,
		\begin{equation*}
			T^j: (u_n)_n, (s_n)_n,(w_k)_k, (\tau_k)_k \longmapsto (u_{n+m \overline{w}_{2j}})_n, (s_{n + m \overline{w}_{2j}})_n,(w_{k+2j})_k, (\tau_{k+j})_k.
		\end{equation*}
		Let $\omega := (u_n)_n, (s_n)_n,(w_k)_k, (\tau_k)_k \sim \PP$.
		We need to check that the joint conditional distribution of $(u_{n+m \overline{w}_{2j}})_n$, $(s_{n + m \overline{w}_{2j}})_n$, $(w_{k+2j})_k$ and $(\tau_{k+j})_k$ with respect to $\mathcal{C}_j$ is almost surely equal to $\PP$ for all $j$. 
		
		Since $(w_k)_k \sim (\mathcal{G}_{\alpha}\otimes\delta_{1})^{\odot\NN}$ the conditional distribution of $(w_{k+2j})_k$ with respect to $\mathcal{C}_j$ is almost surely equal to $(\mathcal{G}_{\alpha}\otimes\delta_{1})^{\odot\NN}$ for all $j$.
		Since $\tilde\kappa$ and $\tilde\mu$ are supported on $\Gamma^m$ and by independence, the joint conditional distribution of $(u_{n + m i})_n$ and $(s_{n + m i})_n$ with respect to the joint data of $(u_n)_{n\le m i}$, $(s_n)_{n \le m i}$ and $(w_k)_k$ is almost surely equal to $\tilde\kappa^{\odot\NN} \otimes \tilde\mu^{\odot\NN}$ for all $i$.
		As a consequence, the joint conditional distribution of $(u_{n + m \overline{w}_{2j}})_n$ and $(s_{n + m \overline{w}_{2j}})_n$ with respect to the joint data of $(u_n)_{n < m\overline{w}_{2j}}$, $(s_n)_{n < m\overline{w}_{2j}}$ and $(w_k)_k$ is almost surely equal to $\tilde\kappa^{\odot\NN} \otimes \tilde\mu^{\odot\NN}$.
		Hence, the joint conditional distribution of $(u_{n+m \overline{w}_2})_{n \ge 0}$, $(s_{n + m \overline{w}_2})_{n \ge 0}$ and $(w_{k+2j})_{k \ge 0}$ with respect to the joint data of $(u_n)_{n < m\overline{w}_{2j}}$, $(s_n)_{n < m\overline{w}_{2j}}$ and $(w_k)_{k < 2j}$ is almost-surely equal to $\tilde\kappa^{\odot\NN}\otimes \tilde\mu^{\odot\NN} \otimes (\mathcal{G}_{\alpha}\otimes\delta_{1})^{\odot\NN}$.
		We conclude by the fact that $(\tau_k)_k$ is i.i.d, independent of everything else and simply shifted by $T$.	
		
		Let us construct the sequence $(v_k)_k$ as a family of measurable maps $(v_k : \Omega \to \NN)_k$.
		Since $v_{4k+1} = v_{4k+2} = v_{4k+3} = 1$ for all $k$, we only need to define the maps $(v_{4k} : \Omega \to \NN)_k$.
		Write $(\gamma_n)_{n\in\NN} := \bigodot_{k=0}^\infty \widetilde{u}^{w}_{2k} \tilde{s}^{w}_{2k+1}$.
		For all $j \in\NN$, we define
		\begin{align*}
			A_j &:= \left(\gamma^{\check{w}}_0\cdots \gamma^{\check{w}}_{2j} \AA \gamma^{\check{w}}_{2j+1} \AA \gamma^{\check{w}}_{2j+2}\right) = \left(u^{\check{w}}_0 s^{\check{w}}_1 \cdots u^{\check{w}}_{2j} \AA s^{\check{w}}_{2j+1} \AA u^{\check{w}}_{2j+2}\right)\\
			A'_{j} &:= \left(\tau_j < \frac{(1-2\rho)\mathds{1}_{A_{j}}}{\PP\left(A_{j}\mid (u_n)_n, (\widetilde{s}^w_k)_{k \neq 2j+1},(w_k)_{k}, (\tau_k)_{k \neq 2j+1}\right)}\right)
		\end{align*}
		Let $j_0 = \min\{j \mid A'_j\}$, \ie $j_0 : \omega \mapsto \min\{j \in\NN \mid \omega \in A'_{j}\}$.
		For all $k \ge 1$ we set $j_{k} := j_0 \circ T^{\overline{v}_{4k}/2}$ and $v_{4k} := 2j_k +1$.
		We remind that $\overline{v}_{4k} = 3 k + \sum_{k' = 0}^{k-1} v_{4k'}$ so $v_{4k}$ is expressed in terms of $(v_{4k'})_{k' < k}$ and is therefore well defined by induction.
		
		Let us prove that point \eqref{item:zouli-alignemon} holds. 
		Clearly, $A'_j \subset A_j$ for all $j$ so $A_{j_0}$ holds when $j_0 \neq + \infty$.
		Hence, $T^{-\overline{v}_{4k}/2}(A_{j_{k}})$ holds almost surely and for all $k$.
		We remind that, by definition:
		\begin{align*}
			T^{-\overline{v}_{4k}/2}\left(A_{j_{k}}\right) & = \left(\gamma^{\check{w}}_{\overline{v}_{4k}} \cdots \gamma^{\check{w}}_{\overline{v}_{4k} + 2 j_k} \AA \gamma^{\check{w}}_{\overline{v}_{4k} + 2 j_k +1} \AA \gamma^{\check{w}}_{\overline{v}_{4k} + 2 j_k +2}\right)\\
			& = \left(\gamma^{\check{w}}_{\overline{v}_{4k}} \cdots \gamma^{\check{w}}_{\overline{v}_{4k} + v_{4k}-1} \AA \gamma^{\check{w}}_{\overline{v}_{4k+1}} \AA \gamma^{\check{w}}_{\overline{v}_{4k+2}}\right)\\
			& = \left(\gamma^{\check{v}}_{4k} \AA \gamma^{\check{v}}_{4k+1} \AA \gamma^{\check{v}}_{4k+2}\right)
		\end{align*}
		
		Let us now prove that $(v_k)$ satisfies point \eqref{item:hat-v-is-stopping-time}.
		For all $j \in\NN$, we denote by $\mathcal{F}_j$ the $\sigma$-algebra generated by the joint data of $(u_n)_n$, $(\widetilde{s}^w_k)_{k \neq 2j + 1}$, $(w_k)_{k}$ and $(\tau_k)_{k \neq j}$.
		For all $j$, we define
		\begin{equation*}
			P_j := \PP(A_j\mid \mathcal{F}_j) = \PP\left(A_{j}\mid (u_n)_n, (\widetilde{s}^w_k)_{k \neq 2j+1},(w_k)_{k}, (\tau_k)_{k \neq 2j+1}\right).
		\end{equation*}
		The fact that $A_j \in \mathcal{C}_{j+2}$ implies that $P_j$ is $\mathcal{C}_{j+2}$-measurable so $P_j = \PP(A_j \mid \mathcal{C}'_j)$, for  $\mathcal{C}'_j := \mathcal{C}_{j+2} \cap \mathcal{F}_j$. 
		By definition, $\tau_j$ is $\mathcal{C}_{j+2}$-measurable, yielding that $A'_j \in \mathcal{C}_{j+2}$.
		Hence $(j_0 = j) \in \mathcal{C}_{j+2}$ for all $j \in\NN$.
		In other words $j_0$ is a stopping time for the filtration $(\mathcal{C}_{j+2})_{j \in \NN}$ so $T^{j_0+2}(\omega)$ is independent of $\mathcal{C}_{j_0+2}$ and the data of $\left(\widetilde{\gamma}^{\check{v}}_0, \widetilde{\gamma}^{\check{v}}_1, \widetilde{\gamma}^{\check{v}}_2, \widetilde{\gamma}^{\check{v}}_3\right)$ is $\mathcal{C}_{j_0+2}$-measurable.
		Therefore, the sequence $((\widetilde{\gamma}^{\check{v}}_{4k}, \dots, \widetilde{\gamma}^{\check{v}}_{4k+3}))_{k \in\NN}$ is i.i.d.
		Moreover, for all $j \in\NN$, the events $(A'_{j'})_{j' \le j}$ and the data of $\left(\widetilde{\gamma}^{\check{w}}_0, \dots, \widetilde{\gamma}^{\check{w}}_{2j+2}\right)$ are measurable for $\mathcal{F}_{j+1}$, which is independent of $\widetilde{s}^w_{2j + 3} = \widetilde\gamma^w_{2j +3}$.
		Therefore, $\widetilde{\gamma}^{\check{v}}_{3} = \widetilde{s}^w_{2j_0+3}$ is independent of the joint data of $j_0$ and $(\widetilde{\gamma}^{\check{v}}_{0}, \widetilde{\gamma}^{\check{v}}_{1}, \widetilde{\gamma}^{\check{v}}_{2})$.
		This proves point \eqref{item:hat-v-is-stopping-time}
		
		To prove points \eqref{indep-v} and \eqref{law-of-v}, there remains to show that $j_0$ is independent of the joint data of $(u_k)$ and $(w_k)$ and has distribution law $\mathcal{G}_{2\rho}$.
		The joint data of $(u_k)$ and $(w_k)$ is $\bigcap_j \mathcal{F}_j$-measurable and $A'_j \in \bigcap_{j' \neq j} \mathcal{F}_{j'}$. 
		To conclude, it suffices to prove that $\PP(A'_j \mid \mathcal{F}_j) = 1 - 2 \rho$ for all $j$. 
		Let $j \in\NN$ be fixed.
		By construction, ${s}^w_{2j+1}$ is independent of $\mathcal{F}_j$ and has distribution law $\mu$, which is $\rho$-Schottky for $\AA$ by assumption. 
		Therefore, $P_j \ge 1 - 2\rho$.
		By construction, $\tau_j$ is uniformly distributed and independent of $\mathcal{F}_j \ni P_j$ so
		\begin{equation}\notag\label{sbdvkqskl}
			\PP\left(\tau_j < \frac{(1-2\rho)}{P_{j}} \,\middle|\,\mathcal{F}_j\right) = \frac{(1-2\rho)}{P_{j}}.
		\end{equation}
		Moreover, $A_j$ belongs to the $\sigma$-algebra generated by $\widetilde{s}^w_{2j+1}$ and $\mathcal{F}_j$, which is independent of $\tau_j$ by construction.
		Therefore:
		\begin{equation*}
			\PP\left(\tau_j < \frac{(1-2\rho) \mathds{1}_{A_j}}{P_{j}} \,\middle|\,\mathcal{F}_j\right) = \PP(A_j \mid \mathcal{F}_j) \PP\left(\tau_j < \frac{(1-2\rho)}{P_{j}} \,\middle|\,\mathcal{F}_j\right).
		\end{equation*}
		Hence $\PP(A'_j \mid \mathcal{F}_j) = 1-2 \rho$.
	\end{proof}

	\subsection{Pivot algorithm}
	
	To construct the sequence $(p_k)_k$, we will apply the following deterministic algorithm to the random sequences $\left(\gamma^{\hat{v}}_k\right)_{k \in\NN},(\tau_{k})_{k\in\NN}$, where $\left(\gamma^{\hat{v}}_k\right)_{k \in\NN}$ is as constructed in Lemma \ref{lem:first-part} and $(\tau_k)_k \sim \mathcal{U}^{\otimes \NN}$ is independent of all previously defined random variables.
	
	\begin{Def}[Pivoting technique with penalty]\label{def:pivot-alg}
		Let $\Gamma$ be a metric semi-group endowed with a measurable binary relation $\AA$. Let $0< \rho < 1/5$ and let $\mu$ be a probability distribution on $\Gamma$ that is $\rho$-Schottky for $\AA$.
		To all pair of sequences $(\gamma_n)\in \Gamma^{\NN}$ and $(\tau_n) \in [0,1]^\NN$, we associate a family of integers $(p^j_k)_{j,k \in \NN}$, defined by induction on $j$ with the following algorithm.
		For $j = 0$ and for all $k \in \NN$ we set $p^0_k := 1$.
		Let $j \ge 0$ and assume that we have constructed a sequence of odd integers $(p^j_k)_{k \ge 0}$ where $p^j_{2k + 1} = 1$ for all $k$ and an integer $m_j$ satisfying $\overline{p}^j_{2m_j + 1} = 2j+1$ and $p^j_{2k} = 1$ for all $k > m_j$.
		\begin{itemize}
			\item If $\gamma^{p^j}_{2m_j} \AA \gamma^{p^j}_{2j+1} \AA \gamma^{p^j}_{2j+2}$ and:
			\begin{equation*}
				\tau_j < \frac{1-2\rho}{\mu\left\{s \in \Gamma\,\middle|\, \gamma^{p^j}_{2m_j} \AA \, s \, \AA \gamma^{p^j}_{2m_j+2}\right\}},
			\end{equation*}
			we set $m_{j+1} := m_j + 1$ and $p^{j+1}_k := p^{j}_k$ for all $k \in\NN$.
			\item Otherwise, and if there exists $0 \le k < m_j$ such that:
			\begin{equation*}
				\tau_{\frac{\overline{p}^j_{2k+1}-1}{2}} < \frac{(1-3\rho) \mathds{1}\left(\gamma^{p^j}_{2k+1} \AA \gamma_{\overline{p}^j_{2k + 2}}\cdots\gamma_{2j+2}\right)}{\mu\left\{s \in \Gamma\,\middle|\, \gamma^{p^j}_{2k} \AA\, s\, \AA \gamma_{\overline{p}^j_{2k + 2}} \cap s \, \AA \gamma_{\overline{p}^j_{2k + 2}} \cdots \gamma_{2j+2} \right\}},
			\end{equation*}
			we define $m_{j+1}$ to be the largest such $k$, for all $k < 2 m_{j+1}$, we set $p^{j+1}_{k} := p^j_k$, we set $p^{j+1}_{2m_{j+1}} = 2j + 3 - \overline{p}^j_{2m_{j+1}}$ and we set $p^j_{k} := 1$ for all $k > 2m_{j+1}$.
			\item Otherwise, and if there is no such $k$, we set $p^{j+1}_0 = 2j+3$ and $p^{j+1}_k = 1$ for all $k >0$.
		\end{itemize}
	\end{Def}
	
	\begin{Def}\label{def:temps-pivots}
		Let $(p^j_k: \Omega \to \NN)_{j,k \ge 0}$ be as in Definition \ref{def:pivot-alg}.
		Let $\Omega' \subset \Omega$ be the set of entries $(\gamma_n)_n, (\tau_n)_n$ such that the sequence $(p^j_{k})_{j \ge 0}$ is eventually stationary for all $k \ge 0$. 
		For all $k \ge 0$, we define $p_k: \Omega' \to \NN$ by $p_k = \lim_{j} p^j_k$.
		For all $k \ge 0$, we say that $\overline{p}_{2k+1}$ is the $k$-th pivotal time of the sequence $(\gamma_k)_{k \ge 0}$ with weights $(\tau_k)_{k \ge 0}$.
	\end{Def}
	
	We use the following notations through the present section.
	For all $j \ge 0$, we define
	\begin{equation}\label{Aj-def}
		A_{j}:= \left(\tau_j < \frac{(1-2\rho) \mathds{1}\left(\gamma^{p^j}_{2m_j} \AA \gamma^{p^j}_{2m_j + 1} \AA \gamma^{p^j}_{2m_j+2}\right)}{\mu\left\{s \in \Gamma\,\middle|\, \gamma^{p^j}_{2m_j} \AA s \AA \gamma^{p^j}_{2m_j+2}\right\}}\right).
	\end{equation}
	Note that for all $j$, we have $A_j = (m_{j+1} = m_j + 1)$.
	We say that the Pivot algorithm goes forward at step $j$ when $A_j$ is satisfied, otherwise we say that it backtracks.
	For all $j,k \ge 0$, we define
	\begin{equation}\label{Bjk-def}
		B_{j,k}:= \left(k < m_j\right) \cap \left(\tau_{\frac{\overline{p}^j_{2k+1}-1}{2}} < \frac{(1-3\rho) \mathds{1}\left(\gamma^{p^j}_{2k+1} \AA \gamma_{\overline{p}^j_{2k + 2}}\cdots\gamma_{2j+2}\right)}{\mu\left\{s\,\middle|\, \gamma^{p^j}_{2k} \AA s \AA \gamma_{\overline{p}^j_{2k + 2}} \cap s \AA \gamma_{\overline{p}^j_{2k + 2}} \cdots \gamma_{2j+2} \right\}}\right).
	\end{equation}
	Note that for all $j,k$, the event $B_{j,k}$ implies that $m_{j+1} \ge k$ and in fact:
	\begin{equation}\label{backtrack}
		\forall j \ge 0,\; m_{j+1} = \mathds{1}(A_j)(m_j + 1) + \mathds{1}(\Omega\setminus A_j) \max_{k \ge 0} k \mathds{1} (B_{j,k}).
	\end{equation}
	Here $\max k \mathds{1}(B_{j,k})(\omega) = \max(\{k \mid \omega \in B_{j,k}\} \cup\{0\})$.
	
	Let us illustrate what the pivot algorithm does on the first 4 steps of a toy model example.
	Let $\Gamma = \mathrm{F}_6 =\langle a,b,c \rangle$, let $\ell$ be the word length associated to the set of generators $S = \{a, b, c, a^{-1}, b^{-1}, c^{-1}\}$ let $\AA = \{f,g \mid  \ell(fg) = \ell(f) + \ell(g)\}$ and let $\mu$ be the uniform distribution on $S$.
	One can easily check that $\mu$ is $\frac{1}{6}$-Schottky for $\AA$.
	Assume that we have drawn $\gamma_0 = a^9$, $\gamma_1 = a$, $\gamma_2 = cba$, $\gamma_3 = c$, $\gamma_4 = b^{-5}$, $\gamma_5 = b^{-1}$, $\gamma_6 = b^{269}$, $\gamma_7 = c^{-1}$ and $\gamma_8 = a$.
	This draw respects our assumptions that all $\gamma_k$'s are non trivial and all the $\gamma_{2k+1}$'s are in $S$.
	The starting decomposition is:
	\begin{multline*}
		(\gamma_0 \mid  [\gamma_1, \tau_0],\gamma_2, [\gamma_3, \tau_1], \gamma_4, [\gamma_5, \tau_2],\gamma_6, [\gamma_7, \tau_3],\gamma_8,\dots) \\ = (a^{9} \mid  [a, \tau_0], aba, [c, \tau_1], b^{-5}, [b^{-1}, \tau_2], b^{269}, [c^{-1}, \tau_3], a, \dots).
	\end{multline*}
	We mark with brackets the candidate pivotal times, which are simply the oddly indexed times and to each of them, we associate a penalty $\tau_k$.
	We denote by a vertical bar the position $m_n$.
	We start with $m_ 0 = 0$.
	At time $1$, we check whether $\gamma_0 \AA \gamma_1 \AA \gamma_2$, this is true because the word $a^9 b a b a$ is reduced. 
	Moreover, we need to check whether $\tau_0 < 4/5$ because $\mu\{s\mid  \gamma_0 \AA \, s \, \AA \gamma_2\} = \mu\left(S \setminus\{ a^{-1} \}\right) = 5/6$.
	Assume that $\tau_0 < 4/5$.
	Then at time $n = 1$, $m_n = 1$ and:
	\begin{multline*}
		\left(\gamma_0^{p^n},[\gamma_1^{p^n}, \tau_{(\overline{p}^n_3 - 1) / 2}],\dots, \gamma_{2m_n}^{p^n} \,\middle|\, [\gamma_{2 n + 1}, \tau_n], \gamma_{2n + 2}, \dots\right) \\ = \left(a c^9, [b, \tau_0], aba \,\middle|\, [c, \tau_1], b^{-5}, \dots\right).
	\end{multline*}
	Now we check whether $\gamma^{p^1}_{2 m_1} = \gamma_2 \AA \gamma_3 \AA \gamma_4$, this is true because $aba c b^{-5}$ is reduced. Moreover, there is no condition to check on $\tau_1$ because $\mu\{s\mid  \gamma_2 \AA \, s \, \AA \gamma_4\} = \mu (S \setminus\{ a^{-1}, b \}) = 4/6$.
	Then, $m_2 = 2$ and:
	\begin{multline*}
		\left(\gamma_0^{p^n},\dots, \gamma_{2m_n}^{p^n} \,\middle|\, [\gamma_{2 n + 1}, \tau_n], \gamma_{2n + 2}, \dots\right) \\ = \left(a c^9, [b, \tau_0], aba, [c, \tau_1], b^{-5}\,\middle|\, [b^{-1}, \tau_2], b^{269}, \dots\right).
	\end{multline*}
	Here, we do not have $b^{-1} \AA\, b^{269}$ so we need to backtrack.
	First, we check whether $\gamma_{3} \AA \gamma_{4} \gamma_5 \gamma_6$. 
	This holds because $c \AA b^{263}$.
	Moreover, there is no condition to check on $\tau_1$ because $\mu\{s\mid  \gamma_2 \AA \, s \, \AA \gamma_4 \cap s \, \AA \gamma_4 \gamma_5 \gamma_6\} = \mu (S \setminus\{ a^{-1}, b , b^{-1}\}) = 3/6$.
	We need to check conditions on $\tau$ when the alignment condition with the product $\gamma_4 \gamma_5 \gamma_6$ is trivial or redundant with one of the former alignment conditions.
	Then at time $n  = 3$, we have backtracked to $m_n  = 1$ again and:
	\begin{equation*}
		\left(a c^9, [b, \tau_0], abacb^{263} \,\middle|\, [c^{-1}, \tau_3], a,\dots\right).
	\end{equation*}
	We check whether $abacb^{263} \AA\, c^{-1} \AA\, a$, this is true.
	Moreover $b \neq a$ so there is no condition to check on $\tau_3$ and we get:
	\begin{equation*}
		\left(a c^9, [b, \tau_0], abacb^{263}, [c^{-1}, \tau_3], a \,\middle|\, [\gamma_9, \tau_4], \gamma_{10}, \dots\right).
	\end{equation*}
	We then keep going with the same algorithm. 
	Assuming only that all the $\gamma_k$'s are non trivial, we can show by induction on $n \ge 0$ that the words $(\gamma^{p^n}_{0},\gamma^{p^n}_{1}, \dots,\gamma^{p^n}_{2 m_n})$ on the left of the vertical bar are all non-trivial elements of $\mathrm{F}_6$ and their product $\gamma^{p^n}_{0} \cdots \gamma^{p^n}_{2 m_n}$ has length $L(\gamma^{p^n}_{0}) + \cdots + L(\gamma^{p^n}_{2 m_n})$, which is equivalent to saying that the last letter of each word is not the inverse of the first letter of the next one \ie they are aligned. 
	The induction is pretty straightforward in this toy model case, if we move the bar forward, then by assumption, the two words we add are non trivial and aligned by assumption and if we move backward, then the last word is the product of three aligned words, the first two of them being non trivial and aligned with the previous word by induction hypothesis so the product is non trivial and aligned with the previous word, which concludes the induction.
	
	In a group whose geometry is not the one of a tree the construction is not as simple but the same proof by induction yields the technical formula of point \eqref{item:ali-rec} in Theorem \ref{th:ex-piv}, which implies alignment by Lemma \ref{lem:Atilde}.
	
	The reason why we need to check conditions on $\tau_k$ is not to guarantee alignment but to make sure that the Markov chain $m_n$ is a renewal process that is independent of the data of $(\gamma_{2k})_{k \ge 0}$. 
	In fact, show that the event $(m_{n+1})$ is independent the joint data of $(\gamma_k)_{k \le 2n}$ and $(m_{j})_{j \le n}$ \ie the past.

	Let us now give more details on what this means exactly.
	
	We fix a metric semi-group\footnote{We call metric semi-group a metric space endowed with a continuous and associative product map.} $\Gamma$, endow it with its Borel $\sigma$-algebra $\mathcal{A}_{\Gamma}$ and fix a measurable binary relation $\AA$. 
	We fix a constant $0 < \rho < 1/5$ and a probability distribution $\mu$ on $\Gamma$ that is $\rho$-Schottky for $\AA$.
	Let $\Omega:= \Gamma^\NN \times [0,1]^\NN$ and let $\mathcal{A}_{\Omega}$ be the infinite product of the Borel $\sigma$-algebras on $\Gamma$ and $[0,1]$.
	For all $n \in\NN$, we write $\gamma_n: \Omega \to \Gamma$ for the map $((\gamma'_k)_k, (\tau'_k )_k) \mapsto \gamma'_n$ and $\tau_n: \Omega \to \Gamma$ for the map $((\gamma'_k)_k, (\tau'_k )_k) \mapsto \tau'_n$.
	Then for all $n \in\NN$, the maps $\gamma_n$ and $\tau_n$ are measurable.
	
	The algorithm described in Definition~\ref{def:pivot-alg} associates a family $(p_k^j)_{j,k \ge 0}$ to all couple of sequences $((\gamma_k)_k, (\tau_k )_k) \in \Omega$ so it describes a family of functions $(p_k^j: \Omega \to \NN)_{j,k \ge 0}$.
	
	We use the following notations. 
	Given $k \ge 0$, given $F(X_1, \dots, X_k)$ a formula and given $f_1, \dots, f_k$ a family of functions on $\Omega$, we write $\left(F(f_1, \dots, f_k)\right)$ for the set $\left\{\omega \in\Omega\,\middle|\, F(f_1(\omega), \dots, f_k(\omega))\right\}$. 
	Given $A \subset \Omega$, we say that the formula $F(f_1, \dots, f_k)$ is satisfied on $A$ when $A \subset (F(f_1, \dots, f_k))$.
	Given $f: \Omega \to \NN$ and $T: \Omega \to \Omega$, we write $T^f: \Omega \to \Omega \;;\, \omega \mapsto T^{f(\omega)}(\omega)$. 
	Note that if $T$ and $f$ are measurable then $T^f$ is.
	
	For now $\Omega$ is nothing more that a metric space endowed with its Borel $\sigma$-algebra, in particular We do not yet have a notion of almost-sure properties on $\Omega$.
	Let:
	\begin{equation}\label{T-omega-def}
		T_\Omega:\, \Omega \to \Omega ; \, (\gamma_n)_{n \ge 0}, (\tau_n)_{n\ge 0} \mapsto (\gamma_{n+2})_{n \ge 0}, (\tau_{n+1})_{n\ge 0}.
	\end{equation}
	That way $\tau_k  = \tau_0 \circ T_\Omega^k$, $\gamma_{2k} = \gamma_0 \circ T_\Omega^k$ and $\gamma_{2k + 1} = \gamma_1 \circ T_\Omega^k$ for all $k \ge 0$.	
	Given two measurable maps $\phi : \Omega \to \NN$ and $\theta: \Omega \to \Omega$, we write $\theta^\phi : \omega \mapsto \theta^{\phi(\omega)}(\omega)$.
	This map is measurable on every level set of $\phi$.
	There are countably many such level sets so $\theta^\phi$ is measurable.
	
	For all $j \in \NN$, we denote by $\mathcal{P}_j$ the $\sigma$-algebra generated by $(\gamma_k)_{k \le 2j}$ and $(\tau_k)_{k < j}$.
	Clearly, $(\mathcal{P}_j)_{j \ge 0}$ is a filtration and by induction, the map $m_j : \Omega \to \NN$ is $\mathcal{P}_j$-measurable for all $j$.
	Indeed, at all time $j$, the data of $m_{j+1}$ is determined by the data of $(\gamma_{k})_{k \le 2j + 2}$, $(\tau_k)_{k \le j}$ and $(p^j_k)_{k \ge 0}$. 
	The data of $(\gamma_{k})_{k \le 2j + 2}$ and $(\tau_k)_{k \le j}$ is $\mathcal{P}_{j+1}$-measurable by definition and $(p^j_k)_{k \ge 0}$ is determined by the data of $(m_{j'})_{j' \le j}$, which we may assume to be $\mathcal{P}_j$-measurable by induction. 
	It follows that $m_{j+1}$ is $\mathcal{P}_{j+1}$-measurable.
	
	The following Lemma does not say anything about probabilities, its role is to break down in details what the pivot algorithm does and introduce useful notations.
	
	\begin{Lem}\label{lem:descri-pivot}
		With the notations of Definition \ref{def:temps-pivots} the functions $(p^j_k: \Omega \to \NN)_{j,k \ge 0}$, set $\Omega'$ and the functions $(p^j_{k}: \Omega' \to \NN)_{j \ge 0}$ are measurable.
		Let $(m_{j})_{j \ge 0}$ be as in Definition \ref{def:pivot-alg}.
		Then $\Omega' = (\lim_j m_j = +\infty)$ and on $\Omega'$,
		\begin{equation}\label{l-k}
			\forall k \ge 0,\;\frac{\overline{p}_{2k+1} - 1}{2} = \max\{j \in \NN\mid m_j \le k\}.
		\end{equation}
		More precisely,
		\begin{equation}\label{l-j-k}
			\forall j, k \ge 0,\;\frac{\overline{p}^j_{2k+1} - 1}{2} = \max\{j' \le j\mid  m_{j'} \le k\} + (k - m_j)^+.
		\end{equation}
		Then $\Omega'$ is $T_\Omega^{\frac{p_0 + 1}{2}}$ invariant.
		Moreover, $(T_\Omega^{\frac{p_0 + 1}{2}})^k = T_\Omega^{\frac{\overline{p}_{2k}}{2}}$ for all $k \ge 0$ and $p_{2k + 2i} = p_{2i} \circ T_\Omega^{\frac{\overline{p}_{2k}}{2}}$ for all $i,k \ge 0$. 
		Furthermore, for $A_j$ as in \eqref{Aj-def}, there exists $S \subset \Omega$ measurable such that,
		\begin{equation}\label{p-zero}
			\forall j, k \ge 0,\; (\overline{p}_{2k+1} = 2j + 1) = (m_j = k) \cap A_j \cap T_\Omega^{-j-1}(S). 
		\end{equation}
		More precisely, there exists a non-increasing sequence of measurable sets $(S^{j'})_{j' \ge 0}$ such that
		\begin{equation}\label{p-j'}
			\forall i, j, k, \;(\overline{p}^{i + j}_{2k+1} = 2j + 1) = (m_j = k) \cap A_j \cap T_\Omega^{-j-1}(S^{i}). 
		\end{equation}
	\end{Lem}
	
	\begin{proof}
		For $(B_{j,k})$ as defined in \eqref{Bjk-def}, we set
		\begin{equation*}
			{S}:= \bigcap_{j = 0}^{+\infty} \left(A_j \cup \bigcup_{k = 0}^{+ \infty} B_{j,k}\right).
		\end{equation*}
		For all $j,k \in\NN$, \eqref{backtrack} yields
		\begin{equation}\label{m-j+1}
			( m_{j+1}  = k ) = A_j \cap (m_j = k-1) \sqcup B_{j,k} \setminus A_j \setminus \bigcup_{k' > k} B_{j,k'}  \sqcup (k = 0) \setminus A_j \setminus \bigcup_{k'} B_{j,k'}.
		\end{equation}
		
		Let us prove \eqref{l-j-k}.		
		First note that $m_{j + 1} \le m_{j} + 1$ for all $j \ge 0$. 
		Let us show by induction on $j \ge 0$ that
		\begin{equation}\label{futur-p}
			\overline{p}^j_{2m_j + 1} = 2j + 1 \quad\text{and} \quad \forall k > 2m_j,\; p^j_k = 1.
		\end{equation}  
		This is true for $j = 0$ because $\overline{p}^0_{1} = p^0_0 = 1$.
		Assume that $\overline{p}^j_{2m_j + 1} = 2j+1$ for an integer $j \ge 0$. 
		It holds on $A_j$ that $\overline{p}^{j+1}_{2m_{j+1} + 1} = \overline{p}^j_{2m_j + 3} = \overline{p}^j_{2m_j + 3} = \overline{p}^j_{2m_j + 1} + p^j_{2m_j + 1} + p^j_{2m_j + 2} = 2j + 3$. Outside of $A_j$, we explicitly set $p^{j+1}_k = 1$ for all $k > 2m_{j+1}$ moreover, $\overline{p}^{j+1}_{2m_{j+1} + 1} = \overline{p}^{j+1}_{2m_{j+1}} + {p}^{j+1}_{2m_{j+1}} = 2j+3$.  
		
		Note also that:
		\begin{equation}\label{past-p}
			\forall j \ge 0, \forall k < 2 m_{j+1},\; p^{j+1}_k = p^j_k.
		\end{equation}
		A notable consequence of \eqref{futur-p} and \eqref{past-p} is that $p^j_{2k+1} = 1$ for all $j$ and all $k$.
		
		Let $j \ge k \ge 0$.
		On the set $(k \le m_j)$, we set $l^j_k := \max\{j' \le j \mid  m_{j'} \le k\}$. 
		Then $m_{l^j_k} = k$ so $\overline{p}^{l^j_k}_{2k+1} = 2l^j_k + 1$ by \eqref{past-p}.
		When $l^j_k = j$, \eqref{l-j-k} holds, otherwise $m_{j'+1} > k$ for all $j' \in\{l^j_k, \dots, j-1\}$ so \eqref{past-p} yields that $p^j_{k'} = p^{l^j_k}_{k'}$ for all $k' \le 2k$.
		Therefore $\overline{p}^j_{2k+1} = 2l^j_k + 1$, which proves that \eqref{l-j-k} holds on the set $(k \le m_j)$. 
		On the set $(k > m_j)$, we have by definition $\overline{p}^j_{2k + 1} = 2j + 1 + 2(k - m_j)$, which concludes the proof of \eqref{l-j-k}.
		
		By \eqref{past-p} and \eqref{futur-p}, the map $j \mapsto \overline{p}^j_k$ is non-decreasing for all $k$ so it has a (maybe infinite) limit as $j \to +\infty$ that we denote by $\overline{p}_k$. 
		Taking the limit $j \to +\infty$ in \eqref{l-j-k} yields \eqref{l-k}. 
		The facts that $p^j_{2k+1} =1$ for all $k$, that $p^j_{2k} = 1$ for all $k > m_j$ and \eqref{l-j-k} guarantee that the data of $(p^j_{k})_{k \ge 0}$ is fully determined by the data of $(m_{j'})_{j' \le j}$. 
		
		Now let us prove \eqref{p-zero}.
		By \eqref{l-k}, we know that $(p_0 = 2j+1)$ is equal to $(m_{j} = 0) \cap \left(\forall j' > j, m_{j'} \ge m_j\right)$.
		Moreover, the set $(\forall j' > j, m_{j'} > m_j)$ is the intersection of $A_j$ with the set
		\begin{equation*}
			S_j := \bigcap_{j' = j+1}^{\infty} (m_{j'+1} \ge m_{j+1}) = \bigcap_{j' = j+1}^{+\infty} \left(A_{j'} \cup \bigcup_{k = 1}^{+ \infty} B_{j',k} \cap (k \ge m_{j+1})\right).
		\end{equation*}
		We claim that $S_j \cap A_j = T_\Omega^{-j-1}(S) \cap A_j$.
		First note that \eqref{futur-p} implies that $\gamma^{p^{j}}_{2m_j + 2 + k} = \gamma_{2j+2+k}$ for all $k \ge 0$.
		Moreover, $A_j$ implies that $\gamma^{p^{j+1}}_{2m_{j+1} + k} = \gamma_{2j+2+k}$ for all $k \ge 0$. 
		For all $j' \ge 0$, set
		\begin{equation*}
			S^{j'}:= \bigcap_{j'' = 0}^{j'-1} \left(A_{j''} \cup \bigcup_{k = 0}^{+ \infty} B_{j'',k}\right)
		\end{equation*} 
		and for all $j' > j$, set
		\begin{equation*}
			S_j^{j'} := \bigcap_{j'' = j + 1}^{j'-1} \left(A_{j''} \cup \bigcup_{k = 1}^{+ \infty} B_{j'',k} \cap (k \ge m_{j+1})\right).
		\end{equation*}
		Note that $S^{j'}_j \subset S^{j' + 1}_{j}$ for all $j \le j'$.
		For all $j < j'$, the set $S_j^{j'} \cap A_j$ corresponds to the sequences $(\gamma, \tau) \in \Omega$ such that $m_{j''} > m_j$ for all $j < j'' \le j'$ but the sets $S_{j}^{j'}$ and $S^{j'}$ can not be described in terms of the sequence $(m_k)$ alone.
		We claim that $S_j^{j'} \cap A_j = T_\Omega^{-j-1}(S^{j'- j - 1}) \cap A_j$ for all $j' > j$ and that on $S_j^{j'} \cap A_j$, the following holds:
		\begin{gather}
			m_{j' + 1} = m_{j' - j} \circ T_\Omega^{j + 1} + m_{j + 1} \label{m-j'}\\
			\forall k \ge 0,\; p^{j'}_{2k + 2m_{j+1}} = p^{j' - j - 1}_{2k} \circ T_\Omega^{j+1}\label{p-k+m}
		\end{gather}
		Let us focus on the case $j' = j + 1$ first.
		Note that $S^0 = A_0$ and $S_{j}^{j+1} = A_{j+1}$.
		Moreover, on $A_j$, the identification $(\gamma^{p^j}_{2m_{j+1} + k}) = \gamma_{2j + 2 + k} = \gamma_k \circ T_\Omega^{j+1}$ holds so $A_{j+1} \cap A_j = T_\Omega^{-j-1}(A_0) \cap A_j$.
		We deduce that $S_j^{j + 1} \cap A_j = T_\Omega^{-j-1}(S^{0}) \cap A_j$.
		Moreover, for $j' = j + 1$, \eqref{p-k+m} follows from \eqref{m-j'} and both hold trivially.
		Indeed, on $S_j^{j + 1} \cap A_j = A_{j+1}\cap A_j$ we have $m_{j + 2} = m_{j} + 2 = m_{j+1} + 1$ and $m_1 = 1$ on $A_0$.
		
		Let $i, j \ge 1$ be fixed and assume for all $j' \le j + i$, that $S_j^{j'} \cap A_j = T_\Omega^{-j-1}(S^{j'- j - 1}) \cap A_j$ and that \eqref{m-j'} and \eqref{p-k+m} hold on $S_j^{j'} \cap A_j$.
		We want to show that the same holds for $j' = j + i + 1$. 
		By \eqref{futur-p} and by induction on $i$, $\overline{p}^{i + j}_{2m_{j+1}} = 2j+2$ on $S_j^{j+i} \cap A_j$, combined with \eqref{p-k+m} for $j' = j + i$, this yields that $\gamma^{p^{i + j}}_{2k + 2m_{j+1}} = \gamma^{p^{i - 1}}_{2k} \circ T_\Omega^{j+1}$ for all $k \ge 0$. 
		Therefore, $B_{j'+1, k} \cap (k > m_j) \cap A_j = T_\Omega^{-j-1} (B_{j'-j, k-m_j}) \cap (k > m_j) \cap A_j$ for all $k \ge 0$ and $A_{j' + 1} \cap A_j = T_\Omega^{-j-1}(A_{j'-j}) \cap A_j$.
		It follows that $S_{j}^{j + i + 1} \cap A_j = T_\Omega^{- j -1}(S^{i}) \cap A_j$ and by \eqref{backtrack}, \eqref{m-j'} holds for $j' := j + i + 1$ and \eqref{p-k+m} follows by construction $p$.
		
		By induction, $S_j^{j'} \cap A_j = T_\Omega^{-j-1}(S^{j'- j - 1}) \cap A_j$ for all $j' > j$ and \eqref{m-j'} and \eqref{p-k+m} hold on this set.
		Taking the limit $j' \to +\infty$, we obtain \eqref{p-zero}.
		Combined with \eqref{futur-p}
		
		From \eqref{p-k+m}, we deduce that:
		\begin{equation}\label{p-j-k+i}
			p^{j'}_{2k + 2i} = p^{j' - \frac{\overline{p}^j_{2k}}{2}}_{2i} \circ T_{\Omega}^{\frac{\overline{p}^{j'}_{2k}}{2}}
		\end{equation}
		on the set $\left(m_j \ge k\right)$ for all $i,j,k\ge 0$.
		Indeed, \eqref{p-j-k+i} is trivial for $k = 0$ and on $\left(m_{j'} \ge k\right)$, let $j < j'$ be maximal and such that $m_{j+1} = k$, then $ \frac{\overline{p}^j_{2k}}{2} = j+1$ by \eqref{l-j-k}. 
		By taking the sum of \eqref{p-j-k+i} on $i$, we deduce that:
		\begin{equation}
			\overline{p}^{j'}_{2k + 2i} = \overline{p}^{j' - \frac{\overline{p}^j_{2k}}{2}}_{2i} \circ T_{\Omega}^{\frac{\overline{p}^{j'}_{2k}}{2}} + \overline{p}^{j'}_{2k}
		\end{equation}
		on the set $\left(m_j \ge k\right)$ for all $i, j, k \ge 0$.
		Therefore, on $\Omega' = \left(\forall k, \overline{p}_k < +\infty\right)$, for all $i,k \ge 0$ and for all $j' \ge \frac{\overline{p}_{2k}}{2}$ we have $\overline{p}^{j'}_{2i} \circ T_{\Omega}^{\frac{\overline{p}_{2k}}{2}} \le \overline{p}_{2k + 2i} < +\infty$.
		Hence $T_{\Omega}^{\frac{\overline{p}_{2k}}{2}}(\Omega') \subset \Omega'$ so $\Omega'$ is $T_\Omega^{p_0 + p_1 \over 2}$-invariant.
	\end{proof}

	With these notations, we can properly study the probabilistic behaviour of the sequence $(m_j)$.
	We remind that for $T_\Omega$ as in \eqref{T-omega-def} and for all $j$, the $\sigma$-algebra $T_\Omega^{j + 1 *} \mathcal{A}_\Omega$ is generated by the data of $(\gamma_k)_{k \ge 2j + 2}$ and $(\tau_k)_{k \ge j + 1}$.
	Therefore $\mathcal{P}_j \vee T_\Omega^{j + 1 *} \mathcal{A}_\Omega$ is generated dy the data of $(\gamma_k)_{k \neq 2j+1}$ and $(\tau_k)_{k \neq j}$.
	
	\begin{Lem}\label{lem:loi-de-m}
		Let $\PP$ be a probability distribution on $\Omega$ such that $(\gamma_{2 k})_{k \in\NN}$ is independent of the joint data of $(\gamma_{2 k + 1})_{k \in \NN}$ and $(\tau_k)_{k \in\NN}$ and $((\gamma_{2 k + 1})_{k \in \NN}, (\tau_k)_{k \in \NN}) \sim \mu^{\otimes \NN}\otimes\mathcal{U}^{\otimes \NN}$.
		Then for all $j, k \in \NN$,
		\begin{gather}
			\PP(m_{j+1} = m_j + 1 \mid m_0, \dots, m_j) = (1 - 2 \rho)\label{estim-forward} \\ 
			\PP\left( m_{j+1} < m_j - k \,\middle|\, m_0, \dots, m_j \right) = 2\rho \left(\frac{\rho}{1 - 2 \rho}\right)^k \mathds{1}(k < m_j). \label{estim-backtrack}
		\end{gather}
	\end{Lem}
	
	\begin{proof}
		Let $A_j$ and $B_j$ be as defined in \eqref{Aj-def} and \eqref{Bjk-def}.
		First note that by construction of $\PP$, the conditional distribution of $\gamma_{2j+1}, \tau_j$ with respect to $\mathcal{P}_j \vee T_{\Omega}^{j+1*}\mathcal{A}_\Omega$ is constant equal to $\mu \otimes \mathcal{U}$. 
		Moreover the random variable $\gamma^{p^j}_{2m_j}$ is $\mathcal{P}_j$-measurable and $\gamma^{p^j}_{2m_j+1}$ is equal to $\gamma_{2j+1}$. 
		Therefore, 
		\begin{equation*}
			\PP\left(\gamma^{p^j}_{2m_j} \AA \gamma^{p^j}_{2m_j+1} \AA \gamma^{p^j}_{2m_j+2}\,\middle|\,\mathcal{P}_j, T_{\Omega}^{j+1*}\mathcal{A}_\Omega\right) = \mu\left\{s \in \Gamma\,\middle|\, \gamma^{p^j}_{2m_j} \AA \, s \, \AA \gamma^{p^j}_{2m_j+2}\right\}.
		\end{equation*}
		Moreover, by independence, for all $\mathcal{P}_j \vee (T_{\Omega}^{j+1})^*\mathcal{A}_\Omega \vee \langle \gamma_{2j +1} \rangle$-measurable map $t_j: \Omega \to [0,1]$,
		\begin{equation*}
			\PP\left(\tau_j < t_j\,\middle|\,\mathcal{P}_j, T_{\Omega}^{j+1*}\mathcal{A}_\Omega, \gamma_{2j+1}\right) = t_j.
		\end{equation*}
		Let $t_j:= \frac{(1-2\rho)\mathds{1}(\gamma^{p^j}_{2m_j} \AA \gamma^{p^j}_{2m_j+1} \AA \gamma^{p^j}_{2m_j+2})}{\mu\left\{s \in \Gamma\,\middle|\, \gamma^{p^j}_{2m_j} \AA s \AA \gamma^{p^j}_{2m_j+2}\right\}}$, the Schottky property for $\mu$ guarantees that $t_j \le 1$.
		Note that $A_j = (\tau_j < t_j)$
		Hence,
		\begin{equation}\notag
			\PP(A_j \mid \mathcal{P}_j, T_{\Omega}^{j+1*}\mathcal{A}_\Omega) = \EE(t_j \mid \mathcal{P}_j, T_{\Omega}^{j+1*} \mathcal{A}_\Omega) = 1-2\rho.
		\end{equation}
		This yields \eqref{estim-forward}.
		
		Now we prove \eqref{estim-backtrack}.
		Let $j \ge 0$.
		For all $k\ge 0$ let $l_k^j = \max\{j' \le j\mid  m_{j'} \le k\}$.
		Note that $l_k$ is measurable for the algebra generated by $(m_0, \dots, m_j)$ and by \eqref{l-j-k} $l_k^j = \frac{\overline{p}^j_{2k+1} - 1}{2}$ for all $k \le m_j$.
		By \eqref{p-j'}, we know that for all fixed $l$,
		\begin{equation*}
			(l_k^j = l)  = (m_l = k) \cap A_l \cap T_{\Omega}^{-l-1}(S^{j-l})
		\end{equation*}
		Since $(m_l = k) \in \mathcal{P}_l$, the distribution of $(\gamma_{2l+1}, \tau_l)$ with respect to the event $(l_k^j = l)$ is almost surely equal to its distribution with respect to $A_l$. 
		Let us define
		\begin{equation*}
			C_{l,k}^j := \left(\tau_l < \frac{(1-3\rho)\mathds{1}\left(\gamma_{2l+1} \AA \gamma_{2l + 2}\cdots \gamma_{2j+2}\right) \mathds{1}\left(\gamma^{p^l}_{2k} \AA \gamma_{2l+1} \AA \gamma_{2l+2}\right)}{\mu\{s\mid  s\, \AA \gamma_{2l + 2}\cdots \gamma_{2j+2} \cap \gamma^{p^l}_{2k} \AA \, s\, \AA \gamma_{2l+2}\}}\right).
		\end{equation*}
		We claim that $C_{l,k}^j \cap (l_k^j = l) \subset A_l$. 
		It is obvious for the indicator function part because $\gamma_{2k}^{p^j} = \gamma_{2m_l}^{p^l}$ when $(l_k^j = l)$ by \eqref{past-p} and \eqref{futur-p} and $\gamma_{2l+1}, \gamma_{2l+2} = \gamma^{p^l}_{2m_l + 1}, \gamma^{p^l}_{2m_l + 2}$, so there remains to show that:
		\begin{equation*}
			\mu\left\{s \,\middle|\,  s\, \AA \gamma_{2l + 2}\cdots \gamma_{2j+2} \cap \gamma^{p^j}_{2k} \AA \, s\, \AA \gamma_{2l+2}\right\} \ge \frac{1-3\rho}{1-2\rho} \mu\left\{s\,\middle|\, \gamma^{p^j}_{2k} \AA s \AA \gamma_{2l+2}\right\}.
		\end{equation*}
		Let $A := \left\{s \,\middle|\, \gamma^{p^j}_{2m_j} \AA s \AA \gamma^{p^j}_{2m_j+2}\right\}$ and $B := \{s\mid  s\, \AA \gamma_{2l + 2}\cdots \gamma_{2j+2}\}$. 
		By the Schottky property, $\mu(A) \ge 1 - 2 \rho$ and $\mu(B) \ge 1 - \rho$.
		Hence $\mu (A \cap B) \ge \mu(A) - \rho = \mu(A) \left( 1 - \frac{\rho}{\mu(A)} \right) \ge \mu(A) \frac{1 - 3 \rho}{1-2\rho}$.
		This proves that $C_{l,k}^j \cap (l_k = l) \subset A_l$.
		Moreover $\PP(C_{l,k}\mid \mathcal{P}_l \vee (T_{\Omega}^{l+1})^*\mathcal{A}_\Omega) = 1-3\rho$ because $\gamma_{2l+1} \sim \mu$.
		Therefore, $\PP(C_{l,k}^j\mid \mathcal{P}_l, (T_{\Omega}^{l+1})^*\mathcal{A}_\Omega, A_l) = \frac{1-3\rho}{1-2\rho}$.
		Furthermore, by \eqref{p-k+m} the data of $m_0, \dots, m_l$ is $\mathcal{P}_l$-measurable and $m_{i + l + 1} = m_{i} \circ T_{\Omega}^{l+1} + m_l + 1$ for all $i \ge 0$ on $A_l$.
		Note also that, on the level sets of $\left(m_0, \dots, m_j\right)$ and for all $k < m_j$,
		\begin{equation}
			(m_{j+1} < m_{j} - k) = \Omega \setminus A_j \setminus \bigcup_{k' = m_j - k}^{m_j - 1} C_{l_{k'}, k'}.
		\end{equation}
		Moreover the $(l_{k'})_{0 \le k' < m_j}$ are distinct by definition and, on each level set of the random variable $(m_0, \dots, m_j)$ (seen as a measurable map) and for all $k' \le m_j$, the event $C_{l_{k'}, k'}$ is independent of $\mathcal{P}_{l_k} \vee T_{\Omega}^{l_k+1*}\mathcal{A}_\Omega$ and has conditional probability $\frac{1-3\rho}{1-2\rho}$.
		Therefore, on each level set of $(m_0, \dots, m_j)$ and for all $k \le m_j$, the event $\bigcap_{k' = m_j - k}^{m_j - 1}(\Omega\setminus C_{l_{k'}, k'})$ has conditional probability $\left(\frac{\rho}{1-2\rho}\right)^k$ and is independent of $A_j$.
		Moreover, $\PP(m_{j+1} < m_j - k) = 0$ when $k \ge m_j$.
		From that we deduce \eqref{estim-backtrack}. 
	\end{proof}
	
	We deduce from Lemma \ref{lem:loi-de-m} that $\PP(\Omega') = 1$ when $\rho < 1/5$.
	Let us give more precisions in the following corollary.
	
	\begin{Cor}\label{cor:indep+exp-moment}
		Let $\PP$ be a probability distribution on $\Omega$ such that $(\gamma_{2k})_{k \in\NN}$ is independent of the joint data of $(\gamma_{2k+1})_{k \in\NN}$ and $(\tau_k)_{k \in\NN}$ and $\left((\gamma_{2k+1})_{k \in\NN}, (\tau_k)_{k \in\NN}\right) \sim \mu^{\otimes\NN}\otimes\mathcal{U}^{\otimes\NN}$.
		Then the data of $(m_j)_{j \ge 0}$ is independent of the data of $(\gamma_{2k})_{k \ge 0}$.
		Moreover $m_j \to +\infty$ almost surely (\ie $\PP(\Omega') = 1$) and there exist constants $C, \beta >0$ such that:
		\begin{equation}\label{moment-de-p}
			\forall k,t \in\NN, \PP(\overline{p}_{3} \ge t) \le C e^{-\beta t}.
		\end{equation}
	\end{Cor}
	
	\begin{proof}
		First we prove the independence result. 
		Note that \eqref{estim-forward} and \eqref{estim-backtrack} together with the fact that $m_0 = 0$, determine the distribution law of the sequence $\left(m_j\right)_{j \ge 0}$. 
		Let $\kappa$ be the distribution of $\left(m_j\right)_{j \ge 0}$.
		Let $X = \Gamma^{2\NN}$ and $Y = \Gamma^{2\NN + 1}\times [0,1] ^\NN$.
		Note that $\Omega$ naturally identifies with $X \times Y$.
		Let $\mu_Y = \mu^{\otimes\NN} \otimes \mathcal{U}^{\otimes\NN}$ and let $\PP_x = \delta_x \otimes \mu_Y$ for all $x \in X$. 
		The fact that $(\gamma_{2k})_{k \in\NN}$ is independent of the joint data of $(\gamma_{2k+1})_{k \in\NN}$ and $(\tau_k)_{k \in\NN}$ means that there exists a probability distribution $\mu_X$ on $X$ such that
		\begin{equation*}
			\PP = \mu_X \otimes \mu_Y = \int_{x \in X} \PP_x d\mu_X(x).
		\end{equation*}
		For all $x$, Lemma \ref{lem:loi-de-m} applied to $\PP_x$ yields
		\begin{equation*}
			((\gamma_{2k})_{k \ge 0}, (m_j)_{j\ge 0}) \sim \int_{x \in X} (\delta_x \otimes \kappa) d\mu_X(x) = \mu_X \otimes \kappa.
		\end{equation*}
		This means that the sequences $(\gamma_{2k})_{k \ge 0}$ and $(m_j)_{j \ge 0}$ are independent.
		
		Let us now prove that each $p_k$ has a finite exponential moment.
		Let $x$ be a random variable that takes values in $\{1\}\cup \mathbb{Z}_-$ such that
		\begin{equation*}
			\PP(x = 1) = 1 - 2\rho \quad \text{and} \quad \forall k < 0,\; \PP(x = k) = 2\rho \frac{1-3\rho}{\rho} \left(\frac{\rho}{1-2\rho}\right)^{k}.
		\end{equation*}
		Then $\EE(x) = \frac{(1-2\rho)(1-5\rho)}{1-3\rho} > 0$ and $x$ has finite exponential moment. 
		It follows from Lemma \ref{lem:ldev-classique}, that $m_j \to +\infty$ almost surely.
		Moreover there exist constants $C, \beta > 0$ such that
		\begin{equation*}
			\forall j \ge 0, \PP\left(m_j \le 1 \right) \le C e^{-\beta j}.
		\end{equation*}
		Hence by taking the sum, we get
		\begin{equation*}
			\forall t \ge 0, \PP\left(\exists j \ge t,\; m_j \le 1 \right) \le \frac{C}{1-e^{-\beta}} e^{-\beta t}.
		\end{equation*}
		Then by \eqref{l-k},
		\begin{equation*}
			\forall t \ge 0, \PP\left(\overline{p}_3 \ge 2j+1\right) \le \frac{C}{1-e^{-\beta}} e^{-\beta t}.
		\end{equation*}
		This yields \eqref{moment-de-p}.
	\end{proof}

	\begin{Lem}[Renewal Property]\label{lem:markov-prop}
		Let $\eta$ be any probability distribution on $\Gamma$ and let $\PP:= \left(\eta \otimes \mu\right)^{\odot \NN} \otimes \mathcal{U}^{\otimes\NN}$.
		Then, for all $k \ge 0$, the data of $(\widetilde{\gamma}^p_{k'})_{k' \le 2k}$ is independent of the data of $(\widetilde{\gamma}^p_{k'})_{k' \ge 2k+2}$ and
		\begin{equation}\label{schottky-conditioned}
			\forall A \in \mathcal{A}_\Gamma, \; \PP\left(\gamma^p_{2k+1} \in A\mid (\widetilde{\gamma}^p_{k'})_{k' \neq 2k+1}\right) = \frac{\mu \{s \in A\mid  \gamma^p_{2k} \AA \,s \,\AA \gamma_{\overline{p}_{2k+2}}\}}{\mu \{s \in \Gamma\mid  \gamma^p_{2k} \AA \, s \, \AA \gamma_{\overline{p}_{2k+2}}\}}.
		\end{equation}
		Therefore $\left(\widetilde{\gamma}^p_{2k + 2}\right)_{k \ge 0}$ is i.i.d. and independent of $\widetilde{\gamma}^p_{0}$. 
	\end{Lem}

	\begin{proof}
		Let $S$ be as in Lemma \ref{lem:descri-pivot}.
		For all $k \in \NN$, we write $\mathcal{P}'_k$ for the $\sigma$-algebra generated by the joint data of $\overline{p}_{2k+1}$ and $\left(\widetilde{\gamma}^p_{k'}\right)_{k'\le 2k}$. 
		We remind that for all $j$, we write $\mathcal{P}_j$ for the $\sigma$-algebra generated by $(\gamma_{n})_{n \le 2j}$ and $(\tau_k)_{k < j}$.
		
		First we show that for all $k \ge 0$, the conditional distribution of $T_{\Omega}^\frac{\overline{p}_{2k+2}}{2}(\omega)$ with respect to $\mathcal{P}'_k$ is almost surely equal to the normalized restriction of $\PP$ to $S$, defined by the formula $\frac{\mathds{1}(S)}{\PP(S)} \PP$.
		Let $l,k \ge 0$. By \eqref{p-zero}
		\begin{equation*}
			(\overline{p}_{2k+1} = 2l+1) = (m_l = k) \cap A_l \cap T_\Omega^{-l-1}(S).
		\end{equation*}
		The event $(m_l = k)$ only depends on $\mathcal{P}_{l}$, which is independent of $(T_\Omega^{l+1})^*  \mathcal{A}_{\Omega}$.
		Moreover, $\PP(A_l \mid  \mathcal{P}_{l} \vee (T_\Omega^{-l-1})^*  \mathcal{A}_{\Omega}) = 1 - 2 \rho$ so $A_l$ is independent of the joint data of $\mathcal{P}_{l}$ and $(T_\Omega^{-l-1})^*  \mathcal{A}_{\Omega}$.
		Reciprocally $(T_{\Omega}^{l+1})^*  \mathcal{A}_{\Omega}$ is independent of the joint data of $\mathcal{P}_{l}$ and $A_l$.
		Therefore the conditional distribution of $T_\Omega^{l+1}$ with respect to the $(\overline{p}_{2k+1} = 2l+1)$ and $\mathcal{P}_\frac{\overline{p}_{2k+1}-1}{2}$ is almost surely equal to the conditional distribution of $T_\Omega^{l+1}$ with respect to $T_\Omega^{-l-1}(S)$ and equals $\frac{\mathds{1}(S)}{{T_{\Omega}^{l+1}}_*\PP(S)} {T_{\Omega}^{l+1}}_*\PP$.
		Moreover, $\frac{\mathds{1}(S)}{{T_{\Omega}^{l+1}}_*\PP(S)} {T_{\Omega}^{l+1}}_*\PP = \frac{\mathds{1}(S)}{\PP(S)} \PP$ because $\PP$ is $T$-invariant.
		We conclude by integrating over all possible values of $l \ge 0$.
		
		By \eqref{p-zero} (and since $\gamma_{2l+1}$ is independent of $\mathcal{P}_{l} \vee T_\Omega^{l+1 }\mathcal{A}_{\Omega}$ by construction of $\PP$), the conditional distribution of $\gamma_{2l+1}$ with respect to $(\overline{p}_{2k+1} = 2l+1)$ and $\mathcal{P}_{l} \vee T_\Omega^{l + 1*} \mathcal{A}_{\Omega}$ is almost surely equal to its conditional distribution with respect to $A_l$, which is characterized by \eqref{schottky-conditioned}.
		We conclude by integrating over all possible values of $l$ to prove \eqref{schottky-conditioned}.
		
		By Lemma \ref{lem:descri-pivot}, $p_{2k + 2 + i} = p_i \circ T_{\Omega}^{\frac{\overline{p}_{2k+2}}{2}}$ for all $k, i \ge 0$.
		From that we deduce that $\widetilde{\gamma}^p_{2k+2 + i} = \widetilde{\gamma}^p_{2k+2 + i} \circ T_{\Omega}^{\frac{\overline{p}_{2k+2}}{2}}$.
		Therefore, the data of $(\widetilde{\gamma}^p_{2k+2+i})_{i \ge 0}$ is independent of $\mathcal{P}'_k$ has distribution $(\bigotimes_{i \ge 0} \widetilde{\gamma}^p_{i})_* ( \PP {\mathds{1}(S)}/{\PP(S)} )$.
		This distribution does not depend on $k$ so the sequence $(\widetilde{\gamma}^p_{2k+2})_{k \ge 0}$ is identically distributed. 
		Moreover, for all $k \ge 0$, the data of $(\widetilde{\gamma}^p_{2k})_{k' \le k}$ is $\mathcal{P}'_k$-measurable so it is independent of $\widetilde{\gamma}^p_{2k+2}$. 
		Therefore, the random variables $(\widetilde{\gamma}^p_{2k})_{k \ge 0}$ are globally independent.
	\end{proof}

	\begin{Lem}\label{lem:ali-pivot}
		With the Notations of Definitions \ref{def:pivot-alg} and \ref{def:temps-pivots}.
		The alignment $\gamma^p_{2k+1} \AA \gamma_{\overline{p}_{2k+2}}$ holds on $\Omega'$ for all $k \ge 0$ and there exists a measurable function defined on $\Omega'$ and valued in the space of finite families of odd integers $1 = c_1^k < c_2^k \dots < c_{j_k}^k = p_{2k+2}$ such that for all $1 \le i < j_k$, 
		\begin{equation}\label{atilde-lemme}
			\gamma_{\overline{p}_{2k+2}} \cdots \gamma_{\overline{p}_{2k+2} + c_i^k - 1} \AA \gamma_{\overline{p}_{2k+2} + c_i^k} \AA \gamma_{\overline{p}_{2k+2} + c_i^k + 1} \cdots \gamma_{\overline{p}_{2k+2} + c_{i+1}^k - 1}.
		\end{equation}
	\end{Lem}
	
	\begin{proof}
		For all $ \ge 0$, let $l_k:= \frac{\overline{p}_{2k+1} - 1}{2} = \max \{j\mid  m_j \le k\}$ and $j_k:= \# \{j > l_k \mid  m_{j} = k + 1\}$.
		These random integers are well defined and measurable on $\Omega'$.
		Moreover $j_k > 0$ because $m_{l_k + 1} = k+1$
		We define by induction $b_1^k := l_k + 1$ and for $i \ge 1$,
		\begin{equation*}
			b_{i+1}^k := \min\{j > b_i^k \mid  m_{j} = k + 1\}.
		\end{equation*}
		Let $c_i^k = 2b_i^k + 1 - \overline{p}_{2k+2}$ for all $1 \le i \le j_k$.
		We get $\gamma^p_{2k+ 1} \AA \gamma_{\overline{p}_{2k+2}}$ because $A_{l_k}$ holds (\ie $\omega \in A_{l_k(\omega)}$ for all $\omega \in \Omega'$). 
		
		Let $1 \le i < j_k$.
		It follows from the maximality of $l_k$ that $m_{b_i^k + 1} \ge k + 1 = m_{b_i^k}$. 
		Therefore $A_{b_i^k}$ holds so $\gamma^{p^{b_i^k}}_{2k+2} \AA \gamma_{2 b_i^k + 1}$.
		Then by \eqref{l-j-k}, and because $b_{i}^{k} \ge l_k$, we have $\overline{p}^{b_i^k}_{2k+2} = \overline{p}_{2k+2}$ so $\gamma^{p^{b_i^k}}_{2k+2} = \gamma_{\overline{p}_{2k+2}} \cdots \gamma_{2b_i^k} = \gamma_{\overline{p}_{2k+2}} \cdots \gamma_{\overline{p}_{2k+2} + c_i^k - 1}$. 
		This proves the first alignment condition in \eqref{atilde-lemme}.
		
		The fact that $b_{i+1}^k - 1 > l_k$ yields that $m_{b_{i+1}^k - 1} > k$ and $m_j$ does not have steps of size $0$ so necessarily, $k + 1 = m_{b_{i+1}^{k}} < m_{b_{i+1}^k - 1}$.
		Hence, $B_{b_{i + 1}^k - 1, k+1}$ holds so $\gamma^{p^{b_{i}^k-1}}_{2k+3} \AA \gamma_{\overline{p}^{b_{i+1}^k-1}_{2k+4}} \cdots \gamma_{2b_i^k}$. 
		Moreover, $m_j > k + 1$ for all $b_{i}^k < j < b_{i+1}^k$ so $\overline{p}^{b_{i+1}^k - 1}_{2k+3} = 2b_{i}^k + 1$ by \eqref{l-j-k}. 
		Therefore $\gamma_{\overline{p}^{b_{i+1}^k-1}_{2k+4}} \cdots \gamma_{2b_i^k} = \gamma_{\overline{p}_{2k+2} + c_i^k + 1} \cdots \gamma_{\overline{p}_{2k+2} + c_{i+1}^k - 1}$.
		This proves the second alignment condition in \eqref{atilde-lemme}.
	\end{proof}

	Now the main Theorem follows directly from the above Lemmas.
	
	\begin{proof}[Proof of Theorem \ref{th:ex-piv}]
		Let $(\hat{v}_k)_{k \ge 0}$ be as in Lemma \ref{lem:first-part}.
		Let $\left(\tau_k\right)_{k \ge 0} \sim \mathcal{U}^{\otimes \NN}$ be independent of $\left( \widetilde{\gamma}^v_k \right)_{k \ge 0}$.
		Let $\eta$ be the distribution law of $\gamma^{\hat{v}}_0$. 
		Then by Lemma \ref{lem:first-part}, $(\gamma^{\hat{v}}_k)_{k \ge 0} \sim \left(\eta \otimes \mu\right)^{\odot \NN}$.
		Let $(\tau_k)_{k \ge 0} \in[0,1]^\NN$ be uniformly distributed, globally independent and independent of the data of $\widetilde(\gamma^{\hat{v}}_k)_{k \ge 0}$.
		Let $\left(p_k\right)_{k \ge 0}$ be the length of the pivotal blocks associated to the sequence $\left(\gamma^{\hat{v}}_k\right)_{k \ge 0}$ with weights $(\tau_k)_{k \ge 0}$.
		Points \eqref{indep-v} to \eqref{item:hat-v-is-stopping-time} follow directly from Lemma \ref{lem:first-part}. 
		Let us now prove points \eqref{item:tail} to \eqref{item:tjr-piv}.
		
		We know by Corollary \ref{cor:indep+exp-moment} that the data of $(p_k)_{k \ge 0}$ is independent of $\left(\gamma^{\hat{v}}_{2k}\right)_{k \ge 0}$.
		Moreover, the joint conditional distribution of $\left(\widetilde{\gamma}^{\hat{v}}_{2k}\right)_{k \ge 0}$, $(v_k)_{k \ge 0}$, $\left(w_k\right)_{k \ge 0}$ and $(u_k)_{k \ge 0}$ with respect to the joint data of $\left(\gamma^{\hat{v}}_{k}\right)_{k \ge 0}$ and $(\tau_k)_{k \ge 0}$ only depends on $\left(\gamma^{\hat{v}}_{2k}\right)_{k \ge 0}$.
		Moreover, the sequence $(p_k)_k$ is given by a measurable function of $\left(\gamma^{\hat{v}}_{k}\right)_{k \ge 0}$ and $(\tau_k)_{k \ge 0}$. 
		Therefore the joint data of $\left(\widetilde{\gamma}^{\hat{v}}_{2k}\right)_{k \ge 0}$, $(v_k)_{k \ge 0}$ and $\left(w_k\right)_{k \ge 0}$ is independent of $\left(p_k\right)_{k \ge 0}$.
		This proves point \eqref{item:tail}.
		
		Let us now prove points \eqref{item:indep-p} and \eqref{item:indep}. 
		For all $k$, we have $p_{2k+1} = 1$ by construction and $\overline{p}_{2k + 1}$ is odd by \eqref{l-k} in Lemma \ref{lem:descri-pivot}. 
		Moreover $p_0 = \overline{p}_{1}$ which is odd and for $k \ge 1$, $p_{2k} = \overline{p}_{2k + 1} - \overline{p}_{2k - 1} -1$, which is also odd.
		If follows that all the $p_k$'s are odd and therefore positive.
		Point \eqref{item:indep} in Theorem \ref{th:ex-piv} directly follows from Lemma \ref{lem:markov-prop}.
		Taking the push-forward by $L^{\otimes\NN}$, we get that $(p_{2k+2})_k$ is i.i.d. and independent of $p_0$.
		
		Point \eqref{item:exp-p}, follows directly from \eqref{moment-de-p} in Corollary \ref{cor:indep+exp-moment} because $\overline{p}_3 = p_0 + p_1 + p_2$ with $p_1 = 1$.
		
		Points \eqref{item:ali} and \eqref{item:ali-rec} follow directly from Lemma \ref{lem:ali-pivot}.
		
		Point \eqref{item:tjr-piv} follows from \eqref{schottky-conditioned} in Lemma \ref{lem:markov-prop} and from the fact that the conditional distribution of $\widetilde{\gamma}^{\check{p}}_{2k+1}$ with respect to $\left(\gamma^{\hat{v}}_{k}\right)_{k \ge 0}$ and $(\tau_k)_{k \ge 0}$ only depends on $\gamma^{\check{p}}_{2k + 1}$.
	\end{proof}

	\section{limit theorems for a product of random matrices}\label{results}

	The present section is dedicated to the demonstration of Theorems \ref{th:escspeed}, \ref{th:cvspeed}, \ref{th:eigenspace} and \ref{th:slln} announced in the introduction of the present paper.
	For that purpose, we use Theorem \ref{th:pivot} as a black box. 
	
	In the following discussion $\nu$ denotes a probability distribution over $\End(E)$ that is $(\alpha, \eps, m)$-contracting in the sense of Definition \ref{def:quant-contraction} for some constants $\alpha, \eps > 0$ and $m \in \NN$.

	\subsection{Preliminaries about one-sided large deviations}

	Let us state basic results about sequences that satisfy large deviations inequalities.
	
	\begin{Def}[One-sided large deviations]
		Let $(\overline{x}_n)\in\RR^\NN$ be a random sequence and $\lambda \in\RR \cup \{+\infty\}$. 
		We say that $(\overline{x}_n)$ satisfies large deviations inequalities below the speed $\lambda$ if
		\begin{equation}\label{def-ldev}
			\forall \alpha < \lambda, \; \exists C, \beta > 0, \; \forall n,\;\PP(\overline{x}_n \le \alpha n) \le C e^{-\beta n}.
		\end{equation}
	\end{Def}
	
	Note that \eqref{def-ldev} admits the following equivalent reformulations:
	\begin{gather*}
		\forall \alpha < \lambda, \; \exists C, \beta > 0, \; \forall n_0,\;\PP(\exists n \ge n_0,\, \overline{x}_n \le \alpha n) \le C e^{-\beta n}, \\
		\forall \alpha < \lambda, \; \exists n_0, \beta > 0, \; \forall n \ge n_0,\;\PP(\overline{x}_n \le \alpha n) \le e^{-\beta n}, \\
		\forall \alpha < \lambda, \; \limsup_{n_0 \to +\infty} \frac{1}{n}\log\PP(\overline{x}_n \le \alpha n) < 0, \\
		\forall \alpha < \lambda, \; \exists \beta > 0, \; \EE\left(\exp(\beta\max\{n\mid \overline{x}_n \le \alpha n\})\right) < +\infty.
	\end{gather*}

	\begin{Lem}[Classical large deviations inequalities]\label{lem:ldev-classique}
		Let $(x_n) \in \RR^\NN$ be a random i.i.d. sequence and let $\beta_0 > 0$ be a constant. 
		Assume that $\EE(e^{- \beta_0 x_0}) < +\infty$. 
		Then the random sequence $(\overline{x}_n)_n$ satisfies large deviations inequalities below the speed $\EE(x_0) \in \RR \cup \{+\infty\}$.
	\end{Lem}
	
	\begin{proof}
		The first step is to prove the following:
		\begin{equation}\label{mon-conv}
			\forall \alpha < \EE(x_0),\; \exists \beta > 0,\; \EE(e^{-\beta x_0}) < e^{-\beta \alpha}.
		\end{equation}
		This follows from the monotonous convergence Theorem. 
		For all $t \in\RR$, the map $\beta \mapsto \frac{1 - e^{-\beta t}}{\beta}$ is non-increasing on $(0, + \infty)$ and $\lim_{\beta \to 0} \frac{1 - e^{-\beta t}}{\beta} = t$.
		Moreover $\frac{1 - e^{-\beta x_0}}{\beta}$ is in $\ELL^1$, therefore $\EE\left(\frac{1 - e^{-\beta x_0}}{\beta}\right)$ is non-increasing and admits $\EE(x_0)$ as its limit.
		Hence, for all $\alpha < \EE(x_0)$ there exists $\beta >0$ such that $\EE\left(\frac{1 - e^{-\beta x_0}}{\beta}\right) > \alpha$, and by linearity $\EE(e^{-\beta x_0}) < 1 - \beta \alpha \le e^{-\beta \alpha}$.
		This proves \eqref{mon-conv}.
		
		Let $\alpha < \EE(x_0)$ and let $\beta >0$ be as in \eqref{mon-conv}. 
		By independence, $\EE(e^{-\beta \overline{x}_n}) = \EE(e^{-\beta x_0})^n$ for all $n$.
		Then by Markov inequality,
		\begin{equation*}
			\PP\left(\overline{x}_n \le \alpha n\right) \le \frac{\EE(e^{-\beta \overline{x}_n})}{e^{\beta \alpha n}} = \left(\frac{\EE(e^{-\beta x_0})}{e^{\beta \alpha}}\right)^n. \qedhere
		\end{equation*}
	\end{proof}
	
	\begin{Lem}\label{lem:ldev-shift}
		Let $(\overline{x}_n)_{n \ge 0}\in\RR^\NN$ and  $(l_n)_{n \ge 0} \in \NN_{\ge 0}^\NN$ be random sequences and let $\lambda \le + \infty$ and $C_0, \beta_0 > 0$ be constants.
		Assume that $(\overline{x}_n)_n$ satisfies large deviations inequalities below the speed $\lambda$ and that $\EE(e^{\beta_0 l_n}) \le C_0$ for all $n$.
		Then $(\overline{x}_{n} - l_n)_n$ satisfies large deviations inequalities below the speed $\lambda$.
	\end{Lem}
	
	\begin{proof}
		Let $\alpha < \alpha' < \lambda$. 
		By assumption, there exist constants $C', \beta' > 0$ such that $\PP(\overline{x}_n\le \alpha' n) \le C' e^{-\beta' n}$ for all $n$.
		Moreover, by Markov's inequality, $\PP\left(l_n \ge (\alpha' -\alpha) n\right) \le C_0 e^{-\beta_0(\alpha'-\alpha n)}$ for all $n$.
		Therefore,
		\begin{equation*}
			\forall n, \; \PP\left(\overline{x}_{n} - l_n \le \alpha n\right) \le C_0 e^{-\beta_0(\alpha'-\alpha n)} + C' e^{-\beta' n}. \qedhere
		\end{equation*}
	\end{proof}
	
	\begin{Lem}\label{lem:compo-ldev}
		Let $(\overline{x}_n)_{n \ge 0}\in\RR^\NN$ and  $(\overline{w}_n)_{n \ge 0} \in \NN_{\ge 0}^\NN$ be random sequences and let $0 \le \delta, \lambda \le + \infty$ be constants.
		Assume that $(\overline{x}_n)_n$ satisfies large deviations inequalities below the speed $\lambda$ and that $(\overline{w}_n)_n$ satisfies large deviations inequalities below the speed $\delta$. 
		Then $\left(\overline{x}_{\overline{w}_n}\right)_n$ satisfies large deviations inequalities below the speed $\lambda\delta$.
	\end{Lem}
	
	\begin{proof}
		Let $\alpha < \lambda \delta$ and let $0 < \alpha_0 < \lambda$ and $0 < \alpha_1 < \delta$ be such that $\alpha \le \alpha_0 \alpha_1$.
		Let $C_0, \beta_0, C_1, \beta_1 > 0$ be such that $\PP(\overline{x}_n \le \alpha_0 n) \le C_0e^{-\beta_0 n}$ and $\PP(\overline{w}_n \le \alpha_1 n) \le C_1 e^{-\beta_1 n}$ for all $n \ge 0$. 
		Then for all $n \ge 0$, 
		\begin{align*}
			\PP\left(\overline{x}_{\overline{w}_n} \le \alpha n\right) & \le \PP\left({\overline{w}_n} \le \alpha_1 n \cup \exists k \ge \alpha_1 n, \overline{x}_k \le \alpha_0 k\right) \\
			& \le C_1 e^{-\beta_1 n} + \sum_{k \ge \alpha_1 n} C_0e^{-\beta_0 k} \le C_1 e^{-\beta_1 n} + \frac{C_0}{1- e^{-\beta_0}} e^{-\beta_0 \alpha_1 n}. \qedhere
		\end{align*}
	\end{proof}
	
	\begin{Lem}\label{lem:sumsum}
		Let $(\overline{x}_n)_{n \ge 0}$ be a random sequence of non-negative numbers that satisfies large deviations inequalities above a speed $\lambda$. Let $w$ be a random non-negative integer that has a finite exponential moment. 
		Then without any independence assumption, the random variable $\overline{x}_w$ has a finite exponential moment.
	\end{Lem}
	
	\begin{proof}
		Let $t \ge 0$, we have $(\overline{x}_w \ge t) \subset (\overline{x}_{\lfloor t/(2\lambda)\rfloor} \ge t) \cup (w > \lfloor t/(2\lambda)\rfloor)$. 
		The probability of the first term decreases exponentially fast in $t$ by the large deviations inequalities and the probability decreases exponentially fast in $t$ because $w$ has a finite exponential moment.
	\end{proof}

	\subsection{The bi-lateral pivoting technique}\label{section:l-r-piv}
	
	In this Paragraph, we illustrate the implications of Point \eqref{pivot:schottky} from Theorem \ref{th:pivot}.
	We remind that \eqref{pivot:schottky} says that for all $k \ge 0$, the conditional distribution of $\gamma^p_{2k+1}$ with respect to $(\widetilde{\gamma}^p_{k'})_{k' \neq 2k+1}$ is almost surely $\frac{1}{4}$-Schottky for the binary relation $\AA^\eps$. 
	It is equivalent to saying that for all $k$ and all random variable $g$, given by the image of $(\widetilde{\gamma}^p_{k'})_{k' \neq 2k+1}$ by a measurable function, $\PP(\gamma^p_{2k+1} \AA^\eps g) \ge 3 / 4$ and $\PP(g \AA^\eps \gamma^p_{2k+1}) \ge 3 / 4$ almost surely.
	
	Let us now describe the left, the right and the cyclical pivoting techniques, which respectively allow us to control the behaviour of the lines, of the columns and of the top eigenspace of $\overline{\gamma}_n$.
	We remind that $\mathcal{U}$ denotes the restriction of the Lebesgue measure to the interval $(0,1)$.

	\begin{Lem}[Pivoting technique]\label{lem:l-r-piv}
		Let $\alpha, \eps > 0$, $m \in\NN$ and let $\nu$ be an $(\alpha, \eps, m)$-contracting probability measure on $\End(E)$ that has positive eventual rank.
		Let $(\gamma_n)_{n \ge 0} \sim \nu^{\otimes\NN}$ and let $(p_n)_{n \ge 0}$ be as in Theorem \ref{th:pivot}.
		Let $(\tau_k)_{k \ge 0} \sim \mathcal{U}^{\otimes\NN}$ be independent of the joint data of $\left(\widetilde{\gamma}^p_k\right)_{k \ge 0}$.
		For all non-zero $f \in E^* \cup \End(E)$, there exists a random integer $l^f \sim \mathcal{G}_{1/4}$, that is independent of $\left(\widetilde{\gamma}^p_{2k}\right)_{k \ge 0}$ and such that:
		\begin{equation}\label{f-ali}
			f \gamma^p_0 \cdots \gamma^p_{2 l^f} \AA^\eps \gamma^p_{2 l^f + 1}
		\end{equation}
		Let $h \in E \cup \End(E)$ be non-zero.
		For all $n \in\NN$, let $q_n := \max\{k \ge 0\mid \overline{p}_{2k} \le n\}$.
		Then for all $n \in \NN$, there exist two random integers $r_n^h \sim \mathcal{G}_{1/4}$ and $c_n \sim \mathcal{G}_{1/2}$, that are independent of $\left(\widetilde{\gamma}^p_{2k}\right)_{k \ge 0}$ and such that:
		\begin{gather}\label{ali-h}
			\gamma^p_{2 q_n - 2r_n^h - 1} \AA^{\eps} \gamma_{\overline{p}_{2q_n - 2r_n^h}} \cdots \gamma_{n-1} h \quad \text{or} \quad r_n^h \ge q_n,
			\\
			\begin{cases}
				\gamma^p_{2 q_n - 2c_n - 1} \AA^{\eps} \gamma_{\overline{p}_{2q_n - 2c_n}} \cdots \gamma_{n-1} \gamma^p_0 \cdots \gamma^p_{2 c_n}\\
				\text{and }\gamma_{\overline{p}_{2q_n - 2c_n - 1}} \cdots \gamma_{n-1} \gamma^p_0 \cdots \gamma^p_{2 c_n} \AA^\eps \gamma^p_{2c_n+1}
			\end{cases}  \quad \text{or} \quad 2c_n + 1 \ge q_n,
		\end{gather}
		Moreover, there exist constants $C, \beta >0$ such that $\EE\left(e^{\beta \overline{p}_{2l^f}}\right) \le C$ for all $f$, $\EE\left(e^{\beta (n- \overline{p}_{2q_n - 2r_n^h})}\right) \le C$, for all $n$ and all $h$ and $\EE\left(e^{\beta (n- \overline{p}_{2q_n - 2c_n - 2} + \overline{p}_{2c_n})}\right) \le C$ for all $n$.
	\end{Lem}
	
	\begin{proof}
		Let $f \in E^* \cup \End(E)$. 
		We construct $l^f$ as:
		\begin{equation*}
			l^f:= \min\left\{k \in\NN\,\middle|\, \tau_k < \frac{3}{4}\frac{\mathds{1}\left(f \gamma^p_0 \cdots \gamma^p_{2 l^f} \AA^\eps \gamma^p_{2 l^f + 1}\right)}{\PP\left(f \gamma^p_0 \cdots \gamma^p_{2 l^f} \AA^\eps \gamma^p_{2 l^f + 1}\,\middle|\,\left(\widetilde{\gamma}^p_{2k}\right)_{k \ge 0}\right)}\right\}.
		\end{equation*}
		Then \eqref{f-ali} holds.
		By point \eqref{pivot:schottky} in Theorem \ref{th:pivot}, $\PP(f \gamma^p_0 \cdots \gamma^p_{2 k} \AA^\eps \gamma^p_{2 k + 1}\mid (\widetilde{\gamma}^p_{k'})_{k'\neq 2k+1}) \ge 3/4$ for all $k$, so $\PP(l^f = k \mid (\widetilde{\gamma}^p_{k'})_{k'\neq 2k+1}, l^f \ge k) = \frac{3}{4}$.
		Therefore, the conditional distribution of $l^f$ with respect to $\left(\widetilde{\gamma}^p_{2k}\right)_{k \ge 0}$ is $\mathcal{G}_{1/4}$.
		If follows from Lemmas \ref{lem:power-compare} and \ref{lem:sumexp} that the random variable $\overline{p}_{2l^f}$ has a finite exponential moment. 
		Let us give more details: 
		\begin{itemize}
			\item The $p_k$'s are independent and each has the distribution of $p_0$, $p_1$ or $p_2$ so $\PP(p_k > t\mid(p_k')_{k' < k}, l^f) \le \max_{0 \le i \le 2} \PP(p_i > t) =: \kappa(t, + \infty)$ for all $t \ge 0$.
			\item Let $\eta$ be the law of $2l^f$. Then $\PP(\overline{p}_{2l^f}) \le \kappa^{*\eta}(t, + \infty)$ for all $t \ge 0$ by Lemma \ref{lem:power-compare}. 
			\item 
			Since the $p_i$'s all have a finite exponential moment, so does $\kappa$ and by Lemma \ref{lem:sumexp}, so does $\kappa^{*\eta}$.
		\end{itemize}
		
		Let $h \in E \cup\End(E)$ and let $n \in\NN$.
		Let us construct $r_n^h$. 
		Since $\left(p_{2k+1}\right)_{k \ge 0}$ is non-random, the data of $q_n$ is measurable for the $\sigma$-algebra generated by $\left(\widetilde{\gamma}^p_{2k}\right)_{k \ge 0}$.
		For all $k \ge 0$, let us define the event:
		\begin{multline*}
			A^h_{n,k} := \left(k \ge q_n\right) \cap \left(\tau_k < \frac{3}{4}\right) 
			\\ 
			\cup \left(k < q_n\right) \cap \left(\tau_k < \frac{3}{4}\frac{\mathds{1}\left(\gamma^p_{2 q_n - 2r_n^h - 1} \AA^{\eps} \gamma_{\overline{p}_{2q_n - 2r_n^h}} \cdots \gamma_{n-1} h\right)}{\PP\left(\gamma^p_{2 q_n - 2 k - 1} \AA^{\eps} \gamma_{\overline{p}_{2q_n - 2k}} \cdots \gamma_{n-1} h\,\middle|\,\left(\widetilde{\gamma}^p_{2j}\right)_{j \ge 0}\right)}\right)
		\end{multline*}
		Then by point \eqref{pivot:schottky} in Theorem \ref{th:pivot}, $\PP\left(A^h_{n,k} \,\middle|\, q_n, \left(\widetilde{\gamma}^p_{j}\right)_{j \neq 2q_n - 2k - 1}, (\tau_j)_{j \neq k}\right) = \frac{3}{4}$. 
		It follows that the $A^h_{n,k}$ are independent and their joint data is independent of $(\widetilde{\gamma}^p_{2j})_{j \ge 0}$.
		Let $r_n^h = \min\{k \ge 0\mid  A_{n,k}^h\}$. 
		Then the conditional distribution of $r_n^h$ with respect to $\left(\widetilde{\gamma}^p_{2k}\right)_{k \ge 0}$ is equal to $\mathcal{G}_{1/4}$.
		
		First, we claim that $\PP\left(n - \overline{p}_{2q_n - 2 k} \ge t\right) \le t \PP(\overline{p}_{2k + 2} \ge t)$ for all $k, n, t \ge 0$. 
		Indeed,
		\begin{align*}
			\PP\left(n - \overline{p}_{2q_n - 2 k} = t\right) & = \sum_{i \ge 0}\PP(\overline{p}_{2i} = n - j \cap \overline{p}_{2i + k} \le n \cap \overline{p}_{2i + 2 k + 2} > n)\\
			& \le \sum_{i \ge 0}\PP(\overline{p}_{2i} = n - j \cap \overline{p}_{2i + 2 k + 2} > n) 
		\end{align*}
		Hence, by independence of the sequence $(p_{i})_{i \ge 0}$, the events $(\overline{p}_{2i} = n - -t)$ and $(p_{2i} + \cdots + {p}_{2i + 2 k + 1} > j)$ are independent for all $i,t,n$.
		Note also that the sequence $\overline{p}_{2i}$ is almost surely increasing, therefore $\sum_{i \ge 0}\PP(\overline{p}_{2i} = n - t) \le 1$ for all $n,t$. 
		Therefore, 
		\begin{align*}
			\PP\left(n - \overline{p}_{2q_n - 2 k} = t\right) & \le \sum_{i \ge 0}\PP(\overline{p}_{2i} = n - j) \PP(p_{2i} + \cdots + {p}_{2i + 2 k + 1} > j) \\
			& \le \max_{i} \PP(p_{2i} + \cdots + {p}_{2i + 2 k + 1} > t)
		\end{align*}
		Let $C, \beta >0$ be such that $\EE(e^{\beta p_{2k}}) \le C$ for all $k \ge 0$.
		We remind that $p_{2k+1} = m$ almost surely and for all $k$ by \eqref{pivot:schottky}.
		Then, by independence, $\EE(e^{\beta(p_{2i} + \cdots + {p}_{2i + 2 k + 1})}) \le (e^{\beta m} C)^k$ for all $i, k \ge 0$ and Markov's inequality yields
		\begin{equation*}
			\forall i, j, k \ge 0,\,\PP\left(p_{2i} + \cdots + {p}_{2i + 2 k + 1} \ge j + k \frac{\beta m+\log(2C)}{\beta}\right) \le 2^{-k} e^{-\beta j}.
		\end{equation*}
		As a consequence, for $C' := m + \frac{\log(2C)}{\beta}$,
		\begin{equation*}
			\PP\left(n - \overline{p}_{2q_n - 2 k} \ge k C'\right) \le \frac{2^{-k}}{1-e^{-\beta}}
		\end{equation*}
		Finally, by independence of $r_n^h$ with $\left(p_i\right)_{i \ge 0}$,
		\begin{align*}
			\left(n - \overline{p}_{2q_n - 2 r_n^h} \ge t\right) & = \sum_{k \ge 0} \PP\left(n - \overline{p}_{2q_n - 2 k} \ge t\right) \PP(r_n^h = k) \\
			& \le \PP(r_n^h \ge  \lfloor t/C'\rfloor) + \PP\left(n - \overline{p}_{2q_n - 2\lfloor t/C'\rfloor} \ge t\right) \\
			& \le 4^{-\lfloor t/C'\rfloor} + 2^{-\lfloor t/C'\rfloor}. 
		\end{align*}
		
		Let us now construct $c_n$.
		We use the same technique as for the construction of $l^f$ and $r_n^h$.
		For all $n, k \ge 0$, we define the event
		\begin{multline*}
			S'_{n,k} =  (2 k + 1 \le q_n) \cup (2 k + 1 < q_n) \cap \left(\gamma^p_{2 q_n - 2k - 1} \AA^{\eps} \gamma_{\overline{p}_{2q_n - 2k}} \cdots \gamma_{n-1} \gamma^p_0 \cdots \gamma^p_{2 k}\right) \\ \cap \left(\gamma_{\overline{p}_{2q_n - 2k - 1}} \cdots \gamma_{n-1} \gamma^p_0 \cdots \gamma^p_{2 k} \AA^\eps \gamma^p_{2 k + 1}\right).
		\end{multline*}
		By \eqref{pivot:schottky}, it holds for all $n, k$ that
		\begin{equation*}
			\PP\left(S'_{n,k}\,\middle|\,\left(\widetilde{\gamma}^p_{k'}\right)_{k' \notin \{2k +1, 2q_n - 2k -1\}}, 2k+1 < q_n\right) \ge \frac{1}{2}
		\end{equation*}
		For all $n, k$ we define the event
		\begin{equation*}
			S_{n,k} := \left(\tau_k < \frac{\mathds{1}(S'_{n,k})}{2\PP\left(S'_{n,k}\,\middle|\, q_n, \left(\widetilde{\gamma}^p_{k'}\right)_{k' \notin \{2k +1, 2q_n - 2k -1\}} \right)}\right).
		\end{equation*}
		Let $c_n:= \min\{k \ge 0\mid S_{n,k}\}$ for all $n$.
		That way, $\PP(S_{n,k}) = \frac{1}{2}$ for all $n,k$. 
		Moreover, the $(S_{n,k})_{k \ge 0}$ are independent and their joint data is independent of $\left(\widetilde{\gamma}^p_{2j}\right)_{j\ge 0}$.
		Therefore $c_n \sim \mathcal{G}_{1/2}$ and $c_n$ is independent of $\left(\widetilde{\gamma}^p_{2j}\right)_{j\ge 0}$.
		With a similar argument as for as for $r_n^h$,
		\begin{equation*}
			\forall t \ge 0,\;\left(n - \overline{p}_{2q_n - 2 c_n - 2} \ge t\right) \le 2^{1 - \lfloor t/C'\rfloor} + 2^{-\lfloor t/C'\rfloor}.
		\end{equation*}
		so $n - \overline{p}_{2q_n - 2 c_n - 2}$ has a bounded exponential moment that does not depend on $n$ and it follows from Lemmas \ref{lem:power-compare} and \ref{lem:sumexp}, $\overline{p}_{2c_n}$ also has a bounded exponential moment that does not depend on $n$. 
		The same holds for their sum, which concludes the proof.
	\end{proof}

	\subsection{Quantitative estimates for the coefficients and spectral radius}\label{section:slln}

	In the present section, we use Lemma \ref{lem:l-r-piv} and mainly point \eqref{pivot:indep-tail} from Theorem \ref{th:pivot} to prove Theorem \ref{th:slln}, from which Theorem \ref{th:intro-lln} follows.
	We follow the plan of the proof sketched in the discussion following the statement of Theorem \ref{th:pivot}.
	
	Let us remind that we write $N: g \mapsto \log\|g\|\|g^{-1}\|$.
	The following result is a more specific version of \eqref{dom-coef} in Theorem \ref{th:slln}.
	The arguments developed in the present section rely on point \eqref{pivot:indep-tail} from Theorem \ref{th:pivot}. 
	We remind that given a measure $\nu$ supported on $\GL(E)$ that is $(\alpha, \eps, m)$-contracting for some $\alpha, \eps> 0$ and $m\in \NN$, we may assume by definition that there exists a measure $\tilde\mu < \nu^{\otimes m} / \alpha$ that is $1/6$-Schottky for $\AA^\eps$.
	By taking a normalized restriction, we may moreover assume that $\tilde\mu$ satisfies $\tilde\mu(N^{-1}[0, B]^m) = 1$ for a large enough constant $B$.

	Note that for a fixed $B$, and fixed values of $(\alpha, \eps, m)$, the space of all such probability measures is weak-$*$ open by Proposition \ref{prop:weak-bounded}.
	This hints towards Wasserstein $\ELL^p$ continuity of some quantities but that is not the purpose of the present paper. 
	 
	\begin{Th}\label{th:coef-control}
		Let $\nu$ be a strongly irreducible and proximal probability distribution on $\GL(E)$.
		Let $\alpha, \eps > 0$ and $m \in \NN$ be as in Lemma \ref{lem:contraction}.
		Let $\tilde\mu < \nu^{\otimes m} / \alpha$ be $1/6$-Schottky for $\AA^\eps$ and let 
		Let $S \subset \GL(E)$ be such that $\tilde\mu(S^m) = 1$ and assume that $B := \max N(S)$ is finite.
		Let $(\gamma_n)_{n \ge 0} \sim \nu^{\otimes\NN}$.
		Let $f \in E^* \setminus\{0\}$ and let $v \in E \setminus \{0\}$.
		There exist constants $C, \beta > 0$ that only depend on $\alpha$ and $m$
		and a constant $D$ that only depends on $\eps$ and $B$, such that for all $f \in E^*\setminus\{0\}$, all $v \in E \setminus\{0\}$ and all $n \ge 0$,
		\begin{equation}\label{coef-control}
			\PP\left(-\log\left(\frac{|f \overline \gamma_n v|}{\|f\|\|\overline{\gamma}_n\|\|v\|}\right) \ge t\right) \le C e^{-\beta n} + \sum_{k \ge 1} C e^{-\beta k} k N_*\nu(t/k - D, +\infty).
		\end{equation}
	\end{Th}
	
	\begin{proof}
		Let $(p_n)_{n \ge 0}$ and be as in Theorem \ref{th:pivot}.
		Assume that the measure $\tilde\mu$ has compact support, as constructed in Corollary \ref{cor:schottky}. 
		Let $l_f$ be as in Lemma \ref{lem:l-r-piv} and for all $n$, let $q_n$ and $r_n^v$ be as in Lemma \ref{lem:l-r-piv}.
		Then by Lemma \ref{lem:l-r-piv}, there exist constants $C, \beta > 0$ such that for all $k$ and all $n$,
		\begin{equation}\label{moment-I}
			\PP\left(n - \overline{p}_{2q_n - 2r_n^v} + \overline{p}_{2l^f + 1} = k\right) \le C e^{-\beta k}.
		\end{equation}
		
		The first issue is that the joint data of $r_n^v$ and $l^f$ is not independent of $\left(\widetilde{\gamma}^p_{2k}\right)_{k \ge 0}$ in general.
		However, there exist constants $C, \beta > 0$ such that for all $n$,
		\begin{equation}
			\PP\left(l^f \ge \lfloor q_n / 2 \rfloor \text{ or } r_n^v - 1 \ge \lceil q_n / 2 \rceil \right) \le C e^{-\beta n}
		\end{equation}
		Moreover, on the intersection of the level sets of $\left(\widetilde{\gamma}^p_{2k}\right)_{k \ge 0}$ with the event $S_n^{f,v}:= (l^f < \lfloor q_n / 2 \rfloor) \cap \left(r_n^v + 1 < \lceil q_n / 2 \rceil \right)$, the random integer $l^f$ satisfies a condition that only depends on the data of $\left(\gamma^p_{2k + 1}\right)_{k < \lfloor q_n / 2 \rfloor}$ and $r_n^v$ satisfies a condition that only depends on the data of $\left(\gamma^p_{2k + 1}\right)_{k > \lfloor q_n / 2 \rfloor}$. 
		Moreover, the data of $(\gamma^p_{2k + 1})_{k < \lfloor q_n / 2 \rfloor}$ and $(\gamma^p_{2k + 1})_{k > \lfloor q_n / 2 \rfloor}$ are almost surely independent with respect to the $\sigma$-algebra $\langle \widetilde{\gamma}^p_{2k}\rangle_{k \ge 0}^\sigma$, generated by $(\widetilde{\gamma}^p_{2k})_{k \ge 0}$, therefore, we may assume that $l^f$ and $r_n^v$ are almost surely independent with respect to $\langle \widetilde{\gamma}^p_{2k}\rangle_{k \ge 0}^\sigma$.
		Moreover, they are both independent of said $\sigma$-algebra by lemma \ref{lem:l-r-piv} so their joint data is independent of $\langle \widetilde{\gamma}^p_{2k}\rangle_{k \ge 0}^\sigma$ on the set $S_n^{v,f}$.
		
		Let us now give an upper bound on $-\log\left(\frac{|f \overline \gamma_n v|}{\|f\|\|\overline{\gamma}_n\|\|v\|}\right)$. 
		By Definition of $l^f$ and $r_n^v$, we have $f \gamma^p_0 \cdots \gamma^p_{2l^f} \AA^\eps \gamma_{2l^f + 1}$ and $\gamma^p_{2q_n - 2r_n -1} \AA^\eps \gamma_{\overline{p}_{2q_n - 2 r_n}} \cdots \gamma_{n-1} v$.
		Therefore, on the set $S_n^{f,v}$, by point \eqref{pivot:herali} in Theorem \ref{th:pivot},
		\begin{equation*}
			f \gamma_0 \cdots \gamma_{\overline{p}_{2l^f + 1} - 1} \AA^\frac{\eps}{2} \gamma_{\overline{p}_{2l^f}} \cdots \gamma_{n-1} v 
			\quad \text{and} \quad 
			\gamma_{\overline{p}_{2l^f}} \cdots \gamma_{\overline{p}_{2q_n - 2 r_n} - 1} \AA^\frac{\eps}{2} \gamma_{\overline{p}_{2q_n - 2 r_n}} \cdots \gamma_{n-1} v
		\end{equation*}
		Therefore,
		\begin{equation*}
			|f \overline{\gamma}_n v| \ge \left\|f \gamma_0 \cdots \gamma_{\overline{p}_{2l^f + 1} - 1}\right\| \frac{\eps}{2} \left\| \gamma_{\overline{p}_{2l^f}} \cdots \gamma_{\overline{p}_{2q_n - 2 r_n} - 1} \right\| \frac{\eps}{2} \left\|\gamma_{\overline{p}_{2q_n - 2 r_n}} \cdots \gamma_{n-1} v\right\|
		\end{equation*}
		Thus,
		\begin{equation*}
			\frac{|f \overline{\gamma}_n v|}{\|f\| \|\overline{\gamma}_n\| \|v\|} \ge \frac{\eps^2}{4} \frac{\left\|f \gamma_0 \cdots \gamma_{\overline{p}_{2l^f + 1} - 1}\right\|}{\left\|f\right\| \left\|\gamma_0 \cdots \gamma_{\overline{p}_{2l^f + 1} - 1}\right\|} \frac{\left\|\gamma_{\overline{p}_{2q_n - 2 r_n}} \cdots \gamma_{n-1} v\right\|}{\left\|\gamma_{\overline{p}_{2q_n - 2 r_n}} \cdots \gamma_{n-1}\right\| \left\|v\right\|}
		\end{equation*}
		Using the inequality $\left\|f g\right\| \left\|g^{-1}\right\| \ge \left\|f\right\|$, valid for all $f,g$, we obtain
		\begin{equation}\label{xwcn}
			-\log\left(\frac{|f \overline{\gamma}_n v|}{\|f\| \|\overline{\gamma}_n\| \|v\|}\right) \le |\log(\eps^2/4)| + \sum_{k = 0}^{\overline{p}_{2l^f + 1} - 1} N(\gamma_k) + \sum_{k = \overline{p}_{2q_n - 2 r_n^v}}^{n - 1} N(\gamma_k).
		\end{equation}
		Let us write $I_{n}^{f,v}$ for the set $\{0, \dots, \overline{p}_{2l^f + 1} - 1\} \cup \{\overline{p}_{2q_n - 2 r_n^v}, \dots, n-1\} \cap \{0,\dots, n\}$.
		Note that $\# I_{n}^{f,v} = n - \overline{p}_{2q_n - 2r_n^v} + \overline{p}_{2l^f + 1} \ge 1$.
		
		Let us now give a probabilistic estimate for $\sum_{k \in I_{n}^{f,v}} N(\gamma_k) + |\log(\eps^2/4)|$.
		Note that $N^{-1}(t, + \infty) \cap S  = \emptyset$ for all $t \ge B$. 
		Hence, by point \eqref{pivot:indep-tail} in Theorem \ref{th:pivot} and making the non-restrictive assumption $\alpha \le 1/2$, we obtain
		\begin{equation}\label{indep-tail}
			\forall t \ge B, \, \forall k,\, \PP\left(N(\gamma_k) > t \,\middle|\, \left(p_j\right)_{j \ge 0}\right) \le \PP(N(g_k) > t) \le 2N_*\nu(t, + \infty).
		\end{equation}
		Note that for all $k$,
		\begin{itemize}
			\item on the measurable set $\mathrm{odd}_p(k) = \left(\exists j, \overline{p}_{2j + 1} \le k < \overline{p}_{2j + 2}\right)$, we have $\gamma_k \in S$ by point \eqref{pivot:schottky} in Theorem \ref{th:pivot} so $\PP\left(N(\gamma_k) > B \,\middle|\, \left(p_j\right)_{j \ge 0}\right) = 0$
			\item on the set $\mathrm{even}_p(k) = \left(\exists j, \overline{p}_{2j} \le k < \overline{p}_{2j + 1}\right)$, the data of $l^f$ and $r_n^v$ is independent of $\gamma_k$ so \eqref{indep-tail} yields
			\begin{equation}\label{jsdhg}
				\forall t \ge B, \, \forall k \le n,\, \PP\left(N(\gamma_k) > t \,\middle|\, \left(p_j\right)_{j \ge 0}, l^f, r_n^v\right) \le 2N_*\nu(t, + \infty).
			\end{equation}
		\end{itemize}
		Note moreover that $N_*\nu(t, + \infty) \le N_*\nu(t- B, +\infty)$ for all $t \ge B$ and $N_*\nu(t - B, +\infty) = 1$ when $t < B$.
		Therefore, for all $t \ge 0$, \eqref{xwcn} yields
		\begin{align*}
			\PP\left(- \log\left(\frac{|f \overline{\gamma}_n v|}{\|f\| \|\overline{\gamma}_n\| \|v\|}\right) > t \cap S_{n}^{f,v}\right) 
			\le \PP\left(\sum_{k \in I_{n}^{f,v}} N(\gamma_k) > t - |\log(\eps^2/4)| \right) &\\
			\le \PP\left(\exists k \in I_n^{f,v}, \, N(\gamma_k) > \frac{t - |\log(\eps^2/4)|}{\# I_n^{f,v}}\right)&\\
			\le \sum_{k  = 1}^{\infty} \PP\left(\# I_n^{f,v} = k \right) 2k N_*\nu\left( \frac{t - |\log(\eps^2/4)|}{k} - B, + \infty\right)&.
		\end{align*}
		By \eqref{moment-I}, there exist constants $C, \beta$ such that $2\PP\left(\# I_n^{f,v} = k\right) \le C e^{-\beta k}$. Without loss of generality, we may also assume that $\PP\left(S_n^{f,v}\right) \ge 1- Ce^{-\beta n}$.
		To conclude, set $D := B + |\log(\eps^2/4)|$ and we obtain \eqref{coef-control}.
	\end{proof}

	We get a stronger result for the spectral radius because it is pinched by the singular values, \ie $s_d \le \rho_1 \le s_1$. 
	This allows us to have a uniform control over the spectral radius at all times, even when $2 c_n + 1 \ge q_n$.
	As a consequence, we get rid of the term $C e^{-\beta n}$ which does no decay to $0$ as $t \to + \infty$.
	That way, the distance between $\lambda_{1}(\overline\gamma_n)$ and $\log\|\overline\gamma_n\|$ is bounded in $\ELL^1$
	It implies that $(\lambda_{1}(\overline\gamma_n)/n)$ converges to $\lambda_1(\nu)$ almost surely and in $\ELL^1$, which is worth mentioning.

	\begin{Th}\label{th:rad-control}
		Let $\nu$ be a strongly irreducible and proximal probability distribution on $\GL(E)$.
		Let $(\gamma_n)_{n \ge 0}\sim \nu^{\otimes\NN}$.
		There exist constants $C, D, \beta > 0$ that depend only on $\alpha$ and $m$ and such that:
		\begin{equation}\label{rad-control}
			\PP\left(-\log\left(\frac{\rho_1(\gamma_n)}{\|\overline{\gamma}_n\|}\right) > t\right) \le \sum_{k = 1}^\infty C e^{-\beta k} k N_*\nu(t/k - D, +\infty).
		\end{equation}
	\end{Th}
	
	\begin{proof}
		Let $(p_n)_{n \ge 0}$ be as in in Theorem \ref{th:pivot}.
		Let $S \subset \GL(E)$ be compact and as in point \eqref{pivot:indep-tail} in Theorem \ref{th:pivot}.
		For all $n$, let $q_n$ and $c_n$ be as in Lemma \ref{lem:l-r-piv}.
		Then by Lemma \ref{lem:l-r-piv}, there exist constants $C, \beta > 0$ such that for all $k$ and all $n$,
		\begin{equation}
			\PP\left(n - \overline{p}_{2q_n - 2c_n} + \overline{p}_{2c_n} = k\right) \le C e^{-\beta k}.
		\end{equation}
		Let us now give an upper bound on $-\log\left(\frac{\rho_1(\gamma_n)}{\|\overline{\gamma}_n\|}\right)$.
		Set $S_n := \left(2 c_n + 1 < q_n\right)$.
		By definition of $c_n$ and by point \eqref{pivot:herali} in Theorem \ref{th:pivot}, it holds that $g_n \AA^\frac{\eps}{4} g_n$ on $S_n$, with
		\begin{equation*}
			g_n := \gamma_{\overline{p}_{2c_n + 1}} \cdots \gamma_{n-1} \gamma_0 \cdots \gamma_{\overline{p}_{2c_n + 1} - 1}.
		\end{equation*}
		Let us give more details.
		Write for all $n$:
		\begin{gather*}
			f_n := \gamma_{\overline{p}_{2q_n - 2c_n - 1}} \cdots \gamma_{n-1} \gamma^p_0 \cdots \gamma^p_{2 c_n},
			\\
			h_n := \gamma_{\overline{p}_{2q_n - 2c_n}} \cdots \gamma_{n-1} \gamma_0^p \cdots \gamma^p_{2 c_n},
			\\
			z_n := \gamma^p_{2 c_n + 1} \cdots \gamma^p_{2q_n - 2c_n - 2}.
		\end{gather*}	
		That way, $g_n = z_n f_n$.
		We have $\gamma^p_{2q_n - 2c_n - 1} \AA^\eps h_n$ and $\sigma\left(\gamma^p_{2q_n - 2c_n - 1}\right) \le \eps^6/48$ by point \eqref{pivot:schottky} in Theorem \ref{th:pivot}. 
		Hence $\sigma(f_n) < \eps^4/48$ by Lemma \ref{lem:c-prod}. 
		Moreover, $f_n \AA^\eps \gamma^p_{2 c_n + 1}$.
		Applying point \eqref{pivot:herali} from Theorem \ref{th:pivot} to both aforementioned conditions yields $f_n \AA^{\frac{\eps}{2}} z_n \AA^{\frac{\eps}{2}} f_n$.
		Moreover $\gamma^p_{2 c_n + 1} \AA^{\frac{\eps}{2}} \gamma^p_{2 c_n + 2} \cdots \gamma^p_{2q_n - 2c_n - 2}$ so $\sigma(z_n) \le \eps^4/12$ by Lemma \ref{lem:c-prod}.
		Then $z_n f_n \AA^{\frac{\eps}{4}} z_n f_n$ by Lemma \ref{lem:herali} applied to $z_n$ and $f_n$.
		
		Then by \eqref{ev-sv} in Lemma \ref{lem:eigen-align},
		\begin{equation*}
			\rho_1(g_n) \ge \frac{\eps}{4} \|g_n\| \ge \frac{\eps^2}{8} \frac{\|z_n\|}{\|f_n\|}
		\end{equation*}
		Moreover, $g_n$ is conjugated to $\overline{\gamma}_n$ so $\rho_1(g_n) = \rho_1(\overline{\gamma}_n)$. 
		Therefore, on $S_n$,
		\begin{align*}
			-\log\left(\frac{\rho_1(\overline{\gamma}_n)}{\|\overline{\gamma}_n\|}\right) & \le |\log(\eps^2/8)| + N(f_n) 
			\\
			& \le |\log(\eps^2/8)| + \sum_{k = 0}^{\overline{p}_{2c_n  +1} - 1} N(\gamma_k) + \sum_{k = \overline{p}_{2q_n - 2c_n - 1}}^{n-1} N(\gamma_k).
		\end{align*}
		It holds outside of $S_n$ that $-\log\left(\frac{\rho_1(\overline{\gamma}_n)}{\|\overline{\gamma}_n\|}\right) \le N(\overline{\gamma}_n) \le \sum_{k = 0}^{n-1} N(\gamma_k)$. 
		Set $J_n:= \{0, \dots , n-1\}$ on $S_n^c$ and $J_n:= \{0, \dots , \overline{p}_{2c_n  +1} - 1\} \cup \{\overline{p}_{2q_n - 2c_n - 1}, \dots, n-1\}$ on $S_n$. 
		Then, for all $t \ge 0$,
		\begin{align*}
			\PP\left(-\log\left(\frac{\rho_1(\gamma_n)}{\|\overline{\gamma}_n\|}\right) > t\right) 
			& \le \PP\left(\sum_{j \in J_n} N(\gamma_j) > t - |\log(\eps^2/8)|\right)
			\\
			& \le \sum_{k  = 0}^n \PP\left(\max_{j \in J_n} N(\gamma_j) > \frac{t - |\log(\eps^2/8)|}{k} \cap \# J_n = k\right)
		\end{align*}
		By Lemma \ref{lem:l-r-piv}, there exist constants $C, \beta > 0$ such that $2 \PP(\# J_n = k) \le C e^{-\beta k}$ for all $k \le n$.
		Moreover, with the same reasoning as in the proof of \eqref{jsdhg} in Theorem \ref{th:coef-control}, there exists a constant $B$ such that:
		\begin{equation*}
			\forall k \le n,\; \PP\left(N(\gamma_k) > t \mid  J_n\right) \le 2N(t-B, + \infty).
		\end{equation*}
		From the above, we deduce \eqref{rad-control} for $D = B + |\log(\eps^2/8)|$, with the same reasoning as in the proof of Theorem \ref{th:coef-control}.
	\end{proof}

	\begin{proof}[Proof of Theorem \ref{th:slln}]
		To derive a proof of Theorem \ref{th:slln} from Theorems \ref{th:coef-control} and \ref{th:rad-control}, we simply need to prove that for all distribution $\eta \neq \delta_{0}$ on $\RR_{\ge 0}$ and for all $C, D, \beta > 0$, there exist constants $C', \beta '> 0$ such that
		\begin{equation}\label{remove-D}
			\forall t \ge 0,\; \sum_{k \ge 1} C e^{-\beta k} k \eta(t/k - D, +\infty) \le \sum_{k \ge 0} C' e^{-\beta' k} \eta(t/k, +\infty).
		\end{equation}
		Note that, if $N_* \nu = \delta_0$, then $\nu$ is supported on $\KK^\times\mathrm{O}(E)$ and therefore $\rho_1 = \rho_2 = \dots =\rho_d$ on $\Gamma_\nu$. 
		This would contradict the proximality so $N_*\nu \neq \delta_0$.
		
		Set $\eta := N_* \nu$ and let $C, \beta > 0$.
		First let $C_1, \beta_1 > 0$ be such that $ C e^{-\beta k} k \le C_1 e^{-\beta_1 k}$ for all $k$.
		Moreover,  $t / k - D \ge t / (2k)$ for all $t \ge 2D k$.
		Hence, for all $t > 0$:
		\begin{align*}
			\sum_{k \ge 0} C e^{-\beta k} k \eta(t/k - D, +\infty) 
			& \le \sum_{1 \le k \le t /(2d)} C_1 e^{-\beta_1 k} \eta(t/(2k), +\infty) + \sum_{k > t/(2d)}  C_1 e^{-\beta_1 k}
			\\
			&\le \sum_{k  = 1}^{+ \infty} 2C_1 e^{-\beta_1 k} \eta(t/k, +\infty) + \frac{C_1}{1 - e^{- \beta_1}} e^{-\beta_1 \lceil t / (2D) \rceil}
		\end{align*}
		Let $\delta > 0$ be such that $\eta(\delta, + \infty) > 0$ and let $\alpha = \eta(\delta, + \infty)$.
		Then, for all $C_2, \beta _2 > 0$,
		\begin{align*}
			\sum_{k = 1}^{\infty}  C_2 e^{-\beta_2 k} \eta(t/k, +\infty) \ge C_2 e^{- \beta_2 \lfloor t / \delta \rfloor} \alpha
		\end{align*}
		Then for all $t \ge 0$, we have $C_2 e^{- \beta_2 \lfloor t / \delta \rfloor} \alpha \ge \frac{C_1}{1 - e^{- \beta_1}} e^{-\beta_1 \lceil t / (2D) \rceil}$ with $\beta_2 := \delta \beta _1 / (2D)$ and $C_2 := \frac{C_1 e^{\beta_1 + \beta_2}}{\alpha(1 - e^{- \beta_1})}$.
		Without loss of generality, we may assume that $2D \ge 1$ so $\beta_2 \le \beta_1$ and $C_2 \ge C_1$.
		Therefore, for all $t \ge 0$,
		\begin{equation*}
			\sum_{k \ge 0} C e^{-\beta k} k \eta(t/k - D, +\infty) \le \sum_{k = 1}^{\infty}  3C_2 e^{-\beta_2 k} \eta(t/k, +\infty). 
		\end{equation*}
		Taking $C' := 3C_2$ and $\beta' = \beta_2$ yields \eqref{remove-D}. 
	\end{proof}


	\subsection{Law of large numbers and large deviations inequalities for the pivot extraction}\label{section-extract}
	
	In this paragraph, we use point \eqref{pivot:renewal} from Theorem \ref{th:pivot} to prove intermediate Lemmas for the proof of Theorem \ref{th:escspeed}.
	
	\begin{Lem}[Escape speed and large deviations inequalities for the aligned extraction] \label{lem:pre-speed}
		Let $\alpha, \eps > 0$, $m \in\NN$ and let $\nu$ be an $(\alpha, \eps, m)$-contracting probability measure on $\End(E)$ that has positive eventual rank.
		Let $(\gamma_n)_{n \ge 0} \sim \nu^{\otimes\NN}$ and let $(p_n)_{n \ge 0}$ be as in Theorem \ref{th:pivot}.
		Then there exists a constant $\lambda' \in  (0, + \infty]$ such that for all $i_0 \ge 1$, and almost surely:
		\begin{equation}
			\lim_{j \to +\infty} \frac{-\log\sigma\left(\gamma^p_{2i_0} \cdots \gamma^p_{2j}\right)}{j - i_0} = \lambda'.
		\end{equation}
		Moreover, for all $\alpha < \lambda'$, there exist constant $C, \beta > 0$ such that for all $i \le j$,
		\begin{equation}
			\PP\left(-\log\sigma\left(\gamma^p_{2i} \cdots \gamma^p_{2j}\right) \le \alpha(j-i)\right) \le C e^{-\beta (j-i)}. \label{ldev-extract}
		\end{equation}
	\end{Lem}
	
	\begin{proof}.
		First we claim that for all $i \le j < k$, the following sub-additivity inequality holds
		\begin{equation}\label{sub-add}
			\log\sigma\left(\gamma^p_{2i} \cdots \gamma^p_{2k}\right) \le \log\sigma\left(\gamma^p_{2i} \cdots \gamma^p_{2j}\right) + \log\sigma\left(\gamma^p_{2j + 2} \cdots \gamma^p_{2k}\right) + \log(\eps^2 / 12).
		\end{equation}
		This is a consequence of Lemma \ref{lem:triple-ali} with $f = \gamma^p_{2i} \cdots \gamma^p_{2j}$, $g = \gamma^p_{2j+1}$ and $h = \gamma^p_{2j + 2} \cdots \gamma^p_{2k}$.
		By point \eqref{pivot:schottky} in Theorem \ref{th:pivot}, $\sigma(g) \le \eps^6/48$ and by point \eqref{pivot:herali}, $f \AA^{\frac{\eps}{2}} g \AA^{\frac{\eps}{2}} h$ so $\sigma(fgh) \le \sigma(f)\frac{\eps^6}{48}\sigma(h)\frac{4}{\eps^4}$, which proves \eqref{sub-add}.
		
		A consequence of point \eqref{pivot:renewal} in Theorem \ref{th:pivot} is that the law of $(\gamma^p_{k + 2})_{k \ge 0}$ is ergodic for the shift map $T: \Gamma^\NN\to\Gamma^\NN\,;\, (g_{k})_{k \ge 0} \mapsto (g_{k + 2})_{k \ge 0}$.
		
		For all $n \ge 0$, let $\phi_{n}:= -\log\sigma(\gamma^p_2 \cdots \gamma^p_{2n + 2})$.
		Then for all $m,n \ge 0$ by \eqref{sub-add} with $i = 1$, $j = m + 1$ and $k = m + n + 1$, we have $\phi_{m + n} \ge \phi_m + \phi_{n} \circ T^m$. 
		Note also that $\phi_n \ge 0$ for all $n$
		Therefore, by Kingman's sub-additive ergodic Theorem \cite{K68}, there exists a constant $\lambda' = \sup \frac{\EE(\phi_n)}{n} = \lim \frac{\EE(\phi_n)}{n} \in [0 , +\infty]$ such that $\frac{\phi_n}{n} \to \lambda'$ almost surely. 
		Not also that $\lambda'> 0$ because $\phi_1 \ge -\log(\eps^2 / 12)$ by \eqref{sub-add}.
		
		Let us now prove \eqref{ldev-extract}. 
		Let $\alpha < \lambda'$ and let $n \in\NN$ be such that $\alpha < \lambda_n = \frac{\EE(-\log\sigma(\gamma^p_2 \cdots \gamma^p_{2 n}))}{n}$. 
		For all $k, r \ge 0$, and by \eqref{sub-add},
		\begin{equation}\label{hgfjfgk}
			-\log\sigma\left(\gamma^p_2 \cdots \gamma^p_{2k n + r}\right) \le \sum_{j = 0} ^{k - 1} -\log\sigma(\gamma^p_{2 j n + 2} \cdots \gamma^p_{2(j+1) n}).
		\end{equation}
		Moreover, for all $i$, the data of $\left(\gamma_{k}\right)_{k \le 2i}$ is independent of the data of $\left(\gamma_{k}\right)_{k \ge 2i + 2}$ so the sequence $\left(-\log\sigma(\gamma^p_{2 j n + 2} \cdots \gamma^p_{2(j+1) n})\right)_{j \ge 0}$ is i.i.d. hence, by Lemma \ref{lem:ldev-classique}, for all $\alpha' < \EE(-\log\sigma(\gamma^p_2 \cdots \gamma^p_{2 n}))$ there exist constants $C, \beta >0$ such that:
		\begin{equation*}
			\forall k\in\NN,\; \PP\left(\sum_{j = 0} ^{k - 1} -\log\sigma(\gamma^p_{2 j n + 2} \cdots \gamma^p_{2(j+1) n}) \le \alpha' k\right) \le C e^{-\beta k}.
		\end{equation*}
		Then by \eqref{hgfjfgk},
		\begin{equation*}
			\forall m \in \NN,\; \PP\left(-\log\sigma(\gamma^p_{2} \cdots \gamma^p_{2m}) \le \alpha' \lfloor m/ n\rfloor \right) \le C e^{-\beta \lfloor m/n \rfloor} \le C e^{\beta} e^{-\frac{\beta}{n}m}.
		\end{equation*}
		Moreover, we may assume that $n \alpha < \alpha'$ hence, there exists $m_0$ such that $\alpha m \le \alpha' \lfloor m/ n\rfloor$ for all $m \ge m_0$, which concludes the proof.
	\end{proof}
	
	Lemma \ref{lem:pre-speed} tells us that a subsequence of $(-\log\sigma(\overline\gamma_n))_n$ escapes linearly to $+\infty$ with large deviations inequalities below the escape speed.
	In dimension $2$, we can conclude using the sub-additivity of $N = -\log\sigma$.
	In higher dimension, we use the lower bound $-\log\sigma(\overline\gamma^p_{2q_n - 2 r_n^{h} - 1}) \le -\log\sigma(\overline\gamma_n)$ for a good choice of $h$ to get the large deviations inequalities below the escape speed.
	
	\begin{Lem}\label{lem:ldev-bloc-central}
		Let $\alpha, \eps > 0$, $m \in\NN$ and let $\nu$ be an $(\alpha, \eps, m)$-contracting probability measure on $\End(E)$ that has positive eventual rank.
		Let $(\gamma_n)_{n \ge 0} \sim \nu^{\otimes\NN}$ and $(p_n)_{n \ge 0}$ be as in Theorem \ref{th:pivot}.
		Let $\lambda'$ be as in Lemma \ref{lem:pre-speed}.
		Let $f \in E^* \cup \End(E)$ and $h \in E \cup \End(E)$.
		Let $l^f$, $(q_n)_{n \ge 0}$ and $(r_n^h)_{n \ge 0}$ be as in Lemma \ref{lem:l-r-piv}.
		Set $\lambda:= \frac{\lambda'}{m + \EE(p_2)}$.
		Then the random sequences $(-\log\sigma(\gamma^p_{2 l^f + 2} \cdots \gamma^p_{2q_n - 2 r_n^h - 2}))_{n \ge 0}$ and $(-\log\sigma(\gamma^p_{2 c_n + 2} \cdots \gamma^p_{2q_n - 2 c_n - 2}))_{n \ge 0}$ both satisfy large deviations inequalities below the speed $\lambda$.
	\end{Lem}
	
	\begin{proof}
		First let us show that $(q_n)_{n \ge 0}$ satisfies large deviations inequalities below the speed $\frac{1}{m + \EE(p_2)}$.
		Let $\alpha > m + \EE(p_2)$. 
		The sequence $(\overline{p}_{2k})_{k \ge 0}$ is non-decreasing so for all $n$,
		\begin{equation*}
			\PP(q_n \le n / \alpha) = \PP\left(\overline{p}_{2\lceil n/\alpha \rceil} > n\right) \le \PP\left(\overline{p}_{2\lceil n/\alpha \rceil} > \alpha \lfloor n/\alpha \rfloor\right)
		\end{equation*}
		Moreover, the sequence $(p_{2k+2} + m)_{k \ge 0}$ is i.i.d and $p_2$ has a finite exponential moment by Lemma \ref{lem:ldev-classique}, the sequence $\left(\overline{p}_{2k + 2} - p_0 - m\right)$ satisfies large deviations inequalities below the speed $-m - \EE(p_2)$ and by Lemma \ref{lem:ldev-shift}, the sequence $\left(\overline{p}_{2k + 2}\right)$ also does.
		Therefore, there exits constants $C, \beta >0$ such that $\PP\left(\overline{p}_{2 k + 2} \ge \alpha k\right) \le C e^{-\beta k}$ for all $k$.
		Let $C, \beta$ be such constants and $n \in\NN$.
		For $k := \lfloor n/\alpha \rfloor$, it holds that
		\begin{equation*}
			\PP(q_n \le n / \alpha) \le  C e^{-\beta \lfloor n/\alpha \rfloor}.
		\end{equation*}
		Therefore $(q_n)_{n \ge 0}$ satisfies large deviations inequalities below the speed $\lambda'':= \frac{1}{m + \EE(p_2)}$. 
		By Lemma \ref{lem:ldev-shift}, the sequences $\left(q_n - r_n^h - 1\right)_{n \ge 0}$ and $\left(q_n - c_n - 1\right)_{n \ge 0}$ also do.
		
		Now let $\alpha < \lambda = \lambda'\lambda''$.
		Let $0 < \delta < \lambda''/2$ and $0 < \alpha_1 \le \lambda'$ be such that $\alpha < (\lambda'' - 2\delta)\alpha_1$.  
		The large deviations inequalities imply that there exist constants $C_0, \beta_0$ such that the following assertions hold:
		\begin{gather}
			\forall n,\; \PP\left(q_n - r_n^h - 1 \le (\lambda'' - \delta) n\right) \le C_0 e^{-\beta_0 n}, \label{qn-rn}\\
			\forall n,\; \PP\left(q_n - c_n - 1 \le (\lambda''-\delta) n\right) \le C_0 e^{-\beta_0 n}.\label{qn-cn}
		\end{gather}
		Since $c_n$ and $l^f$ both follow a geometric law, we may also assume the following:
		\begin{gather}
			\forall n,\; \PP\left(l^f + 1 \ge \delta n\right) \le C_0 e^{-\beta_0 n}, \label{ln}\\
			\forall n,\; \PP\left(c_n + 1 \ge \delta n\right) \le C_0 e^{-\beta_0 n}. \label{cn}
		\end{gather}
		By Lemma \ref{lem:pre-speed}, there exists constants $C_1, \beta_1 > 0$ such that:
		\begin{equation}
			\forall 1 \le i \le j,\; \PP\left(-\log\sigma\left(\gamma^p_{2i} \cdots \gamma^p_{2j}\right) \le \alpha_1 (i-j)\right) \le C_1e^{-\beta_1(j-i)}.\label{i-j}
		\end{equation}
		Let $C_0, \beta_0 > 0$ and $C_1, \beta_1$ be such that \eqref{qn-rn}, \eqref{qn-cn}, \eqref{ln}, \eqref{cn} and \eqref{i-j} hold.
		We have:
		\begin{multline}\label{decomp-l-r}
			\left(-\log\sigma(\gamma^p_{2 l^f + 2} \cdots \gamma^p_{2q_n - 2 r_n^h - 2}) \le \alpha n\right) \subset \left(q_n - r_n^h - 1 \le (\lambda''-\delta) n\right) \\ \cup \left(l^f + 1 \ge \delta n\right) \cup  \left(\exists i \le \delta n,\, \exists j \ge \lambda'' - \delta n ,\; -\log\sigma\left(\gamma^p_{2i} \cdots \gamma^p_{2j}\right) \le \alpha_1 (i-j)\right)
		\end{multline}
		Then by \eqref{decomp-l-r}, 
		\begin{align}
			\PP\left(-\log\sigma(\gamma^p_{2 l^f + 2} \cdots \gamma^p_{2q_n - 2 r_n^h - 2}) \le \alpha n\right) & \le 2C_0e^{-\beta_0} + \sum_{i \le \delta n, j \ge (\lambda''-\delta) n} C_1e^{-\beta_1(j-i)} \notag\\
			\text{idem}& \le 2C_0e^{-\beta_0} + \sum_{i \le \delta n} \frac{C_1}{1-e^{-\beta_1}} e^{-\beta_1((\lambda''-\delta)n-i)} \notag \\
			\text{idem}& \le 2C_0e^{-\beta_0} + \frac{C_1}{(1-e^{-\beta_1})^2}e^{-\beta_1(\lambda''- 2\delta)n}.\label{ccl-ldev}
		\end{align}
		Therefore, the random sequence $(-\log\sigma(\gamma^p_{2 l^f + 2} \cdots \gamma^p_{2q_n - 2 r_n^h - 2}))_n$ satisfies large deviations inequalities below the speed $\lambda$. 
		The exact same computation holds when we consider the sequence $(-\log\sigma(\gamma^p_{2 c_n + 2} \cdots \gamma^p_{2q_n - 2 c_n - 2}))_n$ instead. 
	\end{proof}

	\subsection{Law of large numbers for the singular gap}

	This paragraph is dedicated to the proof of Theorem \ref{th:escspeed}.	
	With the framework developed in the previous paragraphs the large deviation inequalities below the speed $\lambda$ of Lemma \ref{lem:ldev-bloc-central} come naturally.
	However, we can not deduce almost sure convergence from one sided large deviations estimates.
	For that, we use the following Lemma.
	
	\begin{Lem}\label{lem:left-ali-tjrs}
		Let $\nu$ be a probability distribution on $\End(E)$ with positive eventual rank.
		Let $\alpha, \eps > 0$ and $m \in \NN$ be such that $\nu$ is $(\alpha, \eps, m)$-contracting.
		Let $(\gamma_n) \sim \nu^{\otimes\NN}$. 
		There exist constants $l_0 \in \NN$, and $0 < \beta < 1$ such that:
		\begin{equation}\label{ghi}
			\forall f \in E \cup \End(E) \setminus\{0\},\;\PP\left(\forall n\ge l_0,\; f \AA^{\frac{\eps}{2}} \overline{\gamma}_n \right) > \beta.
		\end{equation}
	\end{Lem}
	
	\begin{proof}
		Let $(p_k)$ be as in Theorem \ref{th:pivot}.
		Let $(\gamma_{-m}, \dots,\gamma_{-1}) \sim \nu^{\otimes m}$ be independent of the joint data of $(\gamma_n)_{n \ge 0}$ and $(p_k)_{k \ge 0}$.
		For all $n \ge 0$, we write $q_n:= \max\{k \ge 0\mid \overline{p}_{2k} \le n\}$.
		Then, by Lemma \ref{lem:l-r-piv}, for all $n$, there exists a random integer $r_n = r_n^{\mathrm{Id}} \sim \mathcal{G}_{1/4}$ such that $\gamma^p_{2q_n - 2r_n - 1} \AA^{\eps} \gamma_{\overline{p}_{2q_n - 2r_n}} \cdots \gamma_{n - 1}$ or $r_n \ge q_n$.
		
		Let us apply \eqref{pivot:herali} in Theorem \ref{th:pivot} to $i, j, k = 0, 2, 2q_n - 2r_n$, and $h := \gamma_{\overline{p}_{2q_n - 2r_n}} \cdots \gamma_{n - 1}$.
		Then for all $n$ such that $r_n < q_n$, it holds that $\gamma^p_{0} \gamma^p_{1} \AA^{\frac{\eps}{2}} \gamma_{\overline{p}_2} \cdots \gamma_{n - 1}$. 
		By Lemma \ref{lem:ldev-bloc-central}, $\PP(q_n \le r_n + 1)$ decreases exponentially fast so there exists an integer $l_0$ such that $\PP\left(\exists n \ge l_0 - m, q_n \le r_n + 1\right) \le \frac{1}{2}$.
		Set
		\begin{equation*}
			A_1 := \left(\forall n \ge l_0, \gamma^p_{0} \gamma^p_{1} \AA^{\frac{\eps}{2}} \gamma_{\overline{p}_2} \cdots \gamma_{n - m - 1}\right).
		\end{equation*}
		The event $A_1$ is independent of $(\gamma_{-m}, \dots, \gamma_{-1})$ and $\PP(A_1) \ge 1/2$ by construction of $l_0$.
		
		Let $f \in E \cup \End(E)$.
		By Corollary \ref{cor:schottky} $\nu^{*m}\ge \alpha \mu$ and since $\mu$ is $\frac{1}{4}$-Schottky,
		\begin{align*}
			\PP\left(f \AA^\eps \gamma_{-m}\cdots\gamma_{-1} \AA^\eps \gamma^p_0\gamma^p_1\,\middle|\,(\widetilde{\gamma}^p_k)_{k \ge 0}\right) & = \nu^{* m}\left\{s\,\middle|\,f \AA^\eps \, s \, \AA^\eps \gamma^p_0\gamma^p_1\right\} \\
			& \ge \alpha \mu\left\{s\,\middle|\,f \AA^\eps \, s \, \AA^\eps \gamma^p_0\gamma^p_1\right\} \ge \alpha/2.
		\end{align*}
		Moreover $A_1$ only depends on $(\widetilde{\gamma}^p_k)_{k \ge 0}$, \ie $A_1 \in \langle\widetilde{\gamma}^p_k\rangle^\sigma_{k \ge 0}$.
		Therefore,
		\begin{equation*}
			\PP\left(\forall n \ge l_0,\; f \AA^\eps \gamma_{-m}\cdots\gamma_{-1} \AA^\eps \gamma^p_0\gamma^p_1 \AA^{\frac{\eps}{2}} \gamma_{\overline{p}_2} \cdots \gamma_{n - m - 1}\right) \ge \alpha/4.
		\end{equation*}
		Furthermore, $\gamma^p_0 \AA^\frac{\eps}{2} \gamma^p_1$ so $\sigma(\gamma^p_0 \gamma^p_1) \le \eps^4/12$ and by construction, $\sigma(\gamma_{-m}\cdots \gamma_{-1}) \le \eps^2/12$. 
		Hence, by Lemma \ref{lem:herali} applied twice, the alignment $f\AA^{\frac{\eps}{2}}\gamma_{-m}\cdots\gamma_{n-m -1}$ holds as soon as $f \AA^\eps \gamma_{-m}\cdots\gamma_{-1} \AA^\eps \gamma^p_0\gamma^p_1 \AA^{\frac{\eps}{2}} \gamma_{\overline{p}_2} \cdots \gamma_{n - m - 1}$. 
		It follows that
		\begin{equation*}
			\PP\left(\forall n \ge l_0,\; f \AA^\frac{\eps}{2} \gamma_{-m}\cdots\gamma_{n-m-1}\right) \ge \alpha/4
		\end{equation*}
		and $(\gamma_{n-m})_{n \ge 0} \sim \nu^{\otimes\NN}$, which yields \eqref{ghi} for $\beta = \frac{\alpha}{4}$.
	\end{proof}
	
	Note that Lemma \ref{lem:left-ali-tjrs} is not a direct consequence of Lemma \ref{lem:l-r-piv} alone since we may have $f \gamma_0^p \cdots \gamma^p_{2l^f} = 0$ making the alignment condition meaningless.
	A reformulation of Lemma \ref{lem:left-ali-tjrs} is that the law of $\overline{\gamma}_n$ is eventually $(1 - \beta)$-Schottky for $\AA^\eps$ (\ie it is for all $n \ge l_0$).
	To get the Schottky property on both sides, we need to apply Lemma \ref{lem:left-ali-tjrs} to both $\nu$ and $^t\nu$ but we get the same $\beta$ and $l_0$, because the notion of $(\alpha, \eps, m)$-contraction is symmetric.
	In the invertible case, we moreover know that for all $\rho > 0$, there exists $\eps > 0$ such that the law of $\overline{\gamma}_n$ is eventually $\rho$-Schottky for $\AA^\eps$ this result is due to Guivarc'h and Raugi.
	By \eqref{regxi} in Corollary \eqref{cor:regxi}, it is sufficient to assume that $\zeta^{C, \beta}_\nu(|\log(\eps)|, +\infty) < \rho$ for the law of $\overline{\gamma}_n$ to be eventually $\rho$-Schottky for $\AA^\eps$, yielding a quantitative version of Guivarc'h and Raugi's result.

	\begin{proof}[Proof of Theorem \ref{th:escspeed}]
		Let $(p_k)_k$ be as in Theorem \ref{th:pivot}.
		Let $\lambda$ be as in Lemma \ref{lem:ldev-bloc-central}.
		We want to show that $\lambda_{1,2}(\nu) = \lambda$ satisfies Theorem \ref{th:escspeed}.
		We start with the proof of \eqref{ldsqz}.
		For all $n \ge 0$, let $q_n$, and $r_n = r_n^\mathrm{Id}$ be as in Lemma \ref{lem:l-r-piv}.
		We claim that:
		\begin{equation}
			\sigma(\overline{\gamma}_n) \le \sigma\left(\gamma^p_{0} \cdots \gamma^p_{2q_n - 2 r_n - 2}\right).
		\end{equation}
		Indeed, $\gamma^p_{2q_n - 2 r_n - 2} \AA^\eps \gamma_{\overline{p}_{2q_n - 2 r_n}}\cdots \gamma_{n-1}$ so point \eqref{pivot:herali} in Theorem \ref{th:pivot} yields $\gamma^p_{0} \cdots \allowbreak \gamma^p_{2q_n - 2 r_n - 2} \AA^{\frac{\eps}{2}} \gamma_{\overline{p}_{2q_n - 2 r_n - 1}}\cdot \gamma_{n-1}$.
		Then, by Lemma \ref{lem:c-prod} $\sigma(\gamma_{\overline{p}_{2q_n - 2 r_n - 1}} \allowbreak\cdots \gamma_{n-1}) \allowbreak \le \eps^4/48$ and applying Lemma \ref{lem:c-prod} again yields $\sigma(\overline{\gamma}_n) \le (\eps^2/12) \allowbreak \sigma(\gamma^p_{0} \cdots\allowbreak \gamma^p_{2q_n - 2 r_n - 2})$.
		
		Thus, by Lemma \ref{lem:ldev-bloc-central}, the sequence $\left(-\log\sigma(\overline{\gamma}_n)\right)$ satisfies large deviations inequalities below the speed $\lambda$.
		Then, by Borel Cantelli's Lemma, $\liminf_n \frac{-1}{n} \log \sigma(\overline{\gamma}_n) \ge \lambda$ almost surely. 
		When $\lambda = +\infty$ this suffices to show almost sure convergence.
		Otherwise, we need to prove that $\limsup_n \frac{-\log\sigma(\overline{\gamma}_n)}{n} \le \lambda$ almost surely.
		Let us do it by contraposition. 
		
		By Lemma \ref{lem:pre-speed}, $\frac{-\log\sigma(\gamma^p_2 \cdots \gamma^p_{2 k})}{k} \to \lambda' = ({\EE(p_2) + m})\lambda$ almost surely and in expectation.
		For all $n \ge 0$, we define $k_n:= \min\{k \ge 0 \mid  \overline{p}_{2k} > n\} = q_n + 1$.
		We have $k_n \to +\infty$ almost surely and $\overline{p}_{2k_n - 2} \le n \le \overline{p}_{2 k_n}$.
		Hence, by the law of large numbers, $\frac{k_n}{n} \to \frac{1}{\EE(p_2) + m}$ almost surely.
		Then by Lemma \ref{lem:pre-speed}, $\frac{-\log\sigma(\gamma^p_2 \cdots \gamma^p_{2 k_n})}{n} \to \lambda$.
		In other words, for all $\lambda'' > \lambda$,
		\begin{equation}
			\lim_{n_0 \to + \infty}\PP\left(\exists n \ge n_0,\; -\log\sigma(\gamma^p_2 \cdots \gamma^p_{2 k_n}) > \lambda'' n\right) = 0.
		\end{equation}
		Let $\lambda'' > \lambda$ and assume by contraposition that for all $n_0$,
		\begin{equation*}
			\PP\left(\exists n \ge n_0,\; -\log\sigma(\overline{\gamma}_n) > \lambda'' n\right) \ge \delta.
		\end{equation*}
		Let $l_0, \beta > 0$ be as in Lemma \ref{lem:left-ali-tjrs}.
		Assume that \eqref{ghi} holds for both $\nu$ and for the push-forward $^t\nu$ of $\nu$ by the transposition map.
		Let $l_1$ be such that $\PP(p_0 + p_1 \ge l_1) \le \beta^2\delta / 3$
		and let $n_0$ be such that:
		\begin{equation}\label{enzero}
			\PP\left(\exists n \ge n_0,\; -\log\sigma(\gamma^p_2 \cdots \gamma^p_{2 k_n}) > \lambda'' (n - 2l_0 - l_1) - \log(32\eps^{-4}) \right) < \beta^2\delta / 3.
		\end{equation}
		Such a $n_0$ exists because $\lim_{n}\frac{-\log\sigma(\gamma^p_2 \cdots \gamma^p_{2 k_n}) + C}{n - C} = \lambda$ almost surely and for all constant $C$.
		Let:
		\begin{equation*}
			n_1 = \min\left\{ n \ge n_0,\; -\log\sigma(\gamma_{l_0 + l_1} \cdots {\gamma}_{l_0 + l_1 + n - 1}) > \lambda'' n\right\}.
		\end{equation*}
		That way, $\PP(n_1 \le + \infty) \ge \delta$.
		Let $g_1 = \gamma_{l_0 + l_1} \cdots {\gamma}_{l_0 + l_1 + n_1 - 1}$.
		Let us show that:
		\begin{equation}\label{ali-limsup}
			\PP\left(\forall i < l_1,\, \forall j \ge n_1 + l_1 + 2l_0,\,\gamma_i \cdots \gamma_{l_1 + l_0 - 1} \AA^\frac{\eps}{2} g_1 \AA^\frac{\eps}{2} {\gamma}_{l_0 + l_1 + n} \cdots \gamma_{j-1} \right) \ge \delta \beta^2.
		\end{equation}
		First note that the data of $(\gamma_k)_{k < l_0 + l_1}$ is independent of $g_1$ and $({}^t\gamma_k)_{k < l_0 + l_1} \sim {}^t\nu^{\otimes l_0 +l_1}$.
		Then, by \eqref{ghi} applied to ${}^t \nu$,
		\begin{equation}\label{gauche}
			\PP\left(\forall i' \ge l_0, \, {}^t g_1 \AA^\frac{\eps}{2} {}^t\gamma_{l_0 + l_1 - 1} \cdots {}^t\gamma_{l_0 + l_1 - i'} \,\middle|\, n_1, g_1\right) \ge \beta.
		\end{equation}
		Moreover $l_0 + l_1 + n_1$ is a stopping time with respect to the filtration associated to $(\gamma_k)_k$ (\ie $(l_0 + l_1 + n_1 = n)\in\langle\gamma_k\rangle^\sigma_{k < n}$ for all $n$). 
		Therefore, the conditional distribution of $(\gamma_{k + n_1 + l_0 + l_1})_{k \ge 0}$ with respect to $(\gamma_k)_{k < l_0 + l_1 + n_1}$ is $\nu^{\otimes\NN}$ and therefore,
		\begin{equation}\label{droite}
			\PP\left(\forall j' \ge l_0, \, g_1 \AA^\frac{\eps}{2} \gamma_{l_0 + l_1 + n_1} \cdots \gamma_{l_0 + l_1 + n_1 + j' - 1} \,\middle|\,(\gamma_k)_{k < l_0 + l_1 + n_1} \right) \ge \beta.
		\end{equation}
		Then we compose \eqref{droite} with \eqref{gauche} and with the fact that $n_1$ exists with probability at least $\delta$ and we get \eqref{ali-limsup}.
		Now note that \eqref{ali-limsup} combined with the fact that $\overline{p}_2 < l_1$ outside of an event of probability $1-\beta^2\delta/3$ implies that:
		\begin{equation*}
			\PP\left(\gamma_{\overline{p}_2} \cdots \gamma_{l_1 + l_0 - 1} \AA^\frac{\eps}{2} g_1 \AA^\frac{\eps}{2} {\gamma}_{l_0 + l_1 + n_1} \cdots \gamma_{\overline{p}_{2k_{2l_0 + l_1 + n_1}} - 1} \right) \ge 2\delta \beta^2 / 3.
		\end{equation*}
		Then by Lemma \ref{lem:triple-ali},
		\begin{equation*}
			\PP\left(\sigma(\gamma_{\overline{p}_2} \cdots \gamma_{\overline{p}_{2k_{2l_0 + l_1 + n_1}} - 1}) \le 32\eps^{-4} \sigma(g_1) \right) \ge 2\delta \beta^2/3.
		\end{equation*}
		Moreover, $-\log\sigma(g_1) > \lambda'' n_1$ and $n_1 \ge n_0$ so
		\begin{equation*}
			\PP\left(\exists n \ge n_0,\; -\log\sigma(\gamma^p_2 \cdots \gamma^p_{2 k_n}) > \lambda'' (n - 2l_0 - l_1) - \log(32\eps^{-4}) \right) \ge 2 \beta^2 \delta / 3.
		\end{equation*}
		This contradicts the definition of $n_0$ given by \eqref{enzero}. 
		Finally, by contraposition, $-\log\sigma(\overline\gamma_n) / n \to \lambda$ almost surely.
	\end{proof}

	\subsection{Contraction property}\label{section:contraction}

	The present section is dedicated to the proof of Theorems \ref{th:cvspeed} and \ref{th:eigenspace}.
	We also give a more precise statement of Theorem \ref{th:cvspeed} in Theorem \ref{th:limit-uv}.
	Most technicities of the present section come from the fact that we work with measures that are not supported on $\GL(E)$.
	For example we need to show that we can construct a map $\ell^\infty: \gamma \mapsto \ell^\infty(\gamma)$ that is equivariant and does not depend on $\nu$ since $\lim \eigen^+(\overline\gamma_n)$ is not equivariant and $\lim [\overline\gamma_n v]$ depends on the choice of $v \notin \essker(\nu)$
	To deal with this first technical point, we introduce the following notion of limit line.
	
	\begin{Def}[The space of contracting sequences]\label{def:l-infty}
		Let $E$ be a Euclidean, Hermitian or ultrametric vector space. We define the set $\Omega'(E)$ of contracting sequences as follows:
		\begin{equation*}
			\Omega'(E):= \left\{ \gamma \in \End(E)^\NN\,\middle|\,
			\begin{aligned}
				& \forall n \in\NN,\; \overline{\gamma}_n \neq \{0\} \quad\text{and} \quad \exists \ell^\infty \in \proj(E),\\
				& \forall \eps \in (0,1), \lim_{n}\max_{u \,\in\, U^\eps(\overline{\gamma}_n)} \dist(\ell^\infty, [u]) = 0
			\end{aligned}\right\}.
		\end{equation*}
		We define $T:= \End(E)^\NN \to \End(E)^\NN$ to be the Bernoulli shift and we define $\ell^\infty:\Omega'(E) \to \proj(E)$ to be the unique map such that:
		\begin{equation*}
			\forall \gamma \in \Omega'(E),\; \lim_{n \to +\infty} \max_{u \, \in \, U^\eps(\overline{\gamma}_n) \setminus\{0\}} \dist(\ell^\infty(\gamma), [u]) = 0.
		\end{equation*}
		We define the shift-invariant space of contracting sequences as:
		\begin{equation*}
			\Omega(E):= \left\{ \gamma \in \bigcap_{k = 0}^{+ \infty} T^{-k}\Omega'(E) \,\middle|\, \forall k\in\NN, \gamma_k \ell^\infty(T^{k+1}\gamma) = \ell^\infty(T^k\gamma) \right\}.
		\end{equation*}
	\end{Def}
	
	That way the space $\Omega(E)$ is $T$-invariant and measurable and $\ell^\infty$ is measurable.
	Indeed given $(e_1, \dots, e_d)$ a basis of $E$, for $\eps$ small enough and  for all $g \in\End(E)$ there exists at least one index $1 \le i \le d$ such that $e_i \in V^\eps(g)$.
	So we can define $\ell^\infty(\gamma) = \lim_n [\overline{\gamma} e_{i_n}]$ with $i_n  = \min (\arg \max_j\|\overline{\gamma}_n e_j\|)$.
	Moreover $\ell^\infty$ is $T$-equivariant on $\Omega(E)$ in the sense that $\ell^\infty(\gamma) = \gamma_0 \ell^\infty(T\gamma)$.
	
	Note that $\Omega'(E) \neq \Omega(E)$. 
	For example set $E = \RR^2$, and denote by $\pi_1$ and $\pi_2$ the orthogonal projections onto the first and second coordinates. 
	If $\gamma_0 = \pi_1$ and $\gamma_k = \pi_1 +2\pi_2$ for all $k \ge 1$, then $\overline{\gamma}_n = \pi_1$ for all $n \ge 1$ and $[\gamma_k\cdots \gamma_n] \underset{m \to \infty}{\longrightarrow} [\pi_2]$ for all $k \ge 1$.
	So $\ell^\infty(\gamma)$ is the first coordinate axis and $\ell^\infty(T\gamma)$ is the second coordinate axis.
	Hence $\gamma_0 \ell^\infty(T\gamma) = [0] \neq \ell^\infty(\gamma)$ so $\gamma \in \Omega'(E) \setminus \Omega(E)$.

	\begin{Th}\label{th:limit-uv}
		Let $\nu$ be a probability distribution on $\End(E)$ with positive eventual rank.
		Let $\alpha, \eps > 0$ and $m \in \NN$ be such that $\nu$ is $(\alpha, \eps, m)$-contracting.
		Let $\gamma = (\gamma_n)_{n\in\NN}\sim\nu^{\otimes\NN}$. 
		Then $\gamma \in \Omega(E)$ almost surely. 
		Let $\alpha < \lambda_{1,2}(\nu)$. 
		There exist constants $C,\beta > 0$ such that
		\begin{equation}\label{limit-u}
			\PP\left(\exists u \in U^1(\overline{\gamma}_n)\setminus\{0\},\; \dist([u], \ell^\infty(\gamma))\ge e^{-\alpha n}\right) \le C e^{-\beta n}.
		\end{equation}
		Moreover, for all $v \in E$,
		\begin{equation}\label{limit-v}
			\PP\left(\dist([\overline{\gamma}_n v], \ell^\infty(\gamma))\ge e^{-\alpha n}\mid  \overline{\gamma}_n v \neq 0\right) \le C e^{-\beta n}.
		\end{equation}
	\end{Th}
	
	\begin{proof}
		Let $(p_n)_{n \ge 0}$ be as in Theorem \ref{th:pivot}.
		Let $v \in E \setminus\{0\}$.
		For all $n \in\NN$, let $q_n$ and $r_n^v$ and $r_n^{\mathrm{Id}}$ be as in Lemma \ref{lem:l-r-piv}.
		Let $\alpha < \sigma(\nu)$ and let $C, \beta > 0$ be as in Lemma \ref{lem:ldev-bloc-central}, \ie such that
		\begin{equation}
			\forall n,\,\forall h,\,\PP\left( - \log \sigma\left(\gamma^p_0 \cdots \gamma^p_{2q_n - 2r_n^{h} - 2}\right) \le \alpha n \right) \le C e^{- \beta n}.
		\end{equation}
		For all $n$, set $j'_{n} := \min\{2q_n - 2r_n^{v}, 2q_n - 2r_n^{\mathrm{Id}}\}$.
		Then
		\begin{equation}\label{hgsdghvj}
			\forall n, \;	\PP\left(\exists k \ge n,\; - \log \sigma\left(\overline{\gamma}^p_{j'_{n}}\right) \le \alpha n \right) \le \frac{2C}{1-e^{-\beta}} e^{- \beta n}.
		\end{equation}
		Let $j_n := \min_{k \ge n} j'_{k}$ for all $n$.
		By \eqref{hgsdghvj}, the random integer $n_\alpha := \max\{n\mid  \sigma\left(\overline{\gamma}^p_{j_{n}}\right) \ge e^{-\alpha n}\}$ has a finite exponential moment.
		By point \eqref{pivot:herali} in Theorem \ref{th:pivot}, the alignment $\overline{\gamma}^p_{j_{n}} \AA^\frac{\eps}{2} \gamma_{\overline{p}_{j_{n}}}\cdots \gamma_{k - 1} h$ holds for all $k \ge n$ and $h \in \{v,\mathrm{Id}\}$.
		Then by Lemma \ref{lem:c-prod} and \eqref{cont-prop}, for all $u \in U^1\left(\overline{\gamma}^p_{j_{n}}\right) \setminus \{0\}$ and for all $u' \in [\overline{\gamma}_{k} v] \cup U^\eps\left(\overline{\gamma}_{k}\right)\setminus\{0\}$ for $k \ge n$,
		\begin{equation}\label{hsdjk}
			\dist\left([u], [u']\right) \le \frac{2}{\eps} \sigma\left(\overline{\gamma}^p_{j_{n}}\right) \le \frac{4e^{-\alpha n}}{\eps} + \mathds{1}_{n < n_\alpha}.
		\end{equation}
		Let $u_{n}^{\eps'} \in U^{\eps'}\left(\overline{\gamma}_n\right) \setminus \{0\}$ for all $n \in\NN$ and all $\eps > 0$. 
		Note that \eqref{hsdjk} holds for all $0 < \eps' < \eps$.
		Therefore, $\dist(u_n, u_k) \le 4 e^{-\alpha n} / \eps'$ for all $n_0 \le n \le k$, by \eqref{hsdjk}.
		Therefore $(u_n^{\eps'})_{n \ge 0}$ a Cauchy sequence on the full measure set $(n_\alpha < +\infty)$ so $(\gamma_n)_{n \ge 0} \in \omega'(E)$ on $(n_\alpha < +\infty)$ and therefore $\nu^{\otimes \NN}(\Omega'(E)) = 1$.
		Since \eqref{hsdjk} holds for all $\alpha < \lambda_{1,2}$, the random sequence $\left(-\log\dist\left([u_n^1], \ell^\infty(\Omega)\right)\right)_{n \ge 0}$ satisfies large deviations inequalities below the speed $\lambda_{1,2}$, which implies \eqref{limit-u}.
		
		When $v \notin \essker(\nu)$, the random sequence $\left(-\log\dist\left([\overline{\gamma}_n v], \ell^\infty\right)\right)_{n \ge 0}$ satisfies large deviations inequalities below the speed $\lambda_{1,2}$, which implies \eqref{limit-v}.
		Then, $[\overline{\gamma}_n v] \to \ell^\infty(\gamma)$ and it follows from the same reasoning applied to $(\overline{\gamma}_n)_{n \ge 1}$ that $[\gamma_1 \cdots {\gamma}_n v] \to \ell^\infty(T\gamma)$.
		Hence, $\gamma_0 \ell^\infty(T\gamma) = \ell^\infty(\gamma)$ almost surely so $\nu^{\otimes\NN} (\Omega) = 1$
		By the same argument, for all $\alpha < \lambda_{1,2}$, there exist constants $C', \beta' > 0$ such that:
		\begin{equation*}
			\PP\left(\dist([\overline{\gamma}_n v], \ell^\infty(\gamma))\ge e^{-\alpha n} \cap \overline{\gamma}_n v \neq 0\right) \le C' e^{-\beta' n}.
		\end{equation*}
		By Lemma \ref{lem:descri-ker}, there exists $\delta > 0$ such that $\PP\left(\overline{\gamma}_n v \neq 0 \right) \ge \delta$ for all $n$ and all $v$, which implies \eqref{limit-v}.
	\end{proof}

	Let us now use the $c_n$ of Lemma \ref{lem:l-r-piv} to prove that the dominant eigenspace also converges.
	
	\begin{proof}[Proof of Theorem~\ref{th:eigenspace}]
		The key is to apply Lemma \ref{lem:eigen-align} to a well chosen endomorphism. 
		For that, let $\alpha, \eps > 0$ and $m \in \NN$ be so that $\nu$ is $(\alpha, \eps, m)$-contracting in the sense of Definition \ref{def:quant-contraction}.
		Let $\gamma = (\gamma_n)_{n \ge 0} \sim \nu^{\otimes\NN}$ and $(p_n)_n$ be as in Theorem \ref{th:pivot}.
		For all $n$, let $q_n$ and $c_n$ be as in Lemma \ref{lem:l-r-piv}.
		When $2 c_n + 1 < q_n$, set
		\begin{gather*}
			f_n:= \gamma_{\overline{p}_{2q_n - 2c_n - 1}} \cdots \gamma_{n-1} \gamma^p_0 \cdots \gamma^p_{2 c_n}
			\\
			h_n:= \gamma_{\overline{p}_{2q_n - 2c_n}} \cdots \gamma_{n-1} \gamma_0^p \cdots \gamma^p_{2 c_n}
			\\
			z_n:= \gamma^p_{2 c_n + 1} \cdots \gamma^p_{2q_n - 2c_n - 2}.
		\end{gather*}
		By definition of $c_n$, $\gamma^p_{2q_n - 2c_n - 1} \AA^\eps h_n$ so by point \eqref{pivot:schottky} in Theorem \ref{th:pivot} and Lemma \ref{lem:c-prod}, $\sigma(f_n) \le \eps^4/48$ and by point \eqref{pivot:herali} in Theorem \ref{th:pivot}, $z_n \AA^\frac{\eps}{2} f_n$.
		By definition of $c_n$ again, $f_n \AA^\eps \gamma^p_{2 c_n + 1}$ so by point \eqref{pivot:herali} in Theorem \ref{th:pivot}, $f_n \AA^\frac{\eps}{2} z_n$.
		Moreover, $\gamma^p_{2 c_n + 1} \AA^\frac{\eps}{2} \cdots \gamma^p_{2q_n - 2c_n - 2}$ by Lemma \ref{lem:herali} so $\sigma(z_n) \le \eps^4/12$ by Lemma \ref{lem:c-prod}.
		Then Lemma \ref{lem:herali} yields that $f_n z_n \AA^\frac{\eps}{4} f_n z_n$ and $\sigma(f_n z_n) \le 4 \eps^{-2} \sigma(f_n) \sigma(z_n)$.
		It follows from \eqref{sqz-prox} in Lemma \ref{lem:eigen-align} that
		\begin{equation*}
			\frac{\rho_2}{\rho_1}(\overline{\gamma}_n) = \frac{\rho_2}{\rho_1}(f_n z_n) \le \frac{64 \sigma(f_n z_n)}{\eps^2} \le \sigma\left(\gamma^p_{2 c_n + 2} \cdots \gamma^p_{2q_n - 2c_n - 2}\right) 
		\end{equation*}
		When $2 c_n + 1 \ge q_n$, we still have $\frac{\rho_2}{\rho_1}(\overline{\gamma}_n) \le 1$.
		Then by Lemma \ref{lem:ldev-bloc-central}, for $\lambda_{1,2} = \log(\rho_1 / \rho_2)$, the random sequence $\left(\lambda_{1,2}(\overline{\gamma}_n)\right)_{n \ge 0}$ satisfies large deviations inequalities below the speed $\lambda_{1,2}(\nu)$.
		
		Let us now prove the second part of \eqref{eigenspace}. Namely, we show that the random sequence $(-\log\dist(\eigen^+(\overline{\gamma}_n), \ell^\infty(\gamma))_{n \ge 0}$ satisfies large deviations inequalities below the speed $\lambda_{1,2}(\nu)$.
		Let $n \ge 0$ and assume that $2 c_n + 1 < q_n$.
		If follows from point \eqref{pivot:herali} in Theorem \ref{th:pivot} that $\overline{\gamma}^p_{2q_n - 2c_n - 1} \AA^\frac{\eps}{2} f_n$.
		Then by Lemma~\ref{lem:alipart}, $\overline{\gamma}^p_{2q_n - 2c_n - 1} \AA^\frac{\eps}{4} \allowbreak (f_n z_n)^k$ for all $k \ge 0$. 
		Moreover, for all choice of $e_n \in \eigen^+(f_n z_n) \setminus\{0\}$ and $e'_n \in \eigen^+(^tz_n ^tf_n) \setminus\{0\}$, $[(f_n z_n)^m] \underset{k \to \infty}{\longrightarrow} [e e']$.
		Furthermore, $\AA^\eps$ is bi-homogeneous and closed so $\overline{\gamma}^p_{2q_n - 2c_n - 1} \AA^\frac{\eps}{4} e e'$ and therefore $\overline{\gamma}^p_{2q_n - 2c_n - 1} \AA^\frac{\eps}{4} e$.
		
		Note that $\overline{\gamma}_n$ and $f_n z_n$ are semi-conjugated by $\overline{\gamma}^p_{2q_n - 2c_n - 1}$ so 
		\begin{equation*}
			\eigen^+(\overline{\gamma}_n) = \overline{\gamma}^p_{2q_n - 2c_n - 1} \eigen^+(f_n z_n) = \KK \overline{\gamma}^p_{2q_n - 2c_n - 1} e.
		\end{equation*}
		Therefore, $\eigen^+(\overline{\gamma}_n) \subset U^\frac{\eps}{4}(\overline{\gamma}^p_{2q_n - 2c_n - 1})$.
		Moreover $\overline{\gamma}^p_{2q_n - 2c_n - 1} \AA^\frac{\eps}{2} \gamma^p_{2c_n +1} \cdots \gamma^p_k$ for all $k$ and $\AA^\frac{\eps}{2}$ is closed, so going to the limit yields $\ell^\infty((\gamma_{n + \overline{p}_{2q_n - 2c_n - 1}})_{n \ge 0}) \subset U^\frac{\eps}{4}(\overline{\gamma}^p_{2q_n - 2c_n - 1})$.
		Then by Lemma \ref{lem:cont-prop} and by triangular inequality,
		\begin{equation*}
			\dist(\eigen^+(\overline{\gamma}_n), \ell^\infty) \le \frac{6 \sigma\left(\overline{\gamma}^p_{2q_n - 2c_n - 1}\right)}{\eps} \le \sigma\left(\gamma^p_{2 c_n + 2} \cdots \gamma^p_{2q_n - 2c_n - 2}\right)
		\end{equation*}
		We conclude using Lemma \ref{lem:ldev-bloc-central} again.
	\end{proof}
	
	A direct consequence of Theorem~\ref{th:eigenspace} is that $\liminf_n \frac{1}{n}\log\frac{\rho_1}{\rho_2}(\overline{\gamma}_n) \ge \lambda_{1,2}(\nu)$ almost surely. 
	However, for $d \ge 3$ and without moment condition, we do not know anything about the upper limit. 
	The issue is that the pivoting technique does not tell us much about ${\rho_2}(\overline{\gamma}_n)$ due to the fact that $\bigwedge^2_* \nu$ may not be strongly irreducible. 
	Using a reduction method like for the proof of Corollary \ref{cor:gen-eigen} and the concentration results of \eqref{dom-radius} in Theorem~\ref{th:slln}, one may conjecture that for $\nu$ strongly irreducible, proximal and supported on $\GL(E)$, it moreover holds that $\frac{1}{n}\log\frac{\rho_1}{\rho_2}(\overline{\gamma}_n) \rightharpoonup \lambda_{1,2}(\nu)$ in probability.
	Asking for almost sure convergence seems a bit far fetched in this case.

	\bibliographystyle{alpha}
	\bibliography{biblio-porm.bib}

\end{document}